\documentclass[leqno,11pt]{amsart}
\usepackage[colorlinks,pagebackref,hypertexnames=false]{hyperref}
\hypersetup{
    allcolors = blue,
    pdftitle = {Large mass limits of G2 and Calabi--Yau monopoles: calibrated concentration, Higgs zeros, and abelianization},
    pdfauthor = {Daniel Fadel and Goncalo Oliveira},
    pdfsubject = {Geometric analysis, gauge theory, calibrated geometry, and special holonomy},
    pdfkeywords = {Yang--Mills--Higgs theory, G2-monopoles, Calabi--Yau monopoles, calibrated currents}
}

\usepackage{amsmath,amsthm,amssymb,mathrsfs,mathtools}
\usepackage{esint} 
\usepackage[alphabetic,backrefs]{amsrefs}
\usepackage[T1]{fontenc}
\usepackage[english]{babel}
\usepackage[margin=1in]{geometry}
\usepackage{doi}

\usepackage{xcolor}

\newcounter{dummy}
\usepackage{enumitem}
\makeatletter
\newcommand\myitem[1][]{\item[#1]\refstepcounter{dummy}\def\@currentlabel{#1}}
\makeatother

\usepackage{bookmark}

\definecolor{tocsectioncolor}{RGB}{24,49,83}
\definecolor{tocsubsectioncolor}{RGB}{80,80,80}
\definecolor{tocpagenumcolor}{RGB}{80,80,80}

\makeatletter
\def\l@section{%
  \@tocline{1}{5pt plus 1pt}{0pt}{}%
  {\small\bfseries\color{tocsectioncolor}}%
}
\def\l@subsection{%
  \@tocline{2}{0.5pt plus 0.25pt}{1.5pc}{4.5pc}%
  {\small\color{tocsubsectioncolor}}%
}
\def\@tocpagenum#1{%
  \hss{\mdseries\color{tocpagenumcolor}#1}%
}
\makeatother

\usepackage[nameinlink,noabbrev]{cleveref}
\numberwithin{equation}{section}

\usepackage{comment}

\let\originalleft\left
\let\originalright\right
\renewcommand{\left}{\mathopen{}\mathclose\bgroup\originalleft}
\renewcommand{\right}{\aftergroup\egroup\originalright}

\makeatletter
\renewcommand*{\eqref}[1]{\hyperref[{#1}]{\textup{\tagform@{\ref*{#1}}}}}		
\makeatother

\usepackage{booktabs}

\def\XXint#1#2#3{{\setbox0=\hbox{$#1{#2#3}{\int}$}
    \vcenter{\hbox{$#2#3$}}\kern-.5\wd0}}

\def\YYint#1#2#3{{\setbox0=\hbox{$#1{#2#3}{\int}$}
    \lower1ex\hbox{$#2#3$}\kern-.46\wd0}}
\def\YYYint#1#2#3{{\setbox0=\hbox{$#1{#2#3}{\int}$}
    \lower0.35ex\hbox{$#2#3$}\kern-.48\wd0}}

\def\ZZint#1#2#3{{\setbox0=\hbox{$#1{#2#3}{\int}$}
    \raise1.15ex\hbox{$#2#3$}\kern-.57\wd0}}
\def\ZZZint#1#2#3{{\setbox0=\hbox{$#1{#2#3}{\int}$}
    \raise0.85ex\hbox{$#2#3$}\kern-.53\wd0}}

\DeclareMathOperator\vol{vol}

\newcommand{\cE}{\mathcal{E}}

\usepackage{slashed}

\newcommand{\Ric}{\mathrm{Ric}}

\def\Ric{\mathrm{Ric}}

\def\vol{\mathrm{vol}}

\def\<{\mathopen{}\left<}
\def\>{\right>\mathclose{}}
\def\({\mathopen{}\left(}
\def\){\right)\mathclose{}}

\renewcommand\epsilon{\varepsilon}

\DeclareMathOperator{\inj}{inj}

\newtheorem{theorem}{Theorem}[section]

\newtheorem{introtheorem}{Theorem}

\newtheorem{introcorollary}[introtheorem]{Corollary}

\newtheorem{corollary}[theorem]{Corollary}

\newtheorem{lemma}[theorem]{Lemma}
\newtheorem{proposition}[theorem]{Proposition}
\theoremstyle{definition} \newtheorem{definition}[theorem]{Definition}

\newtheorem{remark}[theorem]{Remark}

\theoremstyle{definition}
\newtheorem*{mainsetting}{Main Setting}

\newcommand{\mainsettingref}{%
    \hyperref[setting:main]{Main Setting}%
}

\begin{document}
	
\author{Daniel Fadel} 
\address[Daniel Fadel]{Instituto de Ci\^{e}ncias Matem\'{a}ticas e de Computa\c{c}\~{a}o, Universidade de S\~ao Paulo, Avenida Trabalhador S\~ao-carlense 400, 13566-590 S\~ao Carlos, SP, Brazil}
\urladdr{\href{https://sites.google.com/view/daniel-fadel-math-homepage/home}{sites.google.com/view/daniel-fadel-math-homepage/home}}
\email{\href{mailto:daniel.fadel@icmc.usp.br}
{daniel.fadel@icmc.usp.br}}

	\author{Gon\c{c}alo Oliveira}
	\address[Gon\c{c}alo Oliveira]{Centro de An\'{a}lise Matem\'{a}tica, Geometria e Sistemas Din\^{a}micos, Departamento de Matem\'{a}tica, Instituto Superior T\'{e}cnico, Universidade de Lisboa, Avenida Rovisco Pais 1, 1049-001 Lisboa, Portugal}
	\urladdr{\href{https://sites.google.com/view/goncalo-oliveira-math-webpage/home}{sites.google.com/view/goncalo-oliveira-math-webpage/home}}
	\email{\href{mailto:galato97@gmail.com}{galato97@gmail.com}}

	\date{August 5, 2026}


\title[Large mass limits of $\mathrm G_2$ and Calabi--Yau monopoles]{Large mass limits of $\mathrm G_2$ and Calabi--Yau monopoles: calibrated concentration, Higgs zeros, and abelianization}

\begin{abstract}
We study large mass monopoles with structure group $\mathrm{SU}(2)$
or $\mathrm{SO}(3)$ on asymptotically conical $\mathrm G_2$-manifolds
and Calabi--Yau $3$-folds, with fixed asymptotic class. After placing
the AC asymptotic theory, the variational compactness theory of
Parise--Pigati--Stern, and Li's singular abelian compactness theory in a
common $\Theta$-monopole framework, we prove that the mass-renormalized
Yang--Mills--Higgs and intermediate energy measures converge to
$8\pi\|T\|$ for a compactly supported calibrated integral
codimension-three cycle $T$. This identifies the two limiting currents
and shows that the variational calibration inequalities are saturated.
Using this common limit as the starting point for a finer analysis, if
$\mathcal S$ is the calibrated support, $\mathcal Z$ the Kuratowski
upper limit of the Higgs zero sets, and $\mathcal C$ the limiting
nonabelian locus, defined as the Kuratowski upper limit of Li's curvature
concentration loci, then
$\mathcal S\subset\mathcal Z\subset\mathcal C=\mathcal S\cup\mathcal O$,
where $\mathcal O$ is precisely the obstruction to effective
codimension-three monotonicity. For the cohomogeneity-one large mass
families on the Bryant--Salamon $\mathrm G_2$-manifolds and the Stenzel
Calabi--Yau $3$-fold, we prove that $\mathcal O=\varnothing$. On $X\setminus\mathcal C$ the sequence
abelianizes; corrected longitudinal curvatures converge smoothly, and
the remaining compactness alternatives are governed by $L^2$-harmonic
$2$-forms.
\end{abstract}

	\maketitle

	{\hypersetup{hidelinks}\tableofcontents}
	
\section{Introduction}\label{sec: intro}

\subsection{Context}

Gauge theory on complete manifolds with special holonomy provides a
natural setting in which calibrated geometry can emerge from large scale
limits of nonlinear elliptic equations
\cites{tian2000gauge,donaldson2009gauge,walpuski2014spin,oliveira2014monopoles,oliveira2016calabi,walpuski2017g2}.
For complete noncompact $\mathrm G_2$-manifolds and Calabi--Yau
$3$-folds, Donaldson and Segal proposed a programme relating
$\mathrm G_2$- and Calabi--Yau monopoles to coassociative and special
Lagrangian geometry, respectively, in analogy with Taubes' correspondence
between Seiberg--Witten theory and pseudoholomorphic curves in dimension
four \cites{donaldson2009gauge,Taubes1996,Taubes1999a,Taubes1999}.

Let $(X^n,g)$ be a complete asymptotically conical manifold endowed with a
closed $(n-3)$-form $\Theta$. Given a principal $G$-bundle $P\to X$, with
$G\in\{\mathrm{SU}(2),\mathrm{SO}(3)\}$, a pair $(\nabla,\Phi)$ consisting
of a connection and a Higgs field is a $\Theta$-monopole if
\[
    F_\nabla\wedge\Theta=*\,\nabla\Phi.
\]
In dimension three this is the Bogomolny equation. It also specializes
to the Calabi--Yau monopole equation on a Calabi--Yau $3$-fold and to
the $\mathrm G_2$-monopole equation on a $\mathrm G_2$-manifold. In the
Calabi--Yau case we additionally impose the primitivity condition
$F_\nabla\wedge\omega^2=0$.

The associated intermediate, or $\Theta$-, energy is
\[
    \mathcal E^\Theta(\nabla,\Phi)
    :=
    \frac12\int_X
    \bigl(|F_\nabla\wedge\Theta|^2+|\nabla\Phi|^2\bigr).
\]
Unlike the full Yang--Mills--Higgs energy
\[
    \mathcal E(\nabla,\Phi)
    :=
    \frac12\int_X
    \bigl(|F_\nabla|^2+|\nabla\Phi|^2\bigr),
\]
which may be infinite even for irreducible higher-dimensional
$\Theta$-monopoles on asymptotically conical manifolds, the
$\Theta$-energy is finite for all currently known irreducible examples
constructed intrinsically in the $\mathrm G_2$ and Calabi--Yau settings.
In the torsion-free settings considered here, $\Theta$-monopoles satisfy
the Yang--Mills--Higgs equations
\[
    \nabla^*\nabla\Phi=0,
    \qquad
    d_\nabla^*F_\nabla=[\nabla\Phi,\Phi].
\]
The large mass regime is characterized by
\[
    \lim_{\rho\to\infty}|\Phi|=m,
    \qquad
    m\to\infty,
\]
together with fixed asymptotic topological data. The Donaldson--Segal
picture predicts that, after normalization by the mass, the energy
concentrates in codimension three along a calibrated cycle. In the
$\mathrm G_2$ case the limiting cycle should be coassociative, and in the
Calabi--Yau case special Lagrangian. Near a smooth point of that cycle,
the monopole should be modelled transversely on a large mass monopole on
$\mathbb R^3$; away from the cycle, the structure group should
abelianize and the limit should be a singular abelian monopole with a
Dirac source along the calibrated support.

The three-dimensional theory initiated in the authors' published work
\cite{fadel2019limit}, developed in corrected and expanded form in
\cite{fadeloliveira2026limitv5}, and continued by the first author in
\cite{fadel2026abelian}, gives a complete description of the degeneration at
finite points. For finite energy $\mathrm{SU}(2)$ monopoles on an
asymptotically conical $3$-manifold, with fixed charge and masses tending
to infinity, the mass-renormalized energy concentrates on a finite set
$\mathcal S$ which agrees with the Kuratowski upper limit $\mathcal Z$ of
the Higgs zero sets. At each point of $\mathcal S$, a complete finite
cluster of mass one Euclidean monopoles carries the full concentration
weight, with no mass-renormalized energy left between the profiles and
the ambient scale
\cite[Theorem~1.1 and Proposition~8.3]{fadeloliveira2026limitv5}. On
$X\setminus\mathcal S$, the fields abelianize exponentially and, after
translation of the Higgs fields by their masses, converge smoothly to a
reducible monopole whose scalar part is the Green potential of the
weighted concentration $0$-cycle and whose singularities are Dirac
singularities with the corresponding total cluster charges
\cite[Theorem~A]{fadel2026abelian}. The remaining charge is precisely the
charge escaping through the AC end
\cite[Corollary~5.2]{fadel2026abelian}. In the natural
three-dimensional analogue of the notation used below, the maximal
abelian region is $X\setminus\mathcal S$, and the corresponding
limiting nonabelian locus satisfies
$\mathcal S=\mathcal Z=\mathcal C^{(3)}$.

Symmetric higher-dimensional examples constructed by the second author
display the analogous geometry. On the Bryant--Salamon
$\mathrm G_2$-manifolds and the Stenzel Calabi--Yau manifold, the zero
section is a compact calibrated codimension-three submanifold, the
$\Theta$-energy concentrates along it, transverse rescalings converge to
monopoles on $\mathbb R^3$, and the sequence converges away from it to an
abelian monopole with a Dirac singularity
\cites{oliveira2014monopoles,oliveira2016calabi}.

\subsection{Relation to previous work and contributions}

The asymptotic foundation is the work of the authors and Ákos Nagy
\cite{fadel2020asymptotic}, hereafter FNO\@. For finite intermediate energy
$\mathrm G_2$-monopoles, FNO proved finite mass under bounded curvature
on nonparabolic manifolds. On an AC $\mathrm G_2$-manifold, assuming
quadratic curvature decay, they obtained sharp asymptotics and a
pseudo-Hermitian--Yang--Mills limit on the link, together with the
topological identity
\[
    \mathcal E^\psi(\nabla,\Phi)
    =
    4\pi m\langle\beta\cup\psi_\infty,[\Sigma]\rangle.
\]
Here $\beta\in H^2(\Sigma;\mathbb R)$ is the asymptotic monopole class.
Li observed that the same arguments apply, with the expected changes, to
Calabi--Yau monopoles \cite{li2025large}. The present work continues the
large mass programme announced in FNO and formulates the required
AC asymptotic theory uniformly for $\Theta$-monopoles.

A second analytic input comes from the first author's AC function theory
for Yang--Mills--Higgs fields \cite{fadel2023asymptotics}, in particular
the decomposition
$\mathscr D^{1,2}(X)=\mathscr D_0^{1,2}(X)\oplus\mathbb R$ and the Green
representation for the Higgs defect. Combined here with the FNO Bochner
estimates, these tools yield the critical Sobolev estimate for
$m-|\Phi|$ and the local largeness estimate needed in the large mass regime to enter
the Parise--Pigati--Stern compactness theory.

Li \cite{li2025large} developed a structure theory for large mass
$\mathrm G_2$ and Calabi--Yau monopoles on AC manifolds. He extracts a
singular abelian monopole with Dirac source carried by a compactly
supported calibrated integral cycle, denoted below by
$T_{\mathrm{Li}}$, and identifies its homology class and mass from the
asymptotic monopole class. He also proves equality of the limiting
mass-renormalized curvature and Higgs gradient measures and introduces
configuration-dependent curvature concentration loci. Li states the
main results for $\mathrm{SU}(2)$ and records the modifications for
$\mathrm{SO}(3)$; in that case the global abelian curvature class can be
half-integral.

Independently, Parise--Pigati--Stern (PPS)
\cite{parise2025nonabelian} proved a variational compactness theorem for
large mass Yang--Mills--Higgs fields. Their charge currents converge to
an integral codimension-three cycle, denoted below by
$T_{\mathrm{PPS}}$. For generalized monopoles they obtain calibration
inequalities and ask when the limiting cycle is calibrated. Their
published theorem is formulated on compact manifolds for the trivial
$\mathrm{SU}(2)$-bundle.

The present paper first places these inputs on a common footing. The
resulting identification furnishes the single calibrated limit from
which we develop the structure theory of the concentration, zero-set,
and abelian regions. The contributions divide into the common
foundation in items~\textup{(i)}--\textup{(ii)} and the principal
structural results in items~\textup{(iii)}--\textup{(iv)}.

\begin{enumerate}[label=\textnormal{(\roman*)}]
\item
\textit{A common AC $\Theta$-monopole framework and the analytic bridge to
PPS.}
We formulate the FNO asymptotic theory for $\mathrm G_2$ and
Calabi--Yau monopoles, prove the required $\mathrm{SO}(3)$ extension,
and keep track of the integral versus half-integral asymptotic charge
lattices. We prove versions of PPS compactness for fixed bundles and for $\mathrm{SO}(3)$, globalize that theorem to AC manifolds, and verify its two
local hypotheses. The Higgs defect estimate described above gives the
local largeness hypothesis; Li's curvature bound at large radius, together
with the fixed class energy identity, gives the local
mass-renormalized energy bound.

\item
\textit{Identification of the common calibrated limit and exact concentration.}
After converting the normalization conventions, we prove
\[
    T_{\mathrm{PPS}}=T_{\mathrm{Li}}=:T,
    \qquad
    \mu=\nu=8\pi\|T\|,
\]
where $\mu$ and $\nu$ are the limiting mass-renormalized
Yang--Mills--Higgs and intermediate energy measures. Thus the PPS
inequalities are saturated, their limiting cycle is calibrated in the
Main Setting, and \cite[Question~1.7]{parise2025nonabelian} is answered
there. Li's formulas then determine $[T]$ and $\|T\|(X)$ from the
asymptotic monopole class.

\item
\textit{A hierarchy of concentration and zero sets.}
Let $\mathcal S=\operatorname{spt}\|T\|$. We introduce the Kuratowski
upper limit $\mathcal Z$ of the Higgs zero sets and the limiting
nonabelian locus $\mathcal C$, defined as the Kuratowski upper limit of
Li's curvature concentration loci. We compare these sets with the usual
subsequential Tian--Uhlenbeck codimension-four concentration sets, and
prove
\[
    \mathcal S\subset\mathcal Z\subset\mathcal C
    =\mathcal S\cup\mathcal O.
\]
Here $\mathcal O$ is a new closed obstruction locus: its complement is
exactly the region where effective codimension-three monotonicity holds
locally. This identity ties the excess of Li's loci over the calibrated
support to the monotonicity problem. A codimension-two Morrey bound for
the positive Yang--Mills--Higgs imbalance excludes $\mathcal O$, and we
verify this bound for every known cohomogeneity-one family, obtaining
$\mathcal S=\mathcal Z=\mathcal C$ in those examples.

\item
\textit{Abelianization and finer compactness away from $\mathcal C$.}
The open set
$\mathcal R_{\mathrm{ab}}:=X\setminus\mathcal C$, called the abelian
region, is the maximal region
on which Li's loci are eventually absent locally. On its compact
subsets, the transverse fields decay exponentially, the full
mass-renormalized energy density tends uniformly to zero, and the
longitudinal curvature satisfies the Hodge equations up to errors that
vanish in the large mass limit. The scalar
limit is the Green potential of $\|T\|$, directly paralleling the
residual Dirac monopole in the three-dimensional theory
\cite[Theorem~A]{fadel2026abelian}. We upgrade Li's corrected
$L^1_{\mathrm{loc}}$ convergence to smooth convergence, prove that the
$L^2$-harmonic corrections are $o(m_i^{1/2})$, and obtain a compactness
dichotomy governed by a global $L^2$-harmonic $2$-form.
We also derive estimates for size, density rates, characteristic radii, and
regularity scales, and isolate the remaining bubbling,
tightness, and nodal set questions.
\end{enumerate}

Items~\textup{(i)}--\textup{(ii)} establish the common analytic,
topological, and measure-theoretic platform. The principal geometric
structure results are items~\textup{(iii)}--\textup{(iv)}: the hierarchy
of concentration and zero sets, the identification of the monotonicity
obstruction, and the abelianization and finer compactness theory on the
maximal open region $\mathcal R_{\mathrm{ab}}$.

\subsection{Main results}

Throughout the paper, our \textit{Main Setting} is the following. We let $(X^n,g,\Theta)$ denote either an AC $\mathrm G_2$-manifold $(X^7,g,\psi)$ or an AC Calabi--Yau $3$-fold $(X^6,g,\operatorname{Re}\Omega)$, and denote by $\Sigma$ the link of the asymptotic cone. We fix $G\in\{\mathrm{SU}(2),\mathrm{SO}(3)\}$ and let
$(\nabla_i,\Phi_i)$ be a sequence of irreducible finite intermediate energy
$\Theta$-monopoles on a fixed principal $G$-bundle $P\to X$, satisfying the
unified asymptotic framework of
\S\ref{subsec: asymptotic_analysis_theta_monopoles}, based on FNO in the
$\mathrm G_2$ case and its Calabi--Yau extension observed in
\cite{li2025large}. In particular, $|F_{\nabla_i}|=O_i(\rho^{-2})$ as $\rho\to\infty$, where the implicit constant may depend on $i$, and $\mathcal E^\Theta(\nabla_i,\Phi_i)<\infty$. Consequently, by
Theorem~\ref{thm:theta_asymptotics}, the masses $m_i:=\lim_{\rho\to\infty}|\Phi_i|$ exist, and each monopole admits an asymptotic abelian model with monopole class
$\beta_i\in H^2(\Sigma;\mathbb R)$. For $G=\mathrm{SU}(2)$ the class $\beta_i$ is integral, whereas for $G=\mathrm{SO}(3)$ the class $2\beta_i$ is integral. We assume that $m_i\to+\infty$ and that $\beta_i=\beta\in H^2(\Sigma;\mathbb R)$ is independent of $i$. Set $k:=\langle\beta\cup\Theta_\infty,[\Sigma]\rangle$. By the topological energy formula \eqref{eq: theta_monopole_energy_formula},
\[
    \mathcal E^\Theta(\nabla_i,\Phi_i)
    =
    \|\nabla_i\Phi_i\|_{L^2}^2
    =
    4\pi k m_i.
\]
Since the monopoles are irreducible and $m_i>0$, this identity forces
$k>0$. The asymptotic monopole class $\beta$ also determines the topological charge
$\operatorname{PD}_X(\delta_\Sigma\beta)\in H_{n-3}(X;\mathbb R)$,
where $\delta_\Sigma:H^2(\Sigma;\mathbb R)\to H_c^3(X;\mathbb R)$ is
the connecting homomorphism for the AC pair $(X,\Sigma)$. By
Theorem~\ref{thm:intro_calibrated_concentration}, this is the homology
class of the limiting calibrated current. We also define the
Kuratowski upper limit of the Higgs zero sets by
\[
    \mathcal Z
    :=
    \bigcap_{N\geqslant1}
    \overline{
        \bigcup_{i\geqslant N}\Phi_i^{-1}(0)
    }.
\]
The following theorem establishes the single calibrated limit and the
quantitative consequences needed for the subsequent structure theory.
Theorem~\ref{thm: finite_mass_and_estimates} gives the critical Sobolev
estimate for the Higgs defect, and Lemma~\ref{lem:largeness} uses it to
verify the PPS local largeness hypothesis. Li's large radius curvature
estimate and the fixed class intermediate energy identity provide the
uniform local mass-renormalized energy bound. Together with the
$\mathrm{SO}(3)$ extensions and the AC version of PPS compactness, these
inputs yield the PPS limiting current and measures. Li's equality of
the limiting mass-renormalized curvature and Higgs gradient measures
then gives saturation of the PPS inequalities and the identification
$T_{\mathrm{PPS}}=T_{\mathrm{Li}}$. The common current, its support, and
the exact limiting measures are the starting data for the comparison of
$\mathcal S$, $\mathcal Z$, $\mathcal C$, and $\mathcal O$ and for the
analysis of the maximal abelian region.

\begin{introtheorem}[Calibrated concentration and approximation by Higgs zeros]
\label{thm:intro_calibrated_concentration}
After passing to a subsequence, the following conclusions hold.

\begin{enumerate}
\item[(i)] The mass-renormalized Yang--Mills--Higgs and
intermediate energy measures
\[
    \mu_i
    :=
    m_i^{-1}
    \bigl(
        |F_{\nabla_i}|^2+|\nabla_i\Phi_i|^2
    \bigr)\,\vol
\]
and
\[
    \nu_i
    :=
    m_i^{-1}
    \bigl(
        |F_{\nabla_i}\wedge\Theta|^2
        +
        |\nabla_i\Phi_i|^2
    \bigr)\,\vol
    =
    2m_i^{-1}|\nabla_i\Phi_i|^2\,\vol
\]
converge weak-* as Radon measures to $\mu$ and $\nu$, respectively. If
$T_{\mathrm{PPS}}$ is the current produced by PPS compactness and
$T_{\mathrm{Li}}$ is Li's limiting calibrated current, then, in the
normalization of the present paper,
\[
    T_{\mathrm{PPS}}
    =
    T_{\mathrm{Li}}
    =:
    T,
    \qquad
    \mu=\nu=8\pi\|T\|.
\]
In particular, $T$ is a compactly supported $\Theta$-calibrated integral
$(n-3)$-cycle. Thus the PPS calibration inequalities are saturated and
\cite[Question~1.7]{parise2025nonabelian} has an affirmative answer for
the sequences in the Main Setting.

\item[(ii)] Li's topological and mass conclusions give
\[
    [T]
    =
    \operatorname{PD}_X(\delta_\Sigma\beta)
    \in
    H_{n-3}(X;\mathbb R),
\]
and
\[
    \|T\|(X)=k,
    \qquad
    \mu(X)=\nu(X)=8\pi k.
\]
Consequently,
\[
    \mathcal S
    :=
    \operatorname{spt}\mu
    =
    \operatorname{spt}\nu
    =
    \operatorname{spt}\|T\|
\]
is compact.

\item[(iii)] Let $\theta_T$ denote the
$\|T\|$-almost-everywhere positive integer multiplicity of $T$, and set
\[
    \ell_T
    :=
    \operatorname*{ess\,inf}_{\|T\|}\theta_T
    \in
    \mathbb Z_{\geqslant1}.
\]
Then $\mathcal{S}$ is detected directly from the approximating sequence as follows:
\[
    x\in\mathcal S
    \quad\Longleftrightarrow\quad
    \liminf_{r\downarrow0}
    \liminf_{i\to\infty}
    \frac{m_i^{-1}}
         {\omega_{n-3}r^{n-3}}
    \int_{B_r(x)}
    \bigl(
        |F_{\nabla_i}|^2+|\nabla_i\Phi_i|^2
    \bigr)\,\vol
    \geqslant8\pi\ell_T,
\]
and the threshold $8\pi\ell_T$ may equivalently be replaced by $8\pi$.
Moreover,
\[
    \mathcal H^{n-3}(\mathcal S)
    \leqslant
    C_n\frac{k}{\ell_T}.
\]

\item[(iv)] Writing $Z_i:=\Phi_i^{-1}(0)$, one has
\[
    \sup_{x\in\mathcal S}d(x,Z_i)
    \longrightarrow0.
\]
Equivalently, every neighbourhood of every point of $\mathcal S$ meets
$Z_i$ for all sufficiently large $i$. In particular,
\[
    \mathcal S\subset\mathcal Z.
\]
\end{enumerate}
\end{introtheorem}

The zeros approaching $\mathcal S$ are forced by nonzero transverse total
charge. More precisely, given $x\in\mathcal S$ and a neighbourhood $V$ of
$x$, one chooses a smooth patch of $T$ of positive measure contained in $V$.
For every sufficiently large $i$, Li's quantitative transverse energy
identity provides a normal fibre over a point of this patch outside the
corresponding exceptional set. The topological argument in Li's proof
produces an outer linking sphere enclosing the monopole clusters at the mass scale in that fibre, whose transverse Chern number equals the positive local multiplicity of $T$. The normal ball bounded by this sphere must therefore contain a zero of $\Phi_i$. Thus $V$ meets the Higgs zero set for every sufficiently large $i$.

Theorem~\ref{thm:intro_calibrated_concentration} follows by combining the
$\mathrm{SO}(3)$ extensions and the verification, in Proposition~\ref{prop:PPS_hypotheses_main_setting}, of the hypotheses {\rm(H1)} and {\rm(H2)} needed to apply the AC version of the PPS compactness theorem stated as
Theorem~\ref{thm:PPS_AC_version}. Saturation, i.e., $\mu=\nu$, the identification with Li's current $T_{\mathrm{PPS}}=T_{\mathrm{Li}}$, compact support, and the conclusions concerning homology and total mass are proved in Proposition~\ref{prop:PPS_saturation_and_Li_identification}. The density characterization and Hausdorff size bound are established in Proposition~\ref{prop:characterization_S},
Lemma~\ref{lem: size_energy_concentration}, and Corollary~\ref{cor: size_S}; and the one-sided Hausdorff convergence of the Higgs zero sets to $\mathcal S$ is Proposition~\ref{prop:S_subset_Z}.

The finite total mass of the calibrated current, together with a
uniform lower growth estimate furnished by calibrated monotonicity,
places a quantitative restriction on how densely the fixed support
$\mathcal S$ can fill a small ball. The precise statement is given in Lemma~\ref{lem:S_not_arbitrarily_dense}, and it roughly states that for sufficiently small $R>0$, if $\mathcal S$ is $\varepsilon$-dense in $B_R(x)$, then
\[
    \epsilon^3
    \gtrsim
     \frac{\ell_TR^n}{\|T\|(B_{2R}(x))}\geqslant \frac{\ell_TR^n}{k}.
\]
\noindent Thus, on any fixed ball of sufficiently small radius, $\mathcal S$
cannot be $\varepsilon_j$-dense for a sequence
$\varepsilon_j\downarrow0$. This does not amount to a porosity estimate
for $\mathcal S$. The Higgs zero sets and curvature concentration loci
require a different argument. Proposition~\ref{prop:Li_cover_Hausdorff_content_off_S}
uses Li's Vitali covers to show that, on compact subsets of
$X\setminus\mathcal S$, their critical unrestricted Hausdorff content
tends to zero, while
Propositions~\ref{prop:Li_loci_density_rate_obstruction} and
\ref{prop:Higgs_zero_density_rate_obstruction} give quantitative
obstructions to rapid density. These estimates refine, but do not rule
out, the phenomenon discussed in \cite[Remark~3.10]{li2025large}.

The principal set-theoretic structure theorem compares the limiting
Higgs zero set with the upper limit of Li's curvature concentration
loci. For the fixed constant
$\Lambda_0$ chosen in Subsection~\ref{subsec:Li_moving_loci}, these loci are defined from Li's weighted curvature potential $\mathcal P_i$ (see equation \eqref{eq:Li_potential_our_convention}) by
\[
    C_{\Lambda_0,i}
    :=
    \left\{
        p\in X:
        \mathcal P_i(p)>\Lambda_0^{-1}
    \right\}.
\]
Li's concentration--decay dichotomy implies, in particular, that Higgs
zeros belong to $C_{\Lambda_0,i}$ for all sufficiently large $i$. Since the loci $C_{\Lambda_0,i}$ need not converge as sets, we consider
their Kuratowski upper limit, which we call the \textbf{limiting
nonabelian locus} of the sequence:
\[
    \mathcal C
    :=
    \bigcap_{N\geqslant1}
    \overline{\bigcup_{i\geqslant N}C_{\Lambda_0,i}}.
\]
Equivalently, $x\in\mathcal C$ if and only if there are indices
$i_j\to\infty$ and points $p_j\in C_{\Lambda_0,i_j}$ such that
$p_j\to x$.

We also use the closed obstruction locus $\mathcal O$ introduced in
Definition~\ref{def:codim3_monotonicity_obstruction}. By definition,
$X\setminus\mathcal O$ is the open union of all open subsets on which
the sequence satisfies effective codimension-three monotonicity. It is the
locus where this property holds locally; on the noncompact manifold $X$, the
union need not itself satisfy the estimate with one set of constants uniform
over all centres. Away from $\mathcal S$, local effective monotonicity is
equivalent to vanishing codimension-three Morrey density uniformly over balls
contained in shrinking neighbourhoods; see
Proposition~\ref{prop:equivalent_effective_monotonicity_off_S}. This
smallness forces Li's weighted curvature potential below the defining threshold. Conversely, Li's decay
estimates, combined with a Bochner mean value argument, give local uniform
decay of the full mass-renormalized energy density on
$\mathcal R_{\mathrm{ab}}$ and hence the reverse inclusion. We obtain:

\begin{introtheorem}[Hierarchy of concentration sets and monotonicity obstruction]
\label{thm:intro_obstruction_hierarchy}
In the Main Setting,
\[
    \mathcal S
    \subset
    \mathcal Z
    \subset
    \mathcal C
    =
    \mathcal S\cup\mathcal O.
\]
Consequently,
\[
\mathcal C\setminus\mathcal S = \mathcal O\setminus\mathcal S,
\qquad
\mathcal R_{\mathrm{ab}} = X\setminus(\mathcal S\cup\mathcal O).
\]
Moreover, let $x\in\mathcal C$, and suppose that
$i_j\to\infty$ and $p_j\in C_{\Lambda_0,i_j}$ satisfy $p_j\to x$.
After relabelling this subsequence, $x$ belongs to the corresponding
codimension-four concentration set. More precisely, for every fixed
$0<r\leqslant r_0$,
\[
    \liminf_{j\to\infty}
    e^{cr^2}r^{4-n}
    \int_{B_r(x)}e_{i_j}
    \geqslant
    \varepsilon_0.
\]
In particular, $\mathcal C=\mathcal S$ if and only if $\mathcal O\subset\mathcal S$. Under this condition,
\[
    \mathcal S=\mathcal Z=\mathcal C,
\]
and, for every compact neighbourhood $K\subset X$ of $\mathcal S$, one has
\[
    d_{\mathcal H}(Z_i\cap K,\mathcal S)
    \longrightarrow0,
    \qquad
    d_{\mathcal H}
    \bigl(\overline{C_{\Lambda_0,i}}\cap K,\mathcal S\bigr)
    \longrightarrow0.
\]
\end{introtheorem}

The inclusions $\mathcal S\subset\mathcal Z\subset\mathcal C$
are proven in Proposition~\ref{prop:S_subset_Z} and
Lemma~\ref{lemma:Z_subset_C_infty}. Proposition~\ref{prop:zero_centered_TU_concentration}
shows that every point of $\mathcal Z$, after passing to a subsequence
realizing it by Higgs zeros, belongs to the codimension-four
concentration set of that subsequence. The stronger statement for every
point of $\mathcal C$, along any subsequence realizing it through
Li's curvature concentration loci, is Corollary~\ref{cor:Cinfty_subsequential_TU}.
Proposition~\ref{prop:mass_renormalized_density_decay_Rab} promotes Li's
transverse decay on $\mathcal R_{\mathrm{ab}}$ to local uniform decay of
the full mass-renormalized energy density. Together with
Proposition~\ref{prop:Cinfty_minus_S_obstructs_positive_codim3}, this gives
the identity $\mathcal C=\mathcal S\cup\mathcal O$ in
Corollary~\ref{cor:nonabelian_locus_obstruction_identity}.

The local Hausdorff convergence under the assumption that $\mathcal O\subset\mathcal S$ is Corollary~\ref{cor:Hausdorff_convergence_unobstructed}. This conclusion is local on fixed compact subsets and does not exclude components of $C_{\Lambda_0,i}$ whose centres escape to infinity. No general nonemptiness or density statement for $\mathcal R_{\mathrm{ab}}$ is asserted in the Main Setting.

Without assuming $\mathcal O\subset\mathcal S$, Li's Vitali covering
lemma nevertheless gives quantitative information about the loci $C_{\Lambda_0,i}$
away from the calibrated support.  Writing $p=n-3$, for every
$K\Subset X\setminus\mathcal S$ and every $0<\kappa<1$,
\[
    \mathcal H^{p+\kappa}_\infty
    (C_{\Lambda_0,i}\cap K)
    =o(m_i^{-\kappa}),
    \qquad
    \mathcal H^p_\infty
    (C_{\Lambda_0,i}\cap K)
    \longrightarrow0,
\]
and the same conclusions hold for $Z_i\cap K$; see
Proposition~\ref{prop:Li_cover_Hausdorff_content_off_S} and
Corollary~\ref{cor:Higgs_zero_Hausdorff_content_off_S}.  If the curvature concentration loci
become dense in a fixed ball disjoint from $\mathcal S$, the
corresponding maximal distance $h_i^{\mathcal C}$ necessarily satisfies
\[
    m_i^{n-3}(h_i^{\mathcal C})^n\longrightarrow+\infty.
\]
For the zero sets, the Higgs defect estimate and ordinary
$\varepsilon$-regularity improve this to
\[
    m_i(h_i^{\mathcal Z})^4\longrightarrow+\infty.
\]
These are estimates for the varying loci $C_{\Lambda_0,i}$. They do not, without a
quantitative common approximation rate, imply a Hausdorff measure bound
for the Kuratowski upper limits $\mathcal C$ or $\mathcal Z$.

By Theorem~\ref{thm:intro_obstruction_hierarchy}, 
$\mathcal R_{\mathrm{ab}}=X\setminus(\mathcal S\cup\mathcal O)$.
If $\mathcal O\subset\mathcal S$ (in particular, if $\mathcal O=\varnothing$), then $\mathcal R_{\mathrm{ab}}=X\setminus\mathcal S$ is open and dense and has full Riemannian measure, since Theorem~\ref{thm:intro_calibrated_concentration} shows that $\mathcal S$ is compact and has finite $(n-3)$-dimensional Hausdorff measure.

Crucially, Theorem~\ref{thm:intro_obstruction_hierarchy} highlights the
importance of understanding $\mathcal O$. A sufficient condition for
$x\notin\mathcal O$ is the codimension-two Morrey bound of
Proposition~\ref{prop:codim2_Morrey_criterion_positive_monotonicity} for
the positive part
\[
    q_i^+
    :=
    \bigl(|F_{\nabla_i}|^2-|\nabla_i\Phi_i|^2\bigr)^+
\]
of the Yang--Mills--Higgs imbalance. Independently of this criterion,
Proposition~\ref{prop:excess_points_determine_characteristic_cores}
associates to every
$x\in\mathcal C\setminus\mathcal S
=\mathcal O\setminus\mathcal S$
points $p_i\to x$ and radii $R_i\downarrow0$ such that Li's weighted
curvature potential on $B_{R_i}(p_i)$ is bounded below, whereas the
Yang--Mills--Higgs energy of these balls is $o(m_i)$. Proposition~\ref{prop:regularity_scale_inside_characteristic_core}
shows that the ordinary regularity scale
$\lambda_i:=\mathfrak r_i(p_i)$ satisfies
$\lambda_i\leqslant R_i$. If
$m_i\lambda_i\to+\infty$, then
Proposition~\ref{prop:characteristic_core_above_mass_scale} shows that
rescaling by $m_i$ gives a flat connection with parallel Higgs field,
while the contribution to the characteristic potential from every ball
$B_{L/m_i}(p_i)$ with fixed $L$ tends to zero. Quantitatively, suitable
annuli between $m_i^{-1}$ and $R_i$ have unbounded scale-invariant
curvature energy in codimension-four.

For the known cohomogeneity-one examples, the Morrey criterion holds on a neighbourhood of every point of $X$, and hence $\mathcal O=\varnothing$.

\begin{introtheorem}[Cohomogeneity-one examples]
\label{thm:intro_homogeneous_examples}
For the families of $\mathrm G_2$-monopoles on the Bryant--Salamon $\mathrm G_2$-manifolds $\Lambda_-^2(S^4)$, $\Lambda_-^2(\mathbb{CP}^2)$ \cite{oliveira2014monopoles}, and the family of Calabi--Yau monopoles on the Stenzel metric in $T^*S^3$ \cite{oliveira2016calabi}, every large mass sequence satisfies
\[
    \mathcal O=\varnothing.
\]
Consequently,
\[
    \mathcal S=\mathcal Z=\mathcal C.
\]
In each of the three families, this common set is the unique compact calibrated codimension-three submanifold and coincides with the Higgs zero set of every member of the family. Moreover, for every compact set $K\subset X$ containing
$\mathcal S$ in its interior,
\[
    d_{\mathcal H}
    \bigl(\overline{C_{\Lambda_0,i}}\cap K,\mathcal S\bigr)
    \longrightarrow0.
\]
\end{introtheorem}

This is proved in Section~\ref{sec:homogeneous_examples}. The proof rests on showing that the leading
terms at the mass scale satisfy the auxiliary transverse Bogomolny equations and
cancel in the positive imbalance. The residual terms are not uniformly
bounded pointwise but satisfy the following tubular estimate, uniformly in the mass,
\[
    q_i^+
    \leqslant
    C\bigl(1+r_N^{-2}+m_i r_N^{-1}\bigr),
\]
where $r_N$ is the distance to the codimension-three calibrated cycle.
Since $N$ has codimension three, this gives the codimension-two Morrey bound
of Proposition~\ref{prop:codim2_Morrey_criterion_positive_monotonicity},
uniformly over all balls in a fixed tubular neighbourhood of $N$.

The abelianization theorem concerns the open set
$\mathcal R_{\mathrm{ab}}=X\setminus\mathcal C$. By
Definition~\ref{def:C_infty}, a point
$x\in\mathcal R_{\mathrm{ab}}$ if and only if there exist an open
neighbourhood $U$ of $x$ and an index $i_0$ such that
\[
    U\cap C_{\Lambda_0,i}=\varnothing
    \qquad
    \text{for every }i\geqslant i_0.
\]
Lemma~\ref{lem:maximal_abelian_clearing_region} upgrades this local
eventual avoidance to uniform avoidance on compact subsets and
characterizes $\mathcal R_{\mathrm{ab}}$ as the largest open subset of
$X$ with this property. Consequently, on every compact subset of
$\mathcal R_{\mathrm{ab}}$, the Higgs field is eventually uniformly
large, while the curvature splits into a longitudinal component along
the Higgs field and a transverse component orthogonal to it. The result
below describes their behaviour in the large mass limit, the quantitative
decay of the full mass-renormalized energy on
$\mathcal R_{\mathrm{ab}}$, and
the resulting regularity scale dichotomy.

To state the result, set
$e_i:=|F_{\nabla_i}|^2+|\nabla_i\Phi_i|^2$ and let
$\mathfrak r_i:X\to(0,r_0]$ denote the ordinary Yang--Mills--Higgs
regularity scale,
\[
    \mathfrak r_i(x)
    :=
    \sup\left\{
        0<r\leqslant r_0:
        s^{4-n}\int_{B_s(x)}e_i<\varepsilon_0
        \text{ for every }0<s\leqslant r
    \right\},
\]
where $\varepsilon_0$ is the constant in
Theorem~\ref{thm: total_epsilon_regularity} and $r_0>0$ is fixed below
the geometric scale.

\begin{introtheorem}[Abelianization and regularity scales on the
maximal open $\mathcal R_{\mathrm{ab}}$]
\label{thm:intro_abelianization_regular_scales}
Let $K\Subset\mathcal R_{\mathrm{ab}}$. Then
\[
    \frac{|\Phi_i|}{m_i}
    \longrightarrow1
    \qquad
    \text{uniformly on }K.
\]
In particular, $|\Phi_i|\geqslant m_i/2$ on $K$ for all sufficiently
large $i$. Set $\Psi_i:=\Phi_i/|\Phi_i|$ and $\widehat\omega_i:=\langle F_{\nabla_i},\Psi_i\rangle$. Relative to the
splitting $\mathfrak g_P|_K=\langle\Psi_i\rangle\oplus\langle\Psi_i\rangle^\perp$, the transverse components $(F_{\nabla_i})^\perp$, $(\nabla_i\Phi_i)^\perp$, and $\nabla_i\Psi_i$ are $O_{C^\ell(K)}(e^{-c_{K,\ell}m_i})$ for every $\ell\geqslant0$ and suitable $c_{K,\ell}>0$. Moreover,
\[
    \sup_K m_i^{-1}e_i\longrightarrow0,
    \qquad
    m_i^{-1/2}
    \|\widehat\omega_i\|_{L^\infty(K)}
    \longrightarrow0,
\] 
and
$d\widehat\omega_i,d^*\widehat\omega_i\to0$ in $C^\ell(K)$ for every
$\ell\geqslant0$. Setting
\[
    u_i
    :=
    \frac{m_i^2-|\Phi_i|^2}{2m_i},
    \qquad
    v_i
    :=
    m_i-|\Phi_i|.
\]
Thus $0\leqslant u_i\leqslant v_i$ and
$\widetilde\Phi_i=-v_i\Psi_i$. One also has
\[
    u_i,\ v_i
    \longrightarrow
    v_T
    :=
    4\pi
    \int_X G(\,\cdot\,,y)\,d\|T\|(y)
\]
uniformly on $K$. In particular,
\[
    \sup_K|\widetilde\Phi_i|
    =
    \sup_K(m_i-|\Phi_i|)
    \leqslant C_K.
\]
For some $r_K>0$,
\[
    \sup_{\substack{x\in K\\0<r<r_K}}
    m_i^{-1}r^{2-n}
    \int_{B_r(x)}e_i
    \longrightarrow0.
\]
In particular, if
$q_i:=|F_{\nabla_i}|^2-|\nabla_i\Phi_i|^2$, then the same conclusion
holds with $e_i$ replaced by $q_i^+$. The ordinary regularity scale also
satisfies the unconditional separation estimate
\[
    \inf_{x\in K}m_i^{1/4}\mathfrak r_i(x)
    \longrightarrow+\infty.
\]
There is a regularity scale dichotomy on $\mathcal R_{\mathrm{ab}}$; exactly one of the following alternatives occurs:

\begin{enumerate}[label=\textnormal{(\alph*)}]
\item If $\sup_i\|\widehat\omega_i\|_{L^2(K)}<\infty$ for every compact
domain $K\Subset\mathcal R_{\mathrm{ab}}$, then the regularity scales are
locally uniformly bounded below. After passing to a subsequence and
applying gauge transformations, the translated configurations
$(\nabla_i,\widetilde\Phi_i)$, where
$\widetilde\Phi_i:=\Phi_i-m_i\Psi_i$, converge smoothly on compact subsets
of $\mathcal R_{\mathrm{ab}}$ to a smooth abelian $\Theta$-monopole.

\item Otherwise, there exist a compact domain
$K\Subset\mathcal R_{\mathrm{ab}}$ and a subsequence, still indexed by
$i$, such that $\mathfrak a_i:=\|\widehat\omega_i\|_{L^2(K)}\to\infty$.
Necessarily,
\[
    \mathfrak a_i=o(m_i^{1/2}).
\]
After passing to a further subsequence,
$\mathfrak a_i^{-1}\widehat\omega_i$ converges smoothly on compact subsets
of $\operatorname{int}K$ to a harmonic $2$-form
$\overline\omega_{\mathrm{har}}$, which is primitive in the Calabi--Yau case. If
$U'\Subset U\Subset\operatorname{int}K$ and
$\inf_U|\overline\omega_{\mathrm{har}}|>0$, then
$\mathfrak r_i(x)\asymp\mathfrak a_i^{-1/2}$ uniformly for $x\in U'$.
At any point $x\in\operatorname{int}K$ satisfying the quantitative
order-$q$ profile condition of
Definition~\ref{def:quantitative_order_q_profile}, one has
$\mathfrak r_i(x)\asymp\mathfrak a_i^{-1/(q+2)}$.
\end{enumerate}

If $\mathcal O\subset\mathcal S$, then
$\mathcal R_{\mathrm{ab}}=X\setminus\mathcal S$, and the preceding
conclusions hold on compact subsets of $X\setminus\mathcal S$.
\end{introtheorem}

The following corollaries combine the local results of Theorem~\ref{thm:intro_abelianization_regular_scales} with Li's corrected
compactness theorem.

\begin{introcorollary}[Smooth convergence of Li's corrected longitudinal
forms on compact subsets of $\mathcal R_{\mathrm{ab}}$]
\label{cor:intro_smooth_Li_corrected}
For each $i$, fix an $L^2$-harmonic correction
$\sigma_i\in \mathscr{H}_{(2)}^2(X)$ furnished by
\cite[Theorem~1.12]{li2025large}, with the $\mathrm{SO}(3)$ extension
described in Remark~\ref{rmk:Li_SO3_variant}. Then, after passing to the subsequence in Li's theorem,
\[
    \frac{1}{4\pi}\widehat\omega_i+\sigma_i
    \longrightarrow
    \widetilde F_\infty
\]
smoothly on compact subsets of $\mathcal R_{\mathrm{ab}}$, where
$\widetilde F_\infty=F_\infty/(2\pi)$ is Li's normalized limiting
abelian curvature form.
\end{introcorollary}

\begin{introcorollary}[Li's harmonic corrections and $\mathcal R_{\mathrm{ab}}$]
\label{cor:intro_harmonic_correction_alternatives}
For every sequence of $L^2$-harmonic corrections $\sigma_i\in \mathscr{H}_{(2)}^2(X)$ furnished by Li's theorem,
\[
   \|\sigma_i\|_{L^2(X)} = o(m_i^{1/2}).
\]
Moreover, exactly one of the following alternatives holds.
\begin{enumerate}
\item  Li's corrections are uniformly bounded in $L^2(X)$:
\[
    \sup_i\|\sigma_i\|_{L^2(X)}<\infty.
\]
Then, after passing to a further subsequence and applying gauge transformations, the translated configurations
\[
    (\nabla_i,\widetilde\Phi_i),
    \qquad
    \widetilde\Phi_i:=\Phi_i-m_i\Psi_i,
\]
converge smoothly on compact subsets of $\mathcal R_{\mathrm{ab}}$ to a
smooth abelian $\Theta$-monopole.

\item The $L^2(X)$ norms have infinite limsup:
\[
    \limsup_{i\to\infty}
    \|\sigma_i\|_{L^2(X)}
    =
    +\infty.
\]
After passing to a further subsequence and relabelling, write
$b_i:=\|\sigma_i\|_{L^2(X)}$, so that $b_i\to+\infty$. Then $b_i=o(m_i^{1/2})$, and there exists a
form $\sigma_\infty\in \mathscr{H}_{(2)}^2(X)$ with $\|\sigma_\infty\|_{L^2(X)}=1$ such that
\[
    \frac{1}{4\pi b_i}\widehat\omega_i
    \longrightarrow
    -\sigma_\infty
\]
smoothly on compact subsets of $\mathcal R_{\mathrm{ab}}$. On every
compact set $K\Subset\mathcal R_{\mathrm{ab}}$, one has
\[
    \inf_{x\in K}m_i^{1/4}\mathfrak r_i(x)
    \longrightarrow
    +\infty,
\]
independently of the alternative for the harmonic corrections. Where
$|\sigma_\infty|$ is bounded below, the sharper comparison
$\mathfrak r_i\asymp b_i^{-1/2}$ holds locally. At every point satisfying
the quantitative order-$q$ profile condition,
\[
    m_i^{1/(2(q+2))}\mathfrak r_i(x)
    \longrightarrow
    +\infty.
\]
\end{enumerate}
In particular, the first alternative always holds when $\mathscr H_{(2)}^2(X)=0$, equivalently when $H_c^2(X;\mathbb R)=0$. If, in addition, $\mathcal O\subset\mathcal S$, then $\mathcal R_{\mathrm{ab}}=X\setminus\mathcal S$, and the resulting smooth
convergence holds on compact subsets of $X\setminus\mathcal S$.
\end{introcorollary}
In the first alternative, the precise comparison between the smooth
abelian limit of the translated configurations and Li's singular
abelian monopole is given in
Remark~\ref{rmk:comparison_with_Li_singular_limit}.

In particular, the last part of Corollary~\ref{cor:intro_harmonic_correction_alternatives} gives an
affirmative answer to \cite[Question~3]{li2025large} under the
hypotheses $H_c^2(X;\mathbb R)=0$ and
$\mathcal O\subset\mathcal S$. By
Theorem~\ref{thm:intro_homogeneous_examples} and
Remark~\ref{rmk:homogeneous_Hc2_vanishing}, the known cohomogeneity-one
families satisfy not only $H_c^2(X;\mathbb R)=0$ but also the stronger
condition $\mathcal O=\varnothing$. These results account, within the present general theory, for
the $C^\infty_{\mathrm{loc}}$ convergence previously proved directly by
the second author for such families in \cites{oliveira2014monopoles,oliveira2016calabi}.

\subsection{Organization}

Section~\ref{sec: foundations} develops the analytic and asymptotic
foundations. We prove that the Higgs defect belongs to
$\mathscr D_0^{1,2}(X)$, establish its critical Sobolev estimate and
capacity, volume and inradius bounds for relative Higgs sublevel
regions, and show that, along a summably growing subsequence, the
normalized Higgs norms converge to $1$ outside a set of zero
$2$-capacity. We also verify the PPS largeness hypothesis along large
mass sequences. We then formulate
the FNO asymptotic theory in the unified $\Theta$-monopole notation,
establish the required $\mathrm{SO}(3)$ variant, and combine Li's
large radius curvature estimate with the fixed class energy identity
to obtain hypothesis {\rm(H1)}.

Section~\ref{sec: concentration} establishes the $\mathrm{SO}(3)$ and AC
versions of PPS compactness and verifies their hypotheses in the Main
Setting. We prove saturation of the PPS inequalities and identify the
PPS current with Li's calibrated current. We then characterize the
common calibrated support by the limiting energy measures and by
pointwise lower energy densities, prove size estimates and quantitative
restrictions on density, relate limiting Higgs zeros to ordinary
codimension-four concentration centred at Higgs zeros, and establish one-sided Hausdorff approximation of the calibrated support
by the Higgs zero sets. The section also introduces the obstruction locus
$\mathcal O$ and a codimension-two Morrey criterion for excluding it.

Section~\ref{sec: Li_concentration} introduces the Kuratowski upper limit
$\mathcal C$ of the curvature concentration loci and the abelian region
$\mathcal R_{\mathrm{ab}}$. We prove
\[
    \mathcal S
    \subset
    \mathcal Z
    \subset
    \mathcal C
    =
    \mathcal S\cup\mathcal O,
\]
and show that the full mass-renormalized energy density converges locally
uniformly to zero on $\mathcal R_{\mathrm{ab}}$.  Li's Vitali covering
lemma yields vanishing critical Hausdorff content for the
curvature concentration loci and Higgs zero sets on compact subsets of
$X\setminus\mathcal S$, together with quantitative obstructions to
rapid density and transfer statements under quantitative rate assumptions for local
Hausdorff limits.  We derive local Hausdorff convergence when
$\mathcal O\subset\mathcal S$. Every point of
$\mathcal C\setminus\mathcal S$ admits sequences $p_i\to x$ and
$R_i\downarrow0$ such that Li's weighted curvature potential on
$B_{R_i}(p_i)$ is bounded below while the mass-renormalized energy tends to
zero. The ordinary regularity scale at $p_i$ is at most $R_i$. When its
ratio to the inverse mass scale $m_i^{-1}$ tends to infinity, the
contribution from each fixed multiple of the ball at the inverse mass scale
$B_{m_i^{-1}}(p_i)$ vanishes, and suitable intervening annuli have
unbounded codimension-four scale-invariant curvature energy. On
$\mathcal R_{\mathrm{ab}}$ we prove quantitative vanishing of the full
mass-renormalized codimension-two Morrey norm and
deduce the universal scale separation
$\inf_K m_i^{1/4}\mathfrak r_i\to+\infty$ on compact subsets.

Section~\ref{sec:homogeneous_examples} treats the known
cohomogeneity-one Bryant--Salamon and Stenzel monopole families. In these cases, good control over the signed imbalance gives the codimension-two Morrey bound;
hence $\mathcal O=\varnothing$ and
$\mathcal S=\mathcal Z=\mathcal C$ for these examples.

Section~\ref{sec: abelianization} proves quantitative abelianization on
compact subsets of $\mathcal R_{\mathrm{ab}}$, derives Hodge equations
for the longitudinal curvature up to errors that vanish in the large mass
limit, and uses the Green
representation to prove locally uniform normalization of the Higgs
fields and identify the scalar residual limit with the Green potential
of the calibrated measure.
Section~\ref{sec:frequency_away_Cinfty} combines these equations with
Li's strong $L^1_{\mathrm{loc}}$ convergence of the corrected forms to
obtain smooth local convergence. We prove that the global harmonic
corrections are $o(m_i^{1/2})$ in $L^2(X)$ and then distinguish the
uniformly $L^2$-bounded alternative from the alternative in which the limsup is infinite;
in the latter case, a subsequence is selected along which the $L^2$ norms
tend to infinity, and the normalized corrections converge to a nonzero
global $L^2$-harmonic form. The section also derives the corresponding
regularity scale estimates and compares these scales with the mass scale.

Section~\ref{sec:resulting_picture} synthesizes the resulting
large mass picture and isolates the remaining problems: effective
codimension-three monotonicity; characteristic radii at points of
$\mathcal C\setminus\mathcal S$; bubbling centred at Higgs zeros and at
points of $\mathcal Z\setminus\mathcal S$; tightness of the loci
$C_{\Lambda_0,i}$; nonemptiness of $\mathcal R_{\mathrm{ab}}$; and
regularity near the nodal set of the limiting harmonic form.

The appendices collect the local mean value, stress-energy,
Bochner--Weitzenb\"ock, annular no-neck, and $\varepsilon$-regularity
tools. They also contain the Euclidean stress-energy identity and the
finite energy rigidity statement.

\subsection*{Acknowledgements}
The authors thank Daniel Stern and Yang Li for helpful exchanges and
comments on earlier drafts. This work was partially supported by FCT/Mobility/1301033342/2024-25
under the FCT Mobility Programme. D.F. was supported by
the MATH-AmSud project \emph{Symmetries in Geometry and Physics}
(SGP 24-MATH-12), CNPq Universal grant 406666/2023-7, the ANR--FAPESP
project \emph{BRIDGES---Brazil--France Interplays in Gauge Theory,
Extremal Structures and Stability} (ANR-21-CE40-0017 and FAPESP
2021/04065-6), and the University of S\~ao Paulo through its
\emph{Programa de Apoio aos Novos Docentes USP}.

\emph{Use of artificial intelligence:} The authors used OpenAI's ChatGPT, including the GPT-5.5 and GPT-5.6
models, as an aid in drafting and polishing portions of the exposition
and for exploratory discussion of preliminary ideas and possible proof
strategies formulated by the authors. All mathematical statements,
arguments, computations, and bibliographic references retained in the
manuscript were independently checked and verified by the authors, who
take full responsibility for the contents of the paper.


\section{Analytic foundations and asymptotic estimates}\label{sec: foundations}

\subsection{Notation and conventions}\label{ss:Notation}

Throughout the paper, $(X^n,g)$ denotes a complete, connected, oriented, noncompact Riemannian manifold without boundary, assumed asymptotically conical (AC) with rate $\nu_0<0$. Thus there exist a compact set $K\subset X$, a closed Riemannian manifold $(\Sigma^{n-1},g_\Sigma)$, and a diffeomorphism
\[
\Upsilon:(1,\infty)_r\times\Sigma\to X\setminus K
\]
such that
\[
|D_C^j(\Upsilon^*g-g_C)|_{g_C}
=
O(r^{\nu_0-j}),
\qquad j\geqslant 0,
\]
where $g_C=dr^2+r^2g_\Sigma$ is the cone metric on $C(\Sigma):=(1,\infty)\times\Sigma$ and $D_C$ its Levi--Civita connection. We equip $\Sigma$ with the unique orientation for which $\Upsilon$ is orientation-preserving when $(1,\infty)\times\Sigma$ is given the product orientation, and denote by
$\operatorname{vol}_\Sigma$ the corresponding Riemannian volume form.

We denote by $\rho:X\to[1,\infty)$ a fixed radius function, namely a smooth extension of the radial coordinate function $r$ along the end. For $R\geqslant 1$ we write
\[
B_R:=\{\rho<R\},
\qquad
\overline B_R:=\{\rho\leqslant R\},
\qquad
\Sigma_R:=\partial\overline B_R.
\]
For all sufficiently large $R$, the set $\overline B_R$ is a compact
smooth $n$-manifold with boundary $\Sigma_R$, and $\Sigma_R$ is
diffeomorphic to $\Sigma$. We orient $\Sigma_R=\partial\overline B_R$ by the boundary orientation; under the identification with $\Sigma$ induced by $\Upsilon$, this agrees with the orientation fixed above.

When we write $B_r(x)$ we mean the open geodesic ball centred at $x\in X$
and of radius $r\in(0,\infty)$. We denote by $d(x,y)$ the Riemannian distance
between two points $x,y\in X$ and by $d_{\mathcal H}(A,B)$ the Hausdorff
distance between compact subsets $A,B\subset X$. For arbitrary subsets
$A,B\subset X$, the notation $A\Subset B$ means that $\overline A$ is
compact and contained in the interior of $B$.

We denote by $\mathcal H^k$ the normalized $k$-dimensional Hausdorff
measure associated with the Riemannian distance $d$, and by $\omega_k$ the
Euclidean volume of the unit ball in $\mathbb R^k$. Thus the normalization
is chosen so that a Euclidean $k$-ball of radius $r$ has
$\mathcal H^k$-measure $\omega_k r^k$. With this convention,
$\mathcal H^n$ agrees with the Riemannian volume measure $\vol=\vol_X$ on $X$,
and, for every smooth embedded $k$-dimensional submanifold $M\subset X$,
the restriction $\mathcal H^k\lfloor M$ agrees with the Riemannian measure
induced on $M$.

Unless otherwise indicated, integrals over subsets of $X$ and over smooth
submanifolds are taken with respect to these measures; the measure symbols
are usually omitted.
Convergence of possibly signed Radon measures is always understood in the
weak-$*$ sense, against functions in $C_c^0(X)$.

We denote by
\[
    \mathscr{H}_{(2)}^2(X)
    :=
    \left\{
        \sigma\in\Omega^2(X):
        d\sigma=0,\quad d^*\sigma=0,\quad
        \|\sigma\|_{L^2(X)}<\infty
    \right\}
\]
the space of $L^2$-harmonic $2$-forms on $X$. The Hodge Laplacian and the Hodge--Dirac operator on $\Omega^k(X)$ are denoted by
\[
    \Delta:=dd^*+d^*d,
    \qquad
    \slashed D:=d+d^*.
\]
Given a principal $G$-bundle $P\to X$, where $G$ is a compact Lie group,
we denote by
\[
    \mathfrak g_P:=P\times_{\mathrm{Ad}}\mathfrak g
\]
the associated adjoint bundle. We equip $\mathfrak g$ with a fixed
$\mathrm{Ad}$-invariant inner product, and use the induced fibrewise
metric on $\mathfrak g_P$, still denoted by $\langle\cdot,\cdot\rangle$.
Together with the Riemannian metric on $X$, this induces the pointwise
inner products and norms on all bundles
$\Lambda^kT^*X\otimes\mathfrak g_P$. 

Unless explicitly stated otherwise, the principal bundles considered in
the paper have structure group
\[
    G\in\{\mathrm{SU}(2),\mathrm{SO}(3)\}.
\]
On $\mathfrak{su}(2)$ we use the $\mathrm{Ad}$-invariant inner product
\begin{equation}\label{eq:inner_prod_convention}
    \langle a,b\rangle:=-2\operatorname{tr}_{\mathbb C^2}(ab),
    \qquad
    a,b\in\mathfrak{su}(2).
\end{equation}
Let $\pi:\mathrm{SU}(2)\longrightarrow\mathrm{SO}(3)$ be the adjoint double cover. On $\mathfrak{so}(3)$ we use the unique $\mathrm{Ad}$-invariant inner product for which $d\pi:\mathfrak{su}(2)\to\mathfrak{so}(3)$ is an isometry. Equivalently,
\begin{equation}\label{eq:SO3_inner_prod_convention}
    \langle A,B\rangle
    :=
    -\frac12\operatorname{tr}_{\mathbb R^3}(AB),
    \qquad
    A,B\in\mathfrak{so}(3).
\end{equation}
For comparison with the compactness results used below, we record that
Parise--Pigati--Stern and Li use, on $\mathfrak{su}(2)$, the convention
\[
    \langle a,b\rangle_{\mathrm{PPS}}
    =
    \langle a,b\rangle_{\mathrm{Li}}
    :=
    -\frac12\operatorname{tr}_{\mathbb C^2}(ab),
    \qquad
    a,b\in\mathfrak{su}(2),
\]
and its transport under $d\pi$ on $\mathfrak{so}(3)$, namely
\[
    \langle A,B\rangle_{\mathrm{PPS}}
    =
    \langle A,B\rangle_{\mathrm{Li}}
    :=
    -\frac18\operatorname{tr}_{\mathbb R^3}(AB).
\]
Consequently, for either structure group,
\begin{equation}
\label{eq:PPS_Li_inner_product_comparison}
    \langle\cdot,\cdot\rangle
    =
    4\langle\cdot,\cdot\rangle_{\mathrm{PPS}}
    =
    4\langle\cdot,\cdot\rangle_{\mathrm{Li}},
\end{equation}
and hence
\begin{equation}
\label{eq:PPS_Li_norm_mass_comparison}
    |\cdot|
    =
    2|\cdot|_{\mathrm{PPS}}
    =
    2|\cdot|_{\mathrm{Li}}.
\end{equation}
In particular, for the same adjoint-valued Higgs field, the masses in the
three conventions satisfy
\begin{equation}
\label{eq:PPS_Li_mass_comparison}
    m
    =
    2m^{\mathrm{PPS}}
    =
    2m^{\mathrm{Li}}.
\end{equation}
Whenever results from
\cites{parise2025nonabelian,li2025large} are stated in the normalization
of the present paper, the corresponding norms, masses, energies, and
measures are understood to have been converted according to
\eqref{eq:PPS_Li_inner_product_comparison}--%
\eqref{eq:PPS_Li_mass_comparison}.

Certain arguments do not depend on the choice of structure group $G$
and will be stated for an arbitrary compact Lie group. This occurs, for
instance, in Subsections~\ref{subsec: YMH_config_finite_mass} and
\ref{subsec: monopoles}, and in the local analytic estimates collected in
Appendices~\ref{app: A}, \ref{app: B}, and \ref{app: C}. Elsewhere the
standing structure group is
$G\in\{\mathrm{SU}(2),\mathrm{SO}(3)\}$.

We let $\Gamma(\mathfrak{g}_P)$ denote the space of smooth sections of the bundle $\mathfrak{g}_P$; its elements are called Higgs fields. Given a Higgs field $\Phi$, we write
\[
Z(\Phi):=\{x\in X:\Phi(x)=0\}
\]
for its zero set. For $G\in\{\mathrm{SU}(2),\mathrm{SO}(3)\}$, the adjoint bundle is an
oriented Euclidean bundle of rank three. Outside $Z(\Phi)$ we use the
orthogonal decomposition
\[
    \mathfrak g_P
    =
    \mathfrak g_P^{\parallel}\oplus\mathfrak g_P^{\perp},
    \qquad
    \mathfrak g_P^{\parallel}=\langle\Phi\rangle,
\]
and accordingly decompose sections $\xi$ as
\[
    \xi=\xi^{\parallel}+\xi^{\perp}.
\]
All identities in the Lie algebra used below are unchanged under the
isometric Lie algebra identification
$\mathfrak{su}(2)\simeq\mathfrak{so}(3)$ fixed above.

We denote by $\mathscr{A}(P)$ the space of smooth connections on $P$. Given $\nabla\in\mathscr{A}(P)$, we denote by $d_\nabla$ the induced exterior
covariant derivative and by
\[
\Delta_\nabla
=
d_\nabla d_\nabla^*
+
d_\nabla^*d_\nabla
\]
the corresponding Hodge Laplacian.  The rough Laplacian is denoted by
$\nabla^*\nabla$.  To avoid confusion with the connection $\nabla$ on the
bundle, we denote by $D$ the Levi--Civita connection of $(X,g)$. Covariant derivatives of adjoint-valued tensor fields are
taken with respect to the coupled Levi--Civita and bundle connection. 

In deriving estimates, we denote by $c,C>0$ generic constants which may
change from line to line unless otherwise stated. Unless explicitly indicated
otherwise, all such constants are allowed to depend only on the fixed
background geometry of $(X,g)$, the asymptotic structure, the form $\Theta$,
and the topology of the bundle under consideration. In particular, they are
independent of the large mass parameter $m$ and of the index $i$ along large mass
sequences.

Given nonnegative quantities $A$ and $B$, we write
\[
A\lesssim B
\]
if there exists a constant $C>0$, depending only on the fixed geometric
background, such that
\[
A\leqslant CB.
\]
More generally, given parameters $c_1,\dots,c_k$, we write
\[
A\lesssim_{c_1,\dots,c_k} B
\]
to mean that
\[
A\leqslant CB
\]
for a constant $C>0$ depending only on
$c_1,\dots,c_k$ and on the fixed background geometry.

When both
\[
A\lesssim_{c_1,\dots,c_k} B
\qquad\text{and}\qquad
B\lesssim_{c_1,\dots,c_k} A
\]
hold, we write
\[
A\sim_{c_1,\dots,c_k} B.
\]
If the dependence on a compact subset $K\subset X$ is relevant, we indicate this explicitly by writing, for instance,
\[
A\lesssim_K B.
\]
We also use the standard Landau notation $O(\cdot)$ and $o(\cdot)$
with the same convention regarding the dependence of the implicit constants.

We write $\mathcal D_k(X)$ for the space of $k$-currents on $X$,
$\mathscr R_{k,\mathrm{loc}}(X)$ for locally (integer) rectifiable $k$-currents, and $\mathscr I_{k,\mathrm{loc}}(X)$ for locally integral $k$-currents, that is, $T\in\mathscr I_{k,\mathrm{loc}}(X)$ if and only if $T\in\mathscr R_{k,\mathrm{loc}}(X)$ and $\partial T\in\mathscr R_{k-1,\mathrm{loc}}(X)$. For a current $T\in\mathcal D_k(X)$, we denote by $\|T\|$ its associated weight (or mass) measure. A closed locally integral current means $T\in\mathscr I_{k,\mathrm{loc}}(X)$ with $\partial T=0$; when $\|T\|(X)<\infty$, we call such a current an integral cycle.

\subsection{The Yang--Mills--Higgs equations, finite mass, and Higgs defect estimates}\label{subsec: YMH_config_finite_mass}

Let $P\to X$ be a principal $G$-bundle, where $G$ is a compact Lie group. Given a
configuration $(\nabla,\Phi)\in\mathscr A(P)\times\Gamma(\mathfrak{g}_P)$, we define the Yang--Mills--Higgs energy density by
\begin{equation}\label{eq: YMH_density}
e(\nabla,\Phi):=|F_\nabla|^2+|\nabla\Phi|^2.
\end{equation}
The corresponding Yang--Mills--Higgs functional over a precompact open set
$U\Subset X$ is
\[
\cE_U(\nabla,\Phi)
:=
\frac{1}{2}\int_U e(\nabla,\Phi).
\]
The Euler--Lagrange equations associated with compactly supported variations
are
\begin{equation}\tag{YMH}\label{eq:YMH_second_order}
\begin{aligned}
\nabla^*\nabla\Phi &=0,
\\
d_\nabla^*F_\nabla &= [\nabla\Phi,\Phi].
\end{aligned}
\end{equation}
We refer to \eqref{eq:YMH_second_order} as the Yang--Mills--Higgs equations.
A configuration $(\nabla,\Phi)$ is called \textbf{irreducible} if
$\nabla\Phi\not\equiv0$, and \textbf{reducible} otherwise. This terminology
concerns the pair $(\nabla,\Phi)$ and does not assert that $\nabla$ is
irreducible as a connection.

We work exclusively on complete asymptotically conical manifolds $X$ with only one end. In this setting, the relevant asymptotic quantity is the mass. 
\begin{definition}[Finite mass]
	A configuration $(\nabla,\Phi)\in\mathscr{A}(P)\times\Gamma(\mathfrak{g}_P)$ is said to have \textbf{finite mass} $m\in[0,\infty)$ if for any $x\in X$ one has
	\begin{equation}\label{eq: finite_mass}
		\lim_{R\to\infty}\sup_{X\setminus B_R(x)} |m-|\Phi|| = 0.
	\end{equation}
\end{definition}
\begin{remark}\label{rmk: finite_mass}
Suppose that $(\nabla,\Phi)$ has finite mass $m\geqslant0$ and satisfies
the first equation in \eqref{eq:YMH_second_order}. Then
\begin{equation}\label{eq: subharmonic}
    \frac12\Delta|\Phi|^2
    =
    \langle\nabla^*\nabla\Phi,\Phi\rangle-|\nabla\Phi|^2
    =
    -|\nabla\Phi|^2
    \leqslant0.
\end{equation}
Thus $|\Phi|^2$ is subharmonic, and the maximum principle
\cite[Proposition~IV.3.3]{Jaffe1980}, together with
\eqref{eq: finite_mass}, gives $|\Phi|\leqslant m$ on $X$. Moreover,
either $|\Phi|\equiv m$ or $|\Phi|<m$ everywhere. In particular, $m=0$
implies $\Phi\equiv0$. If $m>0$, then for every $\delta\in(0,1)$ there is a
sufficiently large compact subset $K_{\Phi,\delta}\subset X$ outside of which
$|\Phi|>\delta m$.
\end{remark}

By Remark~\ref{rmk: finite_mass}, a finite mass irreducible configuration
solving \eqref{eq:YMH_second_order} gives rise to the nonconstant bounded
subharmonic function $|\Phi|^2$. The underlying manifold must therefore
be nonparabolic. We work below on one-ended AC manifolds, which are
nonparabolic in dimensions $n\geqslant3$. The finite energy,
concentration, and residual abelian theories in dimension three are
developed in
\cites{fadel2023asymptotics,fadeloliveira2026limitv5,fadel2026abelian};
the principal applications in this paper occur in dimensions $6$ and
$7$.

From now on, let $(X,g)$ be a complete AC $n$-manifold with only one end, and let $\rho$ be a radius function for $(X,g)$. Let $\mathscr{D}^{1,2}(X)$ denote the space of real-valued functions $f\in W_{\text{loc}}^{1,2}(X)$ for which $df\in L^2(X)$, and let $\mathscr{D}_0^{1,2}(X)$ be the Hilbert space obtained from the completion of the space $C_c^{\infty}(X)$ of smooth compactly supported functions on $X$ with respect to the norm 
\[
\|f\|_{\mathscr{D}_0^{1,2}(X)}:=\|df\|_{L^2(X)}.
\]
By the $L^2$ Sobolev inequality on $(X,g)$
\cite[Theorem~2.6, \S2.2]{van2009regularity}, there is a continuous
embedding
\[
    \mathscr{D}_0^{1,2}(X)
    \hookrightarrow
    L^{\frac{2n}{n-2}}(X).
\]
In particular, by H\"older's inequality on compact subsets, elements of
$\mathscr{D}_0^{1,2}(X)$ belong to $W_{\mathrm{loc}}^{1,2}(X)$, and hence
\[
    \mathscr{D}_0^{1,2}(X)
    \subset
    \mathscr{D}^{1,2}(X).
\]
The following decomposition is proved in
\cite[Lemma~2.27]{fadel2023asymptotics}.
\begin{lemma}\label{lem: mass}
For each $f\in\mathscr{D}^{1,2}(X)$, there exists a unique real number
$m(f)\in\mathbb R$ such that
\[
    f-m(f)\in\mathscr{D}_0^{1,2}(X).
\]
Consequently,
\[
    \mathscr{D}^{1,2}(X)
    =
    \mathscr{D}_0^{1,2}(X)\oplus\mathbb R.
\]
Moreover, $m(f)$ is the unique real number such that
\[
    f-m(f)\in L^{\frac{2n}{n-2}}(X).
\]
\end{lemma}
\begin{remark}\label{rmk: mass_and_asymptotic_value}
Let $f\in\mathscr{D}^{1,2}(X)\cap C^0(X)$ and suppose that
$f\to\ell\in\mathbb R$ uniformly at infinity. Then $m(f)=\ell$.

Indeed, Lemma~\ref{lem: mass} gives
$f-m(f)\in L^{\frac{2n}{n-2}}(X)$. If $m(f)\neq\ell$, the uniform
convergence of $f$ to $\ell$ implies that $|f-m(f)|$ is bounded below by
a positive constant outside a compact set. This contradicts the infinite
volume of the AC end. Hence $m(f)=\ell$.
\end{remark}

We recall the following Poisson representation for harmonic Higgs fields
on one-ended AC $n$-manifolds, $n\geqslant3$, proved in
\cite[Theorem~3.11]{fadel2023asymptotics}.
\begin{theorem}\label{thm: alternative_finite_mass}
Suppose that
$(\nabla,\Phi)\in\mathscr A(P)\times\Gamma(\mathfrak g_P)$ satisfies the
first equation in \eqref{eq:YMH_second_order} and
$\nabla\Phi\in L^2(X)\cap L^{2(n-1)}(X)$. Let $G(x,y)>0$ be the minimal
positive Green's function of $(X^n,g)$. Then
\begin{equation}\label{eq: def_w}
    w(x):=2\int_XG(x,\cdot{})|\nabla\Phi|^2,
    \qquad x\in X,
\end{equation}
defines the unique nonnegative smooth solution of
$\Delta w=2|\nabla\Phi|^2$ that decays uniformly to zero at infinity.
Moreover, there exists $m\in[0,\infty)$ such that
\begin{equation}\label{eq: alt_characterization_m}
    w=m^2-|\Phi|^2.
\end{equation}
In particular, $(\nabla,\Phi)$ has finite mass $m$; that is,
\eqref{eq: finite_mass} holds.
\end{theorem}

The Bochner estimate of Lemma~\ref{lem:bounded_curvature_moser} also gives
the following global integrability and decay statement, by the argument
of \cite[Corollary~3.2]{fadel2020asymptotic}.
\begin{lemma}\label{cor:YMH_Lpbounds}
Let $(\nabla,\Phi)\in\mathscr A(P)\times\Gamma(\mathfrak g_P)$ be a
solution of \eqref{eq:YMH_second_order}. If
$|F_\nabla|\in L^\infty(X)$ and $|\nabla\Phi|\in L^2(X)$, then
$|\nabla\Phi|^2\in L^\infty(X)\cap L^p(X)$ for every
$p\in[1,\infty)$ and decays uniformly to zero along the end.
\end{lemma}

On a one-ended AC $\mathrm G_2$-manifold $X^7$, the authors proved in
\cite[Theorem~4.1]{fadel2020asymptotic} that solutions $(\nabla,\Phi)$ of
\eqref{eq:YMH_second_order} with $|\nabla\Phi|\in L^2(X)$ and bounded
curvature $|F_{\nabla}|\in L^{\infty}(X)$ have finite mass $m$ and admit a
Poisson representation for $m^2-|\Phi|^2$. Combining the more general
Poisson representation result of Theorem~\ref{thm: alternative_finite_mass}
with Lemma~\ref{cor:YMH_Lpbounds}, and crucially Lemma~\ref{lem: mass}, we
obtain the following refinement for solutions of
\eqref{eq:YMH_second_order} on one-ended AC $n$-manifolds. In addition to
finite mass, it identifies the \textbf{Higgs defect} $m-|\Phi|$ as an element of $\mathscr{D}_0^{1,2}(X)$ and gives a critical Sobolev estimate. This estimate is the main input for the large mass defect analysis in
Lemma~\ref{lem:largeness}, which is subsequently used in
Proposition~\ref{prop:PPS_hypotheses_main_setting}.
\begin{theorem}\label{thm: finite_mass_and_estimates}
Let $(\nabla,\Phi)$ be a solution of \eqref{eq:YMH_second_order}, and
assume that $|F_\nabla|\in L^\infty(X)$ and
$|\nabla\Phi|\in L^2(X)$. Then $(\nabla,\Phi)$ has finite mass
$m\geqslant0$, $m-|\Phi|\in\mathscr D_0^{1,2}(X)$, and
    \begin{equation}\label{ineq: mass_Sobolev}
    \|m - |\Phi|\|_{L^{\frac{2n}{n-2}}(X)}\lesssim \|\nabla\Phi\|_{L^2(X)}.
    \end{equation}
\end{theorem}
\begin{proof}
    By Lemma~\ref{cor:YMH_Lpbounds}, we have
    \[
        f:=2|\nabla\Phi|^2\in L^{(2^*)'}(X),
        \qquad
        2^*:=\frac{2n}{n-2},
        \qquad
        (2^*)'=\frac{2n}{n+2}.
    \]
Thus, the linear functional
    \[
        h \longmapsto \int_X fh 
    \]
    is bounded on $\mathscr{D}_0^{1,2}(X)$. Indeed, by H\"older's inequality and the Sobolev
    embedding $\mathscr{D}_0^{1,2}(X)\hookrightarrow L^{2^*}(X)$,
    \[
        \left|\int_X fh \right|
        \leqslant
        \|f\|_{L^{(2^*)'}(X)}
        \|h \|_{L^{2^*}(X)}
        \lesssim
        \|f\|_{L^{(2^*)'}(X)}
        \|h \|_{\mathscr{D}_0^{1,2}(X)}.
    \]
Hence, by the Riesz representation theorem, there exists a unique
    $u\in\mathscr{D}_0^{1,2}(X)$ such that
    \[
        \int_X \langle du,dh \rangle
        =
        \int_X fh ,
        \qquad
        \forall h \in\mathscr{D}_0^{1,2}(X).
    \]
Equivalently,
    \[
        \Delta u = 2|\nabla\Phi|^2
    \]
    weakly on $X$. Since $u\in\mathscr{D}_0^{1,2}(X)\hookrightarrow W^{1,2}_{\mathrm{loc}}(X)$,
    elliptic regularity implies that $u$ is smooth.

    On the other hand, Lemma~\ref{cor:YMH_Lpbounds} gives $|\nabla\Phi|\in L^2(X)\cap L^{2(n-1)}(X)$.
    Therefore, Theorem~\ref{thm: alternative_finite_mass} applies. Thus
    \[
        w(x):=2\int_X G(x,y)|\nabla\Phi(y)|^2\,dy
    \]
    is the unique nonnegative smooth solution of
    \[
        \Delta w=2|\nabla\Phi|^2
    \]
    which decays uniformly to zero at infinity. Moreover, there exists
    $m\in[0,\infty)$ such that
    \[
        w=m^2-|\Phi|^2.
    \]
In particular, $(\nabla,\Phi)$ has finite mass $m$, that is, $|\Phi|\longrightarrow m$ uniformly along the end.

    We claim that $u=w$. First we show that
$w\in\mathscr{D}_0^{1,2}(X)$. We have
    \[
        dw=-d|\Phi|^2=-2\langle \nabla\Phi,\Phi\rangle,
    \]
and the maximum principle recalled in Remark~\ref{rmk: finite_mass} gives
$|\Phi|\leqslant m$. Therefore,
    \[
        |dw|\leqslant 2m|\nabla\Phi|.
    \]
Since $|\nabla\Phi|\in L^2(X)$, it follows that $dw\in L^2(X)$, and hence
$w\in\mathscr{D}^{1,2}(X)$. By Lemma~\ref{lem: mass}, there exists a unique constant
$m(w)\in\mathbb R$ such that
\[
    w-m(w)\in\mathscr{D}_0^{1,2}(X).
\]
Since $w$ decays uniformly to zero at infinity, Remark~\ref{rmk: mass_and_asymptotic_value} gives $m(w)=0$, and hence
$w\in\mathscr{D}_0^{1,2}(X)$, as required.

Now $u,w\in\mathscr{D}_0^{1,2}(X)$ and both satisfy
\[
    \int_X\langle du,dh\rangle
    =
    \int_X\langle dw,dh\rangle
    =
    \int_X fh,
    \qquad
    h\in\mathscr{D}_0^{1,2}(X).
\]
By uniqueness in the Riesz representation theorem, we conclude that
\[
    u=w=m^2-|\Phi|^2.
\]
Finally, since $|\Phi|\longrightarrow m$
    uniformly along the end and $|d|\Phi||\leqslant |\nabla\Phi|$
    by Kato's inequality, we have $|\Phi|\in\mathscr{D}^{1,2}(X)$. By Lemma~\ref{lem: mass} and Remark~\ref{rmk: mass_and_asymptotic_value}, the unique constant $m(|\Phi|)$ such that
    $|\Phi|-m(|\Phi|)\in\mathscr{D}_0^{1,2}(X)$ is characterized by the asymptotic limit
    of $|\Phi|$. Hence
    \[
        m(|\Phi|)=m,
    \]
    and therefore
    \[
        |\Phi|-m\in\mathscr{D}_0^{1,2}(X).
    \]
Applying the Sobolev inequality to $m-|\Phi|\in\mathscr{D}_0^{1,2}(X)$ gives
    \[
        \|m-|\Phi|\|_{L^{2^*}(X)}
        \lesssim
        \|d|\Phi|\|_{L^2(X)}
        \leqslant
        \|\nabla\Phi\|_{L^2(X)},
    \]
    thereby showing \eqref{ineq: mass_Sobolev}. The proof is complete.
\end{proof}

We next record a potential-theoretic form of the Higgs defect estimate.
For a set $E\subset X$, define its homogeneous $2$-capacity by
\begin{equation}\label{eq:def_homogeneous_2_capacity}
\begin{aligned}
    \operatorname{Cap}_2(E)
    :=
    \inf\biggl\{
        \int_X|df|^2:\;&
        f\in\mathscr D_0^{1,2}(X),\\
        &f\geqslant1
        \text{ almost everywhere on an open neighbourhood of }E
    \biggr\}.
\end{aligned}
\end{equation}
We use the standard monotonicity and countable subadditivity of
$\operatorname{Cap}_2$. The latter also follows directly from the
definition: if $f_j$ are admissible for $E_j$, then
\[
    f
    :=
    \min\left\{
        1,
        \left(\sum_j f_j^2\right)^{1/2}
    \right\}
\]
is admissible for $\bigcup_jE_j$ and satisfies
$\int_X|df|^2\leqslant\sum_j\int_X|df_j|^2$.

Suppose that $(\nabla,\Phi)$ has finite mass $m>0$. For
$\delta\in(0,1)$, set
\[
    Z_\delta(\Phi)
    :=
    \{x\in X:|\Phi(x)|\leqslant\delta m\}.
\]
When $(\nabla,\Phi)$ is furthermore a solution of
the first equation in \eqref{eq:YMH_second_order}, the last observation in
Remark~\ref{rmk: finite_mass} shows that the closed set
$Z_\delta(\Phi)$ is bounded and therefore compact in the complete
manifold $(X,g)$.

\begin{proposition}[Capacity and volume of relative Higgs sublevel regions]
\label{prop:Higgs_sublevel_volume}
Let $(\nabla,\Phi)$ satisfy the hypotheses of
Theorem~\ref{thm: finite_mass_and_estimates}, and suppose that its mass
$m$ is positive. Then, for every $\delta\in(0,1)$,
\begin{equation}\label{eq:Higgs_sublevel_capacity}
    \operatorname{Cap}_2\bigl(Z_\delta(\Phi)\bigr)
    \leqslant
    \frac{\|\nabla\Phi\|_{L^2(X)}^2}
         {m^2(1-\delta)^2}.
\end{equation}
Moreover,
\begin{equation}\label{eq: vol_zero_delta_Phi}
    \vol\bigl(Z_\delta(\Phi)\bigr)^{\frac{n-2}{2n}}
    \lesssim
    \operatorname{Cap}_2\bigl(Z_\delta(\Phi)\bigr)^{\frac12}
    \leqslant
    \frac{\|\nabla\Phi\|_{L^2(X)}}
         {m(1-\delta)}.
\end{equation}
In particular, if
\[
    \|\nabla\Phi\|_{L^2(X)}^2
    \leqslant
    \Lambda m,
\]
then
\begin{equation}\label{eq:Higgs_sublevel_volume_linear_growth}
    \operatorname{Cap}_2\bigl(Z_\delta(\Phi)\bigr)
    \leqslant
    \frac{\Lambda}{m(1-\delta)^2},
    \qquad
    \vol\bigl(Z_\delta(\Phi)\bigr)
    \lesssim
    \left(
        \frac{\Lambda}
             {m(1-\delta)^2}
    \right)^{\frac{n}{n-2}}.
\end{equation}
\end{proposition}

\begin{proof}
Fix $\delta<\delta'<1$ and set
\[
    \eta_{\delta'}
    :=
    \min\left\{
        1,
        \frac{m-|\Phi|}{m(1-\delta')}
    \right\}.
\]
Since $m-|\Phi|\in\mathscr D_0^{1,2}(X)$ by
Theorem~\ref{thm: finite_mass_and_estimates}, the Lipschitz truncation
$\eta_{\delta'}$ belongs to $\mathscr D_0^{1,2}(X)$. On
$Z_\delta(\Phi)$ one has
\[
    \frac{m-|\Phi|}{m(1-\delta')}
    \geqslant
    \frac{1-\delta}{1-\delta'}
    >1.
\]
By continuity, $\eta_{\delta'}=1$ on an open neighbourhood of
$Z_\delta(\Phi)$, so it is admissible in
\eqref{eq:def_homogeneous_2_capacity}. Kato's inequality gives
\[
\begin{aligned}
    \operatorname{Cap}_2\bigl(Z_\delta(\Phi)\bigr)
    &\leqslant
    \int_X|d\eta_{\delta'}|^2 \\
    &\leqslant
    \frac{1}{m^2(1-\delta')^2}
    \int_X|d|\Phi||^2 \\
    &\leqslant
    \frac{\|\nabla\Phi\|_{L^2(X)}^2}
         {m^2(1-\delta')^2}.
\end{aligned}
\]
Letting $\delta'\downarrow\delta$ proves
\eqref{eq:Higgs_sublevel_capacity}.

The Sobolev inequality
$\mathscr D_0^{1,2}(X)\hookrightarrow L^{2n/(n-2)}(X)$,
recalled from \cite[Theorem~2.6, \S2.2]{van2009regularity}, implies the
isocapacitary estimate
\[
    \vol(E)^{\frac{n-2}{2n}}
    \lesssim
    \operatorname{Cap}_2(E)^{\frac12}
\]
for every measurable set $E$ of finite volume. Applying this to
$E=Z_\delta(\Phi)$ and using
\eqref{eq:Higgs_sublevel_capacity} proves
\eqref{eq: vol_zero_delta_Phi}. The final assertions follow immediately
from the linear growth bound.
\end{proof}

As an immediate consequence,
\eqref{eq: vol_zero_delta_Phi} controls the inradius of
$Z_\delta(\Phi)$.

\begin{lemma}[Radius estimate for Higgs sublevel balls]
\label{lem: Higgs_field_estimate_I}
Let $(\nabla,\Phi)$ be as in
Proposition~\ref{prop:Higgs_sublevel_volume}, and define
\[
    r_\delta(x)
    :=
    \sup
    \left\{
        r>0:
        \sup_{B_r(x)}|\Phi|
        \leqslant
        \delta m
    \right\},
\]
with the convention that the supremum is zero if the set is empty. Then
\[
    r_\delta(x)^{\frac{n-2}{2}}
    \lesssim
    \frac{\|\nabla\Phi\|_{L^2(X)}}
    {m(1-\delta)}.
\]
In particular, under the linear growth bound
\[
    \|\nabla\Phi\|_{L^2(X)}^2\leqslant\Lambda m,
\]
one has
\[
    r_\delta(x)
    \lesssim
    \left(
        \frac{\Lambda}
        {m(1-\delta)^2}
    \right)^{\frac1{n-2}}.
\]
\end{lemma}

\begin{proof}
If $r_\delta(x)=0$, there is nothing to prove. Let
$0<r<r_\delta(x)$. By the definition of the supremum, there exists
$\rho>r$ such that
\[
    \sup_{B_\rho(x)}|\Phi|
    \leqslant
    \delta m.
\]
Since $B_r(x)\subset B_\rho(x)$, it follows that
$B_r(x)\subset Z_\delta(\Phi)$. The uniform lower volume estimate for balls on the fixed AC manifold
\cite[Corollary~2.8]{van2009regularity} gives
$\vol(B_r(x))\gtrsim r^n$. Hence,
\[
    r^{\frac{n-2}{2}}
    \lesssim
    \vol(B_r(x))^{\frac{n-2}{2n}}
    \leqslant
    \vol\bigl(Z_\delta(\Phi)\bigr)^{\frac{n-2}{2n}}
    \lesssim
    \frac{\|\nabla\Phi\|_{L^2(X)}}{m(1-\delta)},
\]
where the last inequality follows from
\eqref{eq: vol_zero_delta_Phi}. Letting
$r\uparrow r_\delta(x)$ yields
\[
    r_\delta(x)^{\frac{n-2}{2}}
    \lesssim
    \frac{\|\nabla\Phi\|_{L^2(X)}}{m(1-\delta)}.
\]
If $\|\nabla\Phi\|_{L^2(X)}^2\leqslant\Lambda m$, then
\[
    r_\delta(x)^{\frac{n-2}{2}}
    \lesssim
    \frac{\Lambda^{\frac12}}{m^{\frac12}(1-\delta)}.
\]
Raising both sides to the power $\frac{2}{n-2}$ gives
\[
    r_\delta(x)
    \lesssim
    \left(
        \frac{\Lambda}{m(1-\delta)^2}
    \right)^{\frac1{n-2}},
\]
as claimed.
\end{proof}

The preceding estimates concern the relative sublevel regions of a single
finite mass configuration. We now consider a sequence whose masses tend
to infinity. The Sobolev control of the Higgs defect, together with the
harmonicity equation, yields local strong convergence after the natural
mass normalization.

\begin{lemma}[Large mass vanishing of the Higgs defect]
\label{lem:largeness}
Let $(\nabla_i,\Phi_i)$ be a sequence of solutions of \eqref{eq:YMH_second_order}. Assume that $|F_{\nabla_i}|\in L^\infty(X)$, $|\nabla_i\Phi_i|\in L^2(X)$, and let $m_i$ be their finite masses. Further suppose that $m_i\longrightarrow+\infty$, and $ \|\nabla_i\Phi_i\|_{L^2(X)}^2 \leqslant \Lambda m_i$ for a constant $\Lambda>0$ independent of $i$. Then, the sequence $(w_i)$ defined by
\[
    w_i
    :=
    \frac{m_i^2-|\Phi_i|^2}{m_i^{3/2}}
\]
is bounded in $\mathscr{D}_0^{1,2}(X)$ and
\[
    w_i\rightharpoonup0
    \qquad\text{weakly in }\mathscr{D}_0^{1,2}(X).
\]
Moreover, for every compact subset $K\Subset X$ and every $1\leqslant q<\frac{2n}{n-2}$, one has
\[
    w_i\longrightarrow0
    \qquad\text{strongly in }L^q(K),
\]
or equivalently,
\[
    m_i^{-\frac q2}
    \int_K
    (m_i-|\Phi_i|)^q
    \longrightarrow0.
\]
In particular,
\begin{equation}\label{eq: largeness_of_Higgs_field}
    \lim_{i\to\infty}
    \frac1{m_i}
    \int_K
    (m_i-|\Phi_i|)^2
    =
    0.
\end{equation}
\end{lemma}

\begin{proof}
By the maximum principle recalled in
Remark~\ref{rmk: finite_mass},
\[
    0\leqslant|\Phi_i|\leqslant m_i.
\]
Consequently,
\begin{equation}\label{ineq: w_m_equiv}
    m_i^{-1/2}(m_i-|\Phi_i|)
    \leqslant
    w_i
    \leqslant
    2m_i^{-1/2}(m_i-|\Phi_i|).
\end{equation}
Since $|\Phi_i|\to m_i$ uniformly along the end, $w_i$ converges
uniformly to zero at infinity. Moreover,
\[
    dw_i
    =
    -2m_i^{-3/2}
    \langle\nabla_i\Phi_i,\Phi_i\rangle,
\]
and therefore
\[
    |dw_i|
    \leqslant
    2m_i^{-1/2}|\nabla_i\Phi_i|.
\]
It follows that $dw_i\in L^2(X)$ and
\begin{equation}\label{ineq: uniform_bound_L2_grad_w_m}
    \|w_i\|_{\mathscr{D}_0^{1,2}(X)}^2
    =
    \|dw_i\|_{L^2(X)}^2
    \leqslant
    4m_i^{-1}
    \|\nabla_i\Phi_i\|_{L^2(X)}^2
    \leqslant
    4\Lambda.
\end{equation}
Since $w_i$ has finite Dirichlet energy and tends uniformly to zero at
infinity, Lemma~\ref{lem: mass} and
Remark~\ref{rmk: mass_and_asymptotic_value} imply that
\[
    w_i\in\mathscr{D}_0^{1,2}(X).
\]
The Sobolev estimate \eqref{ineq: mass_Sobolev} also gives the uniform
critical bound
\begin{equation}\label{ineq: uniform_L_p_bound}
\begin{split}
    \|w_i\|_{L^{\frac{2n}{n-2}}(X)}
    &\leqslant
    2m_i^{-1/2}
    \|m_i-|\Phi_i|\|_{L^{\frac{2n}{n-2}}(X)}
    \\
    &\lesssim
    m_i^{-1/2}
    \|\nabla_i\Phi_i\|_{L^2(X)}
    \leqslant
    \Lambda^{1/2}.
\end{split}
\end{equation}
We next identify the weak $\mathscr{D}_0^{1,2}(X)$ limit. From
\[
    \Delta|\Phi_i|^2
    =
    -2|\nabla_i\Phi_i|^2
\]
we obtain
\begin{equation}\label{eq: Laplacian_w_m}
    \Delta w_i
    =
    2m_i^{-3/2}|\nabla_i\Phi_i|^2.
\end{equation}
Thus
\begin{equation}\label{ineq: L1_estimate_Laplacian_w_m}
    \|\Delta w_i\|_{L^1(X)}
    =
    2m_i^{-3/2}
    \|\nabla_i\Phi_i\|_{L^2(X)}^2
    \leqslant
    2\Lambda m_i^{-1/2}
    \longrightarrow0.
\end{equation}
The bound \eqref{ineq: uniform_bound_L2_grad_w_m} shows that every
subsequence of $(w_i)$ admits a further subsequence converging weakly in
$\mathscr{D}_0^{1,2}(X)$ to some $w\in\mathscr{D}_0^{1,2}(X)$. For every
$\chi\in C_c^\infty(X)$, integration by parts and
\eqref{ineq: L1_estimate_Laplacian_w_m} give
\[
\begin{split}
    \left|
        \int_X\langle dw_i,d\chi\rangle
    \right|
    &=
    \left|
        \int_X\chi\,\Delta w_i
    \right|
    \\
    &\leqslant
    \|\chi\|_{L^\infty(X)}
    \|\Delta w_i\|_{L^1(X)}
    \longrightarrow0.
\end{split}
\]
Passing to the weak limit yields
\[
    \int_X\langle dw,d\chi\rangle=0
    \qquad
    \text{for every }\chi\in C_c^\infty(X).
\]
Since $C_c^\infty(X)$ is dense in $\mathscr{D}_0^{1,2}(X)$, this identity extends to
every test function in $\mathscr{D}_0^{1,2}(X)$. Testing with $w$ gives
\[
    \|w\|_{\mathscr{D}_0^{1,2}(X)}^2=0,
\]
and therefore $w=0$ in $\mathscr{D}_0^{1,2}(X)$. Since every weakly convergent
subsequence has the same limit, the full sequence satisfies
\[
    w_i\rightharpoonup0
    \qquad\text{weakly in }\mathscr{D}_0^{1,2}(X).
\]
It remains to prove the local strong convergence. Fix
$K\Subset X$, and choose a precompact open set $U$ with smooth boundary
such that
\[
    K\Subset U\Subset X.
\]
The global bounds
\eqref{ineq: uniform_bound_L2_grad_w_m} and
\eqref{ineq: uniform_L_p_bound}, together with H\"older's inequality, imply
that $(w_i)$ is bounded in $W^{1,2}(U)$. By the Rellich theorem, it is
precompact in $L^q(U)$ for every
\[
    1\leqslant q<\frac{2n}{n-2}.
\]
Every strongly convergent subsequence has limit zero, because $w_i$
converges weakly to zero in $\mathscr{D}_0^{1,2}(X)$ and hence weakly in
$L^{\frac{2n}{n-2}}(X)$. The usual subsequence argument therefore gives
\[
    w_i\longrightarrow0
    \qquad\text{strongly in }L^q(K)
\]
for the full sequence.

Finally, the pointwise comparison \eqref{ineq: w_m_equiv} shows that this
is equivalent to
\[
    m_i^{-\frac q2}
    \int_K
    (m_i-|\Phi_i|)^q
    \longrightarrow0.
\]
Taking $q=2$ gives
\eqref{eq: largeness_of_Higgs_field}.
\end{proof}

\begin{remark}\label{rmk:linear_energy_growth}
For the $\Theta$-monopole sequences introduced later in the {\mainsettingref} of Section~\ref{sec: concentration}, the linear energy bound in
Lemma~\ref{lem:largeness} follows from the fixed monopole class assumption.
Indeed, the asymptotic energy identity
\eqref{eq: theta_monopole_energy_formula}, obtained in
Subsection~\ref{subsec: asymptotic_analysis_theta_monopoles}, gives
\[
    \|\nabla_i\Phi_i\|_{L^2(X)}^2
    =
    \mathcal E^\Theta(\nabla_i,\Phi_i)
    =
    4\pi k m_i,
    \qquad
    k:=\langle\beta\cup\Theta_\infty,[\Sigma]\rangle.
\]
Consequently, \eqref{eq: largeness_of_Higgs_field} verifies hypothesis
{\rm(H2)} of Theorem~\ref{thm:PPS_AC_version}, as used in
Proposition~\ref{prop:PPS_hypotheses_main_setting}. The same identity also specializes
Proposition~\ref{prop:Higgs_sublevel_volume} and
Lemma~\ref{lem: Higgs_field_estimate_I}. For every fixed
$\delta\in(0,1)$, writing
$Z_{\delta,i}:=Z_\delta(\Phi_i)$ and denoting by
$r_{\delta,i}(x)$ the corresponding sublevel radius, one obtains
\begin{equation}\label{eq:Main_Setting_Higgs_sublevel_capacity_volume_radius}
\begin{gathered}
    \operatorname{Cap}_2(Z_{\delta,i})
    \leqslant
    \frac{4\pi k}{m_i(1-\delta)^2},\\
    \vol(Z_{\delta,i})
    \lesssim
    \left(
        \frac{k}
             {m_i(1-\delta)^2}
    \right)^{\frac{n}{n-2}},
    \qquad
    \sup_{x\in X}r_{\delta,i}(x)
    \lesssim
    \left(
        \frac{k}
             {m_i(1-\delta)^2}
    \right)^{\frac1{n-2}}.
\end{gathered}
\end{equation}
Thus the relative Higgs sublevel regions have vanishing $2$-capacity,
volume and inradius as $m_i\to\infty$.
\end{remark}

\begin{corollary}[Quasi-everywhere large mass normalization]
\label{cor:quasi_everywhere_Higgs_normalization}
In the {\mainsettingref}, after passing to a further subsequence such that
\[
    \sum_{i=1}^{\infty}\frac1{m_i}<\infty,
\]
the following holds. For every $\delta\in(0,1)$, define the strict
set-theoretic upper limit
\begin{equation}\label{eq:def_recurrent_relative_Higgs_sublevel}
    \mathcal L_\delta
    :=
    \bigcap_{N\geqslant1}
    \bigcup_{i\geqslant N}Z_{\delta,i}.
\end{equation}
Then
\[
    \operatorname{Cap}_2(\mathcal L_\delta)=0.
\]
Consequently, there exists a set $\mathcal N\subset X$ with
$\operatorname{Cap}_2(\mathcal N)=0$ such that
\begin{equation}\label{eq:quasi_everywhere_Higgs_normalization}
    \frac{|\Phi_i(x)|}{m_i}
    \longrightarrow1
    \qquad
    \text{for every }x\in X\setminus\mathcal N.
\end{equation}
Equivalently, the convergence holds $\operatorname{Cap}_2$-quasi-everywhere on $X$. In particular, the strict set-theoretic upper limit of the Higgs zero
sets,
\[
    \bigcap_{N\geqslant1}
    \bigcup_{i\geqslant N}\Phi_i^{-1}(0),
\]
has zero $2$-capacity.
\end{corollary}

\begin{proof}
By monotonicity and countable subadditivity of the capacity, together
with
\eqref{eq:Main_Setting_Higgs_sublevel_capacity_volume_radius},
\[
\begin{aligned}
    \operatorname{Cap}_2(\mathcal L_\delta)
    &\leqslant
    \inf_{N\geqslant1}
    \operatorname{Cap}_2
    \left(
        \bigcup_{i\geqslant N}Z_{\delta,i}
    \right)\\
    &\leqslant
    \frac{4\pi k}{(1-\delta)^2}
    \inf_{N\geqslant1}
    \sum_{i\geqslant N}\frac1{m_i}
    =0.
\end{aligned}
\]
Set
\[
    \mathcal N
    :=
    \bigcup_{j=2}^{\infty}
    \mathcal L_{1-1/j}.
\]
Another application of countable subadditivity gives
$\operatorname{Cap}_2(\mathcal N)=0$. If $x\notin\mathcal N$, then for
every $j\geqslant2$ one has
\[
    |\Phi_i(x)|
    >
    \left(1-\frac1j\right)m_i
\]
for all sufficiently large $i$. Since $|\Phi_i|\leqslant m_i$ by the
maximum principle, letting $j\to\infty$ proves
\eqref{eq:quasi_everywhere_Higgs_normalization}. Finally,
$\Phi_i^{-1}(0)\subset Z_{\delta,i}$ for every $\delta\in(0,1)$, which
gives the last assertion.
\end{proof}

\begin{remark}[Strict versus Kuratowski upper limits]
\label{rmk:capacity_does_not_control_Kuratowski_zero_limit}
The upper limit in
\eqref{eq:def_recurrent_relative_Higgs_sublevel} contains the points
which themselves belong to infinitely many relative sublevel regions.
It does not contain the closures appearing in the Kuratowski upper
limit
\[
    \mathcal Z
    =
    \bigcap_{N\geqslant1}
    \overline{\bigcup_{i\geqslant N}\Phi_i^{-1}(0)}
\]
defined in \eqref{eq:def_limiting_zero_set}. Since a set of zero $2$-capacity may be dense, Corollary~\ref{cor:quasi_everywhere_Higgs_normalization} does not preclude the Higgs zero sets from becoming arbitrarily dense in an open set,
nor does it preclude $\mathcal Z$ from containing such a set. Thus the
capacity conclusion does not by itself resolve the density phenomenon
raised in \cite[Remark~3.10]{li2025large}; it instead shows that
recurrent failure of $|\Phi_i|/m_i$ to approach $1$ at a fixed point
is confined to a set of zero $2$-capacity.

A different consequence of Li's covering estimates is established later
in Proposition~\ref{prop:Li_cover_Hausdorff_content_off_S} and
Corollary~\ref{cor:Higgs_zero_Hausdorff_content_off_S}: on every compact
set $K\Subset X\setminus\mathcal S$, the curvature concentration loci and the
Higgs zero sets have vanishing $(n-3)$-dimensional unrestricted
Hausdorff content. Proposition~\ref{prop:Higgs_zero_density_rate_obstruction}
further gives a quantitative obstruction to rapid density of the
zero sets. As explained in
Remark~\ref{rmk:Li_content_does_not_transfer_to_Cinfty}, these estimates for the varying loci
still do not pass automatically to the Kuratowski upper limits
$\mathcal C$ and $\mathcal Z$.
\end{remark}

\begin{remark}[Subcritical and critical convergence]
\label{rmk: what_is_proved_largeness}
Lemma~\ref{lem:largeness} gives weak convergence in
\[
    \mathscr D_0^{1,2}(X)
\]
and strong local convergence in $L^q$ whenever $q<2n/(n-2)$. It gives
neither strong convergence in the former space nor convergence in the
critical space $L^{2n/(n-2)}$. In particular, the argument does not
exclude concentration of the normalized Higgs energy measures
$m_i^{-1}|\nabla_i\Phi_i|^2\,\vol$. Either critical conclusion would
require additional nonconcentration information.

The quasi-everywhere conclusion of
Corollary~\ref{cor:quasi_everywhere_Higgs_normalization} is instead a
potential-theoretic consequence of the summable capacity bounds along a
further subsequence; it does not strengthen the preceding convergence to
strong convergence in either critical norm.
\end{remark}

\subsection{\texorpdfstring{$\Theta$-monopoles}{Theta-monopoles} and intermediate energies}
\label{subsec: monopoles}

Let $(X^n,g)$ be endowed with a closed $(n-3)$-form
$\Theta\in\Omega^{n-3}(X)$, and let $P\to X$ be a principal $G$-bundle.

\begin{definition}[cf. {\cite[\S1.3]{oliveira2014thesis}}]
A configuration
$(\nabla,\Phi)\in\mathscr A(P)\times\Gamma(\mathfrak g_P)$ is called a
\textbf{$\Theta$-monopole} if
\begin{equation}\label{eq:Theta_monopole}
    F_\nabla\wedge\Theta
    =
    *\nabla\Phi.
\end{equation}
\end{definition}

When $n=3$ and $\Theta=1$, this is the classical Bogomolny monopole
equation. The higher-dimensional special holonomy cases considered in this
paper fit into the same notation as follows.

In the $\mathrm G_2$ case, $X^7$ carries a torsion-free $\mathrm G_2$-structure $\varphi$, inducing the metric and orientation, and we write $\psi:=*\varphi$. The \textbf{$\mathrm G_2$-monopole equation} is obtained from \eqref{eq:Theta_monopole} by taking $\Theta=\psi$,
\[
F_\nabla\wedge\psi
    =
    *\nabla\Phi.
\]
In the Calabi--Yau $3$-fold case, $X^6$ carries a torsion-free
$\mathrm{SU}(3)$-structure $(\omega,\Omega)$, where $\omega$ is the
K\"ahler form and $\Omega$ is the holomorphic volume form. The
\textbf{Calabi--Yau monopole equation} is obtained from \eqref{eq:Theta_monopole} by taking $\Theta=\operatorname{Re}\Omega$,
\[
F_\nabla\wedge\operatorname{Re}\Omega
    =
    *\nabla\Phi
\]
and further imposing
\begin{equation}\label{eq:CY_monopole_primitivity}
F_{\nabla}\wedge\omega^2 = 0.
\end{equation}
The second equation \eqref{eq:CY_monopole_primitivity} is equivalent to the primitivity condition $\Lambda F_\nabla=0$, where $\Lambda=L^*$ is the adjoint of the Lefschetz operator $L(\alpha)=\omega\wedge\alpha$; on $2$-forms, $\Lambda\alpha=\frac12*(\alpha\wedge\omega^2)$. Throughout the paper, when $\Theta=\operatorname{Re}\Omega$, the term $\Theta$-monopole is understood to include this additional equation.

In both special holonomy settings, we also use the auxiliary form
\[
    \Xi\in\Omega^{n-4}(X),
    \qquad
    \Xi=
    \begin{cases}
        \varphi, & \text{in the $\mathrm G_2$ case},\\
        \omega, & \text{in the Calabi--Yau case}.
    \end{cases}
\]
In the $\mathrm G_2$ case, every $\mathfrak g_P$-valued $2$-form $F$
satisfies
\begin{equation}\label{eq:general_calibration_identity}
    |F|^2\operatorname{vol}
    =
    -\langle F\wedge F\rangle\wedge\Xi
    +|F\wedge\Theta|^2\operatorname{vol}.
\end{equation}
In the Calabi--Yau case, the corresponding identity for an arbitrary
$F$ is
\begin{equation}\label{eq:CY_calibration_identity}
    |F|^2\operatorname{vol}
    =
    -\langle F\wedge F\rangle\wedge\omega
    +|F\wedge\operatorname{Re}\Omega|^2\operatorname{vol}
    +|\Lambda F|^2\operatorname{vol}.
\end{equation}
Consequently, \eqref{eq:general_calibration_identity} also holds in the Calabi--Yau case whenever $F$ is primitive, and hence in particular for the curvature of a Calabi--Yau monopole.

The torsion-free special holonomy assumption enters through the closedness
of the forms. In both cases,
\[
    d\Theta=0=d\Xi.
\]
These forms are also calibrations in the sense of Harvey--Lawson: they are closed and have comass one (with respect to the metric they induce). In the $\mathrm G_2$ case, $\varphi$ and $\psi$ calibrate associative $3$-planes and coassociative $4$-planes, respectively. In the Calabi--Yau case, the powers $\omega^k/k!$ calibrate complex $k$-planes, while $\operatorname{Re}\Omega$ calibrates special Lagrangian $3$-planes (of phase zero).

When such manifolds are assumed to be asymptotically conical, we always work with
one AC end. The special holonomy structures are assumed to be asymptotic to
the homogeneous cone structures induced by the nearly K\"ahler structure on
the link in the $\mathrm G_2$ case, and by the Sasaki--Einstein structure
on the link in the Calabi--Yau case. In particular, the homogeneous
asymptotic models $\Theta_C$ and $\Xi_C$ of $\Theta$ and $\Xi$, respectively, have the same algebraic normalization and calibration properties.

The one-ended assumption is natural in the irreducible special holonomy
setting: if a complete Ricci-flat manifold has more than one end, the
Cheeger--Gromoll splitting theorem forces a product splitting, and hence
full holonomy $\mathrm G_2$ or $\mathrm{SU}(3)$ cannot occur. We
nevertheless keep the one-ended AC assumption explicitly, rather than
imposing full holonomy irreducibility at this stage.

The natural variational functional associated with
\eqref{eq:Theta_monopole} is the \textbf{$\Theta$-energy}
\[
\mathcal E_U^\Theta(\nabla,\Phi)
:=
\frac12
\int_U
\left(|F_\nabla\wedge\Theta|^2+|\nabla\Phi|^2\right),
\]
defined over precompact open subsets $U\Subset X$. We write
\[
\mathcal E^\Theta(\nabla,\Phi)
:=
\frac12
\int_X
\left(
|F_\nabla\wedge\Theta|^2+|\nabla\Phi|^2
\right)
\in[0,+\infty]
\]
for the corresponding global quantity. In dimensions $n>3$, we refer
to $\mathcal E^\Theta(\nabla,\Phi)$ as the
\textbf{intermediate energy}; the configuration has finite intermediate
energy when $\mathcal E^\Theta(\nabla,\Phi)<\infty$.

The closedness of $\Theta$, together with the Bianchi identity, gives the
Bogomolny-type identity
\begin{equation}\label{eq:Theta_Bogomolny_identity}
\mathcal E_U^\Theta(\nabla,\Phi)
=
\int_{\partial U}
\langle\Phi,F_\nabla\rangle\wedge\Theta
+
\frac12
\|F_\nabla\wedge\Theta-*\nabla\Phi\|_{L^2(U)}^2
\end{equation}
whenever $U\Subset X$ has smooth boundary. Consequently, $\Theta$-monopoles
minimize $\mathcal E_U^\Theta$ among configurations which agree with them in
a neighbourhood of $\partial U$, and hence are critical points of
$\mathcal E_U^\Theta$ under compactly supported variations.

The Euler--Lagrange equations associated with the $\Theta$-energy are
\begin{equation}\label{eq:Theta_EL}
d_\nabla^*\Pi_\Theta(F_\nabla)
=
[\nabla\Phi,\Phi],
\qquad
\nabla^*\nabla\Phi=0,
\end{equation}
where
\[
\Pi_\Theta(F)
:=
*\bigl(*(F\wedge\Theta)\wedge\Theta\bigr).
\]
Thus every $\Theta$-monopole satisfies \eqref{eq:Theta_EL}. For a general
closed form $\Theta$, the system \eqref{eq:Theta_EL} need not coincide with
\eqref{eq:YMH_second_order}.

In the $\mathrm G_2$ case, however, the Euler--Lagrange equations of $\mathcal E_U^\Theta$ coincide with \eqref{eq:YMH_second_order}. Indeed, by
\eqref{eq:general_calibration_identity}, on every precompact domain
$U\Subset X$, the curvature part of the $\Theta$-energy $\cE_U^{\Theta}(\nabla,\Phi)$ and the full Yang--Mills curvature term of the Yang--Mills--Higgs energy $\cE_U(\nabla,\Phi)$ differ by the Chern--Weil term
\[
\frac12
\int_U
\langle F_\nabla\wedge F_\nabla\rangle\wedge\Xi.
\]
Since $d\Xi=0$, this term has zero first variation under compactly supported
variations, by the Bianchi identity. The Higgs field part of both
$\mathcal E_U^\Theta$ and $\cE_U$ is the usual Dirichlet Higgs energy $\|\nabla\Phi\|_{L^2(U)}^2$. Hence, in this case, \eqref{eq:Theta_EL} is precisely
\eqref{eq:YMH_second_order}.

In the Calabi--Yau case, the $\Theta$-energy alone does not have
\eqref{eq:YMH_second_order} as its Euler--Lagrange system on arbitrary
configurations. Nevertheless, every Calabi--Yau monopole is a solution of
\eqref{eq:YMH_second_order}. Indeed, the preceding algebraic identity
\eqref{eq:CY_calibration_identity} gives, on every precompact domain
$U\Subset X$,
\[
\mathcal E_U(\nabla,\Phi)
=
\mathcal E_U^\Theta(\nabla,\Phi)
+
\frac12\|\Lambda F_\nabla\|_{L^2(U)}^2
-
\frac12\int_U\langle F_\nabla\wedge F_\nabla\rangle\wedge\omega.
\]
The Chern--Weil term has zero first variation under compactly supported
variations, while the first variation of
$\frac12\|\Lambda F_\nabla\|_{L^2(U)}^2$ vanishes at a primitive connection.
Since the monopole equation makes a Calabi--Yau monopole critical for
$\mathcal E_U^\Theta$, it is also critical for the full Yang--Mills--Higgs
energy and therefore satisfies \eqref{eq:YMH_second_order}.

In the classical $3$-dimensional case, the intermediate and full
Yang--Mills--Higgs energies coincide. In higher dimensions, however, they
differ substantially. In particular, irreducible higher-dimensional
$\Theta$-monopoles on asymptotically conical manifolds can have infinite
full Yang--Mills--Higgs energy, whereas the intermediate energy remains
finite (see Remark~\ref{rmk:finite_nu_infinite_mu_i}). This is the main reason why $\mathcal E^\Theta$ is a natural
energy in the noncompact special holonomy setting.

\subsection{Asymptotics of \texorpdfstring{$\Theta$-monopoles}{Theta-monopoles} on AC manifolds}
\label{subsec: asymptotic_analysis_theta_monopoles}

We now recall the asymptotic results for finite intermediate energy
monopoles on asymptotically conical manifolds in a form adapted to the
unified $\Theta$-monopole notation introduced in
\S\ref{subsec: monopoles}.

Throughout the remainder of the paper, $(X^n,g,\Theta)$ denotes one of the two
asymptotically conical special holonomy settings fixed above, i.e.
\begin{itemize}
\item an AC $\mathrm G_2$-manifold, with $\Theta=\psi=*\varphi$ the coassociative $4$-form;

\item an AC Calabi--Yau $3$-fold, with $\Theta=\operatorname{Re}\Omega$ the real part of the holomorphic volume form.
\end{itemize}
In both cases, we also use the associated closed $(n-4)$-form $\Xi$ from
\S\ref{subsec: monopoles}, namely
\[
    \Xi=\varphi
    \quad\text{in the $\mathrm G_2$ case},
    \qquad
    \Xi=\omega
    \quad\text{in the Calabi--Yau case}.
\]
The asymptotic link $\Sigma^{n-1}$ of $(X,g,\Theta)$ is a nearly
K\"ahler $6$-manifold in the $\mathrm G_2$ case and a Sasaki--Einstein
$5$-manifold in the Calabi--Yau case.

The main results of \cite{fadel2020asymptotic} concern the asymptotic geometry
of finite intermediate energy $\mathrm G_2$-monopoles on AC $\mathrm G_2$-manifolds. However, as observed in \cite{li2025large}, the same arguments
apply almost verbatim to Calabi--Yau monopoles on AC Calabi--Yau $3$-folds.
We therefore state the consequences in the unified $\Theta$-monopole
framework.

In what follows, $P\to X$ is a principal $G$-bundle with
$G\in\{\mathrm{SU}(2),\mathrm{SO}(3)\}$, and we let
$(\nabla,\Phi)\in\mathscr A(P)\times\Gamma(\mathfrak g_P)$ be an
irreducible $\Theta$-monopole; that is,
\[
    F_\nabla\wedge\Theta=*\,\nabla\Phi,
    \qquad
    \nabla\Phi\not\equiv0,
\] and, in the Calabi--Yau case, $F_\nabla\wedge\omega^2=0$.
We assume that $(\nabla,\Phi)$ has finite intermediate energy,
\begin{equation*}
    \mathcal E^\Theta(\nabla,\Phi)
    <\infty,
\end{equation*}
and that the curvature decays at least quadratically,
\begin{equation}\label{eq:quadratic_decay_theta}
    |F_\nabla|=O(\rho^{-2})
    \qquad\text{as }\rho\to\infty.
\end{equation}
Under these assumptions, the asymptotic theory gives the following.

\begin{theorem}[{\cite[Main Theorems 1 and 2]{fadel2020asymptotic}}]\label{thm:theta_asymptotics}
Let $(\nabla,\Phi)$ be an irreducible finite intermediate energy
$\Theta$-monopole satisfying \eqref{eq:quadratic_decay_theta}.  Then:
\begin{enumerate}[label=\textnormal{(\alph*)}]
\item There exists $m>0$ such that
\[
    \lim_{\rho\to\infty}|\Phi|=m.
\]
In particular, $(\nabla,\Phi)$ has finite mass and
$\|\Phi\|_{L^\infty(X)}\leqslant m$.

\item The transverse components decay exponentially:
\[
    |[\nabla\Phi,\Phi]|+|[F_\nabla,\Phi]|
    =
    O(e^{-cm\rho})
\]
for some $c>0$.

\item The Higgs field satisfies the sharp decay estimate
\[
    |\nabla\Phi|=O(\rho^{1-n}).
\]

\item There exists a principal $G$-bundle
$P_\infty\to\Sigma$ and a limiting configuration
$(\nabla_\infty,\Phi_\infty)$ on $P_\infty$ such that
\[
    (\nabla,\Phi)|_{\Sigma_R}
    \longrightarrow
    (\nabla_\infty,\Phi_\infty)
\]
uniformly as $R\to\infty$ up to gauge.
Moreover,
\[
    \nabla_\infty\Phi_\infty=0.
\]

\item The limiting connection $\nabla_\infty$ satisfies the
pseudo-Hermitian--Yang--Mills equation on the asymptotic link.  In the
$\mathrm G_2$ case this means that, with respect to the nearly K\"ahler
structure on $\Sigma^6$,
\[
    F_{\nabla_\infty}^{0,2}=0,
    \qquad
    \Lambda F_{\nabla_\infty}=0.
\]
In the Calabi--Yau case, an analogous pseudo-Hermitian--Yang--Mills equation
holds with respect to the Sasaki--Einstein structure on the link $\Sigma^5$.
\end{enumerate}
\end{theorem}

\begin{remark}[The $\mathrm{SO}(3)$ asymptotic variant]
\label{rmk:SO3_asymptotic_variant}
The finite mass conclusion in
Theorem~\ref{thm:theta_asymptotics}\textnormal{(a)} is independent of the
distinction between $\mathrm{SU}(2)$ and $\mathrm{SO}(3)$. Under the
present hypotheses it follows directly from
Theorem~\ref{thm: finite_mass_and_estimates} in
Subsection~\ref{subsec: YMH_config_finite_mass}, which applies to
solutions of \eqref{eq:YMH_second_order} with arbitrary compact structure
group.

The refined asymptotic results of
\cite[Chapters~5--7]{fadel2020asymptotic} are stated for
$\mathrm{SU}(2)$. The structure group dependence in their proofs enters
through the splitting
$\mathfrak g_P=\langle\Phi\rangle\oplus\langle\Phi\rangle^\perp$ on the
nonvanishing locus of $\Phi$, together with the corresponding
$\mathfrak{su}(2)$ bracket identities and norm estimates. With the
inner product \eqref{eq:SO3_inner_prod_convention}, the Lie algebra
isomorphism
$d\pi:\mathfrak{su}(2)\to\mathfrak{so}(3)$ preserves these identities
and estimates. Hence the Bochner inequalities used there hold equally
for the adjoint bundle of a principal $\mathrm{SO}(3)$-bundle. The
subsequent global arguments on the AC end---including the maximum
principle, Hardy inequalities, Moser iteration and weighted integral
estimates---are unchanged, and yield the exponential transverse decay
and the sharp estimate $|\nabla\Phi|=O(\rho^{1-n})$.

Finally, the construction of the limiting configuration uses Uhlenbeck
compactness for compact structure groups, radial gauge and the preceding
decay estimates, and therefore also applies to $\mathrm{SO}(3)$. No
global $\mathrm{SU}(2)$ lift of $P$ is required. The only modification is
the topology of the asymptotic circle reduction and its charge lattice,
as described below; see also
\cite[Remark~1.8]{li2025large}.
\end{remark}

Since $|\Phi_\infty|=m\neq0$ and
$\nabla_\infty\Phi_\infty=0$, the limiting connection reduces to the
stabilizer circle of $\Phi_\infty$ or, equivalently, that of
\[
    \Psi_\infty:=\frac{\Phi_\infty}{m}.
\]
Writing
\[
    T_G:=\operatorname{Stab}_G(\Psi_\infty)
    =
    \begin{cases}
        \mathrm U(1),&G=\mathrm{SU}(2),\\
        \mathrm{SO}(2),&G=\mathrm{SO}(3),
    \end{cases}
\]
the orthogonal complement
\[
    \langle\Psi_\infty\rangle^\perp
    \subset\mathfrak g_{P_\infty}
\]
is a rank-two real bundle over $\Sigma$ on which the stabilizer circle acts by rotations.
We use this action to orient $\langle\Psi_\infty\rangle^\perp$ and regard it as a complex line
bundle, with the sign chosen so that in the $\mathrm{SU}(2)$ case the
resulting class agrees with the eigenline convention used in the present
paper and in \cite{li2025large}. We define the \textbf{monopole class} to be half the Euler class of $\langle\Psi_\infty\rangle^\perp$, viewed in real cohomology:
\[
    \beta
    :=
    \frac12 e(\langle\Psi_\infty\rangle^\perp)
    \in H^2(\Sigma;\mathbb R).
\]
If $G=\mathrm{SU}(2)$ and $L$ is the fundamental eigenline selected by
the preceding convention, then $\langle\Psi_\infty\rangle^\perp\simeq L^{\otimes2}$ and therefore $\beta=c_1(L)\in H^2(\Sigma;\mathbb Z)$. On the other hand, if
$G=\mathrm{SO}(3)$ then $2\beta\in H^2(\Sigma;\mathbb Z)$, but $\beta$ need not be integral.

Another important cohomology class on $\Sigma$, in this asymptotically conical setting, is induced by the calibration $\Theta$. For $R\gg1$, the cohomology class
\[
    [\Theta|_{\Sigma_R}]\in H^{n-3}(\Sigma;\mathbb R)
\]
is independent of $R$. We denote it by
\[
    \Theta_\infty\in H^{n-3}(\Sigma;\mathbb R)
\]
and call it the \textbf{asymptotic cohomology class}.

With the transported inner product convention of \S \ref{ss:Notation}, the Bogomolny
identity and the asymptotic convergence give, for either structure group,
\begin{equation}\label{eq: theta_monopole_energy_formula}
    \mathcal E^\Theta(\nabla,\Phi)
    =
    \lim_{R\to\infty}
    \int_{\Sigma_R}
    \langle\Phi,F_\nabla\rangle\wedge\Theta
    =
    4\pi m\,
    \langle\beta\cup\Theta_\infty,[\Sigma]\rangle.
\end{equation}
In particular, for fixed monopole class $\beta$ and fixed AC geometry,
the intermediate energy grows linearly with the mass.

\begin{remark}[Half-integral asymptotic charge and integral compactly
supported charge]
\label{rmk:SO3_compactly_supported_charge}
In the $\mathrm{SO}(3)$ case, set
$\gamma:=e(\langle\Psi_\infty\rangle^\perp)=2\beta$. Although $\beta$ may be half-integral, the
real compactly supported class $\delta_\Sigma\beta$ lies in the image of
$H_c^3(X;\mathbb Z)$. Indeed,
\[
    \gamma\bmod2
    =
    w_2(\langle\Psi_\infty\rangle^\perp)
    =
    w_2(\mathfrak g_P|_\Sigma),
\]
and the last class extends over $X$. Its connecting image with
$\mathbb Z_2$ coefficients therefore vanishes. By naturality of the
coefficient reduction, the integral class $\delta_\Sigma\gamma$ reduces
to zero modulo two and hence is divisible by two in
$H_c^3(X;\mathbb Z)$. Thus the real class
$\delta_\Sigma\beta=\frac12\delta_\Sigma\gamma$ has an integral lift. No
choice of such a lift is needed below; Li's local linking sphere argument
independently gives integer multiplicity of the limiting current.
\end{remark}

\subsection{Li's large radius curvature estimates}
\label{subsec: Li_curvature_estimate}

The following estimate was proved by Li
\cite[Lemma~2.3]{li2025large} in the AC $\mathrm G_2$ setting, and the same
argument applies to the Calabi--Yau case as discussed there. For convenience,
we reformulate the proof in the unified $\Theta$-monopole framework, making
explicit the geometric inputs needed in both special holonomy cases.

We continue the geometric framework of
\S\ref{subsec: asymptotic_analysis_theta_monopoles}, where $(X^n,g,\Theta)$ is a one-ended asymptotically conical special holonomy manifold and $P\to X$ is a principal $G$-bundle with $G\in\{\mathrm{SU}(2),\mathrm{SO}(3)\}$. Let
$(\nabla,\Phi)\in\mathscr A(P)\times\Gamma(\mathfrak g_P)$ be an
irreducible finite intermediate energy $\Theta$-monopole with quadratically decaying curvature \eqref{eq:quadratic_decay_theta}, as in the asymptotic framework of \S\ref{subsec: asymptotic_analysis_theta_monopoles}.

By the asymptotic theory recalled there, the configuration admits an
asymptotic abelian model determined by a monopole class $\beta\in H^2(\Sigma;\mathbb R)$, where $\beta$ is integral for $G=\mathrm{SU}(2)$ and $2\beta$ is integral for $G=\mathrm{SO}(3)$. In particular, we have
\begin{equation}\label{ineq:abelian_growth_general}
\limsup_{r\to\infty}
r^{4-n}
\int_{\rho\leqslant r}|F_\nabla|^2
\leqslant C_\beta<\infty.
\end{equation}
The first pointwise input that is needed is the algebraic identity
\eqref{eq:general_calibration_identity}, which gives
\[
|F_{\nabla}|^2\operatorname{vol}
=
-\langle F_{\nabla}\wedge F_{\nabla}\rangle\wedge\Xi
+
|F_{\nabla}\wedge\Theta|^2\operatorname{vol}.
\]
As recalled in \S\ref{subsec: monopoles}, this identity holds for every $2$-form in the $\mathrm G_2$ case and for primitive $2$-forms in the Calabi--Yau case. The latter condition is satisfied here by \eqref{eq:CY_monopole_primitivity}.

We also use the induced tangential inequality on the asymptotic cone.
Let $\Xi_C$ denote the homogeneous model of $\Xi$ on the cone
$(C(\Sigma),g_C)$, and let $\partial_\rho$ be the unit radial vector field.
The contraction $\partial_\rho\!\lrcorner\Xi_C$ restricts to the natural form on the link associated with the cone
structure: in the $\mathrm G_2$ case, the fundamental $2$-form of the nearly
K\"ahler link, and in the Calabi--Yau case, the contact $1$-form of the
Sasaki--Einstein link, with the corresponding homogeneous scaling. The
calibration property of these cone forms implies
the pointwise inequality
\begin{equation}\label{ineq:general_boundary_calibration}
-\langle F_{\nabla}\wedge F_{\nabla}\rangle
\wedge
(\partial_\rho\!\lrcorner\Xi_C)
\leqslant
|F_{\nabla}|_{g_C}^2\,dA_{g_C}
\end{equation}
for the tangential restriction of $F_{\nabla}$ to the level sets of $\rho$. This is the boundary algebraic estimate used below on the hypersurfaces $\{\rho=r\}$.

Finally, we use the following elementary consequence of the AC
asymptotics of $\Xi$. Recall that $\nu_0<0$ denotes the AC decay rate of
$(X,g,\Theta)$.

\begin{lemma}[AC primitive for $\Xi$]\label{lem:primitive_AC}
There exist $R_0>0$ and a form
\[
\tau\in\Omega^{n-5}(X)
\]
such that
\[
\mathrm{supp}(\Xi-d\tau)\subset\{\rho\leqslant R_0\},
\]
and, along the AC end,
\[
\tau
=
\frac{1}{n-4}
\rho\,\partial_\rho\!\lrcorner\Xi_C
+
O(\rho^{\nu_0+1}).
\]
\end{lemma}

\begin{proof}
On the cone, $\Xi_C$ is closed and homogeneous of degree $n-4$. Hence
Cartan's formula gives
\[
d(\rho\,\partial_\rho\!\lrcorner\Xi_C)
=
(n-4)\Xi_C.
\]
Since $\Xi$ is closed and asymptotic to $\Xi_C$ with rate $\nu_0<0$, the
standard radial homotopy construction on the end gives a primitive
$\tau_{\mathrm{end}}$ for $\Xi$ satisfying
\[
d\tau_{\mathrm{end}}=\Xi
\]
for $\rho$ sufficiently large, and
\[
\tau_{\mathrm{end}}
=
\frac{1}{n-4}
\rho\,\partial_\rho\!\lrcorner\Xi_C
+
O(\rho^{\nu_0+1}).
\]
Extending $\tau_{\mathrm{end}}$ to $X$ by a cutoff gives the desired form
$\tau$.
\end{proof}

We can now state the large radius estimate.

\begin{lemma}[Li's large radius curvature estimate]\label{lemm: Yang-Li}
There exist constants $R_0>0$ and $C>0$, depending only on the AC geometry
of $(X,g,\Theta)$, the monopole class $\beta$, and the topology of the principal $G$-bundle $P$, such that
\begin{equation}\label{ineq: energy_estimate}
    \int_{\rho\leqslant r}|F_\nabla|^2
    \leqslant
    \mathcal E^\Theta(\nabla,\Phi)
    +
    Cr^{n-4},
    \qquad
    r\geqslant R_0.
\end{equation}
\end{lemma}
\begin{proof}
The argument we shall present is a streamlined version of Li's large radius estimate
\cite[Lemma~2.3]{li2025large}. Set
\[
    I(r):=\int_{\rho\leqslant r}|F_\nabla|^2,
    \qquad
    \mathcal E^\Theta(r)
    := \mathcal E_{\{\rho\leqslant r\}}^\Theta(\nabla,\Phi) =
    \int_{\rho\leqslant r}|F_\nabla\wedge\Theta|^2,
\] where the last equality follows from the $\Theta$-monopole equation.

By Lemma~\ref{lem:primitive_AC}, after increasing $R_0$ if necessary, there
exists a form $\tau\in\Omega^{n-5}(X)$ such that
\[
    \operatorname{supp}(\Xi-d\tau)\subset\{\rho\leqslant R_0\},
\]
and
\[
    \tau
    =
    \frac{1}{n-4}\rho\,\partial_\rho\!\lrcorner\,\Xi_C
    +
    O(\rho^{\nu_0+1})
\]
along the end. Thus, for every $r\geqslant R_0$, the quantity
\[
    \widetilde C
    :=
    -\int_{ \rho\leqslant r }
    \langle F_\nabla\wedge F_\nabla\rangle
    \wedge(\Xi-d\tau)
\]
is independent of $r$. The closed Chern--Weil form
$\langle F_\nabla\wedge F_\nabla\rangle$ represents a fixed real
characteristic class of degree four of $P$: up to the standard normalization,
this is the second Chern class in the $\mathrm{SU}(2)$ case and the first
Pontryagin class in the $\mathrm{SO}(3)$ case. Hence $\widetilde C$ depends
only on the bundle topology and on the fixed compactly supported form
$\Xi-d\tau$. Using the algebraic identity
\eqref{eq:general_calibration_identity} and Stokes' theorem, we obtain
\[
    \widetilde C
    =
    I(r)-\mathcal E^\Theta(r)
    +
    \int_{\rho=r}
    \langle F_\nabla\wedge F_\nabla\rangle\wedge\tau.
\]
Using the asymptotic expansion of $\tau$ together with the boundary
calibration inequality
\eqref{ineq:general_boundary_calibration}, we get
\[
\begin{aligned}
I(r)-\mathcal E^\Theta(r)-\widetilde C
&=
-\int_{\rho=r}
\langle F_\nabla\wedge F_\nabla\rangle\wedge\tau
\\
&\leqslant
\frac{r}{n-4}
\int_{\rho=r}|F_\nabla|^2\,dA_{g_C}
+
Cr^{\nu_0+1}
\int_{\rho=r}|F_\nabla|^2\,dA_{g_C}.
\end{aligned}
\]
Hence
\begin{equation}\label{ineq:large_radius_before_coarea_theta}
I(r)-\mathcal E^\Theta(r)-\widetilde C
\leqslant
\frac{r}{n-4}(1+Cr^{\nu_0})
\int_{\rho=r}|F_\nabla|^2\,dA_{g_C}.
\end{equation}
On the other hand, the coarea formula and the AC asymptotics imply
\[
I'(r)
=
\int_{\rho=r}|F_\nabla|^2|\nabla\rho|^{-1}dA_g
\geqslant
(1-Cr^{\nu_0})
\int_{\rho=r}|F_\nabla|^2dA_{g_C}.
\]
After increasing $R_0$ if necessary, we may assume
$Cr^{\nu_0}<\frac12$ for all $r\geqslant R_0$. Combining this with
\eqref{ineq:large_radius_before_coarea_theta} gives
\[
I(r)-\mathcal E^\Theta(r)-\widetilde C
\leqslant
\frac{r}{n-4}(1+Cr^{\nu_0})I'(r).
\]
We shall now consider
\[
y(r)
:=
I(r)-\mathcal E^\Theta(\nabla,\Phi)-\widetilde C.
\]
Since $\mathcal E^\Theta(r)\leqslant\mathcal E^\Theta(\nabla,\Phi)$, we have
\[
y(r)
\leqslant
\frac{r}{n-4}(1+Cr^{\nu_0})y'(r),
\]
which implies
\begin{equation}\label{ineq:log_y_theta}
\frac{d}{dr}\log y(r)
\geqslant
\frac{n-4}{r}
-
Cr^{\nu_0-1},
\qquad
\text{whenever }y(r)>0.
\end{equation}
Define
\[
A
:=
\sup_{r\geqslant R_0}
\frac{y(r)_+}{r^{n-4}}
\in[0,\infty].
\]
If $A=0$, then $y(r)\leqslant 0$ for all $r\geqslant R_0$, and the estimate follows
immediately.

Assume now that $A>0$. Then there exists a sequence $r_k\geqslant R_0$ such that
\[
\frac{y(r_k)}{r_k^{n-4}}
=
\frac{y(r_k)_+}{r_k^{n-4}}
\to A>0.
\]
In particular, $y(r_k)>0$ for $k\gg1$.  Let $[r_k,R_k)$ denote the maximal
interval on which $y>0$. Integrating \eqref{ineq:log_y_theta} over
$[r_k,r]\subset[r_k,R_k)$ yields
\[
\log\frac{y(r)}{y(r_k)}
\geqslant
(n-4)\log\frac{r}{r_k}
-
C\int_{r_k}^r t^{\nu_0-1}\,dt.
\]
Since $\nu_0<0$, setting $C':=-C/\nu_0>0$, we obtain
\begin{equation}\label{ineq:almost_monotonicity_y_theta}
\frac{y(r)}{r^{n-4}}
\geqslant
\frac{y(r_k)}{r_k^{n-4}}
\exp\!\bigl(C'(r^{\nu_0}-r_k^{\nu_0})\bigr).
\end{equation}
If $R_k<\infty$, continuity and maximality would imply
$y(R_k)=0$, contradicting the positivity given by
\eqref{ineq:almost_monotonicity_y_theta} as $r\uparrow R_k$.
Hence $R_k=\infty$, and
\eqref{ineq:almost_monotonicity_y_theta} holds for all $r\geqslant r_k$. Hence, since $r\geqslant r_k\geqslant R_0$ and $\nu_0<0$, we have
\[
\exp\!\bigl(C'(r^{\nu_0}-r_k^{\nu_0})\bigr)
\geqslant
\exp(-C'R_0^{\nu_0}),
\]
and thus
\[
\limsup_{r\to\infty}
\frac{y(r)}{r^{n-4}}
\geqslant
A\exp(-C'R_0^{\nu_0}).
\]
On the other hand, by
\eqref{ineq:abelian_growth_general},
\[
\limsup_{r\to\infty}
\frac{y(r)}{r^{n-4}}
=
\limsup_{r\to\infty}
r^{4-n}I(r)
\leqslant
C_\beta<\infty.
\]
Thus
\[
0<A\leqslant C_\beta\exp(C'R_0^{\nu_0})<\infty.
\]
It follows that
\[
y(r)\leqslant Ar^{n-4}\leqslant Cr^{n-4}
\]
for every $r\geqslant R_0$. By definition of $y$, this gives
\[
I(r)
\leqslant
\mathcal E^\Theta(\nabla,\Phi)
+
\widetilde C
+
Cr^{n-4}.
\]
Absorbing $\widetilde C$ into the constant concludes the proof.
\end{proof}

\begin{corollary}[A priori local $L^2$-bounds]
\label{cor: local_L2_estimates}
Under the assumptions of Lemma~\ref{lemm: Yang-Li}, possibly after
increasing $R_0$, the following estimates hold.
\begin{enumerate}
\item[(i)] If $\rho(x)\geqslant R_0$, then for every $r\in(0,r_0]$,
\begin{equation}\label{ineq: m_local_L2_bound_far_center}
    r^{4-n}
    \int_{B_r(x)}
    \bigl(|F_\nabla|^2+|\nabla\Phi|^2\bigr)
    \lesssim
    \rho(x)^{4-n}m+1.
\end{equation}

\item[(ii)] If $\rho(x)<R_0$, then for every $r\in(0,r_0]$,
\begin{equation}\label{ineq: m_local_L2_bound_close_center}
    \int_{B_r(x)}
    \bigl(|F_\nabla|^2+|\nabla\Phi|^2\bigr)
    \lesssim
    m+1.
\end{equation}
\end{enumerate}
In particular, if $m\gg_{X,\Theta,P,\beta}1$, then for all $x\in X$ and
all $r\in(0,r_0]$,
\begin{equation}\label{ineq: m_local_L2_bound}
    \int_{B_r(x)}
    \bigl(|F_\nabla|^2+|\nabla\Phi|^2\bigr)
    \lesssim
    m.
\end{equation}
\end{corollary}

\begin{proof}
By the topological energy formula \eqref{eq: theta_monopole_energy_formula}, the fixed monopole class assumption gives
\[
    \mathcal E^\Theta(\nabla,\Phi) = \int_X|\nabla\Phi|^2\lesssim m.
\]
Combining this with Lemma~\ref{lemm: Yang-Li}, we obtain, for all
$R\geqslant R_0$,
\begin{equation}\label{ineq:global_L2_growth_auxiliary}
    \int_{\{\rho\leqslant R\}}
    \bigl(|F_\nabla|^2+|\nabla\Phi|^2\bigr)
    \lesssim
    m+R^{n-4}.
\end{equation}
We first prove (i). By asymptotic conicality, after increasing $R_0$ if
necessary, there is a uniform constant $\alpha\in(0,1)$ such that, for every
$x$ with $\rho(x)\geqslant R_0$, the radius $R:=\alpha\rho(x)$ satisfies $R<\inj_g(x)$ and $B_R(x)\subset \{\rho\leqslant C\rho(x)\}$ for a fixed geometric constant $C$. Applying
\eqref{ineq:global_L2_growth_auxiliary} with $C\rho(x)$ in place of $R$ gives
\[
\begin{aligned}
    R^{4-n}
    \int_{B_R(x)}
    \bigl(|F_\nabla|^2+|\nabla\Phi|^2\bigr)
    &\lesssim
    \rho(x)^{4-n}
    \int_{\{\rho\leqslant C\rho(x)\}}
    \bigl(|F_\nabla|^2+|\nabla\Phi|^2\bigr)  \\
    &\lesssim
    \rho(x)^{4-n}m+1.
\end{aligned}
\]
Choosing $R_0$ still larger if necessary, we may assume $r_0\leqslant R$ for all
$\rho(x)\geqslant R_0$. The standard codimension-four almost monotonicity estimate for solutions of
\eqref{eq:YMH_second_order} (see
Corollary~\ref{cor:codim4_almost_monotonicity}), applied on the AC ball
$B_R(x)$, then gives
\[
    r^{4-n}
    \int_{B_r(x)}
    \bigl(|F_\nabla|^2+|\nabla\Phi|^2\bigr)
    \lesssim
    R^{4-n}
    \int_{B_R(x)}
    \bigl(|F_\nabla|^2+|\nabla\Phi|^2\bigr)
\]
for every $r\in(0,r_0]$. This proves
\eqref{ineq: m_local_L2_bound_far_center}.

We now prove (ii). If $\rho(x)<R_0$ and $r\in(0,r_0]$, then, after
increasing $R_0$ once more if necessary,
\[
    B_r(x)\subset \{\rho\leqslant 2R_0\}.
\]
Hence \eqref{ineq:global_L2_growth_auxiliary}, applied with $R=2R_0$, gives
\[
    \int_{B_r(x)}
    \bigl(|F_\nabla|^2+|\nabla\Phi|^2\bigr)
    \lesssim
    m+1,
\]
which is \eqref{ineq: m_local_L2_bound_close_center}.

Finally, if $m$ is sufficiently large, the additive constant in (ii) is
absorbed into $m$. In case (i), since $r\leqslant r_0$ and
$\rho(x)^{4-n}\lesssim1$ on $\{\rho\geqslant R_0\}$, the estimate
\eqref{ineq: m_local_L2_bound_far_center} also implies
\[
    \int_{B_r(x)}
    \bigl(|F_\nabla|^2+|\nabla\Phi|^2\bigr)
    \lesssim
    m.
\]
This proves \eqref{ineq: m_local_L2_bound}.
\end{proof}


\section{Codimension-three concentration, Higgs zeros, and obstruction to monotonicity}\label{sec: concentration}

We first apply the compactness theorem of Parise--Pigati--Stern in the
AC setting and identify its saturated limit with Li's calibrated current.
We then study the support of this current, its relation to the limiting
Higgs zero set, and the obstruction to effective codimension-three
monotonicity.

\subsection{PPS compactness and calibrated codimension-three concentration}\label{subsec: PPS_calibrated}

We first recall the compactness and calibration statements in the
published version of Parise--Pigati--Stern, translated into the
normalization of the present paper.

\begin{theorem}[PPS compactness and calibration inequality
{\cite[Theorem~1.2 and Corollary~1.8]{parise2025nonabelian}}]
\label{thm:PPS}
Let $M$ be a compact Riemannian $n$-manifold, possibly with boundary, and
let $P=M\times\mathrm{SU}(2)$ be the trivial principal bundle. For
$(\nabla,\Phi)\in\mathscr A(P)\times\Gamma(\mathfrak g_P)$ define the
\emph{charge form}
\begin{equation}
\label{eq:charge_form}
    Z(\nabla,\Phi):=2\langle\nabla\Phi\wedge F_\nabla\rangle
    \in\Omega^3(M),
\end{equation}
and identify it with the dual $(n-3)$-current. Suppose that
$m_i\to+\infty$ and
$(\nabla_i,\Phi_i)\in\mathscr A(P)\times\Gamma(\mathfrak g_P)$ satisfy
\[
    \liminf_{i\to\infty}
    m_i^{-1}\int_M
    \bigl(|F_{\nabla_i}|^2+|\nabla_i\Phi_i|^2\bigr)
    <+\infty
\]
and
\[
    \lim_{i\to\infty}
    m_i^{-1}\int_M(m_i-|\Phi_i|)^2
    =0.
\]
Then, after passing to a subsequence, there are an $(n-3)$-current $T$,
whose restriction to $M^\circ$ is an integral cycle, and a Radon measure
$\mu$ on $M$ such that
\[
    m_i^{-1}Z(\nabla_i,\Phi_i)
    \rightharpoonup^*8\pi T,
\]
\[
    \mu_i
    :=
    m_i^{-1}
    \bigl(|F_{\nabla_i}|^2+|\nabla_i\Phi_i|^2\bigr)\vol
    \rightharpoonup^*\mu,
\]
and
\[
    8\pi\|T\|\leqslant\mu.
\]
If, in addition, $M$ carries a calibration
$\Theta\in\Omega^{n-3}(M)$ and every $(\nabla_i,\Phi_i)$ is a
$\Theta$-monopole, then, after passing to a further subsequence,
\[
    \nu_i
    :=
    m_i^{-1}
    \bigl(|F_{\nabla_i}\wedge\Theta|^2
          +|\nabla_i\Phi_i|^2\bigr)\vol
    =
    2m_i^{-1}|\nabla_i\Phi_i|^2\vol
    \rightharpoonup^*\nu
\]
for a Radon measure $\nu$ satisfying
\[
    \nu\leqslant8\pi\|T\|.
\]
On $M^\circ$, equality holds if and only if $T$ is calibrated by
$\Theta$. In particular, $\nu|_{M^\circ}$ is
$\mathcal H^{n-3}$-rectifiable; if $\mu=\nu$, then on $M^\circ$
\[
    \mu=\nu=8\pi\|T\|,
\]
and $T$ is calibrated.
\end{theorem}

\begin{remark}[Normalization and the published PPS version]
\label{rmk:PPS_published_normalization}
Theorem~\ref{thm:PPS} is the equivalent large mass reformulation of
the parts of \cite[Theorem~1.2 and Corollary~1.8]{parise2025nonabelian}
needed below.
The constants above
have been converted using
\eqref{eq:PPS_Li_inner_product_comparison}--%
\eqref{eq:PPS_Li_mass_comparison}.
\end{remark}

\begin{proposition}[PPS compactness for fixed bundles and $\mathrm{SO}(3)$]
\label{prop:PPS_SO3_variant}
The statement and conclusions of Theorem~\ref{thm:PPS} remain valid for
a fixed principal $G$-bundle $P\to M$, where
$G\in\{\mathrm{SU}(2),\mathrm{SO}(3)\}$. In the $\mathrm{SO}(3)$ case
we use the inner product \eqref{eq:SO3_inner_prod_convention}. The charge
form is still \eqref{eq:charge_form}, and all constants remain unchanged.
\end{proposition}

\begin{proof}
Choose a finite cover of $M$ by contractible coordinate balls and, near
$\partial M$, contractible half-balls. If $G=\mathrm{SU}(2)$, the bundle
is trivial on each member of the cover. If $G=\mathrm{SO}(3)$, its
restriction lifts there to a principal $\mathrm{SU}(2)$-bundle. Through
the isometric Lie algebra identification
$d\pi:\mathfrak{su}(2)\to\mathfrak{so}(3)$, the connection, Higgs field,
energy density, and charge form lift without changing pointwise norms or
constants.

Apply Theorem~\ref{thm:PPS} on a finite refinement and pass to a common
subsequence. The local current limits agree on overlaps because they are
distributional limits of the globally defined charge forms; the measure
limits agree because they are weak limits of the globally defined energy
measures. Hence they patch globally. Integer rectifiability, the cycle
condition in the interior, the calibration inequality, and its equality
case are local properties, so they pass to the patched limits.
\end{proof}

For convenience, we now state a globalization of the
Parise--Pigati--Stern compactness theorem to the asymptotically conical
setting.

\begin{theorem}[AC version of the Parise--Pigati--Stern compactness theorem]
\label{thm:PPS_AC_version}
Let $(X^n,g)$ be a complete asymptotically conical Riemannian manifold
with radius function $\rho$, and let $P\to X$ be a principal $G$-bundle,
where $G\in\{\mathrm{SU}(2),\mathrm{SO}(3)\}$. Let
$(\nabla_i,\Phi_i)\in\mathscr A(P)\times\Gamma(\mathfrak g_P)$ be a
sequence of configurations with finite masses $m_i\to+\infty$, satisfying
the following hypotheses.

\medskip
\noindent
\textit{(H1) (Uniform local mass-renormalized energy bounds.)}
For every $R>0$,
\[
    \sup_i
    m_i^{-1}
    \int_{B_R}
    \bigl(
        |F_{\nabla_i}|^2
        +
        |\nabla_i\Phi_i|^2
    \bigr)\vol
    <\infty.
\]

\medskip
\noindent
\textit{(H2) (Local largeness condition on the Higgs fields.)}
For every $R>0$,
\[
    \lim_{i\to\infty}
    m_i^{-1}
    \int_{B_R}
    (m_i-|\Phi_i|)^2\vol
    =
    0.
\]
Then, after passing to a subsequence, there exist a closed locally integral
$(n-3)$-current
\[
    T\in\mathscr I_{n-3,\mathrm{loc}}(X),
    \qquad
    \partial T=0,
\]
and a Radon measure $\mu$ on $X$ such that
\begin{align*}
    m_i^{-1}Z(\nabla_i,\Phi_i)
    \rightharpoonup^*
    8\pi T
    \quad&\text{in }\mathcal D_{n-3}(X),\\
    \mu_i :=
    m_i^{-1}
    \bigl(
        |F_{\nabla_i}|^2
        +
        |\nabla_i\Phi_i|^2
    \bigr)\vol \rightharpoonup^*\mu
    \quad&\text{as Radon measures on }X,
\end{align*}
and the $\mathcal H^{n-3}$-rectifiable Radon measure $\|T\|$ satisfies
\[
    8\pi\|T\|\leqslant\mu.
\]
Moreover, suppose that $(X,g)$ carries a calibration
$\Theta\in\Omega^{n-3}(X)$ and that each $(\nabla_i,\Phi_i)$ is a
$\Theta$-monopole. Then, after passing to a further subsequence if necessary,
\[
    \nu_i
    :=
    m_i^{-1}
    \bigl(
        |F_{\nabla_i}\wedge\Theta|^2
        +
        |\nabla_i\Phi_i|^2
    \bigr)\vol\rightharpoonup^*\nu ,
\] 
for some $\mathcal H^{n-3}$-rectifiable Radon measure $\nu$ satisfying
\[
    \nu \leqslant 8\pi\|T\|.
\]
Furthermore, $\nu=8\pi\|T\|$ if and only if $T$ is $\Theta$-calibrated.
\end{theorem}

\begin{proof}
Choose an exhaustion of $X$ by compact domains with smooth boundary,
\[
    \overline B_{R_1}\subset \overline B_{R_2}\subset\cdots\subset X,
    \qquad
    R_N\to+\infty,
\]
where the $R_N$ are regular values of the radius function and
$\bigcup_{N\geqslant 1} B_{R_N}=X$.

Fix $N\in\mathbb{N}$. The restrictions of \textup{(H1)} and \textup{(H2)} to
$\overline B_{R_N}$ are precisely the hypotheses needed to apply
Proposition~\ref{prop:PPS_SO3_variant} on the compact manifold with
boundary $\overline B_{R_N}$. Hence, after passing to a subsequence of $(\nabla_i,\Phi_i)$ depending on
$N$, we obtain a rectifiable $(n-3)$-current
\[
    T_{R_N}\in\mathcal D_{n-3}(\overline B_{R_N})
\]
and a Radon measure $\mu_{R_N}$ on $\overline B_{R_N}$ such that
\[
    m_i^{-1}Z(\nabla_i,\Phi_i)\lfloor {\overline B_{R_N}}
    \rightharpoonup^*
    8\pi T_{R_N},
\]
\[
    \mu_i\lfloor {\overline B_{R_N}}
    \rightharpoonup^*
    \mu_{R_N},
\]
and
\[
    8\pi\|T_{R_N}\|
    \leqslant
    \mu_{R_N}.
\]
Moreover, Proposition~\ref{prop:PPS_SO3_variant} gives that $T_{R_N}$ restricts to an integral cycle in the interior $B_{R_N}$. By a diagonal argument, we may choose a single subsequence of $(\nabla_i,\Phi_i)$ for which these convergences hold simultaneously for every $N$. By uniqueness of distributional limits and Radon measure limits, the resulting currents and measures are compatible under restriction. More precisely, if $M\geqslant N$, then
\[
    T_{R_M}\lfloor {B_{R_N}}
    =
    T_{R_N}\lfloor {B_{R_N}},
    \qquad
    \mu_{R_M}\lfloor {B_{R_N}}
    =
    \mu_{R_N}\lfloor {B_{R_N}}.
\]
Thus the $T_{R_N}$ glue to a current $T\in\mathcal D_{n-3}(X)$, and the measures $\mu_{R_N}$ glue to a Radon measure $\mu$ on $X$. The glued objects satisfy
\[
    m_i^{-1}Z(\nabla_i,\Phi_i)
    \rightharpoonup^*
    8\pi T
    \quad\text{in }\mathcal D_{n-3}(X),
\]
\[
    \mu_i\rightharpoonup^*\mu
    \quad\text{as Radon measures on }X,
\]
and
\[
    8\pi\|T\|\leqslant\mu.
\]
We now verify the remaining properties of the current $T$. Since every compact subset of $X$ is contained in some $B_{R_N}$, and $T$ agrees there with $T_{R_N}$, the current $T$ is locally integer rectifiable, $T\in\mathscr R_{n-3,\mathrm{loc}}(X)$; in particular, $\|T\|$ is an $\mathcal H^{n-3}$-rectifiable Radon measure. Moreover, $\partial T=0$. Indeed, if $\alpha\in\Omega_c^{n-4}(X)$, choose $N$ so large that $\operatorname{spt}\alpha\Subset B_{R_N}$. Since $T_{R_N}$ is a cycle in the interior $B_{R_N}$, we have
\[
    \partial T(\alpha)
    =
    T(d\alpha)
    =
    T_{R_N}(d\alpha)
    =
    0,
\] hence $\partial T=0$ as claimed. Thus, $T\in\mathscr I_{n-3,\mathrm{loc}}(X)$ is a closed locally integral $(n-3)$-current. 

Now assume that $(X,g)$ carries a calibration $\Theta$ and that each
$(\nabla_i,\Phi_i)$ is a $\Theta$-monopole. Since $\Theta$ has comass one,
the measures $\nu_i$ satisfy local bounds controlled by the local
mass-renormalized Yang--Mills--Higgs energy bounds in (H1). Therefore, after
passing to a further subsequence if necessary,
\[
    \nu_i\rightharpoonup^*\nu
\]
as Radon measures on $X$.

Applying the calibrated part of Proposition~\ref{prop:PPS_SO3_variant} on each compact
domain $\overline B_{R_N}$ gives local limiting measures $\nu_{R_N}$ such
that
\[
    \nu_i\lfloor {\overline B_{R_N}}
    \rightharpoonup^*
    \nu_{R_N},
    \qquad
    \nu_{R_N}\leqslant 8\pi\|T_{R_N}\|,
\]
and each $\nu_{R_N}$ is $\mathcal H^{n-3}$-rectifiable. Again, by uniqueness of Radon measure limits, the measures $\nu_{R_N}$ are compatible under restriction and glue to the global limit $\nu$ which is $\mathcal H^{n-3}$-rectifiable and satisfies
\[
    \nu\leqslant 8\pi\|T\|.
\]
Finally, the characterization of when equality holds is a local property. On each interior $B_{R_N}$, Proposition~\ref{prop:PPS_SO3_variant} gives that $\nu_{R_N}=8\pi\|T_{R_N}\|$ holds if and only if $T_{R_N}$ is calibrated by $\Theta$. Since the restrictions
of the currents are compatible and calibration is a local condition on the orientation of the approximate tangent planes, this is equivalent to the statement that $\nu=8\pi\|T\|$ on $X$ if and only if the glued current $T$ is $\Theta$-calibrated.
\end{proof}

\begin{remark}[On the role of hypotheses (H1) and (H2)]
\label{rmk:H1_H2_general_configurations}
Theorem~\ref{thm:PPS_AC_version} is a compactness statement for arbitrary
sequences of configurations with finite masses satisfying (H1) and (H2). In
particular, no Yang--Mills--Higgs equation is assumed in the statement of the
compactness theorem. The two hypotheses should therefore be regarded as
independent assumptions at this level.

This point is important for (H2). For a general sequence of configurations,
the local largeness condition
\[
\lim_{i\to\infty}
m_i^{-1}
\int_{B_R}(m_i-|\Phi_i|)^2\vol=0
\]
does not follow from the local mass-renormalized energy bound (H1). In the main
setting below, (H2) will be verified using Lemma~\ref{lem:largeness}, but that
lemma uses \eqref{eq:YMH_second_order}, in particular the equation
$\nabla_i^*\nabla_i\Phi_i=0$, together with the linear energy bound. Thus the implication
\[
\text{(H1)}\Longrightarrow\text{(H2)}
\]
is a feature of the critical point/monopole setting considered later, not a
feature of the general Parise--Pigati--Stern compactness theorem.
\end{remark}

We now specialize to the Main Setting of the paper: large mass sequences of finite intermediate energy $\Theta$-monopoles on AC
$\mathrm G_2$-manifolds or AC Calabi--Yau $3$-folds, with structure group $G\in\{\mathrm{SU}(2),\mathrm{SO}(3)\}$, fixed asymptotic monopole class, and the asymptotic framework of
\S\ref{subsec: asymptotic_analysis_theta_monopoles}.


\phantomsection
\label{setting:main}
\begin{mainsetting}

From this point onward, unless otherwise stated,
$(X^n,g,\Theta)$ denotes either an AC $\mathrm G_2$-manifold
$(X^7,g,\psi)$ or an AC Calabi--Yau $3$-fold
$(X^6,g,\operatorname{Re}\Omega)$; see
\S\S\ref{subsec: monopoles}-\ref{subsec: asymptotic_analysis_theta_monopoles}.
We fix $G\in\{\mathrm{SU}(2),\mathrm{SO}(3)\}$ and let $(\nabla_i,\Phi_i)$ be a sequence of irreducible finite
intermediate energy $\Theta$-monopoles on a fixed principal $G$-bundle
$P\to X$, satisfying the asymptotic framework of
\S\ref{subsec: asymptotic_analysis_theta_monopoles}. In particular,
\begin{itemize}
\item $\mathcal E^\Theta(\nabla_i,\Phi_i)<\infty$; and
\item $|F_{\nabla_i}|=O_i(\rho^{-2})$ as $\rho\to\infty$, where the implicit constant is allowed to depend on $i$.
\end{itemize}

Consequently, by Theorem~\ref{thm:theta_asymptotics}, the masses
\[
    m_i:=\lim_{\rho\to\infty}|\Phi_i|
\]
exist, and each configuration admits an asymptotic abelian model with
monopole class
\[
    \beta_i\in H^2(\Sigma;\mathbb R).
\]
For $G=\mathrm{SU}(2)$ the class $\beta_i$ is integral, while for
$G=\mathrm{SO}(3)$ the class $2\beta_i$ is integral. We assume throughout
that
\begin{itemize}
\item $m_i\longrightarrow+\infty$;

\item the monopole class is fixed,
\[
    \beta_i=\beta\in H^2(\Sigma;\mathbb R),
\]
and
\[
    k:=\langle\beta\cup\Theta_\infty,[\Sigma]\rangle\neq0.
\]
\end{itemize}

By the topological energy formula
\eqref{eq: theta_monopole_energy_formula},
\begin{equation}\label{eq: linear_interm_energy_growth}
    \mathcal E^\Theta(\nabla_i,\Phi_i)
    =
    \|\nabla_i\Phi_i\|_{L^2}^2
    =
    4\pi k m_i.
\end{equation}
Thus, the intermediate energy grows linearly with the mass. Since the
monopoles are irreducible and $m_i>0$, the identity
\eqref{eq: linear_interm_energy_growth} implies $k>0$.
\end{mainsetting}

We now verify that the {\mainsettingref} lies in the scope of
Theorem~\ref{thm:PPS_AC_version}.

\begin{proposition}[Verification of the PPS hypotheses in the Main Setting]
\label{prop:PPS_hypotheses_main_setting}
In the {\mainsettingref}, the sequence $(\nabla_i,\Phi_i)$ satisfies hypotheses
\textup{(H1)} and \textup{(H2)} of
Theorem~\ref{thm:PPS_AC_version}. Consequently, after passing to a
subsequence, there exist a closed locally integral $(n-3)$-current
\[
    T\in\mathscr I_{n-3,\mathrm{loc}}(X),
    \qquad
    \partial T=0,
\]
a Radon measure $\mu$ on $X$, and a finite
$\mathcal H^{n-3}$-rectifiable Radon measure $\nu$ on $X$ such that
\[
    m_i^{-1}Z(\nabla_i,\Phi_i)
    \rightharpoonup^*
    8\pi T
    \quad\text{in }\mathcal D_{n-3}(X),
\]
and
\[
    \mu_i:=
    m_i^{-1}
    \bigl(
        |F_{\nabla_i}|^2+|\nabla_i\Phi_i|^2
    \bigr)\,\vol
    \rightharpoonup^*
    \mu,
\]
\[
    \nu_i:=
    m_i^{-1}
    \bigl(
        |F_{\nabla_i}\wedge\Theta|^2+|\nabla_i\Phi_i|^2
    \bigr)\,\vol
    \rightharpoonup^*
    \nu ,
\]
weakly as Radon measures on $X$. These are related by
\[
    \nu\leqslant 8\pi\|T\|\leqslant\mu.
\]
Moreover, the $\nu_i$ satisfy the total mass identity
\[
    \nu_i(X)=8\pi k,
    \qquad
    \nu(X)\leqslant 8\pi k.
\]
Furthermore, $\nu=8\pi\|T\|$ if and only if $T$ is $\Theta$-calibrated.
\end{proposition}

\begin{proof}
Hypothesis \textup{(H1)}, the uniform local mass-renormalized
Yang--Mills--Higgs energy bound, follows from the large radius curvature
estimate of Lemma~\ref{lemm: Yang-Li}, the local estimate of
Corollary~\ref{cor: local_L2_estimates}, and the linear growth identity
\eqref{eq: linear_interm_energy_growth}. Hypothesis \textup{(H2)}, the
local largeness condition for the Higgs fields, is precisely
Lemma~\ref{lem:largeness}.

Since the configurations in the {\mainsettingref} are $\Theta$-monopoles and
$\Theta$ is a calibration, the calibrated part of
Theorem~\ref{thm:PPS_AC_version} applies. This gives the closed locally
integral current $T$, the Radon measures $\mu$ and $\nu$, the stated
convergences, the inequalities
\[
    \nu\leqslant 8\pi\|T\|\leqslant\mu,
\]
the local $\mathcal H^{n-3}$-rectifiability of $\nu$, and the equality
characterization in terms of $\Theta$-calibration.

It remains to record the global finiteness of $\nu$ in the {\mainsettingref}.
Since $(\nabla_i,\Phi_i)$ is a $\Theta$-monopole,
\[
    |F_{\nabla_i}\wedge\Theta|^2=|\nabla_i\Phi_i|^2,
\]
and therefore
\begin{equation}\label{eq: finite_measures_nu_i}
    \nu_i(X)
    =
    2m_i^{-1}\mathcal E^\Theta(\nabla_i,\Phi_i)
    =
    8\pi k
\end{equation}
by the intermediate energy identity
\eqref{eq: linear_interm_energy_growth}. Hence the measures $\nu_i$ have
uniformly bounded total mass. Since
$\nu_i\rightharpoonup^*\nu$ as Radon measures on $X$, it follows that
\[
    \nu(X)\leqslant 8\pi k.
\]
Thus $\nu$ is a finite $\mathcal H^{n-3}$-rectifiable Radon measure.
\end{proof}

\begin{remark}[Finiteness of $\nu$ and infinite total mass of the measures $\mu_i$]
\label{rmk:finite_nu_infinite_mu_i}
The preceding proposition records the finite total mass of $\nu$ in the {\mainsettingref}. Indeed, by equation \eqref{eq: finite_measures_nu_i} the measures $\nu_i$ have uniformly bounded total mass, and the limit satisfies $\nu(X)\leqslant 8\pi k$. At this stage, equality need not follow purely from convergence against compactly supported test functions, since mass could in principle escape to infinity. The equality
$\nu(X)=8\pi k$ will follow later from the saturation $\mu=\nu=8\pi\|T\|$
in Proposition~\ref{prop:PPS_saturation_and_Li_identification} below.

In contrast, each measure
\[
    \mu_i
    =
    m_i^{-1}
    \bigl(
        |F_{\nabla_i}|^2+|\nabla_i\Phi_i|^2
    \bigr)\,\vol
\]
has infinite total mass whenever $k\neq0$. Indeed, this is equivalent to
\[
    \int_X |F_{\nabla_i}|^2=+\infty.
\]
If instead $F_{\nabla_i}\in L^2(X)$, then the quadratic curvature decay and
standard $\varepsilon$-regularity along the AC end would improve the
curvature decay at infinity. Consequently, the induced asymptotic
connection on the link $\Sigma$ would be flat and $k$ would have to vanish, contrary to the standing assumption. Alternatively, a nonflat limiting connection on the link contributes a nonzero $\rho^{-2}$ leading curvature term, whose square is not integrable in dimensions six and seven.
\end{remark}

\begin{remark}[Li's $\mathrm{SO}(3)$ variant]
\label{rmk:Li_SO3_variant}
Li states the principal results for $\mathrm{SU}(2)$ and records the
$\mathrm{SO}(3)$ modifications in
\cite[Remarks~1.8, 1.14 and~3.1]{li2025large}. With the transported metric
\eqref{eq:SO3_inner_prod_convention}, all local analytic estimates,
measure limits, transverse rescalings, and the singular abelian equations
are unchanged. The global abelian curvature class may be half-integral:
only twice Li's normalized Chern form is necessarily the Euler class of
the $\mathrm{SO}(2)$ reduction. The limiting current nevertheless has
integer multiplicity, because the linking spheres used in Li's integrality
argument lie in balls on which the $\mathrm{SO}(3)$ bundle lifts to an
$\mathrm{SU}(2)$ bundle; see the proof of
\cite[Theorem~4.3]{li2025large}. Accordingly, the conclusions of
\cite[Theorems~1.12 and~1.16 and Corollary~4.6]{li2025large} used below
apply to either structure group in the {\mainsettingref}.
\end{remark}

We now show that, in the present $\Theta$-monopole setting, the PPS
inequalities are saturated. We use Li's compactness theorem for large
mass $\mathrm G_2$ and Calabi--Yau monopoles. More precisely, Li constructs in \cite[Theorem~1.12]{li2025large} a
limiting singular abelian monopole whose Dirac-type singular set is encoded
by a calibrated current $T_{\mathrm{Li}}$, and identifies the homology
class of this current as the Poincar\'e dual of the image of the asymptotic monopole class under the connecting homomorphism 
\[
    \delta_\Sigma:H^2(\Sigma;\mathbb R)\longrightarrow H_c^3(X;\mathbb R)
\] in the long exact sequence for the AC pair $(X,\Sigma)$.
Li also proves in \cite[Theorem~1.16]{li2025large} the equality of the limiting mass-renormalized curvature and Higgs gradient measures.

All quantities below are written in the convention of the present paper.
Their comparison with the norms and masses used by Li is governed by
\eqref{eq:PPS_Li_inner_product_comparison}--%
\eqref{eq:PPS_Li_mass_comparison}.

\begin{proposition}[Saturation and identification with Li's calibrated current]
\label{prop:PPS_saturation_and_Li_identification}
In the {\mainsettingref}, after passing to a further subsequence of the subsequence fixed in
Proposition~\ref{prop:PPS_hypotheses_main_setting}, one has
\[
    \mu=\nu=8\pi\|T\|.
\]
In particular, $T$ has finite mass and is a $\Theta$-calibrated integral
$(n-3)$-cycle.

Moreover, after converting Li's limiting measures according to
\eqref{eq:PPS_Li_inner_product_comparison}--%
\eqref{eq:PPS_Li_mass_comparison}, the PPS cycle $T$ and Li's limiting
calibrated current $T_{\mathrm{Li}}$ agree as calibrated currents:
\[
    T=T_{\mathrm{Li}}.
\]
Consequently, $T$ is compactly supported and represents the topological charge determined by the asymptotic monopole class:
\[
    [T]
    =
    [T_{\mathrm{Li}}]
    =
    \operatorname{PD}_X(\delta_\Sigma\beta)
    \in H_{n-3}(X;\mathbb R).
\]
In particular,
\[
    \|T\|(X)=k,
    \qquad
    \mu(X)=\nu(X)=8\pi k,
\]
and
\[
    \mathcal S
    :=
    \operatorname{spt}\mu
    =
    \operatorname{spt}\nu
    =
    \operatorname{spt}\|T\|
    =
    \operatorname{spt}\|T_{\mathrm{Li}}\|,
\]
so the mass-renormalized concentration set is precisely the compact support of
Li's limiting calibrated current.
\end{proposition}

\begin{proof}
Since $(\nabla_i,\Phi_i)$ is a $\Theta$-monopole,
$|F_{\nabla_i}\wedge\Theta|^2=|\nabla_i\Phi_i|^2$. Therefore
\[
    \nu_i
    =
    2m_i^{-1}|\nabla_i\Phi_i|^2\,\vol,
\]
whereas
\[
    \mu_i
    =
    m_i^{-1}
    \bigl(
        |F_{\nabla_i}|^2+|\nabla_i\Phi_i|^2
    \bigr)\,\vol.
\]
Hence
\begin{equation}\label{eq:mu_split_general}
    \mu_i
    =
    \nu_i
    +
    m_i^{-1}
    \bigl(
        |F_{\nabla_i}|^2-|\nabla_i\Phi_i|^2
    \bigr)\,\vol.
\end{equation}
By Li's equality of the limiting mass-renormalized curvature and Higgs gradient
measures \cite[Theorem~1.16]{li2025large}, together with the
$\mathrm{SO}(3)$ variant recorded in
Remark~\ref{rmk:Li_SO3_variant}, the second term on the right of
\eqref{eq:mu_split_general} converges weakly to zero as a signed Radon
measure on $X$.
Thus
\[
    \mu=\nu.
\]
Combining this with
Proposition~\ref{prop:PPS_hypotheses_main_setting}, which gives
\[
    \nu\leqslant8\pi\|T\|\leqslant\mu,
\]
we obtain
\[
    \mu=\nu=8\pi\|T\|.
\]
The equality $\nu=8\pi\|T\|$ is the equality case in the calibrated part of
Theorem~\ref{thm:PPS_AC_version}. Hence $T$ is $\Theta$-calibrated. Since
$\nu$ is finite in the {\mainsettingref}, by
Proposition~\ref{prop:PPS_hypotheses_main_setting}, the equality
$\nu=8\pi\|T\|$ implies that $\|T\|(X)<\infty$. Together with
$\partial T=0$ and $T\in\mathscr I_{n-3,\mathrm{loc}}(X)$, this shows that
$T$ is an integral $(n-3)$-cycle.

It remains to compare $T$ with Li's current. Li's Theorem~1.16 gives, in
Li's normalization,
\[
    \frac{1}{2\pi m_i^{\mathrm{Li}}}
    |F_{\nabla_i}|_{\mathrm{Li}}^2\,\vol
    \rightharpoonup^*
    \|T_{\mathrm{Li}}\|,
    \qquad
    \frac{1}{2\pi m_i^{\mathrm{Li}}}
    |\nabla_i\Phi_i|_{\mathrm{Li}}^2\,\vol
    \rightharpoonup^*
    \|T_{\mathrm{Li}}\|.
\]
Since
\[
    |F_{\nabla_i}|^2=4|F_{\nabla_i}|_{\mathrm{Li}}^2,
    \qquad
    |\nabla_i\Phi_i|^2
    =
    4|\nabla_i\Phi_i|_{\mathrm{Li}}^2,
    \qquad
    m_i=2m_i^{\mathrm{Li}},
\]
we obtain, in the normalization of the present paper,
\[
    m_i^{-1}|F_{\nabla_i}|^2\,\vol
    \rightharpoonup^*
    4\pi\|T_{\mathrm{Li}}\|,
    \qquad
    m_i^{-1}|\nabla_i\Phi_i|^2\,\vol
    \rightharpoonup^*
    4\pi\|T_{\mathrm{Li}}\|.
\]
Consequently,
\[
    \mu_i
    =
    m_i^{-1}
    \bigl(
        |F_{\nabla_i}|^2+|\nabla_i\Phi_i|^2
    \bigr)\,\vol 
    \rightharpoonup^*
    8\pi\|T_{\mathrm{Li}}\|.
\]
On the other hand, the first part of the proof gives
\[
    \mu_i\rightharpoonup^*\mu=8\pi\|T\|.
\]
By uniqueness of the Radon measure limit,
\[
    \|T\|=\|T_{\mathrm{Li}}\|.
\]
Both $T$ and $T_{\mathrm{Li}}$ are calibrated by the same form $\Theta$.
Since their mass measures agree, their underlying rectifiable supports and
multiplicity functions agree almost everywhere with respect to their common
mass measure. Moreover, the calibration fixes the orientation of the
approximate tangent planes. Hence the calibrated currents themselves agree:
\[
    T=T_{\mathrm{Li}}.
\]
We now invoke Li's theorem on the calibrated cycle
\cite[Theorem~4.3]{li2025large}, using
Remark~\ref{rmk:Li_SO3_variant} when $G=\mathrm{SO}(3)$.
This theorem gives that $T_{\mathrm{Li}}$ is a calibrated integral
$(n-3)$-cycle with compact support. Since we have already shown that
$T=T_{\mathrm{Li}}$ as currents, it follows that $T$ is compactly
supported and that
\[
    \operatorname{spt}\|T\|
    =
    \operatorname{spt}\|T_{\mathrm{Li}}\|.
\]
We next invoke the remaining topological and mass conclusions of
\cite[Corollary~4.6, items~(1) and~(3)]{li2025large}, using
Remark~\ref{rmk:Li_SO3_variant} when $G=\mathrm{SO}(3)$. By item~(1), the
homology class of Li's limiting calibrated current is Poincar\'e dual to
the image of the asymptotic monopole class under the connecting
homomorphism $\delta_\Sigma: H^2(\Sigma;\mathbb R)\longrightarrow H_c^3(X;\mathbb R)$. Since $T=T_{\mathrm{Li}}$, it follows that
\[
    [T]
    =
    [T_{\mathrm{Li}}]
    =
    \operatorname{PD}_X(\delta_\Sigma\beta)
    \in
    H_{n-3}(X;\mathbb R).
\]
Moreover, item~(3) of the same corollary gives the total mass of
$T_{\mathrm{Li}}$. In the unified $\Theta$-monopole notation, this reads
\[
    \|T_{\mathrm{Li}}\|(X)
    =
    \int_{T_{\mathrm{Li}}}\Theta
    =
    \left\langle
        \beta\cup\Theta_\infty,
        [\Sigma]
    \right\rangle
    =
    k.
\]
Using again $T=T_{\mathrm{Li}}$, we obtain
\[
    \|T\|(X)=k.
\]
Finally, since $\mu=\nu=8\pi\|T\|$ as Radon measures,
\[
    \mu(X)=\nu(X)=8\pi k.
\]
Moreover,
\[
    \operatorname{spt}\mu
    =
    \operatorname{spt}\nu
    =
    \operatorname{spt}\|T\|.
\]
Together with
\[
    \operatorname{spt}\|T\|
    =
    \operatorname{spt}\|T_{\mathrm{Li}}\|,
\]
proved above, this gives the stated identification of $\mathcal S$ and
completes the proof.
\end{proof}

\begin{remark}[Normalization of the longitudinal forms]
\label{rem:identification_with_Li_current}
The identification of the currents above is consistent with the
longitudinal forms used in Li's corrected compactness theorem.
The corrected forms in
\cite[Theorem~1.12]{li2025large} have uncorrected parts
\[
    -\frac{\operatorname{Tr}(\Phi_iF_{\nabla_i})}
    {4\pi m_i^{\mathrm{Li}}},
    \qquad
    -\frac{\operatorname{Tr}(\Phi_iF_{\nabla_i})}
    {4\pi|\Phi_i|_{\mathrm{Li}}},
\]
together with an $L^2$-harmonic correction $\sigma_i$. In our
normalization,
\[
    \omega_i
    :=
    m_i^{-1}\langle F_{\nabla_i},\Phi_i\rangle
    =
    -\frac{1}{m_i^{\mathrm{Li}}}
    \operatorname{Tr}(\Phi_iF_{\nabla_i}),
\]
and
\[
    \widehat\omega_i
    :=
    \left\langle
        F_{\nabla_i},
        \frac{\Phi_i}{|\Phi_i|}
    \right\rangle
    =
    -\frac{1}{|\Phi_i|_{\mathrm{Li}}}
    \operatorname{Tr}(\Phi_iF_{\nabla_i}).
\]
Thus the uncorrected parts of Li's longitudinal forms are precisely
\[
    \frac{1}{4\pi}\omega_i,
    \qquad
    \frac{1}{4\pi}\widehat\omega_i.
\]
The relation between $\widehat\omega_i/(4\pi)$ and the actual eigenline
Chern curvature is recorded in
\eqref{eq:eigenline_Chern_curvature}. For $G=\mathrm{SO}(3)$ these forms are still globally defined by the same
Lie algebra formula. Their cohomology classes may be half-integral; twice
the corresponding Chern form represents the Euler class of the
$\mathrm{SO}(2)$ reduction. This charge lattice distinction does not alter
the differentiated current, the measure normalization, or the integer
multiplicity of $T_{\mathrm{Li}}$.

We also record that
\[
    d\omega_i
    =
    \frac12m_i^{-1}Z(\nabla_i,\Phi_i).
\]
Since the harmonic correction $\sigma_i$ is closed, it does not alter the
differentiated current. Thus the calibrated current obtained from Li's
corrected abelian limit and the calibrated current obtained from the
Parise--Pigati--Stern charge current compactness theorem agree after the
normalization conversion above.
\end{remark}

The rest of this section studies the common \textbf{calibrated concentration set}
obtained above. In view of
Proposition~\ref{prop:PPS_saturation_and_Li_identification}, we write
\begin{equation}\label{eq: support_S_various_forms}
    \mathcal S
    :=
    \operatorname{spt}\mu
    =
    \operatorname{spt}\nu
    =
    \operatorname{spt}\|T\|
    =
    \operatorname{spt}\|T_{\mathrm{Li}}\|
\end{equation}
for this support.

\subsection{Calibrated concentration and ordinary regularity scales}
\label{subsec:calibrated_concentration_regular_scales}

We first characterize the calibrated concentration set $\mathcal S$ directly in
terms of the finite mass sequence. For $x\in X$ and $r>0$, define the
\textbf{mass-renormalized codimension-three local energy}
\begin{equation}\label{eq: def_theta_i}
    \theta_i(x,r)
    :=
    r^{3-n}\mu_i(B_r(x))
    =
    m_i^{-1}r^{3-n}
    \int_{B_r(x)}
    \bigl(
        |F_{\nabla_i}|^2
        +
        |\nabla_i\Phi_i|^2
    \bigr).
\end{equation}
Since the $(n-3)$-dimensional density is normalized by
$\omega_{n-3}r^{n-3}$, we also set
\[
    \vartheta_i(x,r)
    :=
    \frac{1}{\omega_{n-3}}\theta_i(x,r)
    =
    \frac{1}{\omega_{n-3}}
    m_i^{-1}r^{3-n}
    \int_{B_r(x)}
    \bigl(
        |F_{\nabla_i}|^2
        +
        |\nabla_i\Phi_i|^2
    \bigr).
\]
By Proposition~\ref{prop:PPS_saturation_and_Li_identification}, we have
$\mu=8\pi\|T\|$, and $T$ is a compactly supported
$\Theta$-calibrated $(n-3)$-cycle. It follows from the fundamental theorem
of Harvey--Lawson \cite[Theorem~4.2]{harvey1982calibrated} that $T$ is homologically mass-minimizing. The standard monotonicity formula for mass-minimizing currents then implies that $\|T\|$ or, equivalently, $\mu$, satisfies the codimension-three monotonicity formula. More precisely, since $(X,g)$ has bounded geometry, there exist constants $c>0$ and $r_0>0$, depending only on the fixed background geometry, such that
\begin{equation}\label{ineq: exact_codim_3_monot_for_mu}
    e^{cs^2}s^{3-n}\mu(B_s(x))
    \leqslant
    e^{cr^2}r^{3-n}\mu(B_r(x))
\end{equation}
for every $x\in X$ and every $0<s<r\leqslant r_0$; see also
\cite[Remark~4.5]{li2025large}.

Let $\theta_T$ denote the multiplicity function of the integral current
$T$, defined $\|T\|$-almost everywhere, so that $\|T\|=\theta_T\mathcal{H}^{n-3}\lfloor \mathcal S$, and set
\begin{equation}\label{eq: dfn_of_l}
    \ell_T
    :=
    \operatorname*{ess\,inf}_{\|T\|}
    \theta_T
    \in
    \mathbb Z_{\geqslant 1}.
\end{equation}
The $(n-3)$-densities of $\|T\|$ and $\mu$ will be denoted by
\[
    \vartheta_T(x)
    :=
    \lim_{r\downarrow0}
    \frac{\|T\|(B_r(x))}
         {\omega_{n-3}r^{n-3}}\quad\text{and}\quad
    \vartheta_\mu(x)
    :=
    \lim_{r\downarrow0}
    \frac{\mu(B_r(x))}
         {\omega_{n-3}r^{n-3}}.
\]
\begin{proposition}[Energy density characterization of $\mathcal S$]
\label{prop:characterization_S}
For every $x\in X$, the density $\vartheta_T(x)$ exists, and hence so does $\vartheta_\mu(x)$ with $\vartheta_\mu(x) = 8\pi\vartheta_T(x)$.
Moreover,
\[
    \vartheta_T(x)
    =
    \theta_T(x)
    \in
    \mathbb Z_{\geqslant 1}\quad\text{for $\|T\|$-almost every $x$,}
\] while pointwise on $X$ one has
\[
    \vartheta_T(x)=0
    \quad\text{for }x\notin \mathcal S,
    \qquad
    \vartheta_T(x)\geqslant \ell_T
    \quad\text{for }x\in \mathcal S.
\]
Consequently,
\[
    x\in \mathcal S
    \quad\Longleftrightarrow\quad
    \vartheta_\mu(x)\geqslant8\pi
    \quad\Longleftrightarrow\quad
    \vartheta_\mu(x)\geqslant8\pi \ell_T.
\]
The same set is characterized directly from the approximating sequence:
\[
    x\in \mathcal S
    \quad\Longleftrightarrow\quad
    \liminf_{r\downarrow0}
    \liminf_{i\to\infty}
    \vartheta_i(x,r)
    \geqslant8\pi \ell_T.
\]
Equivalently, the threshold $8\pi\ell_T$ in the preceding condition may be
replaced by $8\pi$. In terms of the unnormalized energies $\theta_i$,
\[
    x\in \mathcal S
    \quad\Longleftrightarrow\quad
    \liminf_{r\downarrow0}
    \liminf_{i\to\infty}
    \theta_i(x,r)
    \geqslant
    8\pi \ell_T\,\omega_{n-3},
\]
and again the threshold may be replaced by
$8\pi\omega_{n-3}$.
\end{proposition}

\begin{proof}
The monotonicity formula \eqref{ineq: exact_codim_3_monot_for_mu}, equivalently its version for $\|T\|$, implies that $\vartheta_T(x)$ exists at every $x\in \mathcal S$. If
$x\notin \mathcal S$, then some neighbourhood of $x$ has zero $\|T\|$-measure, so
$\vartheta_T(x)=0$. Since $\mu=8\pi\|T\|$, it follows immediately that
\[
    \vartheta_\mu(x)
    =
    8\pi\vartheta_T(x)
\]
for every $x\in X$.

Since $T$ is an integral rectifiable current, the density formula gives
\[
    \vartheta_T(x)
    =
    \theta_T(x)
    \in
    \mathbb Z_{\geqslant 1}
\]
for $\|T\|$-almost every $x$. In particular, $\theta_T\geqslant \ell_T$ 
at $\|T\|$-almost every point.

We claim that the same lower bound holds for the pointwise density at every
$x\in \mathcal S$. Let $G\subset \mathcal S$ be a set of full $\|T\|$-measure on which
$\vartheta_T=\theta_T\geqslant \ell_T$. Fix $x\in \mathcal S$, $\varepsilon\in(0,1)$,
and a sufficiently small $r>0$. Since $x\in\operatorname{spt}\|T\|$, the
ball $B_{\varepsilon r}(x)$ has positive $\|T\|$-measure and therefore
contains a point $q\in G$. The inclusion
\[
    B_{(1-\varepsilon)r}(q)
    \subset
    B_r(x)
\]
and calibrated monotonicity centred at $q$ give
\[
    \frac{\|T\|(B_r(x))}
         {\omega_{n-3}r^{n-3}}
    \geqslant
    e^{-cr^2}
    (1-\varepsilon)^{n-3}
    \vartheta_T(q)
    \geqslant
    e^{-cr^2}
    (1-\varepsilon)^{n-3}\ell_T.
\]
Letting $r\downarrow0$ and then $\varepsilon\downarrow0$ yields
\[
    \vartheta_T(x)\geqslant \ell_T.
\]
This proves the pointwise density characterization of $\mathcal S$.

We now compare the limiting density with the approximating energies. By
weak convergence of Radon measures and the Portmanteau theorem,
\[
    \liminf_{i\to\infty}\mu_i(B_r(x)) 
    \geqslant
    \mu(B_r(x))
\]
for every fixed open ball. Hence
\[
    \liminf_{i\to\infty}
    \vartheta_i(x,r)
    \geqslant
    \frac{\mu(B_r(x))}
         {\omega_{n-3}r^{n-3}}.
\]
If $x\in \mathcal S$, letting $r\downarrow0$ and using
$\vartheta_\mu(x)\geqslant8\pi\ell_T$ gives
\[
    \liminf_{r\downarrow0}
    \liminf_{i\to\infty}
    \vartheta_i(x,r)
    \geqslant8\pi\ell_T.
\]
Conversely, suppose that $x\notin \mathcal S$. Since $\mathcal S$ is closed, there exists
$r_x>0$ such that
\[
    \overline{B_{r_x}(x)}\cap \mathcal S=\varnothing.
\]
For every $0<r<r_x$, the Portmanteau theorem for the closed ball gives
\[
    0
    \leqslant
    \limsup_{i\to\infty}\mu_i(B_r(x))
    \leqslant
    \limsup_{i\to\infty}\mu_i(\overline{B_r(x)})
    \leqslant
    \mu(\overline{B_r(x)})
    =
    0.
\]
Thus $\vartheta_i(x,r)\to0$ for every sufficiently small $r$, and hence
\[
    \liminf_{r\downarrow0}
    \liminf_{i\to\infty}
    \vartheta_i(x,r)
    =
    0.
\]
The remaining equivalent formulations follow from
$\ell_T\geqslant1$ and from
$\theta_i=\omega_{n-3}\vartheta_i$.
\end{proof}

The preceding characterization has the following finite scale form:
\begin{align}
\mathcal S
&=
\bigcap_{0<r\leqslant r_0}
\left\{
x\in X:
\liminf_{i\to\infty}
e^{cr^2}\vartheta_i(x,r)>0
\right\}
\label{eq: characterization_of_S}
\\
&=
\bigcap_{0<r\leqslant r_0}
\left\{
x\in X:
\liminf_{i\to\infty}
e^{cr^2}\vartheta_i(x,r)\geqslant8\pi
\right\}
\\
&=
\bigcap_{0<r\leqslant r_0}
\left\{
x\in X:
\liminf_{i\to\infty}
e^{cr^2}\vartheta_i(x,r)\geqslant8\pi \ell_T
\right\}.
\label{eq: energy_concentration_set_S}
\end{align}
Indeed, if $x\in \mathcal S$, then
\eqref{ineq: exact_codim_3_monot_for_mu} and weak lower semicontinuity give,
for every $0<r\leqslant r_0$,
\[
    \liminf_{i\to\infty}
    e^{cr^2}\vartheta_i(x,r)
    \geqslant
    e^{cr^2}
    \frac{\mu(B_r(x))}
         {\omega_{n-3}r^{n-3}}
    \geqslant
    \vartheta_\mu(x)
    \geqslant
    8\pi \ell_T.
\]
If $x\notin \mathcal S$, one may choose $0<r\leqslant r_0$ such that
$\overline{B_r(x)}\cap \mathcal S=\varnothing$, in which case
$\vartheta_i(x,r)\to0$. This proves all three equalities.

We next record a general size estimate for the corresponding concentration
sets at arbitrary Morrey scale. This will be used below to obtain the
Hausdorff measure bound for the calibrated concentration set $\mathcal S$.

\begin{lemma}[Size estimates on energy concentration sets]
\label{lem: size_energy_concentration}
Let $\alpha,\varepsilon>0$, with $\alpha\leqslant n$, and define
\[
    S_{\alpha,\varepsilon}
    :=
    \bigcap_{0<r\leqslant r_0}
    \left\{
    x\in X:
    \liminf_{i\to\infty}
    m_i^{-1}e^{cr^2}r^{\alpha-n}
    \int_{B_r(x)}
    \bigl(
        |F_{\nabla_i}|^2+|\nabla_i\Phi_i|^2
    \bigr)
    \geqslant\varepsilon
    \right\}.
\]
Then $S_{\alpha,\varepsilon}$ is compact and
\[
    \mathcal H^{n-\alpha}(S_{\alpha,\varepsilon})
    \leqslant
    C_{n,\alpha}
    \frac{8\pi k}{\varepsilon}.
\]
\end{lemma}

\begin{proof}
We first prove compactness of $S_{\alpha,\varepsilon}$. Since $(X,g)$ is complete, it is enough to show that $S_{\alpha,\varepsilon}$ is closed and bounded. To prove
boundedness, fix once and for all $\bar r\in(0,r_0]$. If $\rho(x)\geqslant R_0$, Corollary~\ref{cor: local_L2_estimates} gives
\[
    \bar r^{4-n}
    \int_{B_{\bar r}(x)}
    \bigl(
        |F_{\nabla_i}|^2+|\nabla_i\Phi_i|^2
    \bigr)
    \leqslant
    C\bigl(\rho(x)^{4-n}m_i+1\bigr),
\]
where $C$ is independent of $i$ and $x$. Hence
\[
\begin{aligned}
&\liminf_{i\to\infty}
m_i^{-1}e^{c\bar r^2}\bar r^{\alpha-n}
\int_{B_{\bar r}(x)}
\bigl(
    |F_{\nabla_i}|^2+|\nabla_i\Phi_i|^2
\bigr)
\\
&\qquad\leqslant
C e^{c\bar r^2}\bar r^{\alpha-4}
\liminf_{i\to\infty}
\bigl(
    \rho(x)^{4-n}+m_i^{-1}
\bigr)
\\
&\qquad=
C e^{c\bar r^2}\bar r^{\alpha-4}\rho(x)^{4-n}.
\end{aligned}
\]
Since $n>4$, the right-hand side tends to zero as
$\rho(x)\to\infty$. Therefore the defining inequality for
$S_{\alpha,\varepsilon}$ fails at the fixed radius $\bar r$ for every
point sufficiently far out on the AC end. Thus
$S_{\alpha,\varepsilon}$ is bounded.

To prove closedness, let $x_j\in S_{\alpha,\varepsilon}$ and suppose
$x_j\to x$. Fix $0<r\leqslant r_0$ and choose $0<s<r$. For all sufficiently
large $j$, one has $B_s(x_j)\subset B_r(x)$. Hence
\[
\begin{aligned}
&\liminf_{i\to\infty}
m_i^{-1}e^{cs^2}s^{\alpha-n}
\int_{B_r(x)}
\bigl(
|F_{\nabla_i}|^2+|\nabla_i\Phi_i|^2
\bigr)
\\
&\qquad\geqslant 
\liminf_{i\to\infty}
m_i^{-1}e^{cs^2}s^{\alpha-n}
\int_{B_s(x_j)}
\bigl(
|F_{\nabla_i}|^2+|\nabla_i\Phi_i|^2
\bigr)
\geqslant
\varepsilon.
\end{aligned}
\]
Letting $s\uparrow r$ gives the defining inequality at $x$, and therefore
$x\in S_{\alpha,\varepsilon}$.

It remains to prove the Hausdorff measure bound. Since
$S_{\alpha,\varepsilon}$ is compact, choose $R_*\geqslant R_0$ such that
$S_{\alpha,\varepsilon}\subset\{\rho\leqslant R_*\}$. For
$0<\delta\ll_g1$, choose a finite Vitali subcover
$\{B_{4\delta}(x_\lambda)\}$ of $S_{\alpha,\varepsilon}$ such that
$x_\lambda\in S_{\alpha,\varepsilon}$ and the balls
$B_{2\delta}(x_\lambda)$ are pairwise disjoint. Then
\[
\begin{aligned}
    \sum_\lambda \delta^{n-\alpha}
    &\leqslant
    \frac{1}{\varepsilon}
    \sum_\lambda
    \liminf_{i\to\infty}
    m_i^{-1}e^{c\delta^2}
    \int_{B_\delta(x_\lambda)}
    \bigl(
    |F_{\nabla_i}|^2+|\nabla_i\Phi_i|^2
    \bigr)
    \\
    &\leqslant
    \frac{e^{c\delta^2}}{\varepsilon}
    \liminf_{i\to\infty}
    m_i^{-1}
    \int_{\bigcup_\lambda B_\delta(x_\lambda)}
    \bigl(
    |F_{\nabla_i}|^2+|\nabla_i\Phi_i|^2
    \bigr)
    \\
    &\leqslant
    \frac{e^{c\delta^2}}{\varepsilon}
    \liminf_{i\to\infty}
    m_i^{-1}
    \int_{\{\rho\leqslant R_*+C\delta\}}
    \bigl(
    |F_{\nabla_i}|^2+|\nabla_i\Phi_i|^2
    \bigr).
\end{aligned}
\]
By Lemma~\ref{lemm: Yang-Li} and the linear energy identity
\eqref{eq: linear_interm_energy_growth}, the last line is bounded by
\[
    \frac{e^{c\delta^2}}{\varepsilon}
    \liminf_{i\to\infty}
    m_i^{-1}
    \bigl(
        8\pi k\,m_i+C(R_*+\delta)^{n-4}
    \bigr)
    =
    e^{c\delta^2}\frac{8\pi k}{\varepsilon}.
\]
Letting $\delta\downarrow0$ gives
\[
    \mathcal H^{n-\alpha}(S_{\alpha,\varepsilon})
    \leqslant
    C_{n,\alpha}
    \frac{8\pi k}{\varepsilon}.
\]
\end{proof}

As an immediate consequence of
Lemma~\ref{lem: size_energy_concentration} and the characterization
\eqref{eq: energy_concentration_set_S}, we obtain the Hausdorff size bound
for the calibrated concentration set. Its compactness was already obtained
in Proposition~\ref{prop:PPS_saturation_and_Li_identification}.

\begin{corollary}
\label{cor: size_S}
The compact calibrated concentration set $\mathcal S$ satisfies
\[
    \mathcal H^{n-3}(\mathcal S)
    \leqslant
    C_n\frac{k}{\ell_T},
\] where $\ell_T$ is defined as in \eqref{eq: dfn_of_l}.
\end{corollary}

In \cite[Remark~3.10]{li2025large}, Li points out that when
$d<n$, a sequence of closed $d$-rectifiable currents may have uniformly
bounded mass and supports with uniformly bounded $\mathcal H^d$-measure,
while those supports nevertheless become $\varepsilon$-dense in an open
set for arbitrarily small $\varepsilon>0$. In the present setting, this
caveat is relevant to Li's curvature concentration loci
$C_{\Lambda_0,i}$, which we recall later in
\S\ref{subsec:Li_moving_loci}, and to their Kuratowski upper limit
$\mathcal C$, introduced in Definition~\ref{def:C_infty}. It is
also relevant to the Higgs zero sets $\Phi_i^{-1}(0)$ and to their Kuratowski upper limit, the limiting Higgs zero set $\mathcal Z$ defined below in \eqref{eq:def_limiting_zero_set}. Indeed, neither $\mathcal Z$ nor $\mathcal C$ is known at this stage to be rectifiable, let alone to be the support of a calibrated current. Thus, Li's remark does not apply to these limiting sets literally; rather, it illustrates why uniform size estimates for the zero sets and the curvature concentration loci do not by themselves prevent their Kuratowski upper limits from being large or even from containing open regions.

The situation is different for the fixed concentration set
$\mathcal S=\operatorname{spt}\mu=\operatorname{spt}\|T\|$.
By Proposition~\ref{prop:characterization_S}, calibrated monotonicity
propagates the almost-everywhere integer multiplicity of $T$ to a uniform
positive lower density bound at every point of $\mathcal S$. Together with the fixed total mass $\|T\|(X)=k$, this prevents $\mathcal S$ from becoming arbitrarily dense in any open region, as quantified in the following lemma.

\begin{lemma}[Quantitative density bound for the calibrated support]
\label{lem:S_not_arbitrarily_dense}
There exist constants $c_0>0$ and $r_1\in(0,r_0]$, depending only on the
fixed background geometry, such that the following holds. Suppose that
$0<R\leqslant r_1/4$, $0<\varepsilon<R/10$, and that $\mathcal S$ is
$\varepsilon$-dense in $B_R(x)$. Then
\[
    \|T\|(B_{2R}(x))
    \geqslant
    c_0\ell_T R^n\varepsilon^{-3}.
\]
In particular, since $\|T\|(X)=k$,
\[
    \varepsilon^3
    \geqslant
    \frac{c_0 \ell_TR^n}{k}.
\]
Consequently, $\mathcal S$ cannot be $\varepsilon_j$-dense in any fixed ball
$B_R(x)$ with $0<R\leqslant r_1/4$ for a sequence
$\varepsilon_j\downarrow0$.
\end{lemma}

\begin{proof}
By Proposition~\ref{prop:characterization_S}, one has
$\vartheta_T(p)\geqslant \ell_T$ for every $p\in\mathcal S$. Letting the
smaller radius tend to zero in the monotonicity formula for $\|T\|$, we
obtain
\[
    e^{cs^2}s^{3-n}\|T\|(B_s(p))
    \geqslant
    \omega_{n-3}\vartheta_T(p)
    \geqslant
    \omega_{n-3}\ell_T
\]
for every $p\in \mathcal S$ and every $0<s\leqslant r_0$. After decreasing $r_1$
if necessary, there is therefore a constant $a_0>0$, depending only on the
background geometry, such that
\begin{equation}\label{ineq:lower_mass_density_S}
    \|T\|(B_s(p))
    \geqslant
    a_0 \ell_T s^{n-3}
\end{equation}
for every $p\in \mathcal S$ and every $0<s\leqslant r_1$.

Choose a maximal $8\varepsilon$-separated set $\{x_1,\ldots,x_N\}\subset B_R(x)$. By maximality, the balls $B_{8\varepsilon}(x_j)$ cover $B_R(x)$. Since $R\leqslant r_1/4$ and $(X,g)$ has bounded geometry, the uniform small scale volume bounds give
\[
    cR^n
    \leqslant
    \operatorname{vol}(B_R(x))
    \leqslant
    \sum_{j=1}^N
    \operatorname{vol}(B_{8\varepsilon}(x_j))
    \leqslant
    CN\varepsilon^n.
\]
Consequently,
\[
    N
    \geqslant
    b_0
    \left(
        \frac{R}{\varepsilon}
    \right)^n
\]
for a constant $b_0>0$ depending only on the background geometry.

Since $\mathcal S$ is $\varepsilon$-dense in $B_R(x)$, for every $j$
there exists $p_j\in\mathcal S$ such that $d(p_j,x_j)<\varepsilon$.
The points $p_j$ are $6\varepsilon$-separated, so the balls $B_\varepsilon(p_j)$ are pairwise disjoint. Moreover, $\varepsilon<R/10$ implies $B_\varepsilon(p_j)\subset B_{2R}(x)$. Applying \eqref{ineq:lower_mass_density_S} at each $p_j$, we obtain
\[
    \|T\|(B_{2R}(x))
    \geqslant
    \sum_{j=1}^N
    \|T\|(B_\varepsilon(p_j))
    \geqslant
    a_0 \ell_T N\varepsilon^{n-3}
    \geqslant
    c_0 \ell_T R^n\varepsilon^{-3}.
\]
The final assertion follows from $\|T\|(X)=k$.
\end{proof}

\begin{remark}
\label{rmk:S_vs_moving_dense_loci}
Apart from the uniform small scale volume bounds of the fixed background
manifold, the proof of Lemma~\ref{lem:S_not_arbitrarily_dense} uses two
properties of the measure $\|T\|$: its finite total mass and the uniform
lower growth estimate
\[
    \|T\|(B_s(p))
    \geqslant
    a_0\ell_T s^{n-3}
\]
for every $p\in\mathcal S$ and every $0<s\leqslant r_1$. In the present
setting, this estimate follows from calibrated monotonicity together with
the pointwise lower bound for the density of $T$ on $\mathcal S$ obtained
in Proposition~\ref{prop:characterization_S}. Rectifiability and finite
mass alone do not give such a uniform estimate at every point of the
support.

More generally, the same covering argument applies to any finite Radon
measure $\lambda$ on the fixed manifold for which there are constants
$a>0$ and $r_*>0$ such that
\[
    \lambda(B_s(p))
    \geqslant
    a\,s^{n-3}
\]
for every $p\in\operatorname{spt}\lambda$ and every
$0<s\leqslant r_*$. The conclusion is therefore a consequence of finite
total mass and uniform lower growth, rather than of rectifiability by
itself.

Lemma~\ref{lem:S_not_arbitrarily_dense} concerns the fixed calibrated
support $\mathcal S=\operatorname{spt}\|T\|$. It does not rule out the
possibility raised in \cite[Remark~3.10]{li2025large} that the Higgs zero sets or curvature concentration loci become increasingly
dense before passage to the limit: no single finite measure satisfying
the uniform lower growth estimate above is associated with either
sequence of sets.

The lemma gives a quantitative form, in the present setting, of the
distinction noted in \cite[Remark~4.4]{li2025large}. Once the limiting
current is known to be calibrated, its monotonicity and density properties
give stronger control of its fixed support than is available for the
support of a general rectifiable current.
\end{remark}

The set $\mathcal S$ arises from the mass-renormalized codimension-three large
mass limit and should be distinguished from the concentration set
associated with the ordinary, unrenormalized codimension-four
Yang--Mills--Higgs regularity scale. We briefly compare these two
concentration mechanisms.

After increasing $c$ and decreasing $r_0$ if necessary, we henceforth use
the same constants $c>0$ and $r_0>0$ in the calibrated
codimension-three monotonicity formula
\eqref{ineq: exact_codim_3_monot_for_mu}
and in the codimension-four Yang--Mills--Higgs almost-monotonicity formula
of Corollary~\ref{cor:codim4_almost_monotonicity}. Enlarging $c$ and
shrinking $r_0$ preserve both inequalities.

Writing $e_i:=e(\nabla_i,\Phi_i)$, we define the \textbf{codimension-four concentration set}
\[
    \mathcal{B}_{\mathrm{TU}}
    :=
    \bigcap_{0<r\leqslant r_0}
    \left\{
        x\in X:
        \liminf_{i\to\infty}
        e^{cr^2}r^{4-n}
        \int_{B_r(x)}
        e_i
        \geqslant
        \varepsilon_0
    \right\},
\]
where $\varepsilon_0>0$ is the threshold in
Theorem~\ref{thm: total_epsilon_regularity}, and $c$ and $r_0$ are the
constants in Corollary~\ref{cor:codim4_almost_monotonicity}.

This definition may equivalently be formulated, after changing the
threshold by a fixed factor, using only the curvature energy. Indeed, the
$\Theta$-monopole equation gives
\[
    |\nabla_i\Phi_i|
    =
    |F_{\nabla_i}\wedge\Theta|,
\]
and, since $\Theta$ is parallel, there is a constant
$C_\Theta<\infty$, depending only on the fixed background geometry, such
that
\[
    |F_{\nabla_i}|^2
    \leqslant
    e_i
    \leqslant
    C_\Theta |F_{\nabla_i}|^2.
\]
We also use the \textbf{regularity scale}
$\mathfrak r_i:X\to(0,r_0]$ associated with
$(\nabla_i,\Phi_i)$, defined by
\[
    \mathfrak r_i(x)
    :=
    \sup\left\{
        0<r\leqslant r_0:
        s^{4-n}
        \int_{B_s(x)}
        e_i
        <
        \varepsilon_0
        \text{ for every }0<s\leqslant r
    \right\}.
\]
Since $e_i$ is smooth, the quantity
$s^{4-n}\int_{B_s(x)}e_i$ tends to zero as $s\downarrow0$, and hence
$\mathfrak r_i(x)>0$ for every $x\in X$. The small energy condition is
imposed at all smaller radii in order to make $\mathfrak r_i$ a genuine
regularity scale. By the codimension-four monotonicity formula for Yang--Mills--Higgs fields, see Corollary~\ref{cor:codim4_almost_monotonicity}, this definition is
equivalent, after changing $r_0$ and $\varepsilon_0$ by fixed geometric
factors, to checking the corresponding small energy condition at the
outer radius.

The quantity defining $\mathcal{B}_{\mathrm{TU}}$ is the usual
codimension-four scale-invariant energy controlling ordinary
$\varepsilon$-regularity. For sequences with locally uniformly bounded
unrenormalized energy, its concentration set is the classical
Tian--Uhlenbeck bubbling locus. In the present large mass problem,
however, the unrenormalized energy diverges near the codimension-three
concentration set. Accordingly, $\mathcal{B}_{\mathrm{TU}}$ is the failure
set for ordinary codimension-four regularity;
it need not itself have codimension-four.

By contrast, $\mathcal S$ is detected by the mass-renormalized scale
\[
    \int_{B_r(x)}
    e_i
    \sim
    m_i r^{n-3}.
\]

\begin{proposition}
\label{prop:S_subset_TU}
One has
\[
    \mathcal S\subset\mathcal{B}_{\mathrm{TU}}.
\]
More precisely, if $x\in \mathcal S$, then for every fixed
$0<r\leqslant r_0$,
\[
    \liminf_{i\to\infty}
    e^{cr^2}r^{4-n}
    \int_{B_r(x)}
    e_i
    =
    +\infty.
\]
\end{proposition}

\begin{proof}
Let $x\in \mathcal S$ and fix $0<r\leqslant r_0$. By the finite scale
characterization \eqref{eq: energy_concentration_set_S},
\[
    \liminf_{i\to\infty}
    e^{cr^2}\vartheta_i(x,r)
    \geqslant
    8\pi\ell_T.
\]
Using the definition of $\vartheta_i$, we have
\[
\begin{aligned}
    e^{cr^2}r^{4-n}
    \int_{B_r(x)}
    e_i
    &=
    \omega_{n-3}m_i r\,
    e^{cr^2}\vartheta_i(x,r).
\end{aligned}
\]
Since $m_i\to+\infty$ and $r>0$ is fixed, the right-hand side has
lower limit $+\infty$. Thus $x\in\mathcal{B}_{\mathrm{TU}}$.
\end{proof}

\begin{remark}
The reverse inclusion is not implied by the preceding argument. An
ordinary codimension-four bubble whose unrenormalized Yang--Mills--Higgs energy remains of order one on fixed scales would be detected by $\mathcal{B}_{\mathrm{TU}}$, but its contribution to the mass-renormalized measures is suppressed by the factor $m_i^{-1}$ and may therefore disappear in the limit as $i\to\infty$. Thus $\mathcal{B}_{\mathrm{TU}}\setminus \mathcal S$ may, in principle, contain ordinary
codimension-four bubbling which is invisible to the limiting measure
$\mu$. We do not assert here that such bubbling actually occurs for
sequences in the {\mainsettingref}.
\end{remark}

\subsection{Limiting Higgs zeros and concentration centred at Higgs zeros}
\label{subsec:limiting_Higgs_zeros}

We next compare the calibrated support with the limiting Higgs zero
set. In dimension three, the corrected and expanded concentration theory
shows that the blow-up set agrees with the limiting zero set for
monopoles of fixed charge on AC manifolds, and that every concentration
weight is the total charge of a complete finite cluster of mass one
Euclidean monopoles
\cite[Theorem~1.1 and Proposition~8.3]{fadeloliveira2026limitv5}.
Zeros and charge may nevertheless escape through the AC end. In the
higher-dimensional setting, equality need not follow from the available
estimates: the cohomogeneity-one families satisfy it, but the general
theory in principle allows limiting Higgs zeros away from the calibrated support.

Define the Higgs zero sets $Z_i$ and their Kuratowski upper limit, the
\textbf{limiting Higgs zero set} $\mathcal Z$, by
\begin{equation}
\label{eq:def_limiting_zero_set}
    Z_i
    :=
    \Phi_i^{-1}(0),
    \qquad
    \mathcal Z
    :=
    \bigcap_{N\geqslant1}
    \overline{
        \bigcup_{i\geqslant N}Z_i
    }.
\end{equation}
Thus $\mathcal Z$ is closed and consists precisely of those points $x\in X$ for
which there exist indices $i_j\to\infty$ and points
$x_j\in Z_{i_j}$ such that $x_j\to x$.

\begin{proposition}[Limiting Higgs zeros and ordinary regularity scales]
\label{prop:zero_centered_TU_concentration}
Let $x\in\mathcal Z$. Then there exist indices $i_j\to\infty$ and points
$x_j\in Z_{i_j}$ such that $x_j\to x$ and, setting
\[
    \lambda_j
    :=
    \mathfrak r_{i_j}(x_j),
\]
one has
\[
    \lambda_j\longrightarrow0.
\]
For all sufficiently large $j$,
\begin{equation}
\label{eq:zero_centered_threshold}
    \lambda_j^{4-n}
    \int_{B_{\lambda_j}(x_j)}
    e_{i_j}
    =
    \varepsilon_0.
\end{equation}
Moreover, after relabelling the subsequence realizing the zero, $x$
belongs to the associated codimension-four concentration set. Equivalently,
in the relabelled notation, $x\in\mathcal B_{\mathrm{TU}}$. More precisely,
for every fixed $0<r\leqslant r_0$,
\begin{equation}
\label{eq:zero_centered_TU_lower_bound}
    \liminf_{j\to\infty}
    e^{cr^2}r^{4-n}
    \int_{B_r(x)}
    e_{i_j}
    \geqslant
    \varepsilon_0.
\end{equation}
Finally, there is a constant $c_\Theta>0$, depending only on the fixed
parallel form $\Theta$, such that
\[
    \lambda_j^{4-n}
    \int_{B_{\lambda_j}(x_j)}
    |F_{\nabla_{i_j}}|^2
    \geqslant
    c_\Theta\varepsilon_0
\]
for all sufficiently large $j$.
\end{proposition}

\begin{proof}
Choose indices $i_j\to\infty$ and points $x_j\in Z_{i_j}$ such that
$x_j\to x$. We first prove that
$\lambda_j=\mathfrak r_{i_j}(x_j)\to0$. Suppose otherwise. After passing
to a subsequence, there is $\lambda_*>0$ such that
$\lambda_j\geqslant\lambda_*$ for every $j$. Fix
\[
    0<\rho<\frac{\lambda_*}{4}
\]
small enough that the balls $B_{2\rho}(x_j)$ are contained in a fixed
compact subset $K\Subset X$ for all sufficiently large $j$. By the
definition of the regularity scale,
\[
    (2\rho)^{4-n}
    \int_{B_{2\rho}(x_j)}
    e_{i_j}
    <
    \varepsilon_0.
\]
Theorem~\ref{thm: total_epsilon_regularity} therefore gives a constant
$C_\rho<\infty$, independent of $j$, such that
\[
    \sup_{B_\rho(x_j)}
    |\nabla_{i_j}\Phi_{i_j}|
    \leqslant
    C_\rho.
\]
Since $\Phi_{i_j}(x_j)=0$, it follows that
\[
    |\Phi_{i_j}|
    \leqslant
    C_\rho\rho
    \qquad
    \text{on }B_\rho(x_j).
\]
As $m_{i_j}\to\infty$, for all sufficiently large $j$ one has
\[
    m_{i_j}-|\Phi_{i_j}|
    \geqslant
    \frac{m_{i_j}}{2}
    \qquad
    \text{on }B_\rho(x_j).
\]
Using the uniform lower volume bound for balls on the fixed AC manifold,
we obtain
\[
    \frac{1}{m_{i_j}}
    \int_K
    \bigl(m_{i_j}-|\Phi_{i_j}|\bigr)^2
    \geqslant
    \frac{m_{i_j}}{4}
    \operatorname{vol}(B_\rho(x_j))
    \longrightarrow
    +\infty.
\]
This contradicts
\eqref{eq: largeness_of_Higgs_field}. Hence
$\lambda_j\to0$.

For all sufficiently large $j$, one has $\lambda_j<r_0$. Set
\[
    Q_j(s)
    :=
    s^{4-n}
    \int_{B_s(x_j)}
    e_{i_j}.
\]
The function $Q_j$ extends continuously to $s=0$ by $Q_j(0)=0$. By the definition of $\lambda_j$,
\[
    Q_j(s)<\varepsilon_0
    \qquad
    \text{for every }0<s<\lambda_j.
\]
Continuity therefore gives
\[
    Q_j(\lambda_j)\leqslant\varepsilon_0.
\]
If $Q_j(\lambda_j)<\varepsilon_0$, continuity on the compact interval $[0,\lambda_j]$ would give
\[
    \max_{0\leqslant s\leqslant\lambda_j}Q_j(s)
    <
    \varepsilon_0,
\]
and the defining inequality for $\mathfrak r_{i_j}(x_j)$ would persist
on a slightly larger interval. This contradicts the definition of
$\lambda_j$. Thus
\eqref{eq:zero_centered_threshold} holds.

We now prove
\eqref{eq:zero_centered_TU_lower_bound}. Fix
$0<r\leqslant r_0$ and set
\[
    r_j
    :=
    r-d(x_j,x).
\]
For all sufficiently large $j$, one has
\[
    0<\lambda_j<r_j<r_0,
    \qquad
    B_{r_j}(x_j)\subset B_r(x),
\]
and $r_j\to r$. By
Corollary~\ref{cor:codim4_almost_monotonicity},
\[
\begin{aligned}
    e^{cr_j^2}r_j^{4-n}
    \int_{B_{r_j}(x_j)}
    e_{i_j}
    &\geqslant
    e^{c\lambda_j^2}\lambda_j^{4-n}
    \int_{B_{\lambda_j}(x_j)}
    e_{i_j}
    \\
    &=
    e^{c\lambda_j^2}\varepsilon_0.
\end{aligned}
\]
Since $B_{r_j}(x_j)\subset B_r(x)$, multiplying by
\[
    \frac{e^{cr^2}r^{4-n}}
         {e^{cr_j^2}r_j^{4-n}}
    \longrightarrow1
\]
and taking the lower limit proves \eqref{eq:zero_centered_TU_lower_bound}. After relabelling the subsequence, this is precisely the defining condition for $x$ to belong to the corresponding set $\mathcal B_{\mathrm{TU}}$.

Finally, the $\Theta$-monopole equation and the fact that $\Theta$ is
parallel give
\[
    e_{i_j}
    =
    |F_{\nabla_{i_j}}|^2
    +
    |F_{\nabla_{i_j}}\wedge\Theta|^2
    \leqslant
    C_\Theta
    |F_{\nabla_{i_j}}|^2.
\]
The curvature lower bound follows from \eqref{eq:zero_centered_threshold}, with
$c_\Theta:=C_\Theta^{-1}$.
\end{proof}

Next, we prove not only that
\[
    \mathcal S\subset \mathcal Z,
\]
but also the stronger statement that every neighbourhood of every point
of $\mathcal S$ meets $Z_i$ for all sufficiently large $i$. Since $\mathcal S$ is compact,
this implies the one-sided Hausdorff convergence
\[
    \sup_{x\in \mathcal S}d(x,Z_i)
    \longrightarrow
    0.
\]
The key input is the identification $T=T_{\mathrm{Li}}$ from Proposition~\ref{prop:PPS_saturation_and_Li_identification}. Li's transverse analysis applies on the smooth locus of $T_{\mathrm{Li}}$, in every fixed precompact region, outside exceptional subsets of arbitrarily small $\mathcal H^{n-3}$-measure. Given a neighbourhood of a point of $\mathcal S$, we first choose a smooth patch of $T_{\mathrm{Li}}$ of positive $\mathcal H^{n-3}$ measure and a fixed normal tube compactly contained in that neighbourhood. We then fix the parameter controlling the exceptional set. For every sufficiently large $i$, at least one normal fibre over the patch is such that the topological argument in \cite[proof of Theorem~5.5, Step~4]{li2025large} applies to that fibre and produces an outer normal linking sphere contained in a normal $3$-ball. The bundle over this ball admits a local $\mathrm{SU}(2)$ lift, and the lifted eigenline on the outer sphere has first Chern number equal to the positive local multiplicity of $T_{\mathrm{Li}}$. If the Higgs field were nonvanishing on the whole ball, that eigenline would extend over the ball and have zero first Chern number. This yields a zero of $\Phi_i$ in the prescribed neighbourhood for every sufficiently large index.

In the proof below we use input from Li's
concentration--decay analysis, the curvature concentration loci
$C_{\Lambda_0,i}$. These loci are recalled systematically, with the
conversion to the present normalization, in
Subsection~\ref{subsec:Li_moving_loci}. For the present argument, only
two consequences are needed: first, Li's dichotomy implies that, for all
sufficiently large $i$,
\[
    Z_i\subset C_{\Lambda_0,i};
\]
second, the normal linking spheres produced in
\cite[proof of Theorem~5.5, Step~4]{li2025large} may be chosen disjoint
from the $m_i^{-1}$-neighbourhood of $C_{\Lambda_0,i}$. Thus the present
proof uses only these quoted facts from Li's analysis; the fixed upper
limiting locus $\mathcal C$ and the open set $\mathcal R_{\mathrm{ab}}$
are introduced later in
Definition~\ref{def:C_infty}.

\begin{proposition}[Higgs zeros near the calibrated support]
\label{prop:S_subset_Z}
For every $x\in \mathcal S$ and every open neighbourhood $V$ of $x$, there exists
$i_V\in\mathbb N$ such that
\[
    V\cap Z_i\neq\varnothing
\]
for every $i\geqslant i_V$. Consequently,
\[
    \lim_{i\to\infty}
    \sup_{x\in \mathcal S}d(x,Z_i)
    =
    0.
\]
In particular,
\[
    \mathcal S\subset \mathcal Z.
\]
\end{proposition}

\begin{proof}
By
Proposition~\ref{prop:PPS_saturation_and_Li_identification},
\[
    T=T_{\mathrm{Li}},
    \qquad
    \mathcal S
    =
    \operatorname{spt}\|T\|
    =
    \operatorname{spt}\|T_{\mathrm{Li}}\|.
\]
Fix $x\in \mathcal S$ and an open neighbourhood $V$ of $x$. We prove that
$V\cap Z_i\neq\varnothing$ for every sufficiently large $i$.

Let $Q_{\mathrm{sm}}$ denote the smooth locus of
$T_{\mathrm{Li}}$. By
\cite[Theorem~4.3, item~(3)]{li2025large},
the singular set of $T_{\mathrm{Li}}$ has Hausdorff dimension at most
$n-5$, and therefore has zero $\mathcal H^{n-3}$-measure. Since
$x\in\operatorname{spt}\|T_{\mathrm{Li}}\|$, every neighbourhood of $x$
has positive $\|T_{\mathrm{Li}}\|$-measure. As
$T_{\mathrm{Li}}$ is represented on $Q_{\mathrm{sm}}$ by integration
with positive integer multiplicity, it follows that
\[
    \mathcal H^{n-3}(Q_{\mathrm{sm}}\cap V)>0.
\]
We may therefore choose a connected, relatively open subset
\[
    K\Subset Q_{\mathrm{sm}}\cap V
\]
such that
\[
    0<\mathcal H^{n-3}(K)<\infty.
\]
By the constancy theorem, after shrinking $K$ if necessary, the
restriction of $T_{\mathrm{Li}}$ to $K$ has constant positive
multiplicity:
\[
    T_{\mathrm{Li}}\lfloor K
    =
    q\,[K],
    \qquad
    q\in\mathbb Z_{\geqslant1},
\]
where $K$ carries the orientation induced by the calibration.

Let $\mathcal U$ be the tubular neighbourhood of
$Q_{\mathrm{sm}}$ used in Li's transverse analysis, with projection
\[
    \varpi:\mathcal U\longrightarrow Q_{\mathrm{sm}}.
\]
Since $\overline K\Subset Q_{\mathrm{sm}}\cap V$ and $V$ is open, there
exists a positive continuous fibre radius function on $\overline K$,
bounded away from zero, such that the corresponding closed normal tube is
contained in $\mathcal U\cap V$. After decreasing this radius if
necessary, we obtain a tubular neighbourhood $\mathcal U_K$ of $K$
satisfying
\[
    \overline{\mathcal U_K}
    \Subset
    \mathcal U\cap V.
\]
Choose once and for all $r_*<\infty$ so large that
\[
    \overline{\mathcal U_K}
    \subset
    \{\rho<r_*\}.
\]
The region $\{\rho<r_*\}$ is fixed independently of $i$, as required in
Li's quantitative transverse energy identity \cite[Theorem~5.5]{li2025large}. For $0<\delta<\mathcal H^{n-3}(K)$, this gives an exceptional set
\[
    E_i\subset Q_{\mathrm{sm}},
    \qquad
    \mathcal H^{n-3}(E_i)<\delta,
\]
such that the conclusions of that theorem hold on every normal fibre
over $Q_{\mathrm{sm}}\setminus E_i$. Since
$\mathcal H^{n-3}(K)>\delta$, one has
\[
    K\setminus E_i\neq\varnothing.
\]
Choose $y_i\in K\setminus E_i$. Because the same fixed $\delta$ is used for all sufficiently large
indices, this choice does not require passing to a subsequence depending
on the point $x$ or on the neighbourhood $V$.

Using the forward notation recalled before the statement, the topological degree argument
in \cite[proof of Theorem~5.5, Step~4]{li2025large} produces a normal $2$-sphere
\[
    \Gamma_i
    =
    \partial D_i^3
    \subset
    \varpi^{-1}(y_i)\cap\mathcal U_K
\]
which encloses the balls containing the monopole clusters in that normal fibre and
satisfies
\[
    \Gamma_i
    \cap
    B_{m_i^{-1}}\bigl(C_{\Lambda_0,i}\bigr)
    =
    \varnothing.
\]
By Li's dichotomy \cite[Lemma~2.17]{li2025large} (see also Lemma~\ref{lem:Li_dichotomy_away_Cinfty} and the proof of Lemma~\ref{lemma:Z_subset_C_infty}), every zero of $\Phi_i$ belongs to $C_{\Lambda_0,i}$ for all sufficiently large $i$. Hence
\[
    \Gamma_i\cap Z_i=\varnothing.
\]
Because the normal ball $D_i^3$ is contractible, the restricted principal
$G$-bundle admits a lift
\[
    \widetilde P_i\longrightarrow D_i^3
\]
to a principal $\mathrm{SU}(2)$-bundle; when $G=\mathrm{SU}(2)$ we simply
take the original restriction. Through the isometric Lie algebra
identification, $(\nabla_i,\Phi_i)$ lifts to adjoint data
$(\widetilde\nabla_i,\widetilde\Phi_i)$ on $\widetilde P_i$. Since
$\Gamma_i\cap Z_i=\varnothing$, the normalized lifted Higgs field defines
an eigenline bundle
\[
    \widetilde L_i\longrightarrow\Gamma_i
\]
in the standard complex rank-$2$ representation. Li's Chern form computed
from the lift agrees on $\Gamma_i$ with the globally defined
Lie algebra form $F_{\mathrm U(1),i}$. Step~4 of the proof of
\cite[Theorem~5.5]{li2025large} chooses $\Gamma_i$ so that
\[
    \left\langle
        c_1(\widetilde L_i),
        [\Gamma_i]
    \right\rangle
    =
    \int_{\Gamma_i}F_{\mathrm U(1),i}
    =
    \theta_{T_{\mathrm{Li}}}(y_i)
    =
    q
    \neq0.
\]
Here $\Gamma_i$ is oriented as the boundary of $D_i^3$, with normal
orientation compatible with the calibrated orientation of
$T_{\mathrm{Li}}$ and the ambient orientation.

Suppose, for contradiction, that $\Phi_i$ were nonzero throughout
$D_i^3$. Then the normalized lifted Higgs field would define an eigenline
bundle over the whole ball whose restriction to $\Gamma_i$ is
$\widetilde L_i$. Since $H^2(D_i^3;\mathbb Z)=0$, this restriction would
have zero first Chern number, contradicting the displayed identity.
Therefore there exists $z_i\in D_i^3$ such that $\Phi_i(z_i)=0$.
Since
\[
    D_i^3\subset\mathcal U_K\subset V,
\]
we conclude that
\[
    V\cap Z_i\neq\varnothing
\]
for every sufficiently large $i$.

This proves the first assertion. We next derive the uniform one-sided
Hausdorff convergence. Let $\varepsilon>0$. Since $\mathcal S$ is compact, choose
points $x_1,\ldots,x_N\in \mathcal S$ such that
\[
    \mathcal S
    \subset
    \bigcup_{a=1}^N
    B_{\varepsilon/2}(x_a).
\]
By the first part of the proof, for every $a$ there exists $i_a$ such
that
\[
    B_{\varepsilon/2}(x_a)\cap Z_i
    \neq
    \varnothing
\]
whenever $i\geqslant i_a$. Set
\[
    i_\varepsilon
    :=
    \max_{1\leqslant a\leqslant N}i_a.
\]
For $i\geqslant i_\varepsilon$ and $x\in \mathcal S$, choose $a$ with
$x\in B_{\varepsilon/2}(x_a)$ and then choose
$z_{i,a}\in B_{\varepsilon/2}(x_a)\cap Z_i$. It follows that
\[
    d(x,Z_i)
    \leqslant
    d(x,z_{i,a})
    <
    \varepsilon.
\]
Thus
\[
    \sup_{x\in \mathcal S}d(x,Z_i)
    \longrightarrow
    0.
\]
Finally, let $x\in \mathcal S$ and $N\in\mathbb N$. Every neighbourhood of $x$
meets $Z_i$ for some $i\geqslant N$, and therefore
\[
    x
    \in
    \overline{\bigcup_{i\geqslant N}Z_i}.
\]
Since this holds for every $N$,
\[
    x
    \in
    \bigcap_{N\geqslant1}
    \overline{\bigcup_{i\geqslant N}Z_i}
    =
    \mathcal Z.
\]
Hence $\mathcal S\subset \mathcal Z$.
\end{proof}

\begin{remark}[Transverse charge and possible Higgs zeros away from $\mathcal S$]
\label{rmk:charged_and_neutral_higgs_degenerations}
Proposition~\ref{prop:S_subset_Z} gives more than the set-theoretic
inclusion $\mathcal S\subset \mathcal Z$. The zero sets satisfy the one-sided
Hausdorff estimate
\[
    \sup_{x\in \mathcal S}d(x,Z_i)
    \longrightarrow
    0.
\]
Moreover, the zeros used in the proof are forced by a transverse
topological mechanism. On a carefully chosen $y_i$ over the smooth locus of
$T=T_{\mathrm{Li}}$, a normal linking sphere $\Gamma_i$ lies in a normal
$3$-ball over which the bundle has a local $\mathrm{SU}(2)$ lift. If
$\widetilde L_i\to\Gamma_i$ denotes the corresponding lifted eigenline,
then
\[
    \left\langle
        c_1(\widetilde L_i),
        [\Gamma_i]
    \right\rangle
    =
    \int_{\Gamma_i}F_{\mathrm U(1),i}
    =
    \theta_T(y_i).
\]
The precise relation between the longitudinal form normalized by $|\Phi_i|$ and the
Chern curvature is given by
\eqref{eq:eigenline_Chern_curvature}. In the $\mathrm{SO}(3)$ case the
global Chern form class may be half-integral, but its restriction to this
locally lifted linking sphere is the integral class displayed above.
Consequently, the normal $3$-ball bounded by $\Gamma_i$ must contain a zero
of $\Phi_i$. Since the exceptional fibres
have arbitrarily small $\mathcal H^{n-3}$-measure, such normal balls with nonzero transverse charge
occur arbitrarily close to every point of $\mathcal S$, for every
sufficiently large index.

The charge in this statement belongs to the outer linking sphere and
records the total transverse charge enclosed by it. We do not claim that
a transverse Bogomolny bubble is based at every point of $\mathcal S$, nor that a
single bubble has charge equal to the local multiplicity. The total
multiplicity may be distributed among several monopole clusters inside
the normal fibre. Accordingly, the conclusion is that every point of $\mathcal S$
is approached by Higgs zeros forced by nonzero transverse total charge,
not that each individual zero carries a canonically defined charge equal
to $\theta_T$.

The reverse inclusion $\mathcal Z\subset\mathcal S$ is not automatic in
higher dimensions. Since
$\mathcal S=\operatorname{spt}\|T_{\mathrm{Li}}\|$, Li's limiting
singular abelian monopole is smooth on $X\setminus\mathcal S$, whereas
its curvature satisfies
\[
    dF_\infty
    =
    2\pi T_{\mathrm{Li}}
    =
    2\pi T
\]
on $X$ as currents; see
\cite[Theorem~1.12]{li2025large}. Thus $\mathcal S$ is the support of the limiting codimension-three Dirac source. On $\mathcal R_{\mathrm{ab}}$, whenever the first
alternative of Corollary~\ref{cor:harmonic_correction_alternatives} holds, the
relation between this singular abelian monopole and the smooth limit of
the translated configurations is given in Remark~\ref{rmk:comparison_with_Li_singular_limit}. A point of
$\mathcal Z\setminus\mathcal S$ is approached by Higgs zeros of the
sequence but lies in the smooth, source-free region of the limiting
abelian monopole. Such a point carries no codimension-three charge detected by the
limiting current $T$. We do not introduce separate terminology for these
points; throughout the paper they are denoted simply by
$\mathcal Z\setminus\mathcal S$.

For a $\Theta$-monopole, the relation between the PPS charge current and
the intermediate energy measure is exact already at the level of the
sequence. Indeed, for every $\chi\in C_c^\infty(X)$,
\[
    \left\langle
        m_i^{-1}Z(\nabla_i,\Phi_i),
        \chi\Theta
    \right\rangle
    =
    2m_i^{-1}
    \int_X
    \chi|\nabla_i\Phi_i|^2\,\vol
    =
    \int_X\chi\,d\nu_i.
\]
Since
$m_i^{-1}Z(\nabla_i,\Phi_i)\rightharpoonup^*8\pi T$,
$\nu_i\rightharpoonup^*8\pi\|T\|$, and $T$ is
$\Theta$-calibrated, passage to the limit gives
\[
    8\pi T(\chi\Theta)
    =
    8\pi\int_X\chi\,d\|T\|.
\]
Thus the pairing of the charge current with the calibration has no
cancellation at the mass-renormalized scale. This exact identity does not,
however, imply $\mathcal Z\subset\mathcal S$. A zero-centred Green
identity and a shrinking linking sphere involve, respectively, an
$i$-dependent singular kernel and an $i$-dependent shrinking cycle,
whereas weak current convergence is tested against fixed smooth forms. In
particular, a zero detected on one transverse fibre need not contribute
positive $\|T\|$-mass to any fixed neighbourhood of its limiting point.

Proposition~\ref{prop:zero_centered_TU_concentration} gives an additional
analytic distinction between the two possible differences in the
hierarchy
\[
    \mathcal S
    \subset
    \mathcal Z
    \subset
    \mathcal C.
\]
Every point of $\mathcal Z\setminus\mathcal S$ is, after passing to a subsequence realizing it by Higgs zeros, a concentration point for the ordinary codimension-four scale-invariant Yang--Mills--Higgs energy. Thus $\mathcal Z\setminus\mathcal S$ is the natural locus in which
instanton bubbling centred at Higgs zeros and occurring at scales
$o(m_i^{-1})$ could remain invisible to the limiting mass-renormalized measure.
By contrast, if $x\in\mathcal C\setminus\mathcal Z$, then a fixed
neighbourhood of $x$ contains no Higgs zeros for all sufficiently large
indices, although abelianization still fails there. Proposition~\ref{prop:zero_centered_TU_concentration} does not provide the scale separation needed for this instanton interpretation; the remaining problem is discussed in Section~\ref{sec:resulting_picture}.

A converse characterization of $\mathcal S$ in terms of transverse
topology would therefore require quantitative control of the transverse
degree. The relevant invariant on an outer normal sphere is the first
Chern number of the eigenline after a local $\mathrm{SU}(2)$ lift, or
equivalently half the Euler number of the induced $\mathrm{SO}(2)$
reduction in the $\mathrm{SO}(3)$ description. By
\eqref{eq:eigenline_Chern_curvature}, this integer is not, without the
correction term, the flux of the raw longitudinal form
$\langle F_{\nabla_i},\Phi_i/|\Phi_i|\rangle$ alone. A suitable
quantitative slicing principle would have to produce, at each sufficiently
small fixed scale $r$ around a candidate point, charged normal fibres over
a subset of the base whose $\mathcal H^{n-3}$-measure is bounded below by
a fixed multiple of $r^{n-3}$. The topological energy lower bound on each
such fibre, integrated over this family and combined with the exact
calibrated pairing above, would give positive mass-renormalized
codimension-three density and force the point to lie in $\mathcal S$.
The corresponding questions are discussed in
Section~\ref{sec:resulting_picture}.
\end{remark}

The next subsection isolates the analytic obstruction to promoting local
finite mass concentration into effective codimension-three
monotonicity. In Section~\ref{sec: Li_concentration}, Li's
curvature concentration loci will then be used to place the limiting
zero set inside the limiting nonabelian locus and to show
that $\mathcal Z\setminus\mathcal S$ is contained in the set where
effective codimension-three monotonicity fails locally.

\subsection{Obstruction to codimension-three monotonicity}
\label{subsec:codim3_monotonicity_obstruction}

By the results of \S\ref{subsec: PPS_calibrated}, after passing to a
subsequence, the mass-renormalized Yang--Mills--Higgs energy concentrates
along a $\Theta$-calibrated integral $(n-3)$-cycle $T$:
\[
    \mu=\nu=8\pi\|T\|.
\]
As recalled in \eqref{ineq: exact_codim_3_monot_for_mu}, the limiting
measure $\mu$ satisfies codimension-three monotonicity. This should be
contrasted with the usual codimension-four Yang--Mills--Higgs
monotonicity formula. More precisely, multiplying the standard formula
\cite[Theorem~2.1]{afuni2019regularity} by $m_i^{-1}$ gives
\[
    e^{cs^2}s^{4-n}\mu_i(B_s(x))
    \leqslant
    e^{cr^2}r^{4-n}\mu_i(B_r(x))
\]
for every $x\in X$ and every $0<s<r\leqslant r_0$. After increasing $c$
and decreasing $r_0$ if necessary, we use the same constants throughout
this subsection.

A natural question is whether an effective codimension-three
almost-monotonicity estimate, with an error uniform in the mass, already
holds for the approximating sequence. The relevant identity is obtained
from the stress-energy formula. In the {\mainsettingref}, the bounded geometry
of $(X,g)$ allows all geometric constants to be chosen independently of
$i$ and of the center.

There exist constants $r_0>0$ and $c\geqslant0$, depending only on the
fixed background geometry and on $n$, such that, for every $i$, every
$x\in X$, and every $0<r\leqslant r_0$,
\begin{align}
\frac{d}{dr}
\left(
    e^{cr^2}\theta_i(x,r)
\right)
&\geqslant
2m_i^{-1}e^{cr^2}r^{3-n}
\int_{\partial B_r(x)}
\left(
    |\partial_{r_x}\lrcorner F_{\nabla_i}|^2
    +
    |\nabla_{\partial_{r_x}}\Phi_i|^2
\right)
\nonumber\\
&\quad
-
e^{cr^2}r^{2-n}
\int_{B_r(x)}
m_i^{-1}
\left(
    |F_{\nabla_i}|^2
    -
    |\nabla_i\Phi_i|^2
\right),
\label{eq:almost_monotonicity_original}
\end{align}
where $\partial_{r_x}$ is the outward unit normal to $\partial B_r(x)$. This is a
convenient reformulation of
\cite[Theorem~2.1]{afuni2019regularity}; a proof is included in
Appendix~\ref{app: B}.

Set
\[
    q_i
    :=
    |F_{\nabla_i}|^2
    -
    |\nabla_i\Phi_i|^2.
\]
Since the boundary flux in
\eqref{eq:almost_monotonicity_original} is nonnegative, we obtain
\begin{equation}
\label{eq:differential_almost_monotonicity_short}
\frac{d}{dr}
\left(
    e^{cr^2}\theta_i(x,r)
\right)
\geqslant
-
e^{cr^2}m_i^{-1}r^{2-n}
\int_{B_r(x)}q_i
\geqslant
-
e^{cr^2}m_i^{-1}r^{2-n}
\int_{B_r(x)}q_i^+ ,
\end{equation}
where $q_i^+=\max \{ 0, q_i \}$ denotes the nonnegative part of $q_i$. Thus, of the \textbf{signed codimension-two scale imbalance} 
\begin{equation}
\label{eq: imbalance_term}
    m_i^{-1}r^{2-n}
    \int_{B_r(x)}q_i,
\end{equation} 
only its positive part can decrease the weighted codimension-three
energy ratio.

More generally, suppose that for some nonnegative function
$h_x\in L^1((0,r_x))$ one has
\[
    m_i^{-1}r^{2-n}
    \int_{B_r(x)}q_i^+
    \leqslant
    h_x(r)
\]
uniformly in $i$ for $0<r<r_x$. Integrating
\eqref{eq:differential_almost_monotonicity_short} gives
\[
    e^{cs^2}\theta_i(x,s)
    \leqslant
    e^{cr^2}\theta_i(x,r)
    +
    \int_s^r e^{ct^2}h_x(t)\,dt,
    \qquad
    0<s<r<r_x.
\]
Hence a uniform codimension-two Morrey bound for the \textbf{positive imbalance} $q_i^{+}$
produces an effective codimension-three monotonicity inequality.
The difficulty is not the mere presence of the signed imbalance, but the
possible absence of such a quantitative bound for its positive part.

Li's theorem \cite[Theorem~1.16]{li2025large} implies that the
mass-renormalized curvature and Higgs gradient measures converge to the
same limiting measure. Consequently,
\[
    m_i^{-1}q_i\,\operatorname{vol}
    \rightharpoonup
    0
\]
weakly as signed Radon measures. For every $x\in X$ and every radius
$r>0$ such that
\[
    \mu(\partial B_r(x))=0,
\]
it follows that
\[
    m_i^{-1}
    \int_{B_r(x)}q_i
    \longrightarrow
    0.
\]
For each fixed center $x$, the set of radii which are not continuity
radii of $\mu$ is at most countable.

This convergence is only a fixed scale statement. It neither controls
the positive variation $q_i^+$ nor gives a uniform estimate for
\[
    m_i^{-1}r^{2-n}
    \int_{B_r(x)}q_i^+
\]
which is integrable as $r\downarrow0$. Large positive and negative
contributions to $q_i$ may cancel in the signed integral. Thus the equality
of the limiting curvature and Higgs gradient measures does not, by
itself, yield an effective codimension-three monotonicity estimate for
the approximating sequence.

For $\Theta$-monopoles, the imbalance has a more rigid geometric
interpretation. Using the algebraic identity
\eqref{eq:general_calibration_identity},
\[
    |F_{\nabla_i}|^2\operatorname{vol}
    =
    -
    \left\langle
        F_{\nabla_i}\wedge F_{\nabla_i}
    \right\rangle
    \wedge\Xi
    +
    |F_{\nabla_i}\wedge\Theta|^2
    \operatorname{vol},
\]
together with the monopole equation
\[
    F_{\nabla_i}\wedge\Theta
    =
    *
    \nabla_i\Phi_i,
\]
we obtain
\begin{equation}
\label{eq: imbalance_chern_weil}
    q_i\,\operatorname{vol}
    =
    -
    \left\langle
        F_{\nabla_i}\wedge F_{\nabla_i}
    \right\rangle
    \wedge\Xi.
\end{equation}
Thus the signed imbalance is represented by a Chern--Weil density
associated with the special holonomy structure.

We next explain why the closedness of $\Xi$ does not immediately provide
the required codimension-two Morrey control. Since $(X,g)$ has bounded geometry
and $\Xi$ is parallel, there exist uniform constants $r_1>0$ and
$C<\infty$, with $r_1$ smaller than the uniform convexity radius, such
that for every $x\in X$ there is an $(n-5)$-form
\[
    \tau_x\in\Omega^{n-5}(B_{r_1}(x))
\]
satisfying
\[
    d\tau_x=\Xi,
    \qquad
    |\tau_x(y)|
    \leqslant
    C\,d_g(x,y).
\]
Such a primitive is obtained from the standard radial homotopy operator
in geodesic normal coordinates.

By the Bianchi identity and the invariance of the fibrewise inner product,
\[
    d
    \left\langle
        F_{\nabla_i}\wedge F_{\nabla_i}
    \right\rangle
    =
    0.
\]
Consequently, for every $0<r\leqslant r_1$, Stokes' theorem and
\eqref{eq: imbalance_chern_weil} give the exact local transgression
identity
\begin{equation}
\label{eq: local_chern_weil_boundary}
    \int_{B_r(x)}
    q_i\,\operatorname{vol} = -\int_{\partial B_r(x)} \left\langle F_{\nabla_i}\wedge F_{\nabla_i}\right\rangle\wedge\tau_x.
\end{equation}
The pullback of
$\langle F_{\nabla_i}\wedge F_{\nabla_i}\rangle$
to $\partial B_r(x)$ depends only on the tangential component of the
curvature. Hence the right-hand side of
\eqref{eq: local_chern_weil_boundary} is a tangential Chern--Weil
boundary term. By contrast, the nonnegative boundary term in
\eqref{eq:almost_monotonicity_original} controls the radial quantities
\[
    |{\partial_{r_x}}\lrcorner F_{\nabla_i}|^2
    +
    |\nabla_{\partial_{r_x}}\Phi_i|^2.
\]
The transgression identity therefore does not, by itself, estimate the
positive contribution of the Chern--Weil boundary term by the radial
stress-energy flux. Additional quantitative control of the tangential
curvature would be required.

Li's theorem shows that the corresponding bulk Chern--Weil measures
\[
    m_i^{-1}
    \left\langle
        F_{\nabla_i}\wedge F_{\nabla_i}
    \right\rangle
    \wedge\Xi
\]
converge weakly to zero. This weak bulk convergence gives no uniform
control of the boundary terms in
\eqref{eq: local_chern_weil_boundary} across shrinking spheres. In
particular, it does not provide the required codimension-two Morrey bound for the positive
imbalance in \eqref{eq:differential_almost_monotonicity_short}.

We therefore isolate the additional quantitative property needed in the
later analysis.

\begin{definition}[Effective codimension-three monotonicity on an open subset]
\label{def:effective_codim3_monotonicity}
Let $U\subset X$ be open. We say that the sequence satisfies
\textbf{effective codimension-three monotonicity on $U$} if there exist
$A_U\in(0,\infty)$, $i_U\in\mathbb N$, and $\gamma_U\in(0,\infty)$ such that,
for every $i\geqslant i_U$, every geodesic ball $B_r(x)\subset U$, and
every $0<s<r\leqslant r_0$, one has
\begin{equation}
\label{ineq:effective_codim3_monotonicity}
    e^{cs^2}\theta_i(x,s)
    \leqslant
    e^{cr^2}\theta_i(x,r)
    +
    A_Ur^{\gamma_U}.
\end{equation}
Here $c\geqslant0$ and $r_0>0$ are the geometric constants fixed above.
\end{definition}

\begin{definition}[Obstruction to effective codimension-three monotonicity]
\label{def:codim3_monotonicity_obstruction}
Let $\mathcal O$ denote the
\textbf{obstruction locus for effective codimension-three monotonicity},
namely the set of points $x\in X$ which admit no open neighbourhood on
which effective codimension-three monotonicity holds. Equivalently,
\[
    X\setminus\mathcal O
    =
    \bigcup
    \left\{
        U\subset X:
        \begin{array}{l}
        U\text{ is open and effective codimension-three}\\[-2pt]
        \text{monotonicity holds on }U
        \end{array}
    \right\}.
\]
In particular, $X\setminus\mathcal O$ is open and $\mathcal O$ is
closed.
\end{definition}

\begin{lemma}[Control on balls with nearby centres]
\label{lem:effective_monotonicity_nearby_centres}
Suppose that effective codimension-three monotonicity holds on an open set
$U$. Then there is a geometric constant
\[
    C_\ast:=2^{n-3}e^{cr_0^2/4}
\]
such that, for every $i\geqslant i_U$, every ball $B_r(x)\subset U$, every
$y\in B_{r/2}(x)$, and every $0<s\leqslant r/2$, one has
\begin{equation}
\label{ineq:effective_codim3_nearby_centres}
    s^{3-n}m_i^{-1}\int_{B_s(y)}e_i
    \leqslant
    C_\ast r^{3-n}m_i^{-1}\int_{B_r(x)}e_i
    +
    2^{-\gamma_U}A_Ur^{\gamma_U}.
\end{equation}
\end{lemma}

\begin{proof}
Since $y\in B_{r/2}(x)$, one has
$B_{r/2}(y)\subset B_r(x)\subset U$. Applying
\eqref{ineq:effective_codim3_monotonicity} with centre $y$ and outer
radius $r/2$ gives
\[
    e^{cs^2}\theta_i(y,s)
    \leqslant
    e^{cr^2/4}\theta_i(y,r/2)
    +
    A_U(r/2)^{\gamma_U}.
\]
Moreover,
\[
    \theta_i(y,r/2)
    \leqslant
    2^{n-3}r^{3-n}m_i^{-1}\int_{B_r(x)}e_i.
\]
Since $e^{cs^2}\geqslant1$ and $r\leqslant r_0$, the claim follows.
\end{proof}

\begin{remark}[Local and global effective monotonicity]
\label{rmk:local_global_effective_monotonicity}
If effective codimension-three monotonicity holds on $X$, then
$\mathcal O=\varnothing$. Conversely, suppose that
$\mathcal O=\varnothing$ and let $K\subset X$ be compact. Choose finitely
many open sets $U_1,\ldots,U_N$ on which effective monotonicity holds and
which cover $K$. After decreasing the radius, one may choose
$r_K\in(0,\min\{r_0,1\}]$ such that every ball $B_r(x)$ with
$x\in K$ and $0<r\leqslant r_K$ is contained in one of the $U_j$. Set
\[
    i_K:=\max_j i_{U_j},
    \qquad
    A_K:=\max_j A_{U_j},
    \qquad
    \gamma_K:=\min_j\gamma_{U_j}.
\]
Since $r_K\leqslant1$, one has
$A_{U_j}r^{\gamma_{U_j}}\leqslant A_Kr^{\gamma_K}$ for
$0<r\leqslant r_K$. Hence
\eqref{ineq:effective_codim3_monotonicity} holds for every
$i\geqslant i_K$, every $x\in K$, and every
$0<s<r\leqslant r_K$.

On the noncompact AC manifold, these constants may degenerate as the
compact set exhausts $X$. Thus $\mathcal O=\varnothing$ expresses local
effective codimension-three monotonicity everywhere, but does not by itself
imply an estimate on all of $X$ with constants uniform in the centre.
\end{remark}

For later use, we record the consequence of
Definition~\ref{def:effective_codim3_monotonicity} away from the
calibrated support $\mathcal S$. For $x\in X$ and $r>0$, set
\[
    \mathfrak M_i(x,r)
    :=
    \sup_{\substack{B_s(y)\subset B_r(x)\\0<s\leqslant r}}
    s^{3-n}m_i^{-1}
    \int_{B_s(y)}e_i.
\]

\begin{lemma}[Effective monotonicity gives vanishing Morrey density off $\mathcal S$]
\label{lem:effective_monotonicity_gives_vanishing_Morrey_off_S}
If
\[
    x\in X\setminus(\mathcal S\cup\mathcal O),
\]
then
\[
    \lim_{r\downarrow0}
    \limsup_{i\to\infty}
    \mathfrak M_i(x,r)
    =
    0.
\]
\end{lemma}

\begin{proof}
Choose an open set $U\ni x$ on which effective codimension-three
monotonicity holds. Since $x\notin\mathcal S=\operatorname{spt}\mu$, after
shrinking $U$ choose $r_x>0$ such that
\[
    \overline{B_{r_x}(x)}\subset U,
    \qquad
    \overline{B_{r_x}(x)}\cap\mathcal S=\varnothing.
\]
For every $0<r<r_x$, the Portmanteau theorem gives
\[
    r^{3-n}m_i^{-1}\int_{B_r(x)}e_i\longrightarrow0.
\]
By Lemma~\ref{lem:effective_monotonicity_nearby_centres},
\[
    \limsup_{i\to\infty}\mathfrak M_i(x,r/2)
    \leqslant
    2^{-\gamma_U}A_Ur^{\gamma_U}.
\]
Letting $r\downarrow0$ proves the claim.
\end{proof}

\begin{remark}[The borderline case]
\label{rmk:borderline_gamma_zero}
The strict positivity of $\gamma_U$ in
Definition~\ref{def:effective_codim3_monotonicity} is essential in
Lemma~\ref{lem:effective_monotonicity_gives_vanishing_Morrey_off_S}. A
borderline estimate with $\gamma_U=0$ would give only a finite uniform
bound for the codimension-three Morrey quantities on balls with nearby
centres. It would not force
\[
    \lim_{r\downarrow0}
    \limsup_{i\to\infty}
    \mathfrak M_i(x,r)
    =
    0.
\]
Away from $\mathcal S$, this vanishing is precisely the smallness needed
below to force Li's weighted curvature potential below its defining
threshold; see
Proposition~\ref{prop:Cinfty_minus_S_obstructs_positive_codim3}.
\end{remark}

We close the subsection with a sufficient condition for effective
codimension-three monotonicity. In view of
\eqref{eq:differential_almost_monotonicity_short}, it is enough to control
the positive part $q_i^+$ at codimension-two Morrey scale. The coarse
pointwise inequality $q_i^+\leqslant|F_{\nabla_i}|^2$ then gives a
corresponding curvature criterion.

\begin{proposition}[Codimension-two Morrey criterion for effective monotonicity]
\label{prop:codim2_Morrey_criterion_positive_monotonicity}
Let $U\subset X$ be open. Suppose that there exist
$\Lambda_U<\infty$ and $i_U\in\mathbb N$ such that, for every
$i\geqslant i_U$ and every geodesic ball $B_s(y)\subset U$ with
$0<s\leqslant r_0$, one has
\begin{equation}
\label{ineq:positive_imbalance_Morrey_bound}
    m_i^{-1}s^{2-n}
    \int_{B_s(y)}q_i^+
    \leqslant
    \Lambda_U.
\end{equation}
Then effective codimension-three monotonicity holds on $U$, with
$\gamma_U=1$. The same conclusion holds under the stronger curvature
Morrey bound
\begin{equation}
\label{ineq:curvature_Morrey_bound_coarse}
    m_i^{-1}s^{2-n}
    \int_{B_s(y)}|F_{\nabla_i}|^2
    \leqslant
    \Lambda_U.
\end{equation}
\end{proposition}

\begin{proof}
Fix $i\geqslant i_U$ and a ball $B_r(x)\subset U$, with
$0<r\leqslant r_0$. Applying
\eqref{eq:differential_almost_monotonicity_short} with centre $x$ gives,
for almost every $0<t<r$,
\[
    \frac{d}{dt}
    \left(e^{ct^2}\theta_i(x,t)\right)
    \geqslant
    -e^{ct^2}\Lambda_U.
\]
Integrating from $s$ to $r$ yields
\begin{equation}
\label{ineq:integrated_codim3_monotonicity_from_qplus}
    e^{cs^2}\theta_i(x,s)
    \leqslant
    e^{cr^2}\theta_i(x,r)
    +
    e^{cr_0^2}\Lambda_Ur
\end{equation}
for every $0<s<r\leqslant r_0$. This is
\eqref{ineq:effective_codim3_monotonicity} with
$A_U=e^{cr_0^2}\Lambda_U$ and $\gamma_U=1$.
The curvature criterion follows from
$q_i^+\leqslant|F_{\nabla_i}|^2$.
\end{proof}

\begin{corollary}[Local form of the Morrey criterion]
\label{cor:local_codim2_Morrey_criterion_positive_monotonicity}
Let $x\in X$. Suppose there exist
$r_x>0$, $\Lambda_x<\infty$, and $i_x\in\mathbb N$ such that, for every
$i\geqslant i_x$, every $y\in B_{r_x}(x)$, and every $0<s<r_x$, one has
\[
    m_i^{-1}s^{2-n}\int_{B_s(y)}q_i^+\leqslant\Lambda_x.
\]
Then effective codimension-three monotonicity holds on
$B_{r_x/2}(x)$, and hence $x\notin\mathcal O$. The same conclusion
holds if $q_i^+$ is replaced by $|F_{\nabla_i}|^2$.
\end{corollary}

\begin{proof}
Every ball $B_s(y)\subset B_{r_x/2}(x)$ has
$y\in B_{r_x}(x)$ and $s<r_x$. Apply
Proposition~\ref{prop:codim2_Morrey_criterion_positive_monotonicity} with
$U=B_{r_x/2}(x)$.
\end{proof}


\section{Li's curvature concentration loci and the limiting nonabelian locus}
\label{sec: Li_concentration}

Li's curvature concentration loci depend on the index and need not
converge as sets. In this section, we pass to their Kuratowski upper limit,
which we call the limiting nonabelian locus of the sequence, compare it
with the calibrated support $\mathcal S$ and the limiting Higgs zero set
$\mathcal Z$, and identify its excess over $\mathcal S$ by
proving
\[
    \mathcal C\setminus\mathcal S
    =
    \mathcal O\setminus\mathcal S.
\]
We also derive Hausdorff content and density rate estimates for the
loci $C_{\Lambda_0,i}$ and Higgs zero sets away from $\mathcal S$, describe excess
points using Li's characteristic radii, compare those radii with the
ordinary Yang--Mills--Higgs regularity scales, and analyze the weighted
curvature potential when the ratio of the regularity scale to
$m_i^{-1}$ tends to infinity. We then establish the estimates needed
for the abelian analysis on the complementary open $\mathcal R_{\mathrm{ab}}$.

\subsection{Curvature concentration loci and the limiting nonabelian locus}
\label{subsec:Li_moving_loci}
We now give the promised systematic recall of the curvature concentration
loci introduced by Li for each configuration in the sequence. These are
the loci $C_{\Lambda_0,i}$ used as a forward reference in the proof of
Proposition~\ref{prop:S_subset_Z}. By
\eqref{eq:PPS_Li_norm_mass_comparison} and
\eqref{eq:PPS_Li_mass_comparison}, the prefactor in Li's weighted
curvature potential is unchanged when written in our convention:
\[
    \bigl(m_i^{\mathrm{Li}}\bigr)^{-2}
    |F_{\nabla_i}|_{\mathrm{Li}}^2
    =
    m_i^{-2}|F_{\nabla_i}|^2.
\]
The mass scale cutoff, however, becomes
\[
    \bigl(m_i^{\mathrm{Li}}\bigr)^{-1}
    =
    2m_i^{-1}.
\]
Consequently, Li's potential, written in our norms and masses, is
\begin{equation}
\label{eq:exact_translated_Li_potential}
    \mathcal P_i^{\mathrm{Li}}(p)
    =
    m_i^{-2}
    \int_X
    \frac{|F_{\nabla_i}|^2(x)}
    {\max\{d(x,p),2m_i^{-1}\}^{n-2}}
    \,{\vol}(x).
\end{equation}
For notational convenience, we use the quantitatively equivalent potential
\begin{equation}
\label{eq:Li_potential_our_convention}
    \mathcal P_i(p)
    :=
    m_i^{-2}
    \int_X
    \frac{|F_{\nabla_i}|^2(x)}
    {\max\{d(x,p),m_i^{-1}\}^{n-2}}
    \,{\vol}(x).
\end{equation}
The two potentials satisfy the precise comparison
\begin{equation}
\label{eq:Li_potential_convention_comparison}
    \mathcal P_i^{\mathrm{Li}}(p)
    \leqslant
    \mathcal P_i(p)
    \leqslant
    2^{n-2}\mathcal P_i^{\mathrm{Li}}(p)
    \qquad
    \text{for every }p\in X.
\end{equation}
Accordingly, if $C^{\mathrm{Li}}_{\Lambda,i}$ denotes the locus defined
from \eqref{eq:exact_translated_Li_potential}, then
\begin{equation}
\label{eq:Li_loci_convention_comparison}
    C^{\mathrm{Li}}_{\Lambda,i}
    \subset
    C_{\Lambda,i}
    \subset
    C^{\mathrm{Li}}_{2^{n-2}\Lambda,i},
\end{equation}
where
\[
    C_{\Lambda,i}
    :=
    \left\{
        p\in X:
        \mathcal P_i(p)>\Lambda^{-1}
    \right\}.
\]
We therefore fix once and for all a constant $\Lambda_0>0$, enlarged by
the universal factor $2^{n-2}$ if necessary, so that all conclusions of
\cite[Lemma~2.17 and Proposition~2.18]{li2025large} hold for the loci
$C_{\Lambda_0,i}$ defined above. We call \eqref{eq:Li_potential_our_convention}
\textbf{Li's weighted curvature potential} and refer to the sets
$C_{\Lambda_0,i}$ as Li's \textbf{curvature concentration loci}.

The sets $C_{\Lambda_0,i}$ capture the regions where the sequence may remain genuinely nonabelian. After the same fixed conversion of the mass scale radii and a relabeling of the universal constant $\varepsilon_3>0$, Li's concentration-decay
dichotomy \cite[Lemma~2.17]{li2025large} takes the following form in our convention\footnote{All constants
below are understood in this fixed convention.}. For all sufficiently
large $i$, if either
\[
    \int_{B(p,2m_i^{-1})}|F_{\nabla_i}|^2
    \geqslant
    \varepsilon_3 m_i^{4-n},
\]
or
\[
    \min_{B(p,m_i^{-1})}|\Phi_i|
    \leqslant
    \frac{m_i}{2},
\]
then $p\in C_{\Lambda_0,i}$. Equivalently, outside $C_{\Lambda_0,i}$ the corresponding alternatives of Li are excluded, the Higgs field remains uniformly large at the mass scale, and the transverse curvature components satisfy the exponential decay estimates of \cite[Proposition~2.18]{li2025large}. 

We now introduce a fixed closed set associated with the sequence whose complement is the maximal open set on which the curvature concentration loci $C_{\Lambda_0,i}$ are eventually absent locally.

\begin{definition}[Limiting nonabelian locus and abelian region]
\label{def:C_infty}
Define
\[
\mathcal C
:=
\bigcap_{N\geqslant1}
\overline{
\bigcup_{i\geqslant N}
C_{\Lambda_0,i}
},
\qquad
\mathcal R_{\mathrm{ab}}
:=
X\setminus\mathcal C.
\]
Equivalently, $x\in\mathcal C$ if and only if there are indices
$i_j\to\infty$ and points $p_j\in C_{\Lambda_0,i_j}$ with $p_j\to x$.
Equivalently again, $x\in\mathcal R_{\mathrm{ab}}$ if and only if there
are an open neighbourhood $U$ of $x$ and an index $i_0$ such that
\[
U\cap C_{\Lambda_0,i}
=
\varnothing
\qquad
\text{for every } i\geqslant i_0.
\]
We call $\mathcal C$ the \textbf{limiting nonabelian locus} of the
sequence and $\mathcal R_{\mathrm{ab}}$ its \textbf{abelian region}.
\end{definition}
The set $\mathcal R_{\mathrm{ab}}$ is the natural domain for the local abelianization theory developed
below.

\begin{lemma}[Local characterization of $\mathcal R_{\mathrm{ab}}$]
\label{lem:maximal_abelian_clearing_region}
The set $\mathcal R_{\mathrm{ab}}$ is open. For every compact set
$K\Subset\mathcal R_{\mathrm{ab}}$, there exist an open set $U$ and an
index $i_K$ such that
\[
    K\Subset U\Subset\mathcal R_{\mathrm{ab}},
    \qquad
    U\cap C_{\Lambda_0,i}=\varnothing
    \quad
    \text{for every }i\geqslant i_K.
\]
Conversely, suppose that $V\subset X$ is open and that, for every compact
set $K\Subset V$, there exists an index $i_K\in\mathbb N$ such that
\[
    K\cap C_{\Lambda_0,i}=\varnothing
    \qquad
    \text{for every }i\geqslant i_K.
\]
Then $V\subset\mathcal R_{\mathrm{ab}}$. Thus
$\mathcal R_{\mathrm{ab}}$ is the maximal open subset of $X$ on which
Li's curvature concentration loci are eventually absent locally.
\end{lemma}

\begin{proof}
Openness follows from the closedness of $\mathcal C$. Let
$K\Subset\mathcal R_{\mathrm{ab}}$. For every $x\in K$,
Definition~\ref{def:C_infty} gives an open neighbourhood $U_x$ and an
index $i_x$ such that
\[
    U_x\cap C_{\Lambda_0,i}=\varnothing
    \qquad
    \text{for every }i\geqslant i_x.
\]
For each $x\in K$, choose an open neighbourhood $V_x$ such that
\[
    x\in V_x,
    \qquad
    \overline{V_x}\Subset U_x\cap\mathcal R_{\mathrm{ab}}.
\]
Compactness of $K$ gives a finite subcover
$V_{x_1},\ldots,V_{x_N}$. Their union
\[
    U:=\bigcup_{j=1}^N V_{x_j}
\]
satisfies
\[
    K\Subset U\Subset\mathcal R_{\mathrm{ab}}.
\]
Taking the maximum of the corresponding indices $i_{x_j}$ gives
\[
    U\cap C_{\Lambda_0,i}=\varnothing
    \qquad
    \text{for every }i\geqslant \max_j i_{x_j},
\]
and proves the first assertion.

Conversely, let $x\in V$ and choose a compact neighbourhood
$K\Subset V$ of $x$. By hypothesis, $K$ is disjoint from
$C_{\Lambda_0,i}$ for all sufficiently large $i$. Hence a neighbourhood
of $x$ is eventually disjoint from the loci $C_{\Lambda_0,i}$, so
$x\notin\mathcal C$ by Definition~\ref{def:C_infty}. Therefore
$V\subset\mathcal R_{\mathrm{ab}}$.
\end{proof}

\begin{remark}[Kuratowski upper limits and escape to infinity]
\label{rmk:Cinfty_escape_infinity}
The set $\mathcal C$ records precisely those points represented by convergent sequences from the
loci $C_{\Lambda_0,i}$. The estimates available under the {\mainsettingref} do not yield
a tightness estimate for the sequence $C_{\Lambda_0,i}$ that is uniform in $i$ and, in
particular, do not exclude sequences
\[
    p_i\in C_{\Lambda_0,i},
    \qquad
    \rho(p_i)\longrightarrow\infty.
\]
The large radius energy estimates of
Subsection~\ref{subsec: Li_curvature_estimate} control macroscopic energy
growth on the AC end, but do not by themselves rule out shrinking balls selected at the mass or characteristic scales whose
centres escape to infinity. Consequently, no compactness assertion for
$\mathcal C$, and no global Hausdorff convergence of the loci
$C_{\Lambda_0,i}$, is used below. All conclusions on $\mathcal R_{\mathrm{ab}}$ are
local on fixed compact subsets.
\end{remark}

We shall also use a quantitative consequence of Li's Vitali covering
argument. For $d>0$, write
\begin{equation}\label{eq:def_unrestricted_Hausdorff_content}
    \mathcal H^d_\infty(E)
    :=
    \inf\left\{
        \sum_\alpha(\operatorname{diam}U_\alpha)^d:
        E\subset\bigcup_\alpha U_\alpha
    \right\}
\end{equation}
for the unrestricted $d$-dimensional Hausdorff content; the harmless
normalizing dimensional constant is suppressed. We similarly write
$\mathcal H^d_\delta$ when the covering sets are required to have
diameter at most $\delta$. Thus, up to the same fixed normalization,
$\mathcal H^d(E)=\lim_{\delta\downarrow0}\mathcal H^d_\delta(E)$.

\begin{proposition}[Li covers and Hausdorff content decay away from $\mathcal S$]
\label{prop:Li_cover_Hausdorff_content_off_S}
Put $p:=n-3$. Let $K\Subset X\setminus\mathcal S$ and fix
$0<\kappa<1$. For every sufficiently large $i$, there is a finite
family of pairwise disjoint geodesic balls
\[
    \bigl\{B(p_{i,a},s_{i,a})\bigr\}_{a\in A_i}
\]
such that their fivefold dilations cover
$C_{\Lambda_0,i}\cap K$ and
\begin{equation}\label{eq:Li_cover_radius_bounds_off_S}
    m_i^{-1}
    \leqslant
    s_{i,a}
    \leqslant
    C_{K,\kappa}
    m_i^{-\frac{\kappa}{1+\kappa}}
    \qquad
    \text{for every }a\in A_i.
\end{equation}
Moreover,
\begin{equation}\label{eq:Li_cover_sums_off_S}
    \sum_{a\in A_i}s_{i,a}^{p+\kappa}
    =o(m_i^{-\kappa}),
    \qquad
    \sum_{a\in A_i}s_{i,a}^{p}=o(1),
    \qquad
    N_i:=\#A_i=o(m_i^p).
\end{equation}
Consequently,
\begin{equation}\label{eq:Li_loci_Hausdorff_content_off_S}
    \mathcal H^{p+\kappa}_\infty
    \bigl(C_{\Lambda_0,i}\cap K\bigr)
    =o(m_i^{-\kappa}),
    \qquad
    \mathcal H^p_\infty
    \bigl(C_{\Lambda_0,i}\cap K\bigr)
    \longrightarrow0.
\end{equation}
\end{proposition}

\begin{proof}
Choose bounded open sets with compact closure
\[
    K\Subset U\Subset V\Subset X\setminus\mathcal S.
\]
Let $C^{\mathrm{Li}}_{\Lambda,i}$ denote Li's locus formed using the
exact potential \eqref{eq:exact_translated_Li_potential}. By
\eqref{eq:Li_loci_convention_comparison},
\[
    C_{\Lambda_0,i}
    \subset
    C^{\mathrm{Li}}_{\Lambda_*,i},
    \qquad
    \Lambda_*:=2^{n-2}\Lambda_0.
\]
The parameter $\Lambda_*$ is fixed, and
$(m_i^{\mathrm{Li}})^\kappa\gg\Lambda_*$ for all sufficiently large
$i$. We may therefore apply \cite[Lemma~3.7]{li2025large} to
$C^{\mathrm{Li}}_{\Lambda_*,i}\cap U$. After converting the norm and
mass according to
\eqref{eq:PPS_Li_norm_mass_comparison}--%
\eqref{eq:PPS_Li_mass_comparison} and absorbing the resulting fixed
factors into the constants, Li's lemma gives pairwise disjoint balls
$B(p_{i,a},s_{i,a})$ whose fivefold dilations cover
$C_{\Lambda_0,i}\cap K$, with
\eqref{eq:Li_cover_radius_bounds_off_S}, and
\begin{equation}\label{eq:translated_Li_Vitali_energy_bound}
    m_i^{1+\kappa}
    \sum_{a\in A_i}s_{i,a}^{p+\kappa}
    \leqslant
    C_{\Lambda_0,\kappa}
    \int_{U_i'}|F_{\nabla_i}|^2.
\end{equation}
Here
\[
    U_i'=B_{\rho_i}(U),
    \qquad
    \rho_i
    \leqslant
    C_{\Lambda_0,\kappa}
    m_i^{-\frac{\kappa}{1+\kappa}},
\]
so $U_i'\Subset V$ for all sufficiently large $i$.

By Proposition~\ref{prop:PPS_saturation_and_Li_identification},
\[
    \mu_i
    =m_i^{-1}e_i\,\vol
    \rightharpoonup^*
    8\pi\|T\|,
    \qquad
    \operatorname{spt}\|T\|=\mathcal S.
\]
Since $\overline V\cap\mathcal S=\varnothing$, the Portmanteau theorem
applied to the closed set $\overline V$ gives
\[
\begin{aligned}
    0
    &\leqslant
    \limsup_{i\to\infty}
    m_i^{-1}\int_{U_i'}|F_{\nabla_i}|^2 \\
    &\leqslant
    \limsup_{i\to\infty}\mu_i(\overline V)
    \leqslant
    8\pi\|T\|(\overline V)
    =0.
\end{aligned}
\]
Thus
\[
    \int_{U_i'}|F_{\nabla_i}|^2=o(m_i).
\]
Combining this with
\eqref{eq:translated_Li_Vitali_energy_bound} yields
\[
    \sum_{a\in A_i}s_{i,a}^{p+\kappa}
    =o(m_i^{-\kappa}).
\]
Since $s_{i,a}\geqslant m_i^{-1}$,
\[
    \sum_{a\in A_i}s_{i,a}^{p}
    =
    \sum_{a\in A_i}
    s_{i,a}^{p+\kappa}s_{i,a}^{-\kappa}
    \leqslant
    m_i^\kappa
    \sum_{a\in A_i}s_{i,a}^{p+\kappa}
    =o(1).
\]
The same lower radius bound gives
\[
    N_i m_i^{-(p+\kappa)}
    \leqslant
    \sum_{a\in A_i}s_{i,a}^{p+\kappa}
    =o(m_i^{-\kappa}),
\]
and hence $N_i=o(m_i^p)$. This proves
\eqref{eq:Li_cover_sums_off_S}.

Finally, the balls $B(p_{i,a},5s_{i,a})$ cover
$C_{\Lambda_0,i}\cap K$, and therefore
\[
\begin{aligned}
    \mathcal H^{p+\kappa}_\infty
    \bigl(C_{\Lambda_0,i}\cap K\bigr)
    &\leqslant
    10^{p+\kappa}
    \sum_{a\in A_i}s_{i,a}^{p+\kappa},\\
    \mathcal H^{p}_\infty
    \bigl(C_{\Lambda_0,i}\cap K\bigr)
    &\leqslant
    10^{p}
    \sum_{a\in A_i}s_{i,a}^{p}.
\end{aligned}
\]
The asserted content estimates follow from
\eqref{eq:Li_cover_sums_off_S}.
\end{proof}

\begin{proposition}[Transfer to limiting sets under a rate assumption]
\label{prop:Li_cover_rate_dependent_transfer}
Put $p:=n-3$. Let $E$ and $K$ be compact sets satisfying
\[
    E\subset\operatorname{int}K,
    \qquad
    K\Subset X\setminus\mathcal S.
\]
For the covers in
Proposition~\ref{prop:Li_cover_Hausdorff_content_off_S}, set
$N_i:=\#A_i$, and suppose
\begin{equation}\label{eq:one_sided_Hausdorff_error_to_Li_locus}
    \eta_i
    :=
    \sup_{x\in E}
    d\bigl(x,C_{\Lambda_0,i}\cap K\bigr)
    \longrightarrow0.
\end{equation}
Then
\begin{equation}\label{eq:rate_dependent_Hausdorff_measure_transfer}
    \mathcal H^p(E)
    \leqslant
    C_p
    \liminf_{i\to\infty}N_i\eta_i^p.
\end{equation}
In particular,
\begin{enumerate}[label=\textnormal{(\roman*)}]
\item if $\sup_iN_i\eta_i^p<\infty$, then
$\mathcal H^p(E)<\infty$;
\item if $N_i\eta_i^p\to0$, then $\mathcal H^p(E)=0$;
\item if $\eta_i=O(m_i^{-1})$, then $\mathcal H^p(E)=0$.
\end{enumerate}
The same conclusions hold if
$C_{\Lambda_0,i}\cap K$ in
\eqref{eq:one_sided_Hausdorff_error_to_Li_locus} is replaced by
$Z_i\cap K$.
\end{proposition}

\begin{proof}
Pass to a subsequence along which the lower limit in
\eqref{eq:rate_dependent_Hausdorff_measure_transfer} is attained. If
$N_i=0$ for infinitely many $i$, then
\eqref{eq:one_sided_Hausdorff_error_to_Li_locus} cannot hold unless
$E=\varnothing$, in which case there is nothing to prove. We may thus
assume $N_i\geqslant1$. Choose numbers $\epsilon_i>0$ such that
\[
    \epsilon_i\longrightarrow0,
    \qquad
    N_i\epsilon_i^p\longrightarrow0.
\]
For instance, one may take
$\epsilon_i=(iN_i)^{-1/p}$. For every $x\in E$, choose
$z_{i,x}\in C_{\Lambda_0,i}\cap K$ satisfying
\[
    d(x,z_{i,x})
    \leqslant
    \eta_i+\epsilon_i.
\]
Since the fivefold Li balls cover
$C_{\Lambda_0,i}\cap K$, the enlarged balls
\[
    B\bigl(p_{i,a},5s_{i,a}+\eta_i+\epsilon_i\bigr),
    \qquad
    a\in A_i,
\]
cover $E$. Their maximal diameter tends to zero by
\eqref{eq:Li_cover_radius_bounds_off_S} and
\eqref{eq:one_sided_Hausdorff_error_to_Li_locus}. Moreover,
\[
\begin{aligned}
    \sum_{a\in A_i}
    \bigl(10s_{i,a}+2\eta_i+2\epsilon_i\bigr)^p
    &\leqslant
    C_p\sum_{a\in A_i}s_{i,a}^p
    +C_pN_i\eta_i^p
    +C_pN_i\epsilon_i^p\\
    &=o(1)+C_pN_i\eta_i^p+o(1),
\end{aligned}
\]
where Proposition~\ref{prop:Li_cover_Hausdorff_content_off_S} was used
in the last line. Passing to the limit gives
\eqref{eq:rate_dependent_Hausdorff_measure_transfer}. The first two
consequences are immediate, while the third follows from
$N_i=o(m_i^p)$.

If the approximation is by $Z_i\cap K$, then Li's
concentration--decay dichotomy gives
$Z_i\subset C_{\Lambda_0,i}$ for all sufficiently large $i$. The
one-sided distance to $C_{\Lambda_0,i}\cap K$ is therefore no larger
than the one-sided distance to $Z_i\cap K$, and the preceding argument
applies.
\end{proof}

\begin{proposition}[Obstruction to rapid density of the curvature concentration loci]
\label{prop:Li_loci_density_rate_obstruction}
Put $p:=n-3$. Let
\[
    \overline{B_{3R}(x_0)}
    \Subset
    X\setminus\mathcal S,
\]
and define
\begin{equation}\label{eq:def_local_density_scale_Li_loci}
    h_i^{\mathcal C}
    :=
    \sup_{x\in\overline{B_R(x_0)}}
    d\bigl(x,C_{\Lambda_0,i}\cap\overline{B_{2R}(x_0)}\bigr).
\end{equation}
If $h_i^{\mathcal C}\to0$, then
\begin{equation}\label{eq:Li_loci_density_rate_obstruction}
    m_i^{n-3}
    \bigl(h_i^{\mathcal C}\bigr)^n
    \longrightarrow+\infty.
\end{equation}
Equivalently,
\[
    h_i^{\mathcal C}
    \gg
    m_i^{-\frac{n-3}{n}}.
\]
Thus the loci $C_{\Lambda_0,i}$ cannot become
$O(m_i^{-1/2})$-dense on such a ball when $n=6$, nor
$O(m_i^{-4/7})$-dense when $n=7$.
\end{proposition}

\begin{proof}
Apply Proposition~\ref{prop:Li_cover_Hausdorff_content_off_S} to
$K=\overline{B_{2R}(x_0)}$, and retain the notation
$N_i=\#A_i$. We first note that $h_i^{\mathcal C}>0$ for every
sufficiently large $i$. Indeed, if $h_i^{\mathcal C}=0$ along a
subsequence, choose $\epsilon_i>0$ with
$\epsilon_i\to0$ and $N_i\epsilon_i^n\to0$. The balls
$B(p_{i,a},5s_{i,a}+\epsilon_i)$ then cover
$\overline{B_R(x_0)}$. Uniform upper volume bounds for sufficiently
small geodesic balls would give
\[
    \vol(B_R(x_0))
    \leqslant
    C\sum_{a\in A_i}(s_{i,a}+\epsilon_i)^n.
\]
Since $p=n-3$,
\[
    \sum_{a\in A_i}s_{i,a}^n
    \leqslant
    \left(\max_{a\in A_i}s_{i,a}\right)^3
    \sum_{a\in A_i}s_{i,a}^p
    =o(1),
\]
and the term $N_i\epsilon_i^n$ also tends to zero, a contradiction.

We may therefore assume $h_i^{\mathcal C}>0$. For every
$x\in\overline{B_R(x_0)}$, choose
$z_{i,x}\in C_{\Lambda_0,i}\cap\overline{B_{2R}(x_0)}$ with
\[
    d(x,z_{i,x})
    \leqslant
    2h_i^{\mathcal C}.
\]
It follows that the balls
\[
    B\bigl(p_{i,a},5s_{i,a}+2h_i^{\mathcal C}\bigr),
    \qquad a\in A_i,
\]
cover $\overline{B_R(x_0)}$. Hence
\[
\begin{aligned}
    \vol(B_R(x_0))
    &\leqslant
    C\sum_{a\in A_i}
    \bigl(s_{i,a}+h_i^{\mathcal C}\bigr)^n\\
    &\leqslant
    C\sum_{a\in A_i}s_{i,a}^n
    +CN_i\bigl(h_i^{\mathcal C}\bigr)^n\\
    &=o(1)+CN_i\bigl(h_i^{\mathcal C}\bigr)^n.
\end{aligned}
\]
Therefore
\[
    N_i\bigl(h_i^{\mathcal C}\bigr)^n
    \geqslant c_R>0
\]
for all sufficiently large $i$. Since $N_i=o(m_i^p)$ by
\eqref{eq:Li_cover_sums_off_S},
\[
    m_i^p\bigl(h_i^{\mathcal C}\bigr)^n
    =
    \frac{m_i^p}{N_i}
    N_i\bigl(h_i^{\mathcal C}\bigr)^n
    \longrightarrow+\infty.
\]
This proves \eqref{eq:Li_loci_density_rate_obstruction}.
\end{proof}

\begin{remark}[The loci $C_{\Lambda_0,i}$ and their Kuratowski upper limit]
\label{rmk:Li_content_does_not_transfer_to_Cinfty}
Proposition~\ref{prop:Li_cover_Hausdorff_content_off_S} is a statement
about the loci $C_{\Lambda_0,i}$ on compact subsets of
$X\setminus\mathcal S$. It does not imply
\[
    \mathcal H^{n-3}(\mathcal C\setminus\mathcal S)=0
\]
or local finiteness of this measure. Indeed, membership in the
Kuratowski upper limit $\mathcal C$ permits different points to be
realized along different subsequences, whereas
Proposition~\ref{prop:Li_cover_rate_dependent_transfer} requires a
quantitative one-sided Hausdorff approximation by a common sequence of
loci $C_{\Lambda_0,i}$. In particular, the vanishing content in
\eqref{eq:Li_loci_Hausdorff_content_off_S} does not pass to
$\mathcal C$ without additional information on the convergence rate.

The density rate obstruction
\eqref{eq:Li_loci_density_rate_obstruction} likewise remains compatible
with the possibility that the loci $C_{\Lambda_0,i}$ become dense in an open set,
but forces such density to develop much more slowly than the indicated
powers of the mass. This refines, without resolving, the density issue
raised in \cite[Remark~3.10]{li2025large}.
\end{remark}

The next lemma summarizes the consequences of Li's
concentration-decay dichotomy on $\mathcal R_{\mathrm{ab}}$. We use a
slightly strengthened form of the transverse decay estimate: the pointwise
decay away from the curvature concentration loci is the content of
\cite[Proposition~2.18]{li2025large}, while the corresponding derivative
estimates follow from the $\varepsilon$-regularity input
\cite[Proposition~2.6]{li2025large} used in Li's proof, together with
standard local bootstrapping for the Yang--Mills--Higgs system.

\begin{lemma}[Li's dichotomy on $\mathcal R_{\mathrm{ab}}$]
\label{lem:Li_dichotomy_away_Cinfty}
Let $K\Subset \mathcal R_{\mathrm{ab}}$.
Then there exist an open set $U$,
\[
K\Subset U\Subset \mathcal R_{\mathrm{ab}},
\]
an index $i_K\in\mathbb N$, and constants
\[
r_K>0,
\qquad
C_K<\infty,
\qquad
c_K>0,
\]
such that, for every $i\geqslant i_K$, the following hold.

\begin{enumerate}
\item
\[
U\cap C_{\Lambda_0,i}
=
\varnothing.
\]

\item
\[
|\Phi_i|
\geqslant
\frac{m_i}{2}
\qquad
\text{on } U.
\]

\item
\[
B_{2r_K}(x)\Subset U
\qquad
\text{for every } x\in K.
\]

\item
For every integer $\ell\geqslant 0$, there are constants
$C_{K,\ell}<\infty$ and $c_{K,\ell}>0$ such that
\begin{equation}\label{ineq:Li_transverse_curvature_derivative_decay}
    |\nabla_i^\ell (F_{\nabla_i})^\perp|
    \leqslant
    C_{K,\ell}e^{-c_{K,\ell}m_i}
\end{equation}
on $B_{r_K}(x)$, for every $x\in K$. Consequently, by differentiating the $\Theta$-monopole equation and the
large Higgs field splitting, for every $\ell\geqslant 0$,
\begin{equation}\label{ineq:Li_transverse_Higgs_Psi_derivative_decay}
    |\nabla_i^\ell(\nabla_i\Phi_i)^\perp|
    +
    |\nabla_i^{\ell+1}\Psi_i|
    \leqslant
    C_{K,\ell}e^{-c_{K,\ell}m_i}
\end{equation}
on the same region, after changing the constants. Here $\Psi_i:=\Phi_i/|\Phi_i|$, which is well-defined on $U$ by (2). In particular,
\[
|(F_{\nabla_i})^\perp|^2(p)
\leqslant
C_K e^{-c_Km_i}
\]
for every $p\in B_{r_K}(x)$, $x\in K$, and hence
\[
\int_{B_r(x)}
|(F_{\nabla_i})^\perp|^2
\leqslant
C_K r^n e^{-c_Km_i}
\] for every $x\in K$, $0<r<r_K$.

\item
For every $p\in U$, one has the curvature smallness at the mass scale
\[
\int_{B(p,2m_i^{-1})}|F_{\nabla_i}|^2
<
\varepsilon_3 m_i^{4-n}.
\]
\end{enumerate}
\end{lemma}

\begin{proof}
Since $K\Subset \mathcal R_{\mathrm{ab}}$, there exists an open set
\[
K\Subset U\Subset \mathcal R_{\mathrm{ab}}.
\]
By the definition of $\mathcal C$, after increasing $i_K$ if
necessary, one has
\[
U\cap C_{\Lambda_0,i}
=
\varnothing
\]
for every $i\geqslant i_K$.

The large Higgs field lower bound follows from the contrapositive of
\cite[Lemma~2.17]{li2025large}: outside $C_{\Lambda_0,i}$, the small Higgs field
alternative cannot occur. Hence
\[
|\Phi_i|
\geqslant
\frac{m_i}{2}
\]
on $U$ for all sufficiently large $i$.

The same contrapositive also gives the curvature smallness at the mass scale.
Indeed, since $p\in U$ implies $p\notin C_{\Lambda_0,i}$ for all
$i\geqslant i_K$, Lemma~2.17 gives
\[
\int_{B(p,2m_i^{-1})}|F_{\nabla_i}|^2
<
\varepsilon_3 m_i^{4-n}
\]
for every $p\in U$ and every sufficiently large $i$. Increasing $i_K$ once
more if necessary proves item (5).

Choose $r_K>0$ such that
\[
B_{2r_K}(x)\Subset U
\qquad
\text{for every } x\in K.
\]
Then for every $p\in B_{r_K}(x)$, with $x\in K$, we have
\[
B_{r_K}(p)\Subset B_{2r_K}(x)\Subset U
\]
and hence
\[
B_{r_K}(p)\cap C_{\Lambda_0,i}=\varnothing
\]
for all $i\geqslant i_K$. Therefore Li's exponential decay estimate away
from the curvature concentration locus applies uniformly on these balls.
More precisely, \cite[Proposition~2.18]{li2025large}, whose proof uses
\cite[Proposition~2.15]{li2025large} together with the dichotomy of
\cite[Lemma~2.17]{li2025large}, gives
\[
|(F_{\nabla_i})^\perp|^2(p)
\leqslant
C r_K^{-n}e^{-m_ir_K}
\int_{B_{r_K}(p)}
|(F_{\nabla_i})^\perp|^2.
\]
The curvature integral over the fixed compact set $U$ is bounded by
$C_K(m_i+1)$ by the local a priori curvature bounds following from
Lemma~\ref{lemm: Yang-Li}. Thus
\[
|(F_{\nabla_i})^\perp|^2(p)
\leqslant
C_K(m_i+1)e^{-m_ir_K}.
\]
After decreasing the exponential rate and changing the constant, the
polynomial factor is absorbed into the exponential, giving
\[
|(F_{\nabla_i})^\perp|^2(p)
\leqslant
C_K e^{-c_Km_i}.
\]
We now justify the derivative estimates. The curvature smallness at the mass scale
proved above holds at every point of $U$. After rescaling a ball
of radius comparable to $m_i^{-1}$ to unit size, the rescaled Higgs field
$m_i^{-1}\Phi_i$ is uniformly bounded by the maximum principle.
The local $\varepsilon$-regularity estimate
\cite[Proposition~2.6]{li2025large} supplies the initial bounds at the mass scale
for $F_{\nabla_i}$, $\nabla_i\Phi_i$, and $\nabla_iF_{\nabla_i}$.
Standard interior bootstrapping for the differentiated
system \eqref{eq:YMH_second_order} on slightly smaller balls, equivalently
the estimates collected in Appendix~\ref{app: C}, then gives polynomial bounds
for all covariant derivatives. More precisely, after replacing $U$ by a
slightly smaller open set still containing all the balls $B_{r_K}(x)$,
for every integer $N\geqslant0$ there are constants $C_{K,N}<\infty$ and an integer
$M_N\geqslant0$ such that
\[
    \sum_{j=0}^{N}
    \left(
        \|\nabla_i^jF_{\nabla_i}\|_{L^\infty}
        +
        \|\nabla_i^{j+1}\Phi_i\|_{L^\infty}
        +
        \|\nabla_i^{j+1}\Psi_i\|_{L^\infty}
    \right)
    \leqslant
    C_{K,N}m_i^{M_N}
\]
on the relevant smaller region. Here the estimates for $\Psi_i$ follow
from $|\Phi_i|\geqslant m_i/2$ and
\[
    \nabla_i\Psi_i
    =
    |\Phi_i|^{-1}(\nabla_i\Phi_i)^\perp,
\]
by differentiating the quotient and using the preceding polynomial
bounds. In particular, the same type of polynomial $C^N$ bound holds for
$(F_{\nabla_i})^\perp$.

Fix $\ell\geqslant1$ and choose $N>\ell$. Covariant interpolation on a
slightly smaller ball gives
\[
\begin{aligned}
    \|\nabla_i^\ell(F_{\nabla_i})^\perp\|_{L^\infty}
    &\leqslant
    C
    \|(F_{\nabla_i})^\perp\|_{L^\infty}^{\,1-\ell/N}
    \|(F_{\nabla_i})^\perp\|_{C^N}^{\,\ell/N}  \\
    &\leqslant
    C_{K,\ell}
    m_i^{M_{N}\ell/N}
    \exp\!\left(
        -c_K\left(1-\frac{\ell}{N}\right)m_i
    \right).
\end{aligned}
\]
After decreasing the exponential rate and changing the constant, the
polynomial factor is absorbed into the exponential. Together with the
already established case $\ell=0$, this proves
\eqref{ineq:Li_transverse_curvature_derivative_decay} for every fixed
$\ell$.

Since the orthogonal projection acts only in the Lie algebra factor, the
$\Theta$-monopole equation gives the exact identity
\[
    (F_{\nabla_i})^\perp\wedge\Theta
    =
    *(\nabla_i\Phi_i)^\perp.
\]
The form $\Theta$ is parallel in both special holonomy settings considered
here. Differentiating this identity and using
\eqref{ineq:Li_transverse_curvature_derivative_decay} therefore yields
the corresponding higher order exponential estimates for
$(\nabla_i\Phi_i)^\perp$. Finally, differentiating
\[
    \nabla_i\Psi_i
    =
    |\Phi_i|^{-1}(\nabla_i\Phi_i)^\perp
\]
and using $|\Phi_i|\geqslant m_i/2$, the polynomial derivative bounds
above, and the product rule gives the higher order exponential estimates for
$\Psi_i$, after decreasing the exponential rates once more. This proves
\eqref{ineq:Li_transverse_Higgs_Psi_derivative_decay}.

The integral estimate follows by integrating the case $\ell=0$ of
\eqref{ineq:Li_transverse_curvature_derivative_decay} over $B_r(x)$.
\end{proof}

The preceding lemma has an immediate consequence for the limiting zero
set. Since zeros of the Higgs fields trigger the small Higgs field
alternative in
Li's concentration-decay dichotomy, they cannot accumulate in the region
where the loci $C_{\Lambda_0,i}$ are eventually absent.

\subsection{Limiting Higgs zeros and the limiting nonabelian locus}
\label{subsec:zeros_and_Cinfty}

\begin{lemma}[Zeros lie in the limiting nonabelian locus]
\label{lemma:Z_subset_C_infty}
One has
\[
    \mathcal Z\subset \mathcal C.
\]
Consequently,
\[
    \mathcal S\subset \mathcal Z\subset\mathcal C.
\]
\end{lemma}

\begin{proof}
Let $x\in \mathcal Z$. Let $i_{\mathrm{Li}}$ be large enough that the
concentration-decay dichotomy above holds for every
$i\geqslant i_{\mathrm{Li}}$. By definition of $\mathcal Z$, for every
neighbourhood $U$ of $x$ and every $N\geqslant1$, there exist
\[
    i\geqslant\max\{N,i_{\mathrm{Li}}\},
    \qquad
    z_i\in U\cap\Phi_i^{-1}(0).
\]
The small Higgs field alternative applies at $z_i$, since
\[
    \min_{B(z_i,m_i^{-1})}|\Phi_i|
    =
    0
    \leqslant
    \frac{m_i}{2}.
\]
Hence, by \cite[Lemma~2.17]{li2025large},
\[
    z_i\in C_{\Lambda_0,i}.
\]
Since this holds for arbitrarily large indices and for every neighbourhood
$U$ of $x$, we obtain
\[
    x\in
    \bigcap_{N\geqslant1}
    \overline{
    \bigcup_{i\geqslant N}C_{\Lambda_0,i}
    }
    =
    \mathcal C.
\]
Thus $\mathcal Z\subset\mathcal C$. The inclusion $\mathcal S\subset \mathcal Z$ is Proposition~\ref{prop:S_subset_Z}.
\end{proof}

\begin{corollary}[Hausdorff content decay of the Higgs zero sets away from $\mathcal S$]
\label{cor:Higgs_zero_Hausdorff_content_off_S}
Put $p:=n-3$. For every compact set
$K\Subset X\setminus\mathcal S$ and every $0<\kappa<1$,
\begin{equation}\label{eq:Higgs_zero_Hausdorff_content_off_S}
    \mathcal H^{p+\kappa}_\infty(Z_i\cap K)
    =o(m_i^{-\kappa}),
    \qquad
    \mathcal H^p_\infty(Z_i\cap K)
    \longrightarrow0.
\end{equation}
Moreover, the sets $Z_i\cap K$ admit covers by the fivefold Li balls of
Proposition~\ref{prop:Li_cover_Hausdorff_content_off_S}; in particular,
the number of balls is $o(m_i^p)$ and their maximal radius is
$O(m_i^{-\kappa/(1+\kappa)})$.
\end{corollary}

\begin{proof}
For every sufficiently large $i$, each zero of $\Phi_i$ belongs to
$C_{\Lambda_0,i}$ by the small Higgs field alternative in Li's
concentration--decay dichotomy, as used in the proof of
Lemma~\ref{lemma:Z_subset_C_infty}. Hence
\[
    Z_i\cap K
    \subset
    C_{\Lambda_0,i}\cap K,
\]
and all assertions follow from
Proposition~\ref{prop:Li_cover_Hausdorff_content_off_S}.
\end{proof}

We next combine the local Higgs defect estimate
\eqref{eq: largeness_of_Higgs_field} with ordinary
$\varepsilon$-regularity in codimension four. This gives a stronger
obstruction to rapid density of the Higgs zero sets than Li's covers
alone provide.

\begin{lemma}[Energy required to pass from a Higgs zero to $|\Phi|\geqslant m/2$]
\label{lem:zero_to_large_Higgs_transition_energy}
There exist constants
\[
    \varepsilon_4>0,
    \qquad
    M_0>0,
    \qquad
    r_0>0,
\]
depending only on the fixed bounded geometry data and the structure group,
with the following property. Let $(\nabla,\Phi)$ be a solution of
\eqref{eq:YMH_second_order} with positive finite mass $m$, let $z\in X$, and
let $0<r\leqslant r_0$. If
\[
    \Phi(z)=0,
    \qquad
    \sup_{B_{r/4}(z)}|\Phi|\geqslant\frac m2,
    \qquad
    mr\geqslant M_0,
\]
then
\begin{equation}\label{eq:zero_to_large_Higgs_transition_energy}
    r^{4-n}
    \int_{B_r(z)}e(\nabla,\Phi)
    \geqslant
    \varepsilon_4.
\end{equation}
\end{lemma}

\begin{proof}
Let $r_0$ and $\varepsilon_0$ be the constants in
Theorem~\ref{thm: total_epsilon_regularity}. Put
\[
    E(r)
    :=
    r^{4-n}
    \int_{B_r(z)}e(\nabla,\Phi).
\]
After decreasing a constant $\varepsilon_4>0$, depending only on the
universal $\varepsilon$-regularity threshold, the implication
\[
    E(r)<\varepsilon_4
\]
gives, by the case $j=0$ of
\eqref{ineq:estimates_coulomb},
\begin{equation}\label{eq:gradient_Higgs_from_small_codim4_energy}
    \sup_{B_{r/2}(z)}|\nabla\Phi|^2
    \leqslant
    Cr^{-4}E(r).
\end{equation}
Choose $y\in B_{r/4}(z)$ with $|\Phi(y)|\geqslant m/2$, and let
$\gamma$ be a minimizing geodesic from $z$ to $y$. By Kato's inequality,
\[
\begin{aligned}
    \frac m2
    &\leqslant
    |\Phi(y)|-|\Phi(z)|
    \leqslant
    \int_\gamma|d|\Phi||\\
    &\leqslant
    \frac r4
    \sup_{B_{r/4}(z)}|\nabla\Phi|
    \leqslant
    Cr^{-1}E(r)^{1/2},
\end{aligned}
\]
where \eqref{eq:gradient_Higgs_from_small_codim4_energy} was used in the
last inequality. Consequently,
\[
    E(r)\geqslant c m^2r^2.
\]
Choose $M_0$ so large that $cM_0^2>\varepsilon_4$. If $mr\geqslant M_0$,
the last estimate contradicts $E(r)<\varepsilon_4$. Hence
\eqref{eq:zero_to_large_Higgs_transition_energy} holds.
\end{proof}

\begin{proposition}[Obstruction to rapid density of the Higgs zero sets]
\label{prop:Higgs_zero_density_rate_obstruction}
Let
\[
    \overline{B_{4R}(x_0)}
    \Subset
    X\setminus\mathcal S,
\]
and define
\begin{equation}\label{eq:def_local_density_scale_Higgs_zeros}
    h_i^{\mathcal Z}
    :=
    \sup_{x\in\overline{B_R(x_0)}}d(x,Z_i).
\end{equation}
If $h_i^{\mathcal Z}\to0$, then
\begin{equation}\label{eq:Higgs_zero_density_rate_obstruction}
    m_i\bigl(h_i^{\mathcal Z}\bigr)^4
    \longrightarrow+\infty.
\end{equation}
Equivalently,
\[
    h_i^{\mathcal Z}\gg m_i^{-1/4}.
\]
Thus the Higgs zero sets cannot become $O(m_i^{-1/4})$-dense in a fixed
open set disjoint from $\mathcal S$.
\end{proposition}

\begin{proof}
Suppose, to the contrary, that
\eqref{eq:Higgs_zero_density_rate_obstruction} fails. Passing to a
subsequence, there is $M<\infty$ such that
\begin{equation}\label{eq:bounded_zero_density_rate_for_contradiction}
    m_i\bigl(h_i^{\mathcal Z}\bigr)^4
    \leqslant M
\end{equation}
for every $i$. Set
\[
    \ell_i
    :=
    \max\bigl\{h_i^{\mathcal Z},m_i^{-1/4}\bigr\}.
\]
Then
\begin{equation}\label{eq:ell_zero_density_properties}
    \ell_i\longrightarrow0,
    \qquad
    m_i\ell_i\geqslant m_i^{3/4}\longrightarrow+\infty,
    \qquad
    m_i\ell_i^4\leqslant\max\{M,1\}.
\end{equation}
Choose a maximal $24\ell_i$-separated family
\[
    \{x_{i,1},\ldots,x_{i,N_i}\}
    \subset
    \overline{B_R(x_0)}.
\]
The balls $B_{24\ell_i}(x_{i,a})$ cover
$\overline{B_R(x_0)}$, so uniform upper volume bounds imply
\begin{equation}\label{eq:number_separated_points_zero_density}
    N_i\geqslant c_R\ell_i^{-n}.
\end{equation}
Since $Z_i$ is closed and $(X,g)$ is proper, the distance to $Z_i$ is
attained. Choose $z_{i,a}\in Z_i$ such that
\[
    d(x_{i,a},z_{i,a})
    \leqslant
    h_i^{\mathcal Z}
    \leqslant
    \ell_i.
\]
For $a\neq b$,
\[
    d(z_{i,a},z_{i,b})
    \geqslant
    24\ell_i-2\ell_i
    =22\ell_i.
\]
Thus the balls $B_{8\ell_i}(z_{i,a})$ are pairwise disjoint and, for all
sufficiently large $i$, are contained in $B_{2R}(x_0)$.

Define
\[
    A_i
    :=
    \left\{
        x\in B_{2R}(x_0):
        |\Phi_i(x)|\leqslant\frac{m_i}{2}
    \right\}.
\]
On $A_i$, one has
$(m_i-|\Phi_i|)^2\geqslant m_i^2/4$. Hence
Lemma~\ref{lem:largeness}, applied to
$\overline{B_{2R}(x_0)}$, gives
\begin{equation}\label{eq:small_volume_half_Higgs_region}
    \vol(A_i)
    \leqslant
    \frac4{m_i^2}
    \int_{B_{2R}(x_0)}(m_i-|\Phi_i|)^2
    =o(m_i^{-1}).
\end{equation}
Call an index $a$ bad if
$B_{2\ell_i}(z_{i,a})\subset A_i$. The balls
$B_{2\ell_i}(z_{i,a})$ are pairwise disjoint, and uniform lower volume
bounds give
\[
    c\ell_i^n\#\{a:\ a\text{ is bad}\}
    \leqslant
    \vol(A_i).
\]
Together with \eqref{eq:number_separated_points_zero_density}, this
implies
\[
    \frac{\#\{a:\ a\text{ is bad}\}}{N_i}
    \leqslant
    C_R\vol(A_i)
    \longrightarrow0.
\]
Thus at least $N_i/2$ indices are good for all sufficiently large $i$.
For every good $a$, there exists
$y_{i,a}\in B_{2\ell_i}(z_{i,a})$ such that
$|\Phi_i(y_{i,a})|>m_i/2$.

Apply Lemma~\ref{lem:zero_to_large_Higgs_transition_energy} with
\[
    z=z_{i,a},
    \qquad
    r=8\ell_i.
\]
The point $y_{i,a}$ belongs to $B_{r/4}(z_{i,a})$, while
$m_ir\to+\infty$ by \eqref{eq:ell_zero_density_properties}. Hence, for
every good $a$ and all sufficiently large $i$,
\[
    \int_{B_{8\ell_i}(z_{i,a})}e_i
    \geqslant
    c\varepsilon_4\ell_i^{n-4}.
\]
Summing over the pairwise disjoint good balls and using
\eqref{eq:number_separated_points_zero_density}, we obtain
\[
    \int_{B_{2R}(x_0)}e_i
    \geqslant
    c_RN_i\ell_i^{n-4}
    \geqslant
    c_R\ell_i^{-4}.
\]
Therefore
\begin{equation}\label{eq:normalized_energy_lower_bound_from_dense_zeros}
    \mu_i(B_{2R}(x_0))
    =
    m_i^{-1}
    \int_{B_{2R}(x_0)}e_i
    \geqslant
    \frac{c_R}{m_i\ell_i^4}.
\end{equation}
The right-hand side is bounded below by a positive constant by
\eqref{eq:ell_zero_density_properties}. On the other hand,
$\overline{B_{2R}(x_0)}\cap\mathcal S=\varnothing$. The Portmanteau
theorem therefore gives
\[
    0
    \leqslant
    \limsup_{i\to\infty}\mu_i(B_{2R}(x_0))
    \leqslant
    \limsup_{i\to\infty}\mu_i(\overline{B_{2R}(x_0)})
    \leqslant
    8\pi\|T\|(\overline{B_{2R}(x_0)})
    =0.
\]
This contradicts
\eqref{eq:normalized_energy_lower_bound_from_dense_zeros} and proves
the proposition.
\end{proof}

\begin{corollary}[Rate obstruction for local Hausdorff limits with interior]
\label{cor:local_Hausdorff_limit_interior_rate_obstruction}
Let $Z_{\mathrm H}\subset X$ be closed and suppose that $Z_i$ converge
locally in Hausdorff distance to $Z_{\mathrm H}$. If
\[
    \overline{B_{4R}(x_0)}\Subset X\setminus\mathcal S,
    \qquad
    \overline{B_R(x_0)}\subset Z_{\mathrm H},
\]
then the quantities $h_i^{\mathcal Z}$ in
\eqref{eq:def_local_density_scale_Higgs_zeros} satisfy
\[
    h_i^{\mathcal Z}\longrightarrow0,
    \qquad
    m_i\bigl(h_i^{\mathcal Z}\bigr)^4
    \longrightarrow+\infty.
\]
An analogous statement holds for a local Hausdorff limit of the closures
$\overline{C_{\Lambda_0,i}}$, with the weaker rate
\eqref{eq:Li_loci_density_rate_obstruction}.
\end{corollary}

\begin{proof}
Local Hausdorff convergence and the inclusion of
$\overline{B_R(x_0)}$ in the limit give the relevant one-sided distances
tending to zero. Apply
Propositions~\ref{prop:Higgs_zero_density_rate_obstruction} and
\ref{prop:Li_loci_density_rate_obstruction}, respectively.
\end{proof}

\begin{remark}[Limitations of the quantitative density estimates]
\label{rmk:density_rate_estimates_do_not_exclude_density}
The estimates
\eqref{eq:Li_loci_density_rate_obstruction} and
\eqref{eq:Higgs_zero_density_rate_obstruction} do not exclude density of
the curvature concentration loci or the Higgs zero sets in an
open region disjoint from $\mathcal S$. They show that any such density
must develop slowly relative to the mass. The stronger exponent
$m_i^{-1/4}$ for the zero sets uses the exact vanishing of the Higgs
field, the local defect estimate
\eqref{eq: largeness_of_Higgs_field}, and ordinary codimension-four
$\varepsilon$-regularity; it does not extend to arbitrary points of
$C_{\Lambda_0,i}$.

A qualitative exclusion of dense zero sets would follow from a uniform
codimension-three lower density estimate centred at a Higgs zero on all scales
between $m_i^{-1}$ and a fixed radius. Such an estimate is precisely the
type of conclusion that effective codimension-three monotonicity would
propagate from a nonvanishing initial density at a Higgs zero. It is not
presently available on the possible excess locus
$\mathcal O\setminus\mathcal S
=\mathcal C\setminus\mathcal S$.
\end{remark}

\begin{proposition}[Local uniform decay of the mass-renormalized energy density]
\label{prop:mass_renormalized_density_decay_Rab}
For every compact set $K\Subset\mathcal R_{\mathrm{ab}}$, one has
\begin{equation}
\label{eq:mass_renormalized_density_decay_Rab}
    \sup_K m_i^{-1}e_i
    \longrightarrow
    0.
\end{equation}
In particular,
\[
    m_i^{-1/2}\|F_{\nabla_i}\|_{L^\infty(K)}
    +
    m_i^{-1/2}\|\nabla_i\Phi_i\|_{L^\infty(K)}
    \longrightarrow
    0.
\]
\end{proposition}

\begin{proof}
Choose open sets
\[
    K\Subset V\Subset U\Subset\mathcal R_{\mathrm{ab}}
\]
and apply Lemma~\ref{lem:Li_dichotomy_away_Cinfty} on the compact set
$\overline U$. On $U$, decompose
\[
    F_{\nabla_i}
    =
    (F_{\nabla_i})^\parallel+(F_{\nabla_i})^\perp,
    \qquad
    \nabla_i\Phi_i
    =
    (\nabla_i\Phi_i)^\parallel+(\nabla_i\Phi_i)^\perp
\]
with respect to the line generated by
$\Psi_i:=\Phi_i/|\Phi_i|$. The parallel line is abelian. Consequently,
each nonzero bracket occurring in the cubic terms of the Bochner
identities \eqref{eq: gen_bochner_dphi} and
\eqref{eq: gen_bochner_curv} contains at least one transverse factor.
Thus, writing $a_i:=\nabla_i\Phi_i$, bilinearity of the bracket gives
\[
\begin{aligned}
&\left|
\left\langle *[*F_{\nabla_i},a_i],a_i\right\rangle
\right|
+
\left|
\left\langle[a_i,a_i],F_{\nabla_i}\right\rangle
\right|
+
\left|
\sum_{j,k,\ell}
\left\langle[F_{jk},F_{k\ell}],F_{\ell j}\right\rangle
\right|  \\
&\qquad\leqslant
C
\left(
|(F_{\nabla_i})^\perp|
+
|(\nabla_i\Phi_i)^\perp|
\right)e_i.
\end{aligned}
\]
The remaining zeroth-order terms in the Bochner identities are either
nonpositive or bounded by $C_Ue_i$ through the ambient curvature. By
Lemma~\ref{lem:Li_dichotomy_away_Cinfty}(4), the transverse factor on the
right-hand side converges exponentially to zero on $U$. Hence, after
increasing the threshold index,
\begin{equation}
\label{eq:linearized_Bochner_Rab}
    \Delta e_i
    \leqslant
    C_Ue_i
\end{equation}
on $U$, with $C_U$ independent of $i$.

Choose finitely many balls
$B_{4\rho}(x_\alpha)\Subset V$ such that the balls
$B_\rho(x_\alpha)$ cover $K$.
Theorem~\ref{thm: gen_mean_value}, applied with
$d=n$, $a_0=a=0$, and $a_1=C_U$ (the required integral control being the
trivial inclusion of smaller balls in $B_{4\rho}(x_\alpha)$), gives
\[
    \sup_{B_\rho(x_\alpha)}m_i^{-1}e_i
    \leqslant
    C_{U,\rho}
    m_i^{-1}\int_{B_{4\rho}(x_\alpha)}e_i.
\]
By Lemma~\ref{lemma:Z_subset_C_infty},
$\mathcal S\subset\mathcal C$, and hence
$\overline{B_{4\rho}(x_\alpha)}\cap\mathcal S=\varnothing$. Since
$\mu_i=m_i^{-1}e_i\,\operatorname{vol}\rightharpoonup\mu$ and
$\operatorname{spt}\mu=\mathcal S$, the Portmanteau theorem gives
\[
\begin{aligned}
    0
    &\leqslant
    \limsup_{i\to\infty}
    m_i^{-1}\int_{B_{4\rho}(x_\alpha)}e_i  \\
    &\leqslant
    \mu\bigl(\overline{B_{4\rho}(x_\alpha)}\bigr)
    =0.
\end{aligned}
\]
Taking the maximum over the finite cover proves
\eqref{eq:mass_renormalized_density_decay_Rab}.
\end{proof}

\begin{corollary}[Quantitative Morrey decay on the open set $\mathcal R_{\mathrm{ab}}$]
\label{cor:effective_monotonicity_on_Rab}
For every compact set $K\Subset\mathcal R_{\mathrm{ab}}$, there exist
$r_K>0$ and a sequence $\varepsilon_i(K)\to0$ such that, for every
$x\in K$ and every $0<r<r_K$,
\begin{equation}
\label{eq:clearing_region_average_decay}
    m_i^{-1}r^{-n}\int_{B_r(x)}e_i
    \leqslant
    \varepsilon_i(K),
\end{equation}
and hence
\begin{equation}
\label{eq:clearing_region_Morrey_decay}
\begin{aligned}
    m_i^{-1}r^{2-n}\int_{B_r(x)}e_i
    &\leqslant
    \varepsilon_i(K)r^2,\\
    \theta_i(x,r)
    &\leqslant
    \varepsilon_i(K)r^3.
\end{aligned}
\end{equation}
The first estimate remains valid with $e_i$ replaced by any of
$|F_{\nabla_i}|^2$, $|\nabla_i\Phi_i|^2$, or $q_i^+$.

Every point of $\mathcal R_{\mathrm{ab}}$ therefore admits an open
neighbourhood on which effective codimension-three monotonicity holds;
in fact, one may take the error exponent $\gamma_U=3$. Consequently,
\begin{equation}
\label{eq:O3_subset_Cinfty_early}
    \mathcal O\subset\mathcal C.
\end{equation}
\end{corollary}

\begin{proof}
Choose an open set $U$ with
$K\Subset U\Subset\mathcal R_{\mathrm{ab}}$ and then choose $r_K>0$ so
that $B_{r_K}(x)\Subset U$ for every $x\in K$. Set
\[
    \delta_i(U):=
    \sup_{U}m_i^{-1}e_i.
\]
By Proposition~\ref{prop:mass_renormalized_density_decay_Rab},
$\delta_i(U)\to0$. The local volume upper bound gives
\[
    m_i^{-1}\int_{B_r(x)}e_i
    \leqslant
    C_U\delta_i(U)r^n
\]
for $x\in K$ and $0<r<r_K$. Taking
$\varepsilon_i(K):=C_U\delta_i(U)$ proves
\eqref{eq:clearing_region_average_decay} and
\eqref{eq:clearing_region_Morrey_decay}; the componentwise estimates
follow from
\[
    q_i^+\leqslant |F_{\nabla_i}|^2\leqslant e_i,
    \qquad
    |\nabla_i\Phi_i|^2\leqslant e_i.
\]
To prove the final assertion, let $x\in\mathcal R_{\mathrm{ab}}$ and
choose an open set $W$ with
\[
    x\in W,
    \qquad
    \overline W\Subset U\Subset\mathcal R_{\mathrm{ab}}.
\]
Let $r_W>0$ be the radius supplied by
\eqref{eq:clearing_region_average_decay} for the compact set
$\overline W$, and choose an open neighbourhood $V$ of $x$ such that
\[
    \overline V\Subset W,
    \qquad
    2\operatorname{diam}(V)<r_W.
\]
Applying the estimate to $\overline W$, after increasing the threshold
index one has $\varepsilon_i(\overline W)\leqslant1$. Every geodesic ball
$B_r(y)\subset V$ then has $r<r_W$, and for every
$0<s<r\leqslant r_0$,
\[
    e^{cs^2}\theta_i(y,s)
    \leqslant
    e^{cr_0^2}s^3
    \leqslant
    e^{cr^2}\theta_i(y,r)
    +
    e^{cr_0^2}r^3.
\]
Hence effective codimension-three monotonicity holds on $V$ with
$A_V=e^{cr_0^2}$ and $\gamma_V=3$. Therefore
$x\notin\mathcal O$, and
$\mathcal O\subset\mathcal C$ follows.
\end{proof}

\begin{corollary}[Universal regularity scale separation on $\mathcal R_{\mathrm{ab}}$]
\label{cor:regularity_scale_mass_separation_Rab}
For every compact set $K\Subset\mathcal R_{\mathrm{ab}}$, there exist
$r_K>0$, $c_K>0$, and a sequence $\delta_i(K)\to0$ such that
\begin{equation}
\label{eq:quantitative_regularity_scale_Rab}
    \inf_{x\in K}\mathfrak r_i(x)
    \geqslant
    \min\left\{
        r_K,
        c_K\bigl(m_i\delta_i(K)\bigr)^{-1/4}
    \right\}.
\end{equation}
In particular,
\begin{equation}
\label{eq:regularity_scale_mass_separation_Rab}
    \inf_{x\in K}
    m_i^{1/4}\mathfrak r_i(x)
    \longrightarrow
    +\infty.
\end{equation}
\end{corollary}

\begin{proof}
Choose $U$ and $r_K$ as in the proof of
Corollary~\ref{cor:effective_monotonicity_on_Rab}, and set
$\delta_i(K):=\sup_Um_i^{-1}e_i\to0$. There is a geometric constant $C_U$
such that, whenever $x\in K$ and $0<s<r_K$,
\[
    s^{4-n}\int_{B_s(x)}e_i
    \leqslant
    C_Um_i\delta_i(K)s^4.
\]
With the convention that $\delta_i(K)^{-1/4}=+\infty$ when
$\delta_i(K)=0$, set
\[
    R_i
    :=
    \min\left\{
        r_K,
        \left(
            \frac{\varepsilon_0}{2C_Um_i\delta_i(K)}
        \right)^{1/4}
    \right\}.
\]
For all sufficiently large $i$, the defining small energy inequality for
$\mathfrak r_i(x)$ holds at every scale $0<s\leqslant R_i$, uniformly
for $x\in K$. Hence
\[
    \mathfrak r_i(x)\geqslant R_i.
\]
This is \eqref{eq:quantitative_regularity_scale_Rab} with
$c_K=(\varepsilon_0/(2C_U))^{1/4}$. Finally,
\[
    m_i^{1/4}R_i
    =
    \min\left\{
        m_i^{1/4}r_K,
        \left(
            \frac{\varepsilon_0}{2C_U\delta_i(K)}
        \right)^{1/4}
    \right\}
    \longrightarrow+\infty,
\]
which proves \eqref{eq:regularity_scale_mass_separation_Rab}.
\end{proof}

\begin{remark}
\label{rmk:Cinfty_vs_S}
The inclusions $\mathcal S\subset \mathcal Z\subset\mathcal C$ relate three different objects. The set $\mathcal S$ is the support of the limiting
calibrated current obtained from the mass-renormalized energy measures; it is rectifiable of codimension-three. The set $\mathcal Z$ is the limiting zero set of the Higgs fields. The
set $\mathcal C$ is the limiting nonabelian locus, defined as the
Kuratowski upper limit of Li's curvature concentration loci
$C_{\Lambda_0,i}$.

Thus the first inclusion says that codimension-three concentration carrying nonzero transverse charge
forces limiting Higgs zeros, while the second says that every limiting
Higgs zero lies in the locus where Li's large Higgs field and transverse
decay estimates may fail. Proposition~\ref{prop:zero_centered_TU_concentration}
also shows that, after passing to a subsequence along which the zero is realized, every
point of $\mathcal Z$ is an ordinary codimension-four concentration
point.

The set $\mathcal C$ is closed, but the present analysis gives
neither compactness nor an unconditional Hausdorff dimension or
Hausdorff measure bound for it. Proposition~\ref{prop:Li_cover_Hausdorff_content_off_S}
does show that the loci $C_{\Lambda_0,i}$ have vanishing critical
unrestricted Hausdorff content on every compact subset of
$X\setminus\mathcal S$, and
Proposition~\ref{prop:Li_loci_density_rate_obstruction} prevents them
from becoming dense there at the rate
$O(m_i^{-(n-3)/n})$. These estimates for the varying loci do not pass
unconditionally to the Kuratowski upper limit: different points of
$\mathcal C$ may be realized along different subsequences. We therefore
still cannot preclude the possibility that $\mathcal C$ has nonempty
interior, or even that $\mathcal C=X$; see
Remark~\ref{rmk:Li_content_does_not_transfer_to_Cinfty}. Moreover, as
explained in Remark~\ref{rmk:Cinfty_escape_infinity}, the upper limiting
set does not record components whose centers escape to infinity
\cite[Remark~3.10, Remark~4.5, Question~2, and Remark~6.2]{li2025large}.

This should be contrasted with the calibrated concentration set $\mathcal S$.
Because $\mathcal S$ is the support of a $\Theta$-calibrated integral cycle, it has the
geometric structure and density properties of such currents. The possible
difference $\mathcal C\setminus\mathcal S$ consists of points represented by convergent sequences from the curvature concentration loci that are not detected by the
limiting calibrated measure.

Corollary~\ref{cor:effective_monotonicity_on_Rab} already gives
$\mathcal O\subset\mathcal C$. The next subsection proves the
complementary inclusion
$\mathcal C\setminus\mathcal S\subset\mathcal O$ by proving that
Li's weighted curvature potential lies below the defining threshold under
vanishing codimension-three Morrey density. Together these two statements identify
$\mathcal C=\mathcal S\cup\mathcal O$.
\end{remark}

\subsection{Excluding curvature concentration by effective monotonicity}
\label{subsec:positive_monotonicity_clearing}

We now prove the exclusion result announced in
Remark~\ref{rmk:Cinfty_vs_S}. The argument splits Li's weighted curvature
potential into a local and a far contribution. On a small neighbourhood of a
point in $X\setminus(\mathcal S\cup\mathcal O)$, the local contribution is made
arbitrarily small by
Lemma~\ref{lem:effective_monotonicity_gives_vanishing_Morrey_off_S}. Once
this neighbourhood is fixed, the far contribution tends to zero as
$i\to\infty$ by the large radius curvature estimate
Lemma~\ref{lemm: Yang-Li}.

\begin{proposition}[Vanishing codimension-three Morrey density excludes curvature concentration]
\label{prop:Cinfty_minus_S_obstructs_positive_codim3}
Let $x\in X\setminus\mathcal S$. If
\begin{equation}
\label{eq:vanishing_center_stable_Morrey_for_clearing}
    \lim_{r\downarrow0}
    \limsup_{i\to\infty}
    \mathfrak M_i(x,r)
    =
    0,
\end{equation}
then $x\notin\mathcal C$.
\end{proposition}
\begin{proof}
Let $x\in X\setminus\mathcal S$ and assume \eqref{eq:vanishing_center_stable_Morrey_for_clearing}. We shall prove that $x\notin\mathcal C$. 

Fix a small parameter $\eta>0$, to be chosen below. By
\eqref{eq:vanishing_center_stable_Morrey_for_clearing}, after decreasing
$r>0$ we may assume that $\overline{B_r(x)}\cap \mathcal S=\varnothing$ and that, for all sufficiently large $i$,
\begin{equation}\label{eq:small_center_stable_Morrey_for_clearing}
    \mathfrak M_i(x,r)\leqslant\eta.
\end{equation}
In particular, if $y\in B_{r/4}(x)$ and $0<t\leqslant r/2$, then
$B_t(y)\subset B_r(x)$, and hence
\begin{equation}\label{eq:Morrey_bound_balls_centered_y}
    \int_{B_t(y)}e_i
    \leqslant
    \eta\,m_i\,t^{n-3}.
\end{equation}
We estimate Li's weighted curvature potential on $B_{r/4}(x)$. By
definition, $C_{\Lambda_0,i} = \{p\in X:\mathcal P_i(p)>\Lambda_0^{-1}\}$. Thus, it is enough to prove that
$\mathcal P_i<\Lambda_0^{-1}$ on a neighbourhood of $x$ for all
sufficiently large $i$.

Fix $y\in B_{r/4}(x)$ and split
\[
\begin{aligned}
    \mathcal P_i(y)
    &=
    m_i^{-2}
    \int_{B_{r/2}(y)}
    \frac{|F_{\nabla_i}|^2(z)}
    {\max(d(z,y),m_i^{-1})^{n-2}}
    \,{\vol}(z)  \\
    &\quad
    +
    m_i^{-2}
    \int_{X\setminus B_{r/2}(y)}
    \frac{|F_{\nabla_i}|^2(z)}
    {d(z,y)^{n-2}}
    \,{\vol}(z)                                      \\
    &=:
    I_i^{\mathrm{loc}}(y)+I_i^{\mathrm{far}}(y).
\end{aligned}
\]
We first estimate the local term. Set $\rho_0:=m_i^{-1}$ and assume
$i$ is sufficiently large so that $\rho_0<r/8$. The contribution of
$B_{\rho_0}(y)$ is bounded using \eqref{eq:Morrey_bound_balls_centered_y}:
\[
\begin{aligned}
    m_i^{-2}
    \int_{B_{\rho_0}(y)}
    \frac{|F_{\nabla_i}|^2(z)}
    {\max(d(z,y),m_i^{-1})^{n-2}}
    \,{\vol}(z)
    &\leqslant
    m_i^{-2}m_i^{n-2}\int_{B_{\rho_0}(y)}e_i        \\
    &\leqslant
    m_i^{n-4}\eta m_i\rho_0^{\,n-3}
    =
    \eta.
\end{aligned}
\]
For the remaining part of $B_{r/2}(y)$, decompose dyadically. Let
$\rho_j:=2^j\rho_0$ and let $J_i$ be chosen so that
$\rho_{J_i}\leqslant r/2<2\rho_{J_i}$. On the annulus
$A_j(y):=B_{\rho_{j+1}}(y)\setminus B_{\rho_j}(y)$, the kernel is bounded
by $C\rho_j^{2-n}$. Hence, for $0\leqslant j<J_i$,
\[
\begin{aligned}
    m_i^{-2}
    \int_{A_j(y)}
    \frac{|F_{\nabla_i}|^2(z)}
    {d(z,y)^{n-2}}
    \,{\vol}(z)
    &\leqslant
    C m_i^{-2}\rho_j^{2-n}
    \int_{B_{\rho_{j+1}}(y)}e_i                 \\
    &\leqslant
    C m_i^{-2}\rho_j^{2-n}
    \eta m_i\rho_{j+1}^{n-3}                  \\
    &\leqslant
    C\eta\,(m_i\rho_j)^{-1}
    =
    C\eta\,2^{-j}.
\end{aligned}
\]
The possible final annulus between $B_{\rho_{J_i}}(y)$ and $B_{r/2}(y)$ is estimated in the same way. Indeed, $r/4<\rho_{J_i}\leqslant r/2$, and hence $(m_i\rho_{J_i})^{-1}\leqslant 4/(m_ir)=o_i(1)$. Its contribution is therefore at most $C\eta$ for all sufficiently large
$i$. Summing over $j$ gives
\begin{equation}\label{eq:local_Li_potential_clearing_bound}
    \sup_{y\in B_{r/4}(x)} I_i^{\mathrm{loc}}(y)
    \leqslant
    C\eta ,
\end{equation}
where $C$ is independent of $i$, $y$, $r$, and $\eta$.

We next estimate the far term. We use the following consequence of
Lemma~\ref{lemm: Yang-Li}: for every compact set $K\Subset X$ and every
fixed $\rho>0$,
\begin{equation}\label{eq:far_weighted_potential_clearing}
    \sup_{y\in K}
    m_i^{-2}
    \int_{X\setminus B_\rho(y)}
    \frac{|F_{\nabla_i}|^2(z)}
    {d(z,y)^{n-2}}
    \,{\vol}(z)
    \longrightarrow0.
\end{equation}
Indeed, decompose $X\setminus B_\rho(y)$ into dyadic annuli
\[
    B_{2^{j+1}\rho}(y)\setminus B_{2^j\rho}(y),
    \qquad j\geqslant 0.
\]
For $y$ in a fixed compact set $K$, the properties of the radius function
give a constant $C_K$ such that
\[
    B_R(y)
    \subset
    \{\rho\leqslant C_K(1+R)\}
    \qquad
    \text{for every }R>0.
\]
Lemma~\ref{lemm: Yang-Li}, with the bounded part absorbed into the
constant, therefore yields
\[
    \int_{B_R(y)}|F_{\nabla_i}|^2
    \leqslant
    C_K\bigl(m_i+1+R^{n-4}\bigr)
\]
with $C_K$ independent of $i$, $R$, and $y\in K$. Consequently,
\[
\begin{aligned}
& m_i^{-2}
    \int_{X\setminus B_\rho(y)}
    \frac{|F_{\nabla_i}|^2(z)}
    {d(z,y)^{n-2}}
    \,{\vol}(z)                                      \\
&\qquad\leqslant
    C_Km_i^{-2}
    \sum_{j\geqslant0}
    (2^j\rho)^{2-n}
    \bigl(m_i+1+(2^{j+1}\rho)^{n-4}\bigr)                 \\
&\qquad\leqslant
    C_K\left(
        m_i^{-1}\rho^{2-n}
        +m_i^{-2}\rho^{2-n}
        +m_i^{-2}\rho^{-2}
    \right),
\end{aligned}
\]
which tends to zero for fixed $\rho>0$. This proves
\eqref{eq:far_weighted_potential_clearing}.

Applying \eqref{eq:far_weighted_potential_clearing} with
\[
    K=\overline{B_{r/4}(x)}
    \qquad\text{and}\qquad
    \rho=r/2,
\]
we obtain
\begin{equation}\label{eq:far_Li_potential_clearing_bound}
    \sup_{y\in B_{r/4}(x)}I_i^{\mathrm{far}}(y)
    \longrightarrow0.
\end{equation}
Combining \eqref{eq:local_Li_potential_clearing_bound} and
\eqref{eq:far_Li_potential_clearing_bound}, we find
\[
    \sup_{y\in B_{r/4}(x)}\mathcal P_i(y)
    \leqslant
    C\eta+o_i(1).
\]
Choose $\eta>0$ so small that
\[
    C\eta<\frac{1}{2\Lambda_0}.
\]
Then, for all sufficiently large $i$,
\[
    \sup_{y\in B_{r/4}(x)}\mathcal P_i(y)<\Lambda_0^{-1}.
\]
By the definition of $C_{\Lambda_0,i}$, this implies
\[
    B_{r/4}(x)\cap C_{\Lambda_0,i}=\varnothing
\]
for all sufficiently large $i$. Therefore, by
Definition~\ref{def:C_infty}, we have $x\notin\mathcal C$ as we wanted.
\end{proof}

\begin{corollary}[Identification of the Kuratowski upper limit]
\label{cor:obstruction_hierarchy}
\label{cor:nonabelian_locus_obstruction_identity}
In the {\mainsettingref},
\[
    \mathcal S
    \subset
    \mathcal Z
    \subset
    \mathcal C
    =
    \mathcal S\cup\mathcal O.
\]
Equivalently,
\[
    \mathcal C\setminus\mathcal S
    =
    \mathcal O\setminus\mathcal S,
    \qquad
    \mathcal R_{\mathrm{ab}}
    =
    X\setminus(\mathcal S\cup\mathcal O).
\]
Consequently,
\[
    \mathcal C=\mathcal S
    \quad\Longleftrightarrow\quad
    \mathcal O\subset\mathcal S,
\]
and, since $\mathcal S$ is nowhere dense and $\mathcal O$ is closed,
\begin{equation}
\label{eq:interiors_Cinfty_O3}
    \operatorname{int}(\mathcal C)
    =
    \operatorname{int}(\mathcal O).
\end{equation}
\end{corollary}

\begin{proof}
The inclusions
$\mathcal S\subset\mathcal Z\subset\mathcal C$ are
Proposition~\ref{prop:S_subset_Z} and
Lemma~\ref{lemma:Z_subset_C_infty}. Corollary~\ref{cor:effective_monotonicity_on_Rab}
gives
\[
    \mathcal O\subset\mathcal C.
\]
On the other hand, if
$x\notin\mathcal S\cup\mathcal O$, then
Lemma~\ref{lem:effective_monotonicity_gives_vanishing_Morrey_off_S} and
Proposition~\ref{prop:Cinfty_minus_S_obstructs_positive_codim3} give
$x\notin\mathcal C$. Hence
\[
    \mathcal C
    \subset
    \mathcal S\cup\mathcal O,
\]
and therefore
$\mathcal C=\mathcal S\cup\mathcal O$. The identities for the
set differences and $\mathcal R_{\mathrm{ab}}$ follow immediately.

It remains to prove \eqref{eq:interiors_Cinfty_O3}. One inclusion is
immediate from $\mathcal O\subset\mathcal C$. Conversely, let
$U\subset\mathcal C$ be open. Since $\mathcal S$ is nowhere dense,
$U\setminus\mathcal S$ is dense in $U$. The identity
$\mathcal C=\mathcal S\cup\mathcal O$ gives
$U\setminus\mathcal S\subset\mathcal O$, and the closedness of
$\mathcal O$ therefore implies $U\subset\mathcal O$.
\end{proof}

\begin{remark}[The three-dimensional benchmark]
\label{rmk:three_dimensional_benchmark_C_O}
On a three-manifold $r^{3-n}=1$, and therefore
\[
    r\longmapsto
    \vartheta_i(x,r)
    :=
    r^{3-n}m_i^{-1}\int_{B_r(x)}e_i = m_i^{-1}\int_{B_r(x)}e_i
\]
is nondecreasing. Effective codimension-three monotonicity thus holds
globally with zero error for every nonnegative locally integrable energy
density, independently of the field equations, and
$\mathcal O^{(3)}=\varnothing$. For Bogomolny monopoles one has, in
addition, the exact pointwise identity
\[
    |F_{\nabla_i}|^2=|\nabla_i\Phi_i|^2,
\]
so the positive Yang--Mills--Higgs imbalance vanishes identically. Hence
the sufficient Morrey criterion in
Proposition~\ref{prop:codim2_Morrey_criterion_positive_monotonicity}
is also satisfied globally. The corrected concentration theorem gives
$\mathcal S=\mathcal Z$
\cite[Theorem~1.1]{fadeloliveira2026limitv5}, while local clearing of the
mass-renormalized energy and exponential transverse decay hold on every
compact subset of $X\setminus\mathcal S$
\cite[Proposition~2.2 and Theorem~A]{fadel2026abelian}. Consequently, if
$\mathcal C^{(3)}$ denotes the limiting nonabelian locus in the
three-dimensional theory, then
\[
    \mathcal S
    =
    \mathcal Z
    =
    \mathcal C^{(3)},
    \qquad
    \mathcal O^{(3)}=\varnothing,
    \qquad
    \mathcal R_{\mathrm{ab}}^{(3)}=X\setminus\mathcal S.
\]
Corollary~\ref{cor:obstruction_hierarchy} therefore isolates a genuinely
higher-dimensional possibility: an excess of the limiting nonabelian
locus over the support of the limiting calibrated measure, measured
exactly by the
failure of effective codimension-three monotonicity.
\end{remark}

\begin{proposition}[Equivalent local descriptions away from the calibrated support]
\label{prop:equivalent_effective_monotonicity_off_S}
Let $x\in X\setminus\mathcal S$. The following conditions are equivalent.
\begin{enumerate}[label=\textnormal{(\roman*)}]
\item
There exists an open neighbourhood $U$ of $x$ on which effective
codimension-three monotonicity holds; equivalently,
$x\notin\mathcal O$.

\item
There exist $r_x>0$, $C_x<\infty$, $i_x\in\mathbb N$, and
$\gamma_x>0$ such that, for every $i\geqslant i_x$, every
$0<r\leqslant r_x$, every $y\in B_{r/2}(x)$, and every
$0<s\leqslant r/2$, one has
\begin{equation}
\label{ineq:pointwise_nearby_centres_effective_monotonicity}
    s^{3-n}m_i^{-1}\int_{B_s(y)}e_i
    \leqslant
    C_x
    \left(
        r^{3-n}m_i^{-1}\int_{B_r(x)}e_i
        +
        r^{\gamma_x}
    \right).
\end{equation}

\item
The mass-renormalized energies have vanishing codimension-three Morrey
density uniformly over balls contained in shrinking neighbourhoods of
$x$:
\[
    \lim_{r\downarrow0}
    \limsup_{i\to\infty}
    \mathfrak M_i(x,r)
    =
    0.
\]

\item
$x\in\mathcal R_{\mathrm{ab}}=X\setminus\mathcal C$.

\item
There exists a relatively compact open neighbourhood $U$ of $x$ such
that
\begin{equation}
\label{eq:local_uniform_mass_renormalized_density_characterization}
    \sup_Um_i^{-1}e_i
    \longrightarrow0.
\end{equation}
\end{enumerate}
In particular, $\mathcal R_{\mathrm{ab}}$ is the maximal open subset of
$X$ on which the mass-renormalized energy densities converge locally
uniformly to zero.
\end{proposition}

\begin{proof}
Assume (i). Choose $U$ as in (i) and decrease $r_x>0$ so that
$B_{r_x}(x)\subset U$. Lemma~\ref{lem:effective_monotonicity_nearby_centres}
gives (ii).

Assume (ii). Since $x\notin\mathcal S$, decrease $r_x$ so that
$\overline{B_{r_x}(x)}\cap\mathcal S=\varnothing$. For every fixed
$0<r<r_x$, the Portmanteau theorem gives
\[
    r^{3-n}m_i^{-1}\int_{B_r(x)}e_i\longrightarrow0.
\]
Estimate \eqref{ineq:pointwise_nearby_centres_effective_monotonicity}
therefore implies
\[
    \limsup_{i\to\infty}\mathfrak M_i(x,r/2)
    \leqslant
    C_xr^{\gamma_x},
\]
and (iii) follows by letting $r\downarrow0$.

The implication (iii)$\Rightarrow$(iv) is
Proposition~\ref{prop:Cinfty_minus_S_obstructs_positive_codim3}.
If (iv) holds, choose a relatively compact open neighbourhood
$U\Subset\mathcal R_{\mathrm{ab}}$ of $x$ and apply
Proposition~\ref{prop:mass_renormalized_density_decay_Rab} to
$\overline U$; this gives (v).

Finally, assume (v) and choose an open set
$V$ with $x\in V$ and $\overline V\Subset U$. For all sufficiently large
$i$, one has $\sup_Um_i^{-1}e_i\leqslant1$. The local volume estimate
therefore gives
\[
    m_i^{-1}s^{2-n}
    \int_{B_s(y)}|F_{\nabla_i}|^2
    \leqslant
    C_Us^2
\]
for every ball $B_s(y)\subset V$. By
Proposition~\ref{prop:codim2_Morrey_criterion_positive_monotonicity},
effective codimension-three monotonicity holds on $V$, proving (i).
\end{proof}

\begin{remark}[Equivalent formulations away from the calibrated support]
\label{rmk:pointwise_neighbourhood_effective_monotonicity}
Condition~\textup{(ii)} of
Proposition~\ref{prop:equivalent_effective_monotonicity_off_S}
requires an estimate with the outer ball centred at a fixed point $x$,
while allowing the smaller ball to be centred at any point of
$B_{r/2}(x)$. A priori, this is weaker than effective
codimension-three monotonicity on an open neighbourhood of $x$, since it
does not provide a same centre almost-monotonicity formula at nearby
points and at arbitrarily small scales.

Proposition~\ref{prop:equivalent_effective_monotonicity_off_S} shows that
this distinction disappears away from the calibrated support
$\mathcal S$. More precisely, if $x\notin\mathcal S$, then the estimate in
condition~\textup{(ii)} holds if and only if effective
codimension-three monotonicity holds on some open neighbourhood of $x$.
Thus, on $X\setminus\mathcal S$, either formulation characterizes the
complement of the obstruction locus $\mathcal O$. The same region is
also characterized by local uniform convergence
$m_i^{-1}e_i\to0$, as stated in condition~\textup{(v)}.

Combining Corollary~\ref{cor:obstruction_hierarchy} with the equivalence of
conditions~\textup{(iii)}--\textup{(v)} in
Proposition~\ref{prop:equivalent_effective_monotonicity_off_S} gives the
exact descriptions
\begin{equation}
\label{eq:excess_locus_Morrey_characterization}
\begin{aligned}
    \mathcal C\setminus\mathcal S
    &=
    \mathcal O\setminus\mathcal S
    \\
    &=
    \left\{
        x\in X\setminus\mathcal S:
        \lim_{r\downarrow0}
        \limsup_{i\to\infty}
        \mathfrak M_i(x,r)>0
    \right\}.
\end{aligned}
\end{equation}
Equivalently,
\begin{equation}
\label{eq:excess_locus_density_characterization}
\begin{aligned}
    \mathcal C\setminus\mathcal S
    =
    \bigl\{
        x\in X\setminus\mathcal S:
        &\text{there is no relatively compact open neighbourhood }U\ni x
        \\
        &\text{such that }
        \sup_U m_i^{-1}e_i\longrightarrow0
    \bigr\}.
\end{aligned}
\end{equation}
In particular,
since $\mathcal Z\setminus\mathcal S\subset
\mathcal C\setminus\mathcal S$, every limiting Higgs zero outside
the calibrated support is detected by a nonvanishing codimension-three
Morrey defect at shrinking scales.
\end{remark}

\begin{corollary}[Subsequential Tian--Uhlenbeck concentration on
$\mathcal C$]
\label{cor:Cinfty_subsequential_TU}
Let $x\in\mathcal C$, and let $i_j\to\infty$ and
$p_j\in C_{\Lambda_0,i_j}$ be any sequences such that $p_j\to x$.
After relabelling the sequence by $j$, one has
$x\in\mathcal B_{\mathrm{TU}}$ for the relabelled sequence. Equivalently,
for every fixed $0<r\leqslant r_0$,
\begin{equation}
\label{eq:Cinfty_subsequential_TU_lower_bound}
    \liminf_{j\to\infty}
    e^{cr^2}r^{4-n}
    \int_{B_r(x)}e_{i_j}
    \geqslant
    \varepsilon_0.
\end{equation}
In particular, every point of $\mathcal C$ is a
codimension-four concentration point along a subsequence realizing it
through Li's curvature concentration loci.
\end{corollary}

\begin{proof}
If $x\in\mathcal S$, then
Proposition~\ref{prop:S_subset_TU} gives the stronger conclusion that
$x\in\mathcal B_{\mathrm{TU}}$ for the full sequence, and hence for every
subsequence. We may therefore assume that
$x\in\mathcal C\setminus\mathcal S$.

Suppose, to the contrary, that
\eqref{eq:Cinfty_subsequential_TU_lower_bound} fails. Then there exists
$0<r_*\leqslant r_0$ such that
\[
    \liminf_{j\to\infty}
    e^{cr_*^2}r_*^{4-n}
    \int_{B_{r_*}(x)}e_{i_j}
    <
    \varepsilon_0.
\]
Since $x\notin\mathcal S$ and $\mathcal S$ is closed, we may choose
$0<r\leqslant r_*$ such that
$\overline{B_r(x)}\cap\mathcal S=\varnothing$. By the codimension-four
almost-monotonicity formula
\eqref{eq:codim4_almost_monotonicity}, the same strict inequality holds
with $r$ in place of $r_*$. Passing to a further subsequence, still denoted
by $i_j$, we may therefore assume that
\begin{equation}
\label{eq:Cinfty_TU_contradiction_small_energy}
    e^{cr^2}r^{4-n}
    \int_{B_r(x)}e_{i_j}
    <
    \varepsilon_0
\end{equation}
for every $j$.

Since
\[
    r^{4-n}\mathcal E_{B_r(x)}(\nabla_{i_j},\Phi_{i_j})
    =
    \frac12 r^{4-n}\int_{B_r(x)}e_{i_j}
    <
    \varepsilon_0,
\]
Theorem~\ref{thm: total_epsilon_regularity}, with derivative order zero in
\eqref{ineq:estimates_coulomb}, gives a constant $C_r<\infty$,
independent of $j$, such that
\[
    \sup_{B_{r/2}(x)}e_{i_j}
    \leqslant
    C_r.
\]
Because $m_{i_j}\to\infty$, it follows that
\[
    \sup_{B_{r/2}(x)}m_{i_j}^{-1}e_{i_j}
    \longrightarrow0.
\]
Applying condition~\textup{(v)} implies condition~\textup{(iv)} in
Proposition~\ref{prop:equivalent_effective_monotonicity_off_S} to the
relabelled subsequence, we conclude that $x$ does not belong to the limiting nonabelian locus
of the relabelled subsequence. This contradicts
Definition~\ref{def:C_infty}, because
$p_j\in C_{\Lambda_0,i_j}$ and $p_j\to x$. Hence
\eqref{eq:Cinfty_subsequential_TU_lower_bound} holds.
\end{proof}

\begin{corollary}[Local Hausdorff convergence when $\mathcal O\subset\mathcal S$]
\label{cor:Hausdorff_convergence_unobstructed}
Assume that $\mathcal O\subset\mathcal S$. Then $\mathcal S=\mathcal Z=\mathcal C$ and $\mathcal R_{\mathrm{ab}}=X\setminus \mathcal S$. Moreover, for every compact set $K\Subset X$ satisfying $\mathcal S\subset\operatorname{int}K$, one has
\begin{equation}
\label{eq:local_Hausdorff_convergence_unobstructed}
    d_{\mathcal H}(Z_i\cap K,\mathcal S)\longrightarrow0,
    \qquad
    d_{\mathcal H}
    \bigl(\overline{C_{\Lambda_0,i}}\cap K,\mathcal S\bigr)
    \longrightarrow0.
\end{equation}
In particular, every compact set $L\Subset X\setminus \mathcal S$ is eventually
disjoint from both $Z_i$ and $C_{\Lambda_0,i}$.
\end{corollary}

\begin{proof}
The first part is just Corollary~\ref{cor:obstruction_hierarchy}. We now prove the local Hausdorff convergence part. The small Higgs field alternative in Li's concentration-decay dichotomy implies
that, for all sufficiently large $i$,
\begin{equation}
\label{eq:zeros_contained_moving_Li_loci}
    Z_i\subset C_{\Lambda_0,i}.
\end{equation}
By Proposition~\ref{prop:S_subset_Z},
\[
    \sup_{x\in \mathcal S}d(x,Z_i)\longrightarrow0.
\]
Since $\mathcal S\subset\operatorname{int}K$ and $\mathcal S$ is compact, one has
$d(\mathcal S,X\setminus K)>0$. Thus, for all sufficiently large $i$, the sets
$Z_i\cap K$ and $\overline{C_{\Lambda_0,i}}\cap K$ are nonempty
compact sets and
\[
    \sup_{x\in \mathcal S}d(x,Z_i\cap K)\longrightarrow0.
\]
Together with
\eqref{eq:zeros_contained_moving_Li_loci}, it also gives
\[
    \sup_{x\in \mathcal S}
    d\bigl(x,\overline{C_{\Lambda_0,i}}\cap K\bigr)
    \longrightarrow0.
\]
For the reverse estimates, suppose first that there exist
$\varepsilon>0$, a subsequence, and points
\[
    p_i\in\overline{C_{\Lambda_0,i}}\cap K
\]
such that $d(p_i,\mathcal S)\geqslant\varepsilon$. Passing to a further
subsequence, compactness of $K$ gives $p_i\to p\in K$. Choose
$q_i\in C_{\Lambda_0,i}$ with $d(p_i,q_i)<i^{-1}$. Then $q_i\to p$, so
Definition~\ref{def:C_infty} gives
\[
    p\in\mathcal C=\mathcal S,
\]
a contradiction. Hence
\[
    \sup_{p\in\overline{C_{\Lambda_0,i}}\cap K}d(p,\mathcal S)
    \longrightarrow0.
\]
The corresponding estimate for $Z_i\cap K$ follows from
\eqref{eq:zeros_contained_moving_Li_loci}. This proves
\eqref{eq:local_Hausdorff_convergence_unobstructed}.

Finally, if some compact $L\Subset X\setminus \mathcal S$ met
$C_{\Lambda_0,i}$ for infinitely many $i$, compactness of $L$ would produce a convergent realization with limit in
$\mathcal C=\mathcal S$, a contradiction. Thus the curvature concentration loci, and
hence also the zero sets, eventually avoid $L$.
\end{proof}

\begin{remark}[Interpretation]
\label{rmk:Cinfty_obstruction_interpretation}
Corollary~\ref{cor:obstruction_hierarchy} gives the precise sense in
which limiting Higgs zeros away from the support $\mathcal S$ of the
calibrated current can occur only where effective codimension-three
monotonicity fails locally. Proposition~\ref{prop:S_subset_Z} shows that
points of $\mathcal S$ are approached by Higgs zeros forced by nonzero
total transverse Chern charge on normal fibres outside the exceptional
sets in Li's quantitative transverse energy identity: an outer linking
sphere encloses the monopole clusters at scale $m_i^{-1}$ and has Chern
number equal to the local multiplicity of $T$. Every point of
$\mathcal Z\setminus\mathcal S$, and more generally every point of
$\mathcal C\setminus\mathcal S$, lies in $\mathcal O$.

Corollary~\ref{cor:nonabelian_locus_obstruction_identity} gives
\[
    \mathcal C\setminus\mathcal S
    =
    \mathcal O\setminus\mathcal S.
\]
Thus $\mathcal O\setminus\mathcal S$ is exactly the excess of the
limiting nonabelian locus over the calibrated support.
Proposition~\ref{prop:equivalent_effective_monotonicity_off_S} further
shows that, away from $\mathcal S$, the open set definition of effective
monotonicity is equivalent both to the pointwise estimate
\eqref{ineq:pointwise_nearby_centres_effective_monotonicity} and to
vanishing codimension-three Morrey density on all balls contained in
shrinking neighbourhoods.

Proposition~\ref{prop:codim2_Morrey_criterion_positive_monotonicity}
gives a concrete sufficient condition for excluding this obstruction.
It is a codimension-two Morrey bound for the positive part of the signed
imbalance, required uniformly over all balls in a neighbourhood:
\[
    q_i:=|F_{\nabla_i}|^2-|\nabla_i\Phi_i|^2
\]
implies effective codimension-three monotonicity, and hence places the
point outside $\mathcal O$. By the Chern--Weil reformulation
\eqref{eq: imbalance_chern_weil}, this imbalance is the density represented
by $-\langle F_{\nabla_i}\wedge F_{\nabla_i}\rangle\wedge\Xi$.
In particular, the stronger curvature Morrey bound
\[
    m_i^{-1}s^{2-n}
    \int_{B_s(y)}|F_{\nabla_i}|^2
    \leqslant
    \Lambda_x
\]
is a coarse sufficient condition, since $q_i^+\leqslant |F_{\nabla_i}|^2$.

This proposition gives a criterion for excluding membership in $\mathcal C$. It does not give a Hausdorff
dimension or measure estimate for $\mathcal O$. The fixed centre liminf
estimates of Lemma~\ref{lem: size_energy_concentration} do not control
the uniform neighbourhood estimates entering
Definition~\ref{def:effective_codim3_monotonicity}. Likewise, the
vanishing Hausdorff content of the loci $C_{\Lambda_0,i}$ in
Proposition~\ref{prop:Li_cover_Hausdorff_content_off_S} does not pass
automatically to their Kuratowski upper limit. Thus a size estimate for
$\mathcal O$, and hence for $\mathcal C\setminus\mathcal S$, still
requires an additional argument.
\end{remark}

\subsection{Characteristic radii at points of \texorpdfstring{$\mathcal C\setminus\mathcal S$}{C minus S}}
\label{subsec:characteristic_cores_excess}

Corollary~\ref{cor:obstruction_hierarchy} localizes
$\mathcal C\setminus\mathcal S$ but does not describe the relevant
shrinking radii. Li's characteristic radius provides this information.
We show that every $x\in\mathcal C\setminus\mathcal S$ admits sequences
$p_i\to x$ and $R_i\downarrow0$ such that the weighted curvature
potential on $B_{R_i}(p_i)$ is bounded below, although the mass-renormalized
Yang--Mills--Higgs energy vanishes on every fixed neighbourhood whose
closure is disjoint from $\mathcal S$.

For $p\in X$ and $s\geqslant m_i^{-1}$, define the \textbf{truncated weighted
curvature potential}
\begin{equation}
\label{eq:truncated_Li_potential}
    \mathcal P_i(p;s)
    :=
    m_i^{-2}
    \int_X
    \frac{|F_{\nabla_i}|^2(y)}
    {\max\{d(y,p),s\}^{n-2}}
    \,{\vol}(y).
\end{equation}
Thus
\[
    \mathcal P_i(p)
    =
    \mathcal P_i(p;m_i^{-1}).
\]
For fixed $i$ and $p$, the map
\[
    s\longmapsto\mathcal P_i(p;s)
\]
is continuous and nonincreasing on $[m_i^{-1},\infty)$, and tends to zero
as $s\to\infty$. Indeed, continuity and the limit follow by dominated
convergence, using the integrable kernel at $s=m_i^{-1}$. Hence the
following minimum is attained. Following Li's inverse parameter
convention \cite[Definition~2.16]{li2025large}, for $\Lambda>0$ set
\begin{equation}
\label{eq:Li_characteristic_radius}
    r_{\Lambda,i}(p)
    :=
    \min
    \left\{
        s\geqslant m_i^{-1}:
        \mathcal P_i(p;s)\leqslant\Lambda^{-1}
    \right\}.
\end{equation}
Equivalently,
\[
    C_{\Lambda,i}
    =
    \left\{
        p\in X:
        r_{\Lambda,i}(p)>m_i^{-1}
    \right\}
    =
    \left\{
        p\in X:
        \mathcal P_i(p)>\Lambda^{-1}
    \right\}.
\]
Let $\widetilde r_{\Lambda,i}(p)$ denote Li's exact characteristic
radius after translation to our convention; thus its lower endpoint is
$2m_i^{-1}$ and it is defined using
\eqref{eq:exact_translated_Li_potential}. Since the truncated potentials
coincide for every $s\geqslant2m_i^{-1}$, one has
\begin{equation}
\label{eq:characteristic_radii_comparison}
    r_{\Lambda,i}(p)>2m_i^{-1}
    \quad\Longrightarrow\quad
    r_{\Lambda,i}(p)=\widetilde r_{\Lambda,i}(p).
\end{equation}
Thus the estimates for the characteristic radius in
\cite[Lemma~2.19]{li2025large} apply directly whenever the radius in our
convention is larger than $2m_i^{-1}$; the complementary case already
lies at the mass scale.

\begin{proposition}[Weighted curvature concentration near $\mathcal C\setminus\mathcal S$]
\label{prop:excess_points_determine_characteristic_cores}
Let $x\in\mathcal C\setminus \mathcal S$. Choose $r>0$ such that
\[
    \overline{B_{4r}(x)}\cap \mathcal S
    =
    \varnothing.
\]
Then, after passing to a subsequence, there exist points
\[
    p_i\in B_r(x)\cap C_{\Lambda_0,i},
    \qquad
    p_i\longrightarrow x,
\]
and radii
\[
    R_i
    :=
    r_{2\Lambda_0,i}(p_i)
    \longrightarrow
    0
\]
such that
\begin{equation}
\label{ineq:characteristic_core_potential_lower_bound}
    m_i^{-2}
    \int_{B_{R_i}(p_i)}
    \frac{|F_{\nabla_i}|^2(y)}
    {\max\{d(y,p_i),m_i^{-1}\}^{n-2}}
    \,{\vol}(y)
    \geqslant
    \frac{1}{2\Lambda_0}.
\end{equation}
At the same time,
\begin{equation}
\label{ineq:fixed_scale_energy_vanishing_at_excess_point}
    m_i^{-1}
    \int_{B_{4r}(x)}
    \left(
        |F_{\nabla_i}|^2
        +
        |\nabla_i\Phi_i|^2
    \right)
    \longrightarrow
    0.
\end{equation}
Consequently,
\begin{equation}
\label{ineq:core_energy_vanishing}
    m_i^{-1}
    \int_{B_{R_i}(p_i)}
    \left(
        |F_{\nabla_i}|^2
        +
        |\nabla_i\Phi_i|^2
    \right)
    \longrightarrow
    0.
\end{equation}
Thus every point of $\mathcal C\setminus \mathcal S$ determines shrinking
balls which remain visible to Li's weighted curvature
potential but are invisible to the limiting calibrated measure.
\end{proposition}

\begin{proof}
Since $x\in\mathcal C$, the definition of the Kuratowski upper limit
gives, after passing to a subsequence, points $p_i\in B_r(x)\cap C_{\Lambda_0,i}$
such that $p_i\to x$. Set $R_i := r_{2\Lambda_0,i}(p_i)$. We claim that $R_i\to0$. 

Since the points $p_i$ remain in a fixed compact subset of $X$, the compact centre form of the estimate in
\cite[Lemma~2.19, item~(2)]{li2025large} is uniform along the sequence.
To account explicitly for the convention above, consider two cases. If
$R_i\leqslant2m_i^{-1}$, then $R_i\to0$ immediately. Otherwise,
\eqref{eq:characteristic_radii_comparison} gives
\[
    R_i
    =
    \widetilde r_{2\Lambda_0,i}(p_i).
\]
Moreover, the minimality of $R_i>2m_i^{-1}$ implies
\[
    \mathcal P_i^{\mathrm{Li}}(p_i)
    =
    \mathcal P_i(p_i;2m_i^{-1})
    >
    (2\Lambda_0)^{-1},
\]
so $p_i$ belongs to Li's exact locus with parameter $2\Lambda_0$.
Li's Lemma~2.19(2), applied with this fixed parameter, then gives
$R_i\to0$. This proves the claim in both cases.

Since $p_i\in C_{\Lambda_0,i}$,
\[
    \mathcal P_i(p_i;m_i^{-1})
    >
    \Lambda_0^{-1}.
\]
On the other hand, by the definition of $R_i$,
\[
    \mathcal P_i(p_i;R_i)
    \leqslant
    (2\Lambda_0)^{-1}.
\]
Hence
\[
    \mathcal P_i(p_i;m_i^{-1})
    -
    \mathcal P_i(p_i;R_i)
    \geqslant
    (2\Lambda_0)^{-1}.
\]
The difference of the two kernels vanishes outside
$B_{R_i}(p_i)$ and is bounded above on this ball by
$\max\{d(y,p_i),m_i^{-1}\}^{2-n}$. This proves \eqref{ineq:characteristic_core_potential_lower_bound}.

It remains to prove the energy vanishing. Since $\overline{B_{4r}(x)}\cap\operatorname{spt}\mu=\varnothing$,
one has $\mu(\overline{B_{4r}(x)})=0$. The Portmanteau theorem therefore gives
\[
    0\leqslant\limsup_{i\to\infty}\mu_i(B_{4r}(x))
    \leqslant
    \limsup_{i\to\infty}
    \mu_i(\overline{B_{4r}(x)})
    \leqslant
    \mu(\overline{B_{4r}(x)})=0.
\]
This is
\eqref{ineq:fixed_scale_energy_vanishing_at_excess_point}. Since
$p_i\to x$ and $R_i\to0$, one has $B_{R_i}(p_i)\subset B_{4r}(x)$ for all sufficiently large $i$, and \eqref{ineq:core_energy_vanishing} follows.
\end{proof}

\begin{proposition}[Comparison with the ordinary regularity scale]
\label{prop:regularity_scale_inside_characteristic_core}
Let $x\in\mathcal C\setminus\mathcal S$, and let
$p_i\to x$ and
\[
    R_i=r_{2\Lambda_0,i}(p_i)\longrightarrow0
\]
be furnished by
Proposition~\ref{prop:excess_points_determine_characteristic_cores}.
Set
\[
    \lambda_i:=\mathfrak r_i(p_i).
\]
Choose the universal threshold $\varepsilon_0$ in the definition of
$\mathfrak r_i$ sufficiently small, as permitted by the preceding
$\varepsilon$-regularity results.  Then
\begin{equation}
\label{ineq:regularity_scale_below_characteristic_radius}
    \lambda_i\leqslant R_i
\end{equation}
for every sufficiently large $i$.  In particular,
$\lambda_i\to0$.
\end{proposition}

\begin{proof}
Write $\rho_i:=m_i^{-1}$.  We next record a dimensional constant that
will be used to fix the regularity threshold.  Set
\[
    C_n
    :=
    1+2^{n-4}\sum_{j=0}^{\infty}2^{-2j}
    =
    1+\frac{2^{n-2}}{3}.
\]
Decrease $\varepsilon_0$ once and for all so that
\begin{equation}
\label{ineq:epsilon_threshold_characteristic_core}
    C_n\varepsilon_0
    <
    \frac{1}{2\Lambda_0}.
\end{equation}
This choice is compatible with every preceding smallness requirement on
$\varepsilon_0$.

Suppose, to the contrary, that $\lambda_i>R_i$ along a subsequence.
By the definition of the regularity scale,
\begin{equation}
\label{ineq:small_codim4_energy_below_characteristic_radius}
    s^{4-n}
    \int_{B_s(p_i)}e_i
    <
    \varepsilon_0
    \qquad
    \text{for every }0<s\leqslant R_i.
\end{equation}
Recall that $R_i\geqslant\rho_i$ by
\eqref{eq:Li_characteristic_radius}.  Choose $J_i\geqslant0$ so that
\[
    2^{J_i}\rho_i
    \leqslant
    R_i
    <
    2^{J_i+1}\rho_i.
\]
Decompose $B_{R_i}(p_i)$ as the union of $B_{\rho_i}(p_i)$ and the
annuli
\[
    A_{j,i}
    :=
    B_{\min\{2^{j+1}\rho_i,R_i\}}(p_i)
    \setminus
    B_{2^j\rho_i}(p_i),
    \qquad
    0\leqslant j\leqslant J_i,
\]
omitting the last annulus if it is empty.

On the innermost ball, using
\eqref{ineq:small_codim4_energy_below_characteristic_radius} with
$s=\rho_i$, one obtains
\begin{align*}
    &m_i^{-2}
    \int_{B_{\rho_i}(p_i)}
    \frac{|F_{\nabla_i}|^2}
    {\rho_i^{n-2}}
    \,\vol
    \\
    &\qquad\leqslant
    \rho_i^{4-n}
    \int_{B_{\rho_i}(p_i)}e_i
    <
    \varepsilon_0.
\end{align*}
For $A_{j,i}$, put
\[
    t_{j,i}:=\min\{2^{j+1}\rho_i,R_i\}.
\]
Since $t_{j,i}\leqslant R_i$, the same small energy estimate gives
\[
    \int_{B_{t_{j,i}}(p_i)}e_i
    <
    \varepsilon_0t_{j,i}^{n-4}
    \leqslant
    2^{n-4}\varepsilon_0(2^j\rho_i)^{n-4}.
\]
Moreover, on $A_{j,i}$ one has
$d(p_i,\cdot)\geqslant2^j\rho_i$.  Consequently,
\begin{align*}
    &m_i^{-2}
    \int_{A_{j,i}}
    \frac{|F_{\nabla_i}|^2}
    {d(p_i,y)^{n-2}}
    \,\vol(y)
    \\
    &\qquad\leqslant
    m_i^{-2}(2^j\rho_i)^{2-n}
    \int_{B_{t_{j,i}}(p_i)}e_i
    \\
    &\qquad\leqslant
    2^{n-4}\varepsilon_0
    m_i^{-2}(2^j\rho_i)^{-2}
    \\
    &\qquad=
    2^{n-4}\varepsilon_0 2^{-2j}.
\end{align*}
Summing the inner ball and annular estimates yields
\[
    m_i^{-2}
    \int_{B_{R_i}(p_i)}
    \frac{|F_{\nabla_i}|^2(y)}
    {\max\{d(y,p_i),m_i^{-1}\}^{n-2}}
    \,\vol(y)
    <
    C_n\varepsilon_0.
\]
This contradicts
\eqref{ineq:characteristic_core_potential_lower_bound} and
\eqref{ineq:epsilon_threshold_characteristic_core}.  Hence
$\lambda_i\leqslant R_i$ for all sufficiently large $i$.  Since
$R_i\to0$, the final conclusion follows.
\end{proof}

\begin{proposition}[Characteristic radius larger than the scale $m_i^{-1}$]
\label{prop:characteristic_core_above_mass_scale}
Retain the hypotheses and notation of
Proposition~\ref{prop:regularity_scale_inside_characteristic_core}, and
assume, after passage to a further subsequence, that
\begin{equation}
\label{eq:characteristic_regular_scale_above_mass_scale}
    \tau_i
    :=
    m_i\lambda_i
    \longrightarrow
    +\infty.
\end{equation}
Choose an oriented orthonormal frame
$(E_1,\ldots,E_n)$ near $x$, let
$\iota_i:\mathbb R^n\to T_{p_i}X$ be the oriented linear isometry
satisfying $\iota_i(e_a)=E_a(p_i)$, and set
\[
    \delta_i^{(m)}(z)
    :=
    \exp_{p_i}\bigl(m_i^{-1}\iota_i(z)\bigr),
\]
\[
    g_i^{(m)}
    :=
    m_i^2(\delta_i^{(m)})^*g,
    \qquad
    \nabla_i^{(m)}
    :=
    (\delta_i^{(m)})^*\nabla_i,
    \qquad
    \Phi_i^{(m)}
    :=
    m_i^{-1}(\delta_i^{(m)})^*\Phi_i.
\]
Then the following conclusions hold.

\begin{enumerate}[label=\textnormal{(\alph*)}]
\item After passing to a further subsequence and applying gauge
transformations, $(\nabla_i^{(m)},\Phi_i^{(m)})$ converges smoothly on
compact subsets of $\mathbb R^n$ to a flat connection
$\nabla_\infty$ and a parallel Higgs field $\Phi_\infty$ satisfying
$|\Phi_\infty|\leqslant1$.

\item For every fixed $L>1$,
\begin{equation}
\label{ineq:characteristic_potential_outside_mass_ball}
\liminf_{i\to\infty}
    m_i^{-2}
    \int_{B_{R_i}(p_i)\setminus B_{L/m_i}(p_i)}
    \frac{|F_{\nabla_i}|^2(y)}{d(y,p_i)^{n-2}}
    \,\vol(y)
    \geqslant
    \frac{1}{2\Lambda_0}.
\end{equation}

\item For $s>0$, define the annular codimension-four curvature energy
\[
    \mathcal A_i(p_i;s)
    :=
    s^{4-n}
    \int_{B_{2s}(p_i)\setminus B_s(p_i)}
    |F_{\nabla_i}|^2\,\vol.
\]
Then, for every fixed $L>1$,
\begin{equation}
\label{ineq:large_annular_codim4_energy_characteristic_core}
    \liminf_{i\to\infty}
    \sup_{L/m_i\leqslant s<R_i}
    \mathcal A_i(p_i;s)
    \geqslant
    \frac{3}{8\Lambda_0}L^2.
\end{equation}
Consequently, after passing to a subsequence, there are numbers
$L_i\to+\infty$ and scales $s_i$ satisfying
\[
    \frac{L_i}{m_i}
    \leqslant
    s_i
    <
    R_i
\]
such that
\begin{equation}
\label{eq:divergent_annular_codim4_energy_characteristic_core}
    m_is_i\longrightarrow+\infty,
    \qquad
    s_i^{4-n}
    \int_{B_{2s_i}(p_i)\setminus B_{s_i}(p_i)}
    |F_{\nabla_i}|^2\,\vol
    \longrightarrow+\infty.
\end{equation}
\end{enumerate}
\end{proposition}

\begin{proof}
By
Proposition~\ref{prop:regularity_scale_inside_characteristic_core},
\[
    m_iR_i
    \geqslant
    m_i\lambda_i
    =
    \tau_i
    \longrightarrow+\infty.
\]
In particular, for every fixed $L>1$, the ball
$B_{L/m_i}(p_i)$ is contained in $B_{R_i}(p_i)$ for all sufficiently
large $i$.

We first prove part~\textnormal{(a)}.  Fix $R>0$.  By
\eqref{eq:characteristic_regular_scale_above_mass_scale}, one has
$R/m_i<\lambda_i$ for all sufficiently large $i$.  The definition of
$\lambda_i$ therefore gives
\[
    \left(\frac{R}{m_i}\right)^{4-n}
    \int_{B_{R/m_i}(p_i)}e_i
    <
    \varepsilon_0.
\]
After the mass rescaling, this becomes
\begin{equation}
\label{ineq:mass_rescaled_small_energy_characteristic_core}
    R^{4-n}
    \int_{B_R(0)}
    e(\nabla_i^{(m)},\Phi_i^{(m)})
    \,\operatorname{vol}_{g_i^{(m)}}
    <
    \varepsilon_0.
\end{equation}
The rescaled metrics converge smoothly on compact subsets to the
Euclidean metric.  Moreover, the maximum principle gives
$|\Phi_i^{(m)}|\leqslant1$.  Applying
Theorem~\ref{thm: total_epsilon_regularity} on successively larger fixed
balls, followed by local Uhlenbeck gauge fixing and a diagonal argument,
produces smooth convergence modulo gauge on compact subsets to a smooth
solution $(\nabla_\infty,\Phi_\infty)$ of
\eqref{eq:YMH_second_order} on $\mathbb R^n$.

Passing to the limit in the interior estimate associated with
\eqref{ineq:mass_rescaled_small_energy_characteristic_core} gives, for
every $R>0$,
\[
    \sup_{B_{R/2}(0)}
    \bigl(
        |F_{\nabla_\infty}|^2
        +
        |\nabla_\infty\Phi_\infty|^2
    \bigr)
    \leqslant
    C\varepsilon_0R^{-4},
\]
where $C$ is independent of $R$.  Given any $z\in\mathbb R^n$, let
$R\to\infty$ with $z\in B_{R/2}(0)$.  It follows that
$F_{\nabla_\infty}=0$ and
$\nabla_\infty\Phi_\infty=0$ everywhere.  The bound
$|\Phi_\infty|\leqslant1$ follows from the maximum principle bound on the
rescaled sequence.  This proves part~\textnormal{(a)}.

For fixed $L>1$, the change of variables under $\delta_i^{(m)}$ gives
\begin{align}
\label{eq:inner_characteristic_potential_mass_rescaling}
    &m_i^{-2}
    \int_{B_{L/m_i}(p_i)}
    \frac{|F_{\nabla_i}|^2(y)}
    {\max\{d(y,p_i),m_i^{-1}\}^{n-2}}
    \,\vol(y)
    \\
    &\qquad=
    \int_{B_L(0)}
    \frac{|F_{\nabla_i^{(m)}}|_{g_i^{(m)}}^2(z)}
    {\max\{|z|,1\}^{n-2}}
    \,\operatorname{vol}_{g_i^{(m)}}(z).
\end{align}
For fixed $L$, the right-hand side tends to zero by the smooth
convergence in part~\textnormal{(a)}.  Subtracting
\eqref{eq:inner_characteristic_potential_mass_rescaling} from
\eqref{ineq:characteristic_core_potential_lower_bound}, and observing
that $d(y,p_i)\geqslant Lm_i^{-1}>m_i^{-1}$ on the complementary
annulus, proves
\eqref{ineq:characteristic_potential_outside_mass_ball}.

It remains to prove part~\textnormal{(c)}.  Fix $L>1$.  For all
sufficiently large $i$, let $J_i$ be the largest nonnegative integer
such that
\[
    s_{J_i,i}
    :=
    2^{J_i}\frac{L}{m_i}
    <
    R_i,
\]
and set
\[
    s_{j,i}:=2^j\frac{L}{m_i},
    \qquad
    A_{j,i}
    :=
    B_{\min\{2s_{j,i},R_i\}}(p_i)
    \setminus
    B_{s_{j,i}}(p_i),
    \qquad
    0\leqslant j\leqslant J_i.
\]
These annuli cover
$B_{R_i}(p_i)\setminus B_{L/m_i}(p_i)$.  Since
$d(p_i,y)\geqslant s_{j,i}$ on $A_{j,i}$, one has
\begin{align*}
    &m_i^{-2}
    \int_{A_{j,i}}
    d(y,p_i)^{2-n}|F_{\nabla_i}|^2(y)
    \,\vol(y)
    \\
    &\qquad\leqslant
    m_i^{-2}s_{j,i}^{2-n}
    \int_{B_{2s_{j,i}}(p_i)\setminus B_{s_{j,i}}(p_i)}
    |F_{\nabla_i}|^2\,\vol
    \\
    &\qquad=
    \frac{\mathcal A_i(p_i;s_{j,i})}
    {(m_is_{j,i})^2}.
\end{align*}
Summing and using $m_is_{j,i}=2^jL$ gives
\begin{align*}
    &m_i^{-2}
    \int_{B_{R_i}(p_i)\setminus B_{L/m_i}(p_i)}
    d(y,p_i)^{2-n}|F_{\nabla_i}|^2(y)
    \,\vol(y)
    \\
    &\qquad\leqslant
    \left(
        \sup_{L/m_i\leqslant s<R_i}
        \mathcal A_i(p_i;s)
    \right)
    \sum_{j=0}^{J_i}\frac{1}{(2^jL)^2}
    \\
    &\qquad\leqslant
    \frac{4}{3L^2}
    \sup_{L/m_i\leqslant s<R_i}
    \mathcal A_i(p_i;s).
\end{align*}
Taking the lower limit and applying
\eqref{ineq:characteristic_potential_outside_mass_ball} yields
\eqref{ineq:large_annular_codim4_energy_characteristic_core}.

Finally, apply
\eqref{ineq:large_annular_codim4_energy_characteristic_core} with
$L=1,2,3,\ldots$ and use a diagonal subsequence.  This gives numbers
$L_i\to+\infty$ and scales
$s_i\in[L_i/m_i,R_i)$ for which the left-hand side defining
$\mathcal A_i(p_i;s_i)$ tends to infinity.  Since
$m_is_i\geqslant L_i$, one obtains
\eqref{eq:divergent_annular_codim4_energy_characteristic_core}.
\end{proof}

\begin{remark}[Comparison with three-dimensional scale selection]
\label{rmk:comparison_three_dimensional_scale_selection}
The moving centre analyses in the three-dimensional large mass problem
show that every nontrivial finite mass Euclidean monopole profile occurs
at a scale comparable to $m_i^{-1}$ and that a mass one profile does not
bubble again at smaller scales; see
\cite[Lemma~8.7 and Proposition~8.3]{fadeloliveira2026limitv5}.  A related scale
comparison, bubble extraction and no-neck identity for three-dimensional
Yang--Mills--Higgs critical points with Higgs self-interaction are proved
in
\cite[Lemma~5.17 and Propositions~5.19 and~5.27]{cheng2025su2}.

Propositions~\ref{prop:regularity_scale_inside_characteristic_core} and
\ref{prop:characteristic_core_above_mass_scale} identify the part of
that mechanism which remains valid at an arbitrary point of $\mathcal C\setminus\mathcal S$ in
ambient dimensions six and seven.  The ordinary regularity scale is
forced to be no larger than Li's characteristic radius.  However, the higher-dimensional
equations do not presently exclude
$m_i\lambda_i\to+\infty$; instead, the weighted curvature potential lower bound
forces the annular conclusion
\eqref{ineq:large_annular_codim4_energy_characteristic_core}.  Thus, when the ordinary regularity scale is larger than $m_i^{-1}$, the selected ball cannot be described by
a finite collection of annular regions with uniformly bounded
codimension-four scale-invariant curvature energy.
\end{remark}

\begin{remark}[Relation with ordinary Yang--Mills--Higgs concentration]
\label{rmk:characteristic_cores_imply_subsequential_TU}
Proposition~\ref{prop:regularity_scale_inside_characteristic_core}
gives a scale level strengthening of
Corollary~\ref{cor:Cinfty_subsequential_TU} at points of
$\mathcal C\setminus\mathcal S$.  Namely, the collapsing ordinary
regularity scale $\lambda_i=\mathfrak r_i(p_i)$ may be chosen at the same centres $p_i$ selected by Li's characteristic radius and satisfies
\[
    \lambda_i
    \leqslant
    r_{2\Lambda_0,i}(p_i)
    \longrightarrow0.
\]
Thus the ordinary codimension-four loss of regularity occurs inside the
region carrying the fixed positive weighted curvature potential contribution.

This comparison does not imply
$R_i=O(m_i^{-1})$.  Li's general estimate
\cite[Lemma~2.19\textup{(2)}]{li2025large} gives, for centres in a fixed
compact set and a fixed characteristic parameter, an upper bound of
order $m_i^{-1/2}$.  If the stronger inverse mass bound fails in the
sense that $m_i\lambda_i\to+\infty$, then
Proposition~\ref{prop:characteristic_core_above_mass_scale} gives the
more precise annular conclusions
\eqref{ineq:characteristic_potential_outside_mass_ball} and
\eqref{eq:divergent_annular_codim4_energy_characteristic_core}.
\end{remark}

\begin{remark}[Weighted potential versus calibrated mass]
\label{rmk:Li_potential_vs_calibrated_mass}
Proposition~\ref{prop:excess_points_determine_characteristic_cores}
does not contradict the lower bound
\eqref{ineq:characteristic_core_potential_lower_bound}, because Li's
weighted potential and the mass-renormalized energy have different
scaling.

Indeed, suppose that curvature localized at a scale
$\rho_i\geqslant m_i^{-1}$ satisfies
\[
    m_i^{-2}\rho_i^{2-n}
    \int_{B_{\rho_i}(p_i)}
    |F_{\nabla_i}|^2
    \geqslant
    c_0>0.
\]
Then
\[
\begin{aligned}
    m_i^{-1}
    \int_{B_{\rho_i}(p_i)}
    |F_{\nabla_i}|^2
    &=
    \rho_i^{n-2}
    \left(
        m_i^{-1}\rho_i^{2-n}
        \int_{B_{\rho_i}(p_i)}
        |F_{\nabla_i}|^2
    \right)
    \\
    &\geqslant
    c_0m_i\rho_i^{n-2}.
\end{aligned}
\]
Thus such a single scale contribution forces a uniform positive amount
of mass-renormalized energy only when
\[
    \liminf_{i\to\infty}
    m_i\rho_i^{n-2}
    >
    0.
\]
If
\[
    m_i\rho_i^{n-2}
    \longrightarrow
    0,
\]
the weighted curvature potential lower bound alone gives no positive uniform lower
bound for the mass-renormalized energy. Hence the existence of the balls $B_{R_i}(p_i)$ at points of $\mathcal C\setminus \mathcal S$ is not excluded by scaling. At points of $\mathcal Z\setminus\mathcal S$, such balls coexist,
after passage to a subsequence along which the zero is realized, with the Tian--Uhlenbeck scales centred at Higgs zeros of Proposition~\ref{prop:zero_centered_TU_concentration}. Points of $\mathcal C\setminus\mathcal Z$, if they occur, would instead
have a fixed neighbourhood which is eventually free of Higgs zeros. Determining whether
the $\Theta$-monopole equation rules out either phenomenon is a separate
problem, discussed further in Section~\ref{sec:resulting_picture}.
\end{remark}

Having described the possible excess
$\mathcal C\setminus\mathcal S$, we now turn to the analysis on
$\mathcal R_{\mathrm{ab}}=X\setminus\mathcal C$.

\subsection{Codimension-two Higgs monotonicity and stress-energy control with vanishing error}
\label{subsec:codim2_Higgs_monotonicity_Cinfty}

By Lemma~\ref{lem:maximal_abelian_clearing_region}, Li's large Higgs field and
transverse decay estimates are eventually uniform on compact subsets of
the open set $\mathcal R_{\mathrm{ab}}$. The results below are therefore
local regularity and structure statements on the abelian region
$\mathcal R_{\mathrm{ab}}$ and require no global size hypothesis.
Corollary~\ref{cor:effective_monotonicity_on_Rab} already shows that the
full mass-renormalized energy has vanishing codimension-two Morrey norm
there. The localized Higgs energy identity provides, in addition, a
genuine codimension-two almost-monotonicity formula whose error is
exponentially small in the mass.

\begin{proposition}[Codimension-two Higgs monotonicity on
$\mathcal R_{\mathrm{ab}}$]
\label{prop:Higgs_monotonicity_away_from_Cinfty}
Let $K\Subset \mathcal R_{\mathrm{ab}}$. Then there exist
$r_K>0$, a constant $C_0<\infty$ depending only on the bounded geometry
of $(X^n,g)$ (and the dimension $n$), and a sequence $\varepsilon_i(K)\to0$ such that, for every
$x\in K$ and every $0<s<r<r_K$,
\begin{equation}\label{ineq:Higgs_energy_monotonicity_away_Cinfty}
m_i^{-1}s^{2-n}
\int_{B_s(x)}
|\nabla_i\Phi_i|^2
\leqslant
C_0
m_i^{-1}r^{2-n}
\int_{B_r(x)}
|\nabla_i\Phi_i|^2
+
\varepsilon_i(K).
\end{equation}
\end{proposition}

\begin{proof}
For simplicity, write
\[
    H_{x,i}(r)
    :=
    \int_{B_r(x)}
    |\nabla_i\Phi_i|^2.
\]
By the bounded geometry assumptions in the {\mainsettingref}, the constants in
Lemma~\ref{lem:localized_Higgs_energy_identity} may be chosen uniformly
with respect to the center $x\in X$. Thus there exist constants
$r_0>0$ and $c'>0$, depending only on the geometry of $(X^n,g)$ (and the dimension $n$), such that
for every $x\in X$ and every $0<r<r_0$,
\begin{align}
\frac{d}{dr}
\left(
e^{c'r^2}r^{2-n}H_{x,i}(r)
\right)
&\geqslant
-
2e^{c'r^2}r^{1-n}
\left|
\int_{B_r(x)}
r_x
\left\langle
[\partial_{r_x}\lrcorner F_{\nabla_i},\Phi_i],
\nabla_i\Phi_i
\right\rangle
\right|.
\label{ineq:Higgs_monotonicity_weighted}
\end{align}
We now estimate the error term on compact subsets of
$\mathcal R_{\mathrm{ab}}$. Since $[\partial_{r_x}\lrcorner F_{\nabla_i},\Phi_i]
\perp
\Phi_i$, 
it pairs only with the transverse component of $\nabla_i\Phi_i$. Hence
\[
\left\langle
[\partial_{r_x}\lrcorner F_{\nabla_i},\Phi_i],
\nabla_i\Phi_i
\right\rangle
=
\left\langle
[\partial_{r_x}\lrcorner F_{\nabla_i},\Phi_i],
(\nabla_i\Phi_i)^\perp
\right\rangle.
\]
Moreover, using the $\Theta$-monopole equation and the bound
$|\Phi_i|\leqslant m_i$, we have
\[
|(\nabla_i\Phi_i)^\perp|
\leqslant
C |(F_{\nabla_i})^\perp|,
\qquad
|[\partial_{r_x}\lrcorner F_{\nabla_i},\Phi_i]|
\leqslant
C m_i |(F_{\nabla_i})^\perp|.
\]
Therefore
\begin{align}\label{ineq:error_term_transverse}
m_i^{-1}
\left|
\int_{B_r(x)}
r_x
\left\langle
[\partial_{r_x}\lrcorner F_{\nabla_i},\Phi_i],
\nabla_i\Phi_i
\right\rangle
\right|
&\leqslant
C r
\int_{B_r(x)}
|(F_{\nabla_i})^\perp|^2.
\end{align}
Since $K\Subset \mathcal R_{\mathrm{ab}}$, choose an open set
$U$ such that
\[
    K\Subset U\Subset \mathcal R_{\mathrm{ab}}.
\]
Then, for all sufficiently large $i$,
\[
    U\cap C_{\Lambda_0,i}=\varnothing.
\]
Choose $0<r_K<r_0$ so that $B_{r_K}(x)\Subset U$ for every $x\in K$. By Li's transverse decay estimates, there exist
constants $a_K>0$ and $B_K<\infty$ such that
\[
    |(F_{\nabla_i})^\perp|^2(p)
    \leqslant
    B_K m_i e^{-a_Km_i}
\]
for every $p\in U$ and all sufficiently large $i$. Consequently, for
every $x\in K$ and every $0<r<r_K$,
\[
\int_{B_r(x)}
|(F_{\nabla_i})^\perp|^2
\leqslant
B_K r^n m_i e^{-a_Km_i}.
\]
Substituting this into \eqref{ineq:error_term_transverse}, we obtain
\[
m_i^{-1}
\left|
\int_{B_r(x)}
r_x
\left\langle
[\partial_{r_x}\lrcorner F_{\nabla_i},\Phi_i],
\nabla_i\Phi_i
\right\rangle
\right|
\leqslant
B_K r^{n+1}m_i e^{-a_Km_i}.
\]
Multiplying \eqref{ineq:Higgs_monotonicity_weighted} by $m_i^{-1}$ and
using the previous estimate gives, for $0<r<r_K$,
\[
\frac{d}{dr}
\left(
m_i^{-1}e^{c'r^2}r^{2-n}H_{x,i}(r)
\right)
\geqslant
-
B_K r^2 m_i e^{-a_Km_i}.
\]
Integrating from $s$ to $r$ yields
\[
m_i^{-1}e^{c's^2}s^{2-n}H_{x,i}(s)
\leqslant
m_i^{-1}e^{c'r^2}r^{2-n}H_{x,i}(r)
+
B_K r_K^3m_i e^{-a_Km_i}.
\]
Set
\[
    \varepsilon_i(K):=
    B_K r_K^3m_i e^{-a_Km_i}.
\]
Then $\varepsilon_i(K)\to0$. Since $0<s<r<r_K<r_0$, we have
\[
    e^{c'(r^2-s^2)}\leqslant e^{c'r_0^2}.
\]
Thus
\[
m_i^{-1}s^{2-n}H_{x,i}(s)
\leqslant
e^{c'r_0^2}
m_i^{-1}r^{2-n}H_{x,i}(r)
+
\varepsilon_i(K).
\]
Taking
\[
    C_0:=e^{c'r_0^2}
\]
proves \eqref{ineq:Higgs_energy_monotonicity_away_Cinfty}.
\end{proof}

\begin{remark}
The key point is that the error term in the localized Higgs energy
monotonicity formula only involves transverse components. Away from the
limiting nonabelian locus $\mathcal C$, these transverse components decay
exponentially fast with the mass. Thus, unlike the full stress-energy
monotonicity formula, the Higgs energy identity does not contain a
tangential Chern--Weil or Yang--Mills obstruction.

This is naturally a codimension-two estimate, rather than a
codimension-three estimate. Indeed, the Higgs term is the Dirichlet
energy of a section and therefore has the stress-energy scaling of a
$1$-form, leading to the factor
\[
    r^{2-n}
    \int_{B_r(x)}
    |\nabla\Phi|^2.
\]
By contrast, the curvature energy has the natural codimension-four scaling
\[
    r^{4-n}
    \int_{B_r(x)}
    |F_\nabla|^2.
\]
Therefore, the expected codimension-three behaviour of large mass monopoles is
not the intrinsic scaling of either component separately, but rather emerges
from the first-order coupling between $F_\nabla$ and $\nabla\Phi$ in the
$\Theta$-monopole equation. The difficulty in proving an effective
codimension-three monotonicity formula lies precisely in reconciling these
two stress-energy scalings.
\end{remark}

\begin{corollary}[Codimension-three monotonicity with vanishing error on $\mathcal R_{\mathrm{ab}}$]
\label{cor:codim3_obstruction_reduction_away_Cinfty}
Let $K\Subset\mathcal R_{\mathrm{ab}}$. Then there exist
$r_K>0$ and a sequence $\varepsilon_i(K)\to0$ such that, for every
$x\in K$ and every $0<r<r_K$,
\begin{equation}
\label{eq:qplus_Morrey_vanishes_Rab}
    m_i^{-1}r^{2-n}
    \int_{B_r(x)}q_i^+
    \leqslant
    \varepsilon_i(K)r^2.
\end{equation}
Consequently, for almost every $0<r<r_K$,
\begin{equation}
\label{eq:vanishing_error_differential_codim3_Rab}
    \frac{d}{dr}
    \left(
        e^{cr^2}\theta_i(x,r)
    \right)
    \geqslant
    -e^{cr_0^2}\varepsilon_i(K)r^2,
\end{equation}
uniformly for $x\in K$. Equivalently, for every
$0<s<r<r_K$,
\begin{equation}
\label{eq:vanishing_error_integrated_codim3_Rab}
    e^{cs^2}\theta_i(x,s)
    \leqslant
    e^{cr^2}\theta_i(x,r)
    +
    \frac{e^{cr_0^2}}{3}\varepsilon_i(K)(r^3-s^3).
\end{equation}
\end{corollary}

\begin{proof}
Estimate \eqref{eq:qplus_Morrey_vanishes_Rab} follows immediately from
\eqref{eq:clearing_region_Morrey_decay} and
$q_i^+\leqslant e_i$. The stress-energy inequality
\eqref{eq:differential_almost_monotonicity_short} then gives
\eqref{eq:vanishing_error_differential_codim3_Rab}, and integration
from $s$ to $r$ gives
\eqref{eq:vanishing_error_integrated_codim3_Rab}.
\end{proof}

\begin{remark}[Mass-renormalized and unrenormalized longitudinal curvature]
On $\mathcal R_{\mathrm{ab}}$, the entire positive stress-energy
imbalance vanishes at the mass-renormalized codimension-two Morrey scale;
there is therefore no residual codimension-three monotonicity obstruction
on this region. Nevertheless, the unrenormalized longitudinal curvature $\widehat\omega_i = \langle F_{\nabla_i},\Psi_i\rangle$
need not remain locally bounded. It is the only curvature component not
forced to decay exponentially and may grow, albeit only at the
subcritical rate $o(m_i^{1/2})$ on compact subsets. This unrenormalized
component governs ordinary regularity scales and gauge compactness, and is
the object studied in
Sections~\ref{sec: abelianization}
and~\ref{sec:frequency_away_Cinfty}.
\end{remark}


\section{Cohomogeneity-one examples and absence of the obstruction}
\label{sec:homogeneous_examples}

In this section we verify that the known cohomogeneity-one large mass
families satisfy the Morrey criterion of
Proposition~\ref{prop:codim2_Morrey_criterion_positive_monotonicity}.
The relevant statement is not a pointwise bound for the signed imbalance that is uniform in the mass,
imbalance
\[
    q_m:=|F_{\nabla_m}|^2-|\nabla_m\Phi_m|^2.
\]
For a $\Theta$-monopole, the algebraic identity
\eqref{eq:general_calibration_identity} and the monopole equation give
an equivalent Chern--Weil expression,
\begin{equation}
\label{eq:homogeneous_Chern_Weil_imbalance}
    q_m\,\vol
    =
    -\langle F_{\nabla_m}\wedge F_{\nabla_m}\rangle\wedge\Xi,
\end{equation}
where $\Xi=\varphi$ in the $\mathrm G_2$ case and $\Xi=\omega$ in
the Calabi--Yau case. Below we compute the full pointwise norms of
$F_{\nabla_m}$ and $\nabla_m\Phi_m$; equation
\eqref{eq:homogeneous_Chern_Weil_imbalance} provides an independent check
of the signs and normalizations.

Indeed, after the leading transverse Bogomolny terms cancel, the residual
imbalance may still be of order $m^2$ in a tubular neighbourhood of the
zero section of radius comparable to $m^{-1}$. What is
uniform is the codimension-two Morrey control over balls with arbitrary
centres in a fixed tubular neighbourhood. More precisely, if $N$ denotes
the calibrated zero section and $r_N=\operatorname{dist}(\,\cdot\,,N)$, then the explicit homogeneous
formulae give, on a fixed tubular neighbourhood of $N$,
\begin{equation}
\label{eq:homogeneous_tubular_q_bound}
    q_m^+
    \leqslant
    C\bigl(1+r_N^{-2}+m r_N^{-1}\bigr).
\end{equation}
Since $N$ has codimension three, the two singular weights on the right-hand
side have precisely the integrability required by the Morrey criterion.

The examples considered below are the Bryant--Salamon monopoles on
$\Lambda_-^2(S^4)$ and $\Lambda_-^2(\mathbb{CP}^2)$ constructed in
\cite{oliveira2014monopoles}, and the $\mathrm{SU}(2)^2$-invariant Stenzel
Calabi--Yau monopoles constructed in \cite{oliveira2016calabi}. The latter
were subsequently identified by Stein as the unique irreducible,
quadratically decaying, invariant Calabi--Yau monopole family with nonzero
Higgs field on the known cohomogeneity-one AC Calabi--Yau threefolds; see
\cite{stein2023invariant}.

We first isolate the geometric estimate which converts
\eqref{eq:homogeneous_tubular_q_bound} into absence of the obstruction.

\begin{lemma}[Tubular Morrey control excludes the obstruction]
\label{lem:tubular_Morrey_control_excludes_obstruction}
Let $(\nabla_i,\Phi_i)$ be a large mass sequence in the {\mainsettingref}, with
masses $m_i\to+\infty$, and let $N\subset X$ be a compact embedded
submanifold of codimension three. Set $r_N=\operatorname{dist}(\,\cdot\,,N)$.
Assume that there exist a tubular neighbourhood $U$ of $N$ and a constant
$C<\infty$ such that
\begin{equation}
\label{eq:tubular_Morrey_hypothesis}
    q_i^+
    \leqslant
    C\bigl(1+r_N^{-2}+m_i r_N^{-1}\bigr)
    \qquad\text{on }U\setminus N
\end{equation}
for all sufficiently large $i$. Assume also that, for every compact set
$K\Subset X\setminus N$, there is $C_K<\infty$ such that
\begin{equation}
\label{eq:away_Morrey_hypothesis}
    q_i^+\leqslant C_K(1+m_i)
    \qquad\text{on }K
\end{equation}
for all sufficiently large $i$. Then
\[
    \mathcal O=\varnothing.
\]
Consequently,
\[
    \mathcal S=\mathcal Z=\mathcal C.
\]
\end{lemma}

\begin{proof}
Fix a relatively compact tubular neighbourhood
\[
    U_0\Subset U
\]
of $N$. After shrinking $U_0$ if necessary, choose
\[
    0<r_0<
    \frac12
    \operatorname{dist}
    \bigl(
        \overline{U_0},
        X\setminus U
    \bigr).
\]
Since $N$ is compact, $\overline{U_0}$ is covered by finitely many
Fermi coordinate charts whose metric coefficients and Jacobians are
uniformly controlled. It follows that there is $C_0<\infty$ such that,
for every $y\in U_0$, every $0<s<r_0$, and every
$0\leqslant\alpha<3$,
\begin{equation}
\label{eq:tubular_weight_integral}
    \int_{B_s(y)}r_N^{-\alpha}\,\vol
    \leqslant
    C_0s^{n-\alpha}.
\end{equation}
Indeed, the choice of $r_0$ ensures that $B_s(y)\subset U$. If
$\operatorname{dist}(y,N)\geqslant2s$, then
$r_N\geqslant s$ on $B_s(y)$ and the estimate is immediate.
Otherwise, $B_s(y)$ is contained in one of the uniformly controlled
Fermi coordinate regions, with tangential size $O(s)$ and normal size
$O(s)$. Since the normal bundle has rank three,
\[
    \int_0^{Cs}t^{2-\alpha}\,dt
    \asymp
    s^{3-\alpha},
\]
which gives the required normal contribution and proves
\eqref{eq:tubular_weight_integral} with a constant independent of
$y$ and $s$.

Using \eqref{eq:tubular_Morrey_hypothesis} and
\eqref{eq:tubular_weight_integral} with $\alpha=1,2$, we obtain
\begin{align*}
    m_i^{-1}s^{2-n}\int_{B_s(y)}q_i^+\,\vol
    &\leqslant
    C m_i^{-1}s^{2-n}
    \bigl(s^n+s^{n-2}+m_is^{n-1}\bigr)\\
    &\leqslant
    C\bigl(m_i^{-1}s^2+m_i^{-1}+s\bigr).
\end{align*}
This is uniformly bounded for $y\in U_0$, $0<s<r_0$, and all sufficiently
large $i$.

If $x\in X\setminus N$, choose $r_x>0$ with
$\overline{B_{2r_x}(x)}\subset X\setminus N$. From
\eqref{eq:away_Morrey_hypothesis}, for $y\in B_{r_x}(x)$ and
$0<s<r_x$,
\[
    m_i^{-1}s^{2-n}\int_{B_s(y)}q_i^+\,\vol
    \leqslant
    C_xm_i^{-1}(1+m_i)s^2
    \leqslant
    2C_xr_x^2
\]
for all large $i$. Thus, near every point of $X$, the codimension-two Morrey bound
\eqref{ineq:positive_imbalance_Morrey_bound} holds uniformly over all balls in a sufficiently small neighbourhood. By Corollary~\ref{cor:local_codim2_Morrey_criterion_positive_monotonicity}, one has $\mathcal O=\varnothing$, and the last assertion follows from Corollary~\ref{cor:obstruction_hierarchy}.
\end{proof}

\subsection{Monopoles on the Bryant--Salamon \texorpdfstring{$\mathrm G_2$-manifolds}{G2-manifolds}}

Let
\[
    X=\Lambda_-^2(S^4)
    \qquad\text{or}\qquad
    X=\Lambda_-^2(\mathbb{CP}^2),
\]
and let $N$ denote the zero section. In the first case $N=S^4$ and the
homogeneous monopoles have structure group $\mathrm{SU}(2)$. In the second
case $N=\mathbb{CP}^2$ and the homogeneous family constructed in
\cite{oliveira2014monopoles} has structure group $\mathrm{SO}(3)$; see
Remark~\ref{rmk:CP2_structure_group} below.

We use the notation of \cite{oliveira2014monopoles}, with one change which
is important for the calculation. Let $s\geqslant0$ be the orbit
coordinate and set
\[
    \rho(s)=\int_0^s f(t^2)\,dt,
    \qquad
    f(s^2)=(1+s^2)^{-1/4},
    \qquad
    g(s^2)=\sqrt{2}(1+s^2)^{1/4}.
\]
Thus $d\rho=f\,ds$ and $g^2=2f^{-2}$. We denote by $\alpha_m$ the
coefficient which appears in the invariant connection; this is the function
called $a$ in \cite{oliveira2014monopoles}. The Bogomolny profile is instead
\begin{equation}
\label{eq:BS_b_alpha_relation}
    b_m=f^{-2}\alpha_m.
\end{equation}
It is the pair $(b_m,\phi_m)$, not $(\alpha_m,\phi_m)$, which satisfies the
standard spherically symmetric Bogomolny system
\begin{equation}
\label{eq:BS_auxiliary_Bogomolny}
    \frac{d\phi_m}{d\rho}
    =
    \frac{b_m^2-1}{2h^2},
    \qquad
    \frac{db_m}{d\rho}
    =
    2\phi_m b_m.
\end{equation}
Here
\[
    h^2=2s^2f^{-2}
    \quad\text{on }\Lambda_-^2(S^4),
    \qquad
    h^2=s^2f^{-2}
    \quad\text{on }\Lambda_-^2(\mathbb{CP}^2).
\]
We choose the sign of the Higgs field so that
\begin{equation}
\label{eq:BS_profile_bounds}
    0<b_m\leqslant1,
    \qquad
    -m\leqslant\phi_m\leqslant0,
    \qquad
    \lim_{\rho\to\infty}\phi_m=-m.
\end{equation}
These are the standard properties of the solutions of
\eqref{eq:BS_auxiliary_Bogomolny}; see the Appendix of
\cite{oliveira2014monopoles}.

Let $T_1,T_2,T_3$ be the Lie algebra basis used in
\cite{oliveira2014monopoles}, so that
$[T_i,T_j]=2\varepsilon_{ijk}T_k$, and put
\[
    \kappa:=|T_1|^2=|T_2|^2=|T_3|^2.
\]
For the $\mathrm{SU}(2)$ inner product \eqref{eq:inner_prod_convention} and the
usual matrix realization of this basis, $\kappa=4$. Keeping $\kappa$
explicit also covers the fixed normalization used for the
$\mathrm{SO}(3)$ family.

The explicit curvature formulae in Sections~3 and~4 of
\cite{oliveira2014monopoles} give the following common norm identities for
both Bryant--Salamon families:
\begin{align}
\label{eq:BS_full_curvature_norm}
    \kappa^{-1}|F_{\nabla_m}|^2
    &={}
    \frac{1+2\alpha_m^2}{2g^4}
    +
    \frac{(1-\alpha_m^2)^2}{4s^4f^4}
    +
    \frac{(\partial_s\alpha_m)^2}{2s^2f^4},\\
\label{eq:BS_full_Higgs_norm}
    \kappa^{-1}|\nabla_m\Phi_m|^2
    &={}
    \frac{(\partial_s\phi_m)^2}{f^2}
    +
    \frac{2\phi_m^2\alpha_m^2}{s^2f^2}.
\end{align}
For completeness, in the $S^4$ case these identities follow from
\begin{align*}
    F_{\nabla_m}
    ={}&
    \left(\frac12\Omega_1-2(1-\alpha_m^2)\omega^{23}\right)T_1
    +
    \frac{\alpha_m}{2}(\Omega_2T_2+\Omega_3T_3)\\
    &+
    (\partial_s\alpha_m)
    \bigl(ds\wedge\omega^2T_2+ds\wedge\omega^3T_3\bigr),\\
    \nabla_m\Phi_m
    ={}&
    (\partial_s\phi_m)ds\,T_1
    +2\phi_m\alpha_m(\omega^3T_2-\omega^2T_3),
\end{align*}
whereas in the $\mathbb{CP}^2$ case they follow from
\begin{align*}
    F_{\nabla_m}
    ={}&
    \bigl(\Omega_1+2(\alpha_m^2-1)\nu^{12}\bigr)T_1\\
    &+
    \bigl((\partial_s\alpha_m)ds\wedge\nu^1+\alpha_m\Omega_2\bigr)T_2
    +
    \bigl((\partial_s\alpha_m)ds\wedge\nu^2+\alpha_m\Omega_3\bigr)T_3,\\
    \nabla_m\Phi_m
    ={}&
    (\partial_s\phi_m)ds\,T_1
    +2\phi_m\alpha_m(\nu^2T_2-\nu^1T_3).
\end{align*}
The apparently identical expressions
\eqref{eq:BS_full_curvature_norm}--\eqref{eq:BS_full_Higgs_norm} result from
the different normalizations of the horizontal forms in the two metrics.

The equations for the actual connection coefficient $\alpha_m$ are not
\eqref{eq:BS_auxiliary_Bogomolny}. On $\Lambda_-^2(S^4)$ they are
\begin{equation}
\label{eq:BS_S4_alpha_phi_ODE}
    \partial_s\phi_m
    =
    -\frac{1-\alpha_m^2}{2s^2f}
    +\frac{f}{g^2},
    \qquad
    \partial_s\alpha_m
    =
    -s f^4\alpha_m+2f\phi_m\alpha_m,
\end{equation}
while on $\Lambda_-^2(\mathbb{CP}^2)$ they are
\begin{equation}
\label{eq:BS_CP2_alpha_phi_ODE}
    \partial_s\phi_m
    =
    \frac{f}{2s^2}\bigl(f^{-4}\alpha_m^2-1\bigr),
    \qquad
    \partial_s\alpha_m
    =
    -s f^4\alpha_m+2f\phi_m\alpha_m.
\end{equation}
Equations \eqref{eq:BS_auxiliary_Bogomolny} are obtained from these only
after the substitution \eqref{eq:BS_b_alpha_relation} and the change from
$s$ to $\rho$.

Substitution of \eqref{eq:BS_S4_alpha_phi_ODE} into
\eqref{eq:BS_full_curvature_norm}--\eqref{eq:BS_full_Higgs_norm} gives the
exact signed imbalance on $\Lambda_-^2(S^4)$:
\begin{equation}
\label{eq:BS_S4_exact_imbalance}
    \kappa^{-1}q_m
    =
    \frac{1-\alpha_m^2}{s^2f^2g^2}
    -
    \frac{4\alpha_m^2\phi_m}{sfg^2}
    +
    \frac{6\alpha_m^2-1}{2g^4}.
\end{equation}
Likewise, substitution of \eqref{eq:BS_CP2_alpha_phi_ODE} gives on
$\Lambda_-^2(\mathbb{CP}^2)$
\begin{equation}
\label{eq:BS_CP2_exact_imbalance}
    \kappa^{-1}q_m
    =
    \frac{1+2\alpha_m^2}{2g^4}
    +
    \frac{f^4-\alpha_m^4}{4s^2f^4}
    -
    \frac{2f\alpha_m^2\phi_m}{s}
    +
    \frac12f^4\alpha_m^2.
\end{equation}
Thus the quadratic Bogomolny terms cancel, but the cancellation does not
leave a pointwise remainder that is uniform in the mass: the terms linear in
$\phi_m/s$ remain.

\begin{proposition}[Morrey control for the Bryant--Salamon families]
\label{prop:Bryant_Salamon_Morrey_control}
Let $(\nabla_m,\Phi_m)$ denote either Bryant--Salamon family. There exist a
fixed tubular neighbourhood $U$ of the zero section $N$ and a constant
$C<\infty$, independent of $m\geqslant1$, such that
\begin{equation}
\label{eq:BS_tubular_q_bound}
    q_m^+
    \leqslant
    C\bigl(1+r_N^{-2}+m r_N^{-1}\bigr)
    \qquad\text{on }U\setminus N.
\end{equation}
Moreover, for every compact set $K\Subset X\setminus N$, there is
$C_K<\infty$ such that
\begin{equation}
\label{eq:BS_away_q_bound}
    q_m^+\leqslant C_K
    \qquad\text{on }K.
\end{equation}
In particular, the positive imbalance satisfies the codimension-two Morrey hypothesis of
Proposition~\ref{prop:codim2_Morrey_criterion_positive_monotonicity} locally uniformly near every point of $X$.
\end{proposition}

\begin{proof}
By \eqref{eq:BS_b_alpha_relation}, \eqref{eq:BS_profile_bounds}, and
$0<f\leqslant1$,
\[
    0<\alpha_m=f^2b_m\leqslant1,
    \qquad
    |\phi_m|\leqslant m.
\]
On a fixed neighbourhood $0<s<s_0$ of the zero section, the functions
$f$ and $g$ are bounded above and below by positive constants and
$r_N=\rho\asymp s$. Formula \eqref{eq:BS_S4_exact_imbalance} therefore
gives
\[
    q_m^+
    \leqslant
    C\bigl(1+s^{-2}+ms^{-1}\bigr)
\]
in the $S^4$ case. In the $\mathbb{CP}^2$ case,
$\alpha_m=f^2b_m$ implies $\alpha_m^4\leqslant f^8\leqslant f^4$; hence
\eqref{eq:BS_CP2_exact_imbalance} gives the same estimate. This proves
\eqref{eq:BS_tubular_q_bound}.

On a compact subset of $X\setminus N$, the connections in the
translated configurations converge smoothly, with all derivatives, to the
corresponding reducible Dirac-type connection by
\cite[Theorem~1]{oliveira2014monopoles}. Hence $F_{\nabla_m}$ is uniformly
bounded there. The monopole equation then gives
$|\nabla_m\Phi_m|=|F_{\nabla_m}\wedge\psi|\leqslant C|F_{\nabla_m}|$,
so the Higgs derivatives are uniformly bounded as well. This proves
\eqref{eq:BS_away_q_bound}. The final assertion follows from
Lemma~\ref{lem:tubular_Morrey_control_excludes_obstruction}.
\end{proof}

\begin{remark}[The structure group of the $\mathbb{CP}^2$ family]
\label{rmk:CP2_structure_group}
The monopoles on $\Lambda_-^2(\mathbb{CP}^2)$ constructed in
\cite{oliveira2014monopoles} live on a homogeneous
$\mathrm{SO}(3)$-bundle. They therefore fall within our expanded Main
Setting. The transported inner product
\eqref{eq:SO3_inner_prod_convention} is the normalization used when
comparing this family with the general compactness and obstruction theory.
The explicit norm and imbalance calculations above are unchanged.
\end{remark}

For each Bryant--Salamon family,
Lemma~\ref{lem:tubular_Morrey_control_excludes_obstruction} and
Corollary~\ref{cor:obstruction_hierarchy} give
\[
    \mathcal O=\varnothing,
    \qquad
    \mathcal S=\mathcal Z=\mathcal C.
\]
By \cite[Theorem~1]{oliveira2014monopoles}, the zero set of every monopole
in either family is the zero section. Hence
\begin{equation}
\label{eq:S4_homogeneous_identification}
    \mathcal S=\mathcal Z=\mathcal C=S^4
\end{equation}
for the $\mathrm{SU}(2)$ family on $\Lambda_-^2(S^4)$, while
\begin{equation}
\label{eq:CP2_homogeneous_identification}
    \mathcal S=\mathcal Z=\mathcal C=\mathbb{CP}^2
\end{equation}
for the $\mathrm{SO}(3)$ family on
$\Lambda_-^2(\mathbb{CP}^2)$.

\subsection{\texorpdfstring{$\mathrm{SU}(2)^2$-invariant}{SU(2) x SU(2)-invariant} Calabi--Yau monopoles}

We now consider the Stenzel metric on $T^*S^3$. In
\cite{oliveira2016calabi}, the second author constructed a one-parameter
family of invariant $\mathrm{SU}(2)$ Calabi--Yau monopoles, parametrized by
the positive mass. Stein's classification identifies this as the unique
irreducible quadratically decaying invariant monopole family with nonzero
Higgs field on the Stenzel smoothing and excludes such families on
$O(-1)\oplus O(-1)$ and $O(-2,-2)$; see
\cite[Theorems~A and~E and Proposition~4.25]{stein2023invariant}.

We carry out the calculation in the notation of
\cite{oliveira2016calabi}. Put
\[
    R_\pm^2=\frac{r^2\pm\varepsilon^2}{2},
    \qquad
    G=R_+R_-F'(r^2),
\]
and write the Stenzel metric as
\begin{equation}
\label{eq:Stenzel_metric_coefficients}
    g
    =
    u^2dr^2
    +v_1^2(\theta^1)^2
    +v_2^2\bigl((\theta^2)^2+(\theta^3)^2\bigr)
    +v_4^2\bigl((\theta^4)^2+(\theta^5)^2\bigr),
\end{equation}
where
\[
    u^2=\dot G\frac{r}{2R_+R_-},
    \qquad
    v_1^2=\dot G\frac{2R_+R_-}{r},
    \qquad
    v_2^2=G\frac{R_+}{R_-},
    \qquad
    v_4^2=G\frac{R_-}{R_+}.
\]
Here a dot denotes $d/dr$. Let $\rho$ be geodesic distance to the zero
section, so $d\rho=u\,dr$, and set
\begin{equation}
\label{eq:Stenzel_h_definition}
    h^2
    :=
    \varepsilon^{-2}R_+R_-G.
\end{equation}
The Ricci-flat equation is
\[
    2\dot G G^2=rR_+R_-.
\]
The genuine six-dimensional monopole has the form
\begin{equation}
\label{eq:Stenzel_genuine_monopole}
    A_m
    =
    A_c^1+
    \frac{\varepsilon a_m}{2R_+}
    (\theta^4T_2+\theta^5T_3),
    \qquad
    \Phi_m=\phi_mT_1.
\end{equation}
The pair $(a_m,\phi_m)$ satisfies
\begin{equation}
\label{eq:Stenzel_auxiliary_Bogomolny}
    \phi_m'
    =
    \frac{a_m^2-1}{2h^2},
    \qquad
    a_m'=2\phi_m a_m,
\end{equation}
where a prime denotes $d/d\rho$. Thus
\eqref{eq:Stenzel_auxiliary_Bogomolny} is the Bogomolny equation for the
auxiliary three-dimensional metric
$d\rho^2+h^2g_{S^2}$. It is important that the connection in
\eqref{eq:Stenzel_genuine_monopole} contains the additional factor
$\varepsilon/R_+$. The auxiliary fibrewise connection used in the bubbling
argument of \cite{oliveira2016calabi} does not contain this factor.

Choose the sign so that
\begin{equation}
\label{eq:Stenzel_profile_bounds}
    0<a_m\leqslant1,
    \qquad
    -m\leqslant\phi_m\leqslant0,
    \qquad
    \lim_{\rho\to\infty}\phi_m=-m.
\end{equation}
Set
\[
    c_m:=\frac{\varepsilon a_m}{R_+}.
\]
The curvature formula in \cite[Remark~8]{oliveira2016calabi} and the
covariant derivative of the Higgs field are
\begin{align}
\label{eq:Stenzel_curvature_explicit}
    F_{A_m}
    ={}&
    \frac12\bigl((c_m^2-1)\theta^{45}-\theta^{23}\bigr)T_1
    +\frac{c_m}{2}(\theta^{12}T_2+\theta^{13}T_3)\\
    &+
    \frac{\partial_rc_m}{2}
    \bigl(dr\wedge\theta^4T_2+dr\wedge\theta^5T_3\bigr),\nonumber\\
\label{eq:Stenzel_Higgs_explicit}
    \nabla_{A_m}\Phi_m
    ={}&
    \phi_m' d\rho\,T_1
    +c_m\phi_m(\theta^5T_2-\theta^4T_3).
\end{align}
Consequently,
\begin{align}
\label{eq:Stenzel_norm_F_v}
    \kappa^{-1}|F_{A_m}|^2
    ={}&
    \frac{(c_m^2-1)^2}{4v_4^4}
    +
    \frac{1}{4v_2^4}
    +
    \frac{c_m^2}{2v_1^2v_2^2}
    +
    \frac{(\partial_\rho c_m)^2}{2v_4^2},\\
\label{eq:Stenzel_norm_Higgs_v}
    \kappa^{-1}|\nabla_{A_m}\Phi_m|^2
    ={}&
    (\phi_m')^2
    +
    \frac{2c_m^2\phi_m^2}{v_4^2}.
\end{align}
To expose the cancellation, define
\begin{equation}
\label{eq:Stenzel_delta_gamma}
    \delta
    :=
    \frac{R_+^2}{\varepsilon^2}-1,
    \qquad
    \gamma
    :=
    \frac{d}{d\rho}\log R_+.
\end{equation}
The identities
\[
    v_2^2=\frac{\varepsilon^2h^2}{R_-^2},
    \qquad
    v_4^2=\frac{\varepsilon^2h^2}{R_+^2},
    \qquad
    v_1^2v_2^2=\frac{R_+^3R_-}{G},
\]
and
\[
    \partial_\rho c_m
    =
    \frac{\varepsilon}{R_+}(a_m'-\gamma a_m)
\]
transform \eqref{eq:Stenzel_norm_F_v}--\eqref{eq:Stenzel_norm_Higgs_v}
into
\begin{align}
\label{eq:Stenzel_full_curvature_norm}
    \kappa^{-1}|F_{A_m}|^2
    ={}&
    \frac{R_-^4}{4\varepsilon^4h^4}
    +
    \frac{\varepsilon^2a_m^2G}{2R_+^5R_-}
    +
    \frac{(a_m^2-1-\delta)^2}{4h^4}
    +
    \frac{(a_m'-\gamma a_m)^2}{2h^2},\\
\label{eq:Stenzel_full_Higgs_norm}
    \kappa^{-1}|\nabla_{A_m}\Phi_m|^2
    ={}&
    (\phi_m')^2
    +
    \frac{2a_m^2\phi_m^2}{h^2}.
\end{align}
Using \eqref{eq:Stenzel_auxiliary_Bogomolny}, the quadratic Bogomolny
terms cancel and one obtains the exact signed imbalance
\begin{equation}
\label{eq:Stenzel_exact_imbalance}
    \begin{aligned}
    \kappa^{-1}q_m
    ={}&
    \frac{R_-^4}{4\varepsilon^4h^4}
    +
    \frac{\varepsilon^2a_m^2G}{2R_+^5R_-}\\
    &+
    \frac{\delta^2-2\delta(a_m^2-1)}{4h^4}
    +
    \frac{\gamma^2a_m^2-2\gamma a_ma_m'}{2h^2}.
    \end{aligned}
\end{equation}
This formula makes precise why the three-dimensional Bogomolny reduction is
not an exact equality between the vertical six-dimensional curvature energy
and the Higgs energy: the factors $\delta$ and $\gamma$ measure the
difference between the genuine Stenzel connection and the auxiliary normal
Bogomolny connection.

\begin{proposition}[Morrey control for the Stenzel family]
\label{prop:Stenzel_Morrey_control}
Let $(A_m,\Phi_m)$ denote the Stenzel family. There exist a fixed tubular
neighbourhood $U$ of the zero section $S^3$ and a constant $C<\infty$,
independent of $m\geqslant1$, such that
\begin{equation}
\label{eq:Stenzel_tubular_q_bound}
    q_m^+
    \leqslant
    C\bigl(1+r_{S^3}^{-2}+m r_{S^3}^{-1}\bigr)
    \qquad\text{on }U\setminus S^3.
\end{equation}
Moreover, for every compact set $K\Subset T^*S^3\setminus S^3$, there is
$C_K<\infty$ such that
\begin{equation}
\label{eq:Stenzel_away_q_bound}
    q_m^+\leqslant C_K
    \qquad\text{on }K.
\end{equation}
Consequently,
\[
    \mathcal O=\varnothing,
    \qquad
    \mathcal S=\mathcal Z=\mathcal C=S^3.
\]
\end{proposition}

\begin{proof}
The smoothness expansions of the Stenzel metric at the zero section give
\[
    h(\rho)\asymp\rho,
    \qquad
    R_-(\rho)=O(\rho),
    \qquad
    \frac{dR_-}{d\rho}=O(1),
    \qquad
    \frac{G(\rho)}{R_-(\rho)}=O(1).
\]
The remaining estimates in
\eqref{eq:Stenzel_zero_section_expansions} follow directly from the
definitions. Indeed, since
\[
    R_+^2-R_-^2=\varepsilon^2,
\]
one has the exact identity
\[
    \delta
    =
    \frac{R_+^2}{\varepsilon^2}-1
    =
    \frac{R_-^2}{\varepsilon^2},
\]
and therefore $\delta(\rho)=O(\rho^2)$. Differentiating
$R_+^2-R_-^2=\varepsilon^2$ with respect to $\rho$ gives
\[
    R_+\frac{dR_+}{d\rho}
    =
    R_-\frac{dR_-}{d\rho}.
\]
Since $R_+(0)=\varepsilon$ and hence $R_+$ is bounded away from zero
near the zero section,
\[
    \gamma
    =
    \frac{d}{d\rho}\log R_+
    =
    \frac{R_-}{R_+^2}\frac{dR_-}{d\rho}
    =
    O(\rho).
\]
Thus
\begin{equation}
\label{eq:Stenzel_zero_section_expansions}
    h(\rho)\asymp\rho,
    \qquad
    R_-(\rho)=O(\rho),
    \qquad
    \frac{G(\rho)}{R_-(\rho)}=O(1),
    \qquad
    \delta(\rho)=O(\rho^2),
    \qquad
    \gamma(\rho)=O(\rho).
\end{equation}
All constants implicit in these estimates depend only on the fixed
Stenzel metric. From \eqref{eq:Stenzel_profile_bounds} and
\eqref{eq:Stenzel_auxiliary_Bogomolny},
\[
    0<a_m\leqslant1,
    \qquad
    |\phi_m|\leqslant m,
    \qquad
    |a_m'|=2|\phi_ma_m|\leqslant2m.
\]
Substituting these bounds and
\eqref{eq:Stenzel_zero_section_expansions} into
\eqref{eq:Stenzel_exact_imbalance} gives, for $0<\rho<\rho_0$,
\[
    q_m^+
    \leqslant
    C\bigl(1+\rho^{-2}+m\rho^{-1}\bigr).
\]
Since $\rho=r_{S^3}$ in the tubular region, this is
\eqref{eq:Stenzel_tubular_q_bound}.

On a compact subset of $T^*S^3\setminus S^3$, the connections in the
translated configurations converge smoothly, with all derivatives, to the
reducible Dirac Calabi--Yau connection by
\cite[Theorem~1]{oliveira2016calabi}. Hence $F_{A_m}$ is uniformly bounded
there. The monopole equation gives
$|\nabla_{A_m}\Phi_m|=|F_{A_m}\wedge\operatorname{Re}\Omega|
\leqslant C|F_{A_m}|$, so the Higgs derivatives are uniformly bounded as
well. This proves \eqref{eq:Stenzel_away_q_bound}. Lemma~\ref{lem:tubular_Morrey_control_excludes_obstruction} therefore gives
$\mathcal O=\varnothing$, and Corollary~\ref{cor:obstruction_hierarchy} gives $\mathcal S=\mathcal Z=\mathcal C$. 

It remains to identify the Higgs zero set of each member of the family.
The regular boundary conditions in the construction give
\[
    a_m(0)=1,
    \qquad
    \phi_m(0)=0.
\]
Moreover, by
\eqref{eq:Stenzel_profile_bounds} and
\eqref{eq:Stenzel_auxiliary_Bogomolny},
\[
    \phi_m'
    =
    \frac{a_m^2-1}{2h^2}
    \leqslant0
    \qquad
    \text{on }(0,\infty).
\]
We claim that
\[
    \phi_m(\rho)<0
    \qquad
    \text{for every }\rho>0.
\]
Suppose instead that $\phi_m(\rho_0)=0$ for some $\rho_0>0$.
Since $\phi_m$ is nonincreasing and $\phi_m(0)=0$, it follows that
\[
    \phi_m\equiv0
    \qquad
    \text{on }[0,\rho_0].
\]
The equation
\[
    a_m'=2\phi_ma_m
\]
then gives $a_m\equiv1$ on the same interval. At every positive radius
the system \eqref{eq:Stenzel_auxiliary_Bogomolny} is a regular first-order
ODE system. Uniqueness applied at any point of $(0,\rho_0)$ therefore
forces
\[
    a_m\equiv1,
    \qquad
    \phi_m\equiv0
    \qquad
    \text{on }(0,\infty),
\]
contradicting
\[
    \lim_{\rho\to\infty}\phi_m=-m<0.
\]
Hence $\phi_m(\rho)<0$ for every $\rho>0$, and the Higgs field vanishes
precisely on the zero section:
\[
    \Phi_m^{-1}(0)=S^3.
\]
Consequently $\mathcal Z=S^3$, and the identification obtained above becomes
\[
    \mathcal S=\mathcal Z=\mathcal C=S^3.
\]
\end{proof}

\begin{remark}[No pointwise imbalance bound uniform in the mass]
\label{rmk:no_uniform_pointwise_imbalance_homogeneous}
The estimates above deliberately control $q_m^+$ in Morrey scale rather
than $q_m$ in $L^\infty$. The exact formulae
\eqref{eq:BS_S4_exact_imbalance},
\eqref{eq:BS_CP2_exact_imbalance}, and
\eqref{eq:Stenzel_exact_imbalance} show that the leading quadratic
Bogomolny block cancels, but terms of size $m r_N^{-1}$ remain. On the
mass scale $r_N\asymp m^{-1}$ these are compatible with an $m^2$ pointwise
imbalance. This does not affect the codimension-two Morrey bound because
the normal bundle has rank three.
\end{remark}

\begin{corollary}[Cohomogeneity-one examples in the unobstructed case]
\label{cor:homogeneous_examples_dream_case}
For the $\mathrm{SU}(2)$ Bryant--Salamon family on
$\Lambda_-^2(S^4)$, the $\mathrm{SO}(3)$ Bryant--Salamon family on
$\Lambda_-^2(\mathbb{CP}^2)$, and the Stenzel family of
$\mathrm{SU}(2)^2$-invariant Calabi--Yau monopoles, every large mass
sequence satisfies
\[
    \mathcal O=\varnothing.
\]
Consequently,
\[
    \mathcal S=\mathcal Z=\mathcal C,
    \qquad
    \mathcal R_{\mathrm{ab}}=X\setminus \mathcal S,
\]
and the common limiting set is the compact calibrated zero section. For
every member of these three explicit families, the Higgs zero set is the
zero section itself. Moreover, for every compact set $K\Subset X$
containing the zero section in its interior,
\[
    d_{\mathcal H}
    \bigl(\overline{C_{\Lambda_0,i}}\cap K,\mathcal S\bigr)
    \longrightarrow0.
\]
\end{corollary}

For these cohomogeneity-one families, the calibrated support and the
Higgs zero sets are the zero section; Li's curvature concentration loci have no
additional points in the Kuratowski upper limit, and the complement carries the
smooth reducible Dirac-type limit of
\cites{oliveira2014monopoles,oliveira2016calabi}. The key input is the
Morrey estimate for the positive imbalance, rather than a pointwise
estimate uniform in the mass bound.

Remark~\ref{rmk:homogeneous_Hc2_vanishing} gives a second derivation of
the smooth convergence on the complement of the zero section from the
vanishing of $H_c^2(X;\mathbb R)$ and the unobstructedness established
above.

\section{Abelianization on \texorpdfstring{$\mathcal R_{\mathrm{ab}}$}{the abelian region}}
\label{sec: abelianization}

We study the sequence locally on compact domains contained in
$\mathcal R_{\mathrm{ab}}$. By the definition of this set and compactness,
Li's curvature concentration loci are eventually absent from each such
domain. Consequently, the Higgs
field is uniformly large there, while the transverse curvature and all
local derivatives of the transverse variables decay exponentially. This
allows us to separate the resulting longitudinal scalar variables from
their transverse adjoint-valued counterparts and to derive the local
Hodge equations satisfied by the longitudinal curvature.

Corollary~\ref{cor:nonabelian_locus_obstruction_identity}
gives $\mathcal R_{\mathrm{ab}}=X\setminus(\mathcal S\cup\mathcal O)$. No general nonemptiness or density statement for this set is currently
known. If $\mathcal O\subset\mathcal S$, then
$\mathcal R_{\mathrm{ab}}=X\setminus\mathcal S$. Since $\mathcal S$ is
compact and has finite $(n-3)$-dimensional Hausdorff measure,
$\mathcal R_{\mathrm{ab}}$ is then an open dense subset of full
Riemannian measure, whose complement has Hausdorff codimension at least
three. For the analysis carried out in this and the following section,
we assume that $\mathcal R_{\mathrm{ab}}\neq\varnothing$.

\medskip

\noindent\textbf{Convention.}
Throughout Sections~\ref{sec: abelianization}
and~\ref{sec:frequency_away_Cinfty}, the sets denoted by
$K,K',K'',\ldots$ are compact domains contained in
$\mathcal R_{\mathrm{ab}}$, by which we mean compact sets satisfying
\[
K=\overline{\operatorname{int}K}\neq\varnothing.
\]
Moreover, the sets denoted by $U,V,W,\ldots$ are nonempty open sets
compactly contained in $\mathcal R_{\mathrm{ab}}$. No regularity of their
boundaries is assumed.

\subsection{Longitudinal and transverse splitting and the translated Higgs fields}\label{subsec: translated_fields}

We begin by decomposing the monopole fields into longitudinal and
transverse components relative to the normalized Higgs direction. The
lower bound for $|\Phi_i|$ on compact subsets of
$\mathcal R_{\mathrm{ab}}$ makes this splitting uniform there. 

Let $K\Subset\mathcal R_{\mathrm{ab}}$. By
Lemma~\ref{lem:Li_dichotomy_away_Cinfty}, after increasing $i_K$ if
necessary, $|\Phi_i|\geqslant m_i/2$ on $K$ for every
$i\geqslant i_K$. Hence the normalized Higgs direction
$\Psi_i:=\Phi_i/|\Phi_i|$ is well-defined on $K$. It induces the
orthogonal splitting
\[
    \mathfrak{g}_P|_K
    =
    \langle\Psi_i\rangle
    \oplus
    \langle\Psi_i\rangle^\perp.
\]
We denote the corresponding components of an adjoint-valued form $\alpha$ by
$\alpha^\parallel$ and $\alpha^\perp$.

We next record the exact scalar abelian equations associated with this
splitting.

\begin{lemma}[Exact scalar abelian equations]
\label{lem:exact_scalar_abelian_equations}
On $K$, define
\[
    v_i
    :=
    m_i-|\Phi_i|,
    \qquad
    u_i
    :=
    \frac{m_i^2-|\Phi_i|^2}{2m_i}.
\]
Thus $v_i$ is the nonnegative Higgs defect and $u_i$ is its quadratic
counterpart. Also define the real-valued $2$-forms
\[
    \widehat\omega_i
    :=
    \langle F_{\nabla_i},\Psi_i\rangle,
    \qquad
    \omega_i
    :=
    m_i^{-1}\langle F_{\nabla_i},\Phi_i\rangle
    =
    \frac{|\Phi_i|}{m_i}\widehat\omega_i.
\]
Then the $\Theta$-monopole equation implies the exact identities
\begin{equation}
\label{eq:exact_scalar_abelian_equations}
    \widehat\omega_i\wedge\Theta
    =
    -*dv_i,
    \qquad
    \omega_i\wedge\Theta
    =
    -*du_i,
\end{equation} and in the Calabi--Yau case both longitudinal forms are primitive: $\Lambda\omega_i=\Lambda\widehat\omega_i=0$.
Moreover,
\begin{equation}
\label{eq:scalar_amplitude_comparison}
    v_i-u_i
    =
    \frac{(m_i-|\Phi_i|)^2}{2m_i}
    \geqslant0,
\end{equation}
and hence, for every compact $K'\Subset K$,
\[
    \|v_i-u_i\|_{L^1(K')}
    \longrightarrow0.
\]
Along the subsequence of \cite[Theorem~1.12]{li2025large}, both $u_i$
and $v_i$ converge strongly in $L^1_{\mathrm{loc}}$ to
$-2\Phi_\infty^{\mathrm{Li}}$, where
$\Phi_\infty^{\mathrm{Li}}$ denotes Li's scalar Higgs limit.
\end{lemma}

\begin{proof}
Pairing $F_{\nabla_i}\wedge\Theta=*\nabla_i\Phi_i$ with $\Psi_i$ gives
\[
    \widehat\omega_i\wedge\Theta
    =
    *\langle\nabla_i\Phi_i,\Psi_i\rangle
    =
    *d|\Phi_i|
    =
    -*dv_i.
\]
Similarly, pairing with $m_i^{-1}\Phi_i$ gives
\[
\begin{aligned}
    \omega_i\wedge\Theta
    &=
    *m_i^{-1}
    \langle\nabla_i\Phi_i,\Phi_i\rangle  \\
    &=
    -*d\left(
        \frac{m_i^2-|\Phi_i|^2}{2m_i}
    \right)
    =
    -*du_i.
\end{aligned}
\]
This proves \eqref{eq:exact_scalar_abelian_equations}. In the Calabi--Yau case, $\Lambda$ acts on the differential form factor, so $\Lambda F_{\nabla_i}=0$ gives $\Lambda\omega_i=\Lambda\widehat\omega_i=0$. 

A direct
calculation gives
\[
\begin{aligned}
    v_i-u_i
    &=
    (m_i-|\Phi_i|)
    -
    \frac{m_i^2-|\Phi_i|^2}{2m_i}  \\
    &=
    \frac{(m_i-|\Phi_i|)^2}{2m_i},
\end{aligned}
\]
which is \eqref{eq:scalar_amplitude_comparison}. Lemma~\ref{lem:largeness}
then yields
\[
    \|v_i-u_i\|_{L^1(K')}
    =
    \frac{1}{2m_i}
    \int_{K'}
    (m_i-|\Phi_i|)^2
    \longrightarrow0.
\]
To compare with Li's scalar variable, define
\[
    \phi_i^{\mathrm{Li}}
    :=
    \frac{|\Phi_i|_{\mathrm{Li}}^2
    -(m_i^{\mathrm{Li}})^2}
    {2m_i^{\mathrm{Li}}}
    \leqslant0.
\]
Li proves that
$\phi_i^{\mathrm{Li}}\to\Phi_\infty^{\mathrm{Li}}$ strongly in
$L^1_{\mathrm{loc}}$. The norm and mass conversions
\eqref{eq:PPS_Li_norm_mass_comparison} and
\eqref{eq:PPS_Li_mass_comparison} give
\[
    u_i=-2\phi_i^{\mathrm{Li}}.
\]
Therefore \cite[Theorem~1.12]{li2025large} yields
$u_i\to-2\Phi_\infty^{\mathrm{Li}}$ strongly in
$L^1_{\mathrm{loc}}$, and
\eqref{eq:scalar_amplitude_comparison} gives the same limit for $v_i$.
\end{proof}

\begin{proposition}[Scalar Green convergence on $\mathcal R_{\mathrm{ab}}$]
\label{prop:scalar_Green_convergence_Rab}
Define
\begin{equation}\label{eq:def_scalar_Green_limit}
    v_T(x)
    :=
    4\pi
    \int_XG(x,y)\,d\|T\|(y),
    \qquad
    x\in X\setminus\mathcal S.
\end{equation}
Then $v_T$ is smooth, harmonic and strictly positive on
$X\setminus\mathcal S$. For every compact domain
$K\Subset\mathcal R_{\mathrm{ab}}$, the scalar fields of
Lemma~\ref{lem:exact_scalar_abelian_equations} satisfy
\begin{equation}\label{eq:scalar_Green_convergence_Rab}
    u_i\longrightarrow v_T,
    \qquad
    v_i\longrightarrow v_T
    \qquad
    \text{uniformly on }K.
\end{equation}
Consequently,
\begin{equation}\label{eq:uniform_Higgs_normalization_Rab}
    \frac{|\Phi_i|}{m_i}
    \longrightarrow1
    \qquad
    \text{uniformly on }K,
\end{equation}
and there is a constant $C_K<\infty$ such that
\begin{equation}\label{eq:uniform_translated_Higgs_bound_Rab}
    \sup_K(m_i-|\Phi_i|)
    =
    \sup_K|\widetilde\Phi_i|
    \leqslant C_K
\end{equation}
for all sufficiently large $i$.
\end{proposition}

\begin{proof}
The Green representation in
Theorem~\ref{thm: alternative_finite_mass} and the definition
\[
    \nu_i
    =
    2m_i^{-1}|\nabla_i\Phi_i|^2\,\vol
\]
give
\begin{equation}\label{eq:ui_Green_representation}
    u_i(x)
    =
    \frac12
    \int_XG(x,y)\,d\nu_i(y).
\end{equation}
By Proposition~
\ref{prop:PPS_saturation_and_Li_identification},
\[
    \nu_i\rightharpoonup^*8\pi\|T\|,
    \qquad
    \nu_i(X)=8\pi k,
    \qquad
    \|T\|(X)=k.
\]
Choose compact domains
\[
    K\Subset K'\Subset\mathcal R_{\mathrm{ab}}
\]
and $\rho>0$ such that
$B_{2\rho}(x)\Subset K'$ for every $x\in K$. Set
\[
    \varepsilon_i(K')
    :=
    \sup_{K'}m_i^{-1}e_i.
\]
By Proposition~\ref{prop:mass_renormalized_density_decay_Rab},
$\varepsilon_i(K')\to0$. The local estimate for the Green function therefore gives
\begin{equation}\label{eq:local_Green_source_vanishes_Rab}
\begin{aligned}
    \sup_{x\in K}
    \int_{B_\rho(x)}G(x,y)\,d\nu_i(y)
    &\leqslant
    2\varepsilon_i(K')
    \sup_{x\in K}
    \int_{B_\rho(x)}G(x,y)\,\vol(y)\\
    &\leqslant
    C_K\varepsilon_i(K')\rho^2
    \longrightarrow0.
\end{aligned}
\end{equation}
Since $\mathcal S\subset\mathcal C$ and
$K'\Subset\mathcal R_{\mathrm{ab}}$, the limiting measure
$\|T\|$ vanishes on $K'$, so its contribution to the same local region
is zero.

It remains to work away from the diagonal. The Green function is smooth
there and decays uniformly at infinity when its first variable ranges
over $K$. The total masses of $\nu_i$ are uniformly bounded. Hence the
tails of the Green potentials outside a sufficiently large compact set
are uniformly small. On the remaining compact subset away from the
diagonal, the family of functions
\[
    y\longmapsto G(x,y),
    \qquad x\in K,
\]
is compact in $C^0$, and weak convergence of $\nu_i$ gives convergence
of the corresponding integrals uniformly in $x\in K$. Together with
\eqref{eq:local_Green_source_vanishes_Rab}, this proves
\[
    \int_XG(x,y)\,d\nu_i(y)
    \longrightarrow
    8\pi\int_XG(x,y)\,d\|T\|(y)
\]
uniformly for $x\in K$. Equation
\eqref{eq:ui_Green_representation} yields the first convergence in
\eqref{eq:scalar_Green_convergence_Rab}.

In particular, $u_i$ is uniformly bounded on $K$. Since
\[
    1-\frac{|\Phi_i|^2}{m_i^2}
    =
    \frac{2u_i}{m_i},
\]
we obtain \eqref{eq:uniform_Higgs_normalization_Rab}. Moreover,
\[
    m_i-|\Phi_i|
    =
    \frac{2u_i}
         {1+|\Phi_i|/m_i},
\]
so \eqref{eq:uniform_translated_Higgs_bound_Rab} follows. Finally,
\eqref{eq:scalar_amplitude_comparison} and the uniform bound just proved
give
\[
    0\leqslant v_i-u_i
    =
    \frac{(m_i-|\Phi_i|)^2}{2m_i}
    \longrightarrow0
\]
uniformly on $K$, proving the second convergence in
\eqref{eq:scalar_Green_convergence_Rab}.

The Green function is positive and
$\|T\|(X)=k>0$, so $v_T>0$. Since
$\operatorname{spt}\|T\|=\mathcal S$, the defining Green potential is
smooth and harmonic on $X\setminus\mathcal S$.
\end{proof}

\begin{remark}[Dirac character of the scalar Green potential]
\label{rmk:Dirac_character_scalar_Green_potential}
The singularities of $v_T$ along $\mathcal S$ are genuine Dirac-type
singularities. Indeed, by Li's local Dirac model
\cite[Theorem~4.3]{li2025large}, at $\|T\|$-almost every point
$p\in\mathcal S$ the blow-up of the singular abelian limit is the
translation-invariant Dirac model in the normal $3$-plane to the
approximate tangent plane of $T$. Under the normalization fixed in \S\ref{ss:Notation}, Li's scalar field
is $-v_T/2$. His local Dirac model therefore gives the corresponding
blow-up statement for the positive Green potential $v_T$.

In particular, the local behaviour is that of the Newtonian potential in
three normal variables. Thus $v_T$ is locally $L^p$ near the regular part
of $\mathcal S$ for every $p<3$, and in particular locally $L^2$, while
$dv_T$ is not locally square-integrable at points of positive
multiplicity. Hence local $L^2$-integrability of $v_T$ is compatible with
the Dirac singularity and does not imply removability.
\end{remark}

\begin{remark}[Comparison with the three-dimensional residual monopole]
\label{rmk:comparison_three_dimensional_residual_monopole}
The positive Green potential in
Proposition~\ref{prop:scalar_Green_convergence_Rab} is the direct
codimension-three extension of the residual field in the
three-dimensional theory. There the limiting measure has the form
\[
    4\pi\sum_a K_a\delta_{x_a},
\]
and, on $X\setminus\mathcal S$, the positive scalar Higgs defect
$m_i-|\Phi_i|$ converges smoothly to
\[
    4\pi\sum_aK_aG(\,\cdot\,,x_a).
\]
The longitudinal curvature converges to the Hodge dual of the negative
of its differential, and the translated configurations converge smoothly
after gauge to a reducible monopole with Dirac singularities of charges
$K_a$ \cite[Proposition~3.1, Corollary~3.2, and
Theorem~A]{fadel2026abelian}.

For comparison, \cite{fadel2026abelian} denotes the linear and
quadratic scalar defect fields by
\[
    u_i^{(3)}:=m_i-|\Phi_i|,
    \qquad
    v_i^{(3)}:=m_i^{-1}(m_i^2-|\Phi_i|^2).
\]
Consequently,
\[
    v_i=u_i^{(3)},
    \qquad
    2u_i=v_i^{(3)},
    \qquad
    v_T=u^{(3)}=\frac12v^{(3)}.
\]
When the calibrated cycle is replaced by the weighted concentration
$0$-cycle $\|T\|=\sum_aK_a\delta_{x_a}$, both theories therefore use
\[
    v_T
    =
    4\pi\int_XG(\,\cdot\,,y)\,d\|T\|(y),
    \qquad
    \widetilde\Phi_i=-v_i\Psi_i,
\]
and the limiting translated Higgs field is $-v_T\Psi_\infty$, with
longitudinal curvature $-*dv_T$.

There are two important differences. In dimension three the convergence
holds on all of $X\setminus\mathcal S$, because the limiting
nonabelian locus is precisely $\mathcal S$; here it is proved on the
potentially smaller open
set $\mathcal R_{\mathrm{ab}}=X\setminus\mathcal C$. Moreover, the
three-dimensional Bogomolny equation determines the limiting
longitudinal curvature directly from the scalar Green potential, so no
$L^2$-harmonic correction is needed at the curvature level. A flat
abelian holonomy ambiguity for the limiting connection may nevertheless
remain, as described in \cite[Theorem~A]{fadel2026abelian}. In the
present dimensions, Li's $L^2$-harmonic $2$-forms measure the additional
global ambiguity in the longitudinal curvature.
\end{remark}

\begin{corollary}[Subcritical longitudinal curvature on $\mathcal R_{\mathrm{ab}}$]
\label{cor:subcritical_longitudinal_curvature_Rab}
For every compact domain $K\Subset\mathcal R_{\mathrm{ab}}$ and every
$p\in[1,\infty]$, one has
\begin{equation}
\label{eq:subcritical_longitudinal_Lp_Rab}
    m_i^{-1/2}
    \left(
        \|\widehat\omega_i\|_{L^p(K)}
        +
        \|\omega_i\|_{L^p(K)}
    \right)
    \longrightarrow0.
\end{equation}
Moreover, for every
$1\leqslant p<\frac{2n}{n-2}$,
\begin{equation}
\label{eq:longitudinal_normalizations_agree_Rab}
    \|\widehat\omega_i-\omega_i\|_{L^p(K)}
    \longrightarrow0.
\end{equation}
In particular,
\[
    \|\widehat\omega_i\|_{L^2(K)}
    =
    o(m_i^{1/2}).
\]
\end{corollary}

\begin{proof}
Since
\[
    |\widehat\omega_i|
    \leqslant
    |F_{\nabla_i}|,
    \qquad
    |\omega_i|
    =
    \frac{|\Phi_i|}{m_i}|\widehat\omega_i|
    \leqslant
    |F_{\nabla_i}|,
\]
the $L^\infty$ case of
\eqref{eq:subcritical_longitudinal_Lp_Rab} follows from
Proposition~\ref{prop:mass_renormalized_density_decay_Rab}; the remaining
values of $p$ follow because $K$ has finite volume.

The relation
\[
    \widehat\omega_i-\omega_i
    =
    \frac{m_i-|\Phi_i|}{m_i}\widehat\omega_i
\]
gives, for every
$1\leqslant p<\frac{2n}{n-2}$,
\[
\begin{aligned}
    \|\widehat\omega_i-\omega_i\|_{L^p(K)}
    &\leqslant
    \left(
        m_i^{-1/2}
        \|m_i-|\Phi_i|\|_{L^p(K)}
    \right)
    \left(
        m_i^{-1/2}
        \|\widehat\omega_i\|_{L^\infty(K)}
    \right).
\end{aligned}
\]
The first factor tends to zero by Lemma~\ref{lem:largeness}, while the
second tends to zero by
\eqref{eq:subcritical_longitudinal_Lp_Rab}. This proves
\eqref{eq:longitudinal_normalizations_agree_Rab}.
\end{proof}

Define the \textbf{translated Higgs field}
\[
    \widetilde\Phi_i
    :=
    \Phi_i-m_i\Psi_i
    =
    -v_i\Psi_i.
\]
Then $|\widetilde\Phi_i|=v_i=m_i-|\Phi_i|$ is the Higgs defect.
Proposition~\ref{prop:scalar_Green_convergence_Rab} gives the
unconditional local bound
\[
    \sup_K|\widetilde\Phi_i|\leqslant C_K
\]
and identifies its positive scalar defect:
\[
    v_i\longrightarrow v_T
    \qquad
    \text{uniformly on }K.
\]
Compactness of $\widetilde\Phi_i$ as an adjoint-valued field still
requires control of the directions $\Psi_i$, and hence ultimately a
local gauge choice. The mass-rescaled $L^2$ convergence from
Lemma~\ref{lem:largeness},
\[
    m_i^{-1/2}\widetilde\Phi_i\longrightarrow0
    \qquad
    \text{in }L^2_{\mathrm{loc}},
\]
will also be used below as part of the basic large Higgs field identities.

\begin{lemma}[Basic large Higgs field identities]
\label{lem:basic_high_Higgs_identities}
On $K$, for all sufficiently large $i$, one has
\begin{equation}\label{eq:nabla_Psi_identity}
    \nabla_i\Psi_i
    =
    \frac{1}{|\Phi_i|}
    (\nabla_i\Phi_i)^\perp ,
\end{equation}
and
\begin{equation}\label{eq:nabla_Phi_tilde_identity}
    \nabla_i\widetilde\Phi_i
    =
    \left(
    1-\frac{m_i}{|\Phi_i|}
    \right)
    (\nabla_i\Phi_i)^\perp
    +
    (\nabla_i\Phi_i)^\parallel.
\end{equation}
Consequently,
\begin{equation}\label{eq:nabla_Phi_tilde_norm}
    |\nabla_i\widetilde\Phi_i|^2
    =
    \frac{(m_i-|\Phi_i|)^2}{|\Phi_i|^2}
    |(\nabla_i\Phi_i)^\perp|^2
    +
    |(\nabla_i\Phi_i)^\parallel|^2.
\end{equation}
Moreover,
\begin{equation}\label{eq:monopole_translated_error}
    |F_{\nabla_i}\wedge\Theta-*\nabla_i\widetilde\Phi_i|
    \leqslant
    C |(\nabla_i\Phi_i)^\perp|.
\end{equation}
Finally,
\begin{equation}\label{eq:translated_Higgs_mass_scaled_L2}
    m_i^{-1/2}\widetilde\Phi_i
    \longrightarrow 0
    \qquad
    \text{in }L^2(K).
\end{equation}
\end{lemma}

\begin{proof}
Since $\Phi_i=|\Phi_i|\Psi_i$, differentiating gives
\[
    \nabla_i\Phi_i
    =
    d|\Phi_i|\otimes\Psi_i
    +
    |\Phi_i|\nabla_i\Psi_i.
\]
The first term is parallel to $\Psi_i$, while $\nabla_i\Psi_i$ is orthogonal
to $\Psi_i$ because $|\Psi_i|=1$. Therefore
\[
    (\nabla_i\Phi_i)^\parallel=d|\Phi_i|\otimes\Psi_i,
    \qquad
    (\nabla_i\Phi_i)^\perp=|\Phi_i|\nabla_i\Psi_i,
\]
which proves \eqref{eq:nabla_Psi_identity}.

Next,
\[
    \widetilde\Phi_i
    =
    \Phi_i-m_i\Psi_i,
\]
so
\[
    \nabla_i\widetilde\Phi_i
    =
    \nabla_i\Phi_i-m_i\nabla_i\Psi_i.
\]
Using \eqref{eq:nabla_Psi_identity} and decomposing
$\nabla_i\Phi_i=(\nabla_i\Phi_i)^\parallel+(\nabla_i\Phi_i)^\perp$ gives
\eqref{eq:nabla_Phi_tilde_identity}. Since the two terms in
\eqref{eq:nabla_Phi_tilde_identity} are orthogonal, we obtain
\eqref{eq:nabla_Phi_tilde_norm}.

Finally, the $\Theta$-monopole equation gives
\[
    F_{\nabla_i}\wedge\Theta=*\nabla_i\Phi_i.
\]
Therefore
\[
    F_{\nabla_i}\wedge\Theta-*\nabla_i\widetilde\Phi_i
    =
    *(\nabla_i\Phi_i-\nabla_i\widetilde\Phi_i)
    =
    *\left(\frac{m_i}{|\Phi_i|}(\nabla_i\Phi_i)^\perp\right).
\]
Since $|\Phi_i|\geqslant m_i/2$ on $K$, this gives
\[
    |F_{\nabla_i}\wedge\Theta-*\nabla_i\widetilde\Phi_i|
    \leqslant
    C|(\nabla_i\Phi_i)^\perp|,
\]
which proves \eqref{eq:monopole_translated_error}.

It remains to prove \eqref{eq:translated_Higgs_mass_scaled_L2}. Since
\[
    \widetilde\Phi_i
    =
    -v_i\Psi_i
\]
and $|\Psi_i|=1$, we have
\[
    |\widetilde\Phi_i|
    =
    m_i-|\Phi_i|.
\]
Therefore
\[
    \|m_i^{-1/2}\widetilde\Phi_i\|_{L^2(K)}^2
    =
    m_i^{-1}
    \int_K
    (m_i-|\Phi_i|)^2.
\]
By Lemma~\ref{lem:largeness}, the right-hand side converges to zero.
This proves \eqref{eq:translated_Higgs_mass_scaled_L2}.
\end{proof}

\begin{corollary}[Decay of the transverse Higgs direction to all orders]
\label{cor:exp_decay_transverse_Higgs_direction}
Let $K\Subset \mathcal R_{\mathrm{ab}}$. Then, for every integer
$\ell\geqslant 0$, there exist constants $C_{K,\ell}<\infty$ and
$c_{K,\ell}>0$ such that, for all sufficiently large $i$,
\[
    |\nabla_i^\ell (F_{\nabla_i})^\perp|
    +
    |\nabla_i^\ell(\nabla_i\Phi_i)^\perp|
    +
    |\nabla_i^{\ell+1}\Psi_i|
    \leqslant
    C_{K,\ell}e^{-c_{K,\ell}m_i}
\]
on $K$.
\end{corollary}

\begin{proof}
This is the restriction to $K$ of
Lemma~\ref{lem:Li_dichotomy_away_Cinfty}. The estimate for
$(F_{\nabla_i})^\perp$ is
\eqref{ineq:Li_transverse_curvature_derivative_decay}, while the estimates
for $(\nabla_i\Phi_i)^\perp$ and $\Psi_i$ are
\eqref{ineq:Li_transverse_Higgs_Psi_derivative_decay}.
\end{proof}

The preceding identities show that, on compact subsets of
$\mathcal R_{\mathrm{ab}}$, the nonabelian monopole equation reduces to
its abelian longitudinal part up to exponentially small transverse error
terms. Consequently, the first gauge-invariant object to control is the
longitudinal curvature component
\[
    F_{\nabla_i}^{\parallel}
    =
    \langle F_{\nabla_i},\Psi_i\rangle\Psi_i.
\]
Since the transverse components vanish exponentially, any limiting curvature
on $\mathcal R_{\mathrm{ab}}$ must lie in the real line determined by the
large Higgs field direction. Gauge convergence of translated configurations would
require additional local control of this longitudinal curvature. The next
subsection isolates precisely this residual component.

\subsection{Hodge equations for the longitudinal curvature in the large mass limit}
\label{subsec: longitudinal_curvature}

We now explain how the remaining curvature component on
$\mathcal R_{\mathrm{ab}}$ satisfies the Hodge equations up to errors
that vanish in the large mass limit.

By the results of the previous subsection, on compact subsets of
$\mathcal R_{\mathrm{ab}}$ the Higgs fields remain uniformly large,
\[
    |\Phi_i|\geqslant \frac{m_i}{2},
\]
while the transverse components $(F_{\nabla_i})^\perp$ and
$(\nabla_i\Phi_i)^\perp$ decay exponentially fast, in all local
derivatives, with the mass. Thus the covariant derivative of the
normalized Higgs direction
\[
    \Psi_i=\frac{\Phi_i}{|\Phi_i|}
\]
decays exponentially in all local derivatives, and the monopole equation
reduces to its abelian longitudinal part up to errors of the same order.

Moreover,
Corollary~\ref{cor:codim3_obstruction_reduction_away_Cinfty} shows that
the full positive stress-energy imbalance tends to zero at
mass-renormalized codimension-two Morrey scale on
$\mathcal R_{\mathrm{ab}}$. Thus no codimension-three monotonicity
obstruction remains there. The longitudinal curvature is nevertheless
the only component which is not exponentially small before mass
renormalization, and it may fail to be locally bounded. It is this
unrenormalized component, rather than a residual stress-energy error,
which controls the subsequent compactness and regularity scale analysis.

Recall the notation from
Lemma~\ref{lem:exact_scalar_abelian_equations}:
\[
    \omega_i
    =
    m_i^{-1}\langle F_{\nabla_i},\Phi_i\rangle,
    \qquad
    \widehat\omega_i
    =
    \langle F_{\nabla_i},\Psi_i\rangle,
\]
so that
\[
    \omega_i
    =
    \frac{|\Phi_i|}{m_i}\widehat\omega_i.
\]
The form $\omega_i$ is the \textbf{mass-renormalized longitudinal curvature} adapted to the concentration measures and to Li's compactness theorem, while $\widehat\omega_i$ is the \textbf{longitudinal curvature normalized by $|\Phi_i|$} adapted to the large Higgs field region.

There is a related but distinct geometric object. On the complement of the
Higgs zero set, the unit section $\Psi_i$ determines the oriented rank-two
bundle
\[
    \langle\Psi_i\rangle^\perp\subset\mathfrak g_P.
\]
Following Li, write $F_{\mathrm U(1),i}$ for the normalized abelian
Chern form. It is defined globally for either structure group by the
Lie algebra formula
\begin{equation}
\label{eq:eigenline_Chern_curvature}
    F_{\mathrm U(1),i}
    =
    \frac{1}{4\pi}\widehat\omega_i
    -
    \frac{1}{8\pi}
    \left\langle
        \Psi_i,
        [\nabla_i\Psi_i\wedge\nabla_i\Psi_i]
    \right\rangle.
\end{equation}
The form $F_{\mathrm U(1),i}$ is closed. If
$G=\mathrm{SU}(2)$, it represents the first Chern class of the eigenline
selected by the convention in \cite[equation~(22)]{li2025large}. If
$G=\mathrm{SO}(3)$, then
\[
    [2F_{\mathrm U(1),i}]
    =
    e(\langle\Psi_i\rangle^\perp)
    \in H^2(X\setminus Z_i;\mathbb Z),
\]
where $\langle\Psi_i\rangle^\perp$ is oriented by the same stabilizer convention
while $[F_{\mathrm U(1),i}]$ may be only half-integral. Equivalently, on
every ball over which the $\mathrm{SO}(3)$ bundle lifts, the same form is
the first Chern form of the lifted eigenline.

In either case, $\widehat\omega_i/(4\pi)$ is the longitudinal part of the
abelian Chern curvature, but is not itself exactly closed. By
Corollary~\ref{cor:exp_decay_transverse_Higgs_direction}, for every
$K'\Subset\mathcal R_{\mathrm{ab}}$ and every $\ell\geqslant0$,
\[
    \left\|
        F_{\mathrm U(1),i}
        -
        \frac{1}{4\pi}\widehat\omega_i
    \right\|_{C^\ell(K')}
    \leqslant
    C_{K',\ell}e^{-c_{K',\ell}m_i}.
\]
Thus the distinction between the actual Chern curvature and its
longitudinal part is exponentially small on $\mathcal R_{\mathrm{ab}}$, while the
integral lattice of its cohomology class depends on the structure group.

As already noted in
Remark~\ref{rem:identification_with_Li_current}, it is useful to compare
this local analysis with Li's construction of the singular abelian limit. In
\cite[Theorem~1.12]{li2025large}, Li studies the corrected longitudinal
forms whose uncorrected parts are
\[
    -\frac{\operatorname{Tr}(\Phi_iF_{\nabla_i})}
    {4\pi m_i^{\mathrm{Li}}},
    \qquad
    -\frac{\operatorname{Tr}(\Phi_iF_{\nabla_i})}
    {4\pi|\Phi_i|_{\mathrm{Li}}}.
\]
By
\eqref{eq:PPS_Li_inner_product_comparison}--%
\eqref{eq:PPS_Li_mass_comparison}, these are precisely
\[
    \frac{1}{4\pi}\omega_i,
    \qquad
    \frac{1}{4\pi}\widehat\omega_i.
\]
Li proves that, after adding suitable $L^2$-harmonic corrections
$\sigma_i$, with
\[
    \|\sigma_i\|_{L^2(X)}\leqslant Cm_i,
\]
these corrected forms converge strongly in $L^1_{\mathrm{loc}}$ to a
normalized limiting form $\widetilde F_\infty$. Li then defines
\[
    F_\infty
    :=
    2\pi\widetilde F_\infty.
\]
Together with Li's limiting Higgs field, $F_\infty$ satisfies the
corresponding singular abelian $\mathrm G_2$- or Calabi--Yau monopole
equations, with
\[
    dF_\infty
    =
    2\pi T_{\mathrm{Li}}.
\]
In the Calabi--Yau case,
Proposition~\ref{prop:smooth_convergence_Li_corrected_longitudinal}
also gives
\[
    \Lambda F_\infty=0
\]
in the distributional sense. Thus, the corrected forms converge to
$\widetilde F_\infty=F_\infty/(2\pi)$, rather than directly to
$F_\infty$. In the alternative with bounded harmonic corrections, the
precise comparison with the smooth abelian limit of the translated
configurations is recorded in
Remark~\ref{rmk:comparison_with_Li_singular_limit}.

The purpose of the present subsection is different. We do not introduce
Li's global harmonic corrections. Instead, on compact subsets of
$\mathcal R_{\mathrm{ab}}$, where the large Higgs field splitting is
available, we study the uncorrected longitudinal forms of the
approximating sequence itself. We are therefore not re-proving the
extraction of the global singular abelian monopole. Rather, we isolate
the local differential and elliptic behaviour of the actual
longitudinal curvature, which is the component entering the
finite mass stress-energy and regularity scale analysis.

The key point is that the exterior derivatives and codifferentials of
both forms tend to zero on $\mathcal R_{\mathrm{ab}}$, but in different
topologies. For $\omega_i$, closedness and coclosedness hold only in a
weak mass-renormalized sense, because its differential contains terms
controlled using the vanishing of the mass-renormalized energy on compact
subsets of $X\setminus \mathcal S$. For $\widehat\omega_i$, they hold
up to errors tending to zero in all local $C^\ell$ norms, because the
errors involve only transverse components, which decay exponentially in
all local derivatives on $\mathcal R_{\mathrm{ab}}$.

These local Hodge estimates are used in
Section~\ref{sec:frequency_away_Cinfty} both in the regularity scale
analysis of the uncorrected sequence and in the proof of smooth convergence
of Li's corrected longitudinal forms on $\mathcal R_{\mathrm{ab}}$.

\begin{lemma}[Closedness of the longitudinal curvature in the large mass limit]
\label{lem:closedness_longitudinal_curvature_large_mass}
Let $K\Subset \mathcal R_{\mathrm{ab}}$. Then
\[
    d\omega_i\longrightarrow0
\]
locally on $K$ in $L^1$, while for every $\ell\geqslant 0$,
\[
    d\widehat\omega_i\longrightarrow0
\]
locally on $K$ in $C^\ell$.

More precisely, on $K$ one has
\begin{equation}\label{eq:domega_longitudinal_curvature}
    d\omega_i
    =
    |\Phi_i|^{-1}\omega_i\wedge d|\Phi_i|
    +
    \eta_i,
\end{equation}
where
\[
    \eta_i
    :=
    m_i^{-1}
    \left\langle
    F_{\nabla_i}\wedge(\nabla_i\Phi_i)^\perp
    \right\rangle.
\]
Moreover, for every $\ell\geqslant 0$,
\[
    \eta_i\longrightarrow0
\]
locally on $K$ in $C^\ell$, while the first term in
\eqref{eq:domega_longitudinal_curvature} converges to zero in
$L^1_{\mathrm{loc}}(K)$.

Furthermore, the longitudinal curvature forms normalized by $|\Phi_i|$ satisfy
\begin{equation}\label{eq:d_omega_hat}
    d\widehat\omega_i
    =
    \langle F_{\nabla_i}\wedge\nabla_i\Psi_i\rangle.
\end{equation}
\end{lemma}

\begin{proof}
Since $K\Subset \mathcal R_{\mathrm{ab}}$, Lemma~\ref{lem:Li_dichotomy_away_Cinfty}
gives, for all sufficiently large $i$,
\[
    |\Phi_i|\geqslant \frac{m_i}{2}
    \qquad\text{on }K.
\]
Thus $\Psi_i=\Phi_i/|\Phi_i|$ is well-defined on $K$, and
$\Phi_i=|\Phi_i|\Psi_i$. Differentiating gives
\[
    \nabla_i\Phi_i
    =
    d|\Phi_i|\otimes\Psi_i
    +
    |\Phi_i|\nabla_i\Psi_i.
\]
The first term is parallel to $\Psi_i$, while the second is orthogonal to
$\Psi_i$. Hence
\[
    (\nabla_i\Phi_i)^\parallel
    =
    d|\Phi_i|\otimes\Psi_i,
    \qquad
    (\nabla_i\Phi_i)^\perp
    =
    |\Phi_i|\nabla_i\Psi_i.
\]
Using the Bianchi identity $d_{\nabla_i}F_{\nabla_i}=0$, we compute
\[
\begin{aligned}
    d\omega_i
    &=
    m_i^{-1}d\langle F_{\nabla_i},\Phi_i\rangle  \\
    &=
    m_i^{-1}\langle F_{\nabla_i}\wedge\nabla_i\Phi_i\rangle  \\
    &=
    m_i^{-1}
    \left\langle
    F_{\nabla_i}\wedge
    \bigl(d|\Phi_i|\otimes\Psi_i
    +
    |\Phi_i|\nabla_i\Psi_i\bigr)
    \right\rangle.
\end{aligned}
\]
The contribution of the parallel part of $\nabla_i\Phi_i$ is
\[
    m_i^{-1}\langle F_{\nabla_i},\Psi_i\rangle\wedge d|\Phi_i|
    =
    m_i^{-1}\widehat\omega_i\wedge d|\Phi_i|.
\]
Since $\omega_i=(|\Phi_i|/m_i)\widehat\omega_i$, this term may
equivalently be written as
\[
    |\Phi_i|^{-1}\omega_i\wedge d|\Phi_i|.
\]
The remaining contribution is
\[
    m_i^{-1}
    \left\langle
    F_{\nabla_i}\wedge(\nabla_i\Phi_i)^\perp
    \right\rangle.
\]
This proves the decomposition \eqref{eq:domega_longitudinal_curvature}.

Similarly, using again the Bianchi identity,
\[
    d\widehat\omega_i
    =
    d\langle F_{\nabla_i},\Psi_i\rangle
    =
    \langle F_{\nabla_i}\wedge\nabla_i\Psi_i\rangle,
\]
which proves \eqref{eq:d_omega_hat}.

We now estimate the error terms. Since
$(\nabla_i\Phi_i)^\perp\perp\Psi_i$ and
$F_{\nabla_i}^{\parallel}=\widehat\omega_i\otimes\Psi_i$, the
longitudinal component of $F_{\nabla_i}$ does not contribute to
$\eta_i$. Thus
\[
    \eta_i
    =
    m_i^{-1}
    \left\langle
    (F_{\nabla_i})^\perp\wedge(\nabla_i\Phi_i)^\perp
    \right\rangle.
\]
Let $K'\Subset K$. By the higher order transverse decay estimates in
Lemma~\ref{lem:Li_dichotomy_away_Cinfty}, for every $\ell\geqslant 0$,
\[
    \|(F_{\nabla_i})^\perp\|_{C^\ell(K')}
    +
    \|(\nabla_i\Phi_i)^\perp\|_{C^\ell(K')}
    \leqslant
    C_{K',\ell}e^{-c_{K',\ell}m_i}.
\]
Using the product rule, we obtain, after changing constants,
\[
    \|\eta_i\|_{C^\ell(K')}
    \leqslant
    C_{K',\ell}e^{-c_{K',\ell}m_i}.
\]
Hence $\eta_i\to0$ locally in $C^\ell$ for every $\ell\geqslant0$.

We next estimate $d\widehat\omega_i$. Since
$\nabla_i\Psi_i\perp\Psi_i$, only the transverse curvature component
contributes to $\langle F_{\nabla_i}\wedge\nabla_i\Psi_i\rangle$.
Thus
\[
    d\widehat\omega_i
    =
    \left\langle
    (F_{\nabla_i})^\perp\wedge\nabla_i\Psi_i
    \right\rangle.
\]
Again by Lemma~\ref{lem:Li_dichotomy_away_Cinfty}, for every $\ell\geqslant 0$,
\[
    \|(F_{\nabla_i})^\perp\|_{C^\ell(K')}
    +
    \|\nabla_i\Psi_i\|_{C^\ell(K')}
    \leqslant
    C_{K',\ell}e^{-c_{K',\ell}m_i}.
\]
Therefore, by the product rule,
\[
    \|d\widehat\omega_i\|_{C^\ell(K')}
    \leqslant
    C_{K',\ell}e^{-c_{K',\ell}m_i}.
\]
Since $K'\Subset K$ was arbitrary, this proves that $d\widehat\omega_i\to0$ locally in $C^\ell$ for every $\ell\geqslant0$.

It remains to prove that $|\Phi_i|^{-1}\omega_i\wedge d|\Phi_i|\to0$ in
$L^1_{\mathrm{loc}}(K)$. Let $K'\Subset K$. Since $|\Phi_i|^{-1}\leqslant 2m_i^{-1}$ on $K'$ for all sufficiently large
$i$, we have
\[
\begin{aligned}
\left\|
|\Phi_i|^{-1}\omega_i\wedge d|\Phi_i|
\right\|_{L^1(K')}
&\leqslant
C m_i^{-1}
\|\omega_i\|_{L^2(K')}
\|d|\Phi_i|\|_{L^2(K')}.
\end{aligned}
\]
Moreover,
$|\omega_i|=m_i^{-1}|\langle F_{\nabla_i},\Phi_i\rangle|
\leqslant|F_{\nabla_i}|$. Since
$K'\Subset\mathcal R_{\mathrm{ab}}\subset X\setminus\mathcal S$ and
since $\mathcal S$ is the support of the limiting mass-renormalized
energy measure, we have
\[
    m_i^{-1}\int_{K'}|F_{\nabla_i}|^2\longrightarrow0,
    \qquad
    m_i^{-1}\int_{K'}|\nabla_i\Phi_i|^2\longrightarrow0.
\]
It follows that $\|\omega_i\|_{L^2(K')}=o(m_i^{1/2})$ and, by Kato's
inequality,
\[
    \|d|\Phi_i|\|_{L^2(K')}
    \leqslant
    \|\nabla_i\Phi_i\|_{L^2(K')}
    =
    o(m_i^{1/2}).
\]
Therefore $m_i^{-1}\|\omega_i\|_{L^2(K')}\|d|\Phi_i|\|_{L^2(K')}=o(1)$, and $|\Phi_i|^{-1}\omega_i\wedge d|\Phi_i|\to0$ in $L^1(K')$.

Since $K'\Subset K$ was arbitrary, and since $\eta_i\to0$ locally in
$C^\ell$ for every $\ell\geqslant 0$, the decomposition \eqref{eq:domega_longitudinal_curvature} gives $d\omega_i\to0$ locally on $K$ in $L^1$.
\end{proof}

We next analyse the codifferential of the longitudinal component. The
resulting identities are structurally analogous, but now the
Yang--Mills--Higgs equation enters through the divergence term
$d_{\nabla_i}^*F_{\nabla_i}$.

\begin{lemma}[Coclosedness of the longitudinal curvature in the large mass limit]
\label{lem:coclosedness_longitudinal_curvature_large_mass}
Let $K\Subset \mathcal R_{\mathrm{ab}}$. Then
\[
    d^*\omega_i\longrightarrow0
\]
locally on $K$ in $L^1$, while, for every $\ell\geqslant 0$,
\[
    d^*\widehat\omega_i\longrightarrow0
\]
locally on $K$ in $C^\ell$.

More precisely, on $K$ one has
\begin{equation}\label{eq:dstar_omega_decomposition}
    d^*\omega_i
    =
    -m_i^{-1} (\nabla |\Phi_i|)\lrcorner\,\widehat\omega_i
    +
    \zeta_i,
\end{equation}
where
\[
    \zeta_i
    :=
    -m_i^{-1}
    \sum_a
    \left\langle
    (F_{\nabla_i})^\perp(e_a,\cdot),
    (\nabla_i\Phi_i)^\perp_{e_a}
    \right\rangle ,
\]
for any local orthonormal frame $\{e_a\}$. Moreover, for every
$\ell\geqslant 0$,
\[
    \zeta_i\longrightarrow0
\]
locally on $K$ in $C^\ell$, while the first term in
\eqref{eq:dstar_omega_decomposition} converges to zero in
$L^1_{\mathrm{loc}}(K)$.

Furthermore,
\footnote{Here $\langle\cdot,\cdot\rangle_{\mathrm{contr}}$
denotes the contraction of the $T^*X$-indices using the Riemannian
metric together with the fibrewise inner product on
$\mathfrak g_P$. Thus, if
$F_{\nabla_i}=F_{ab}\,dx^a\wedge dx^b$ and
$\nabla_i\Psi_i=(\nabla_c\Psi_i)\,dx^c$, then
$\langle F_{\nabla_i},\nabla_i\Psi_i\rangle_{\mathrm{contr}}$
denotes the $1$-form with components
$g^{ac}\langle F_{ab},\nabla_c\Psi_i\rangle\,dx^b$.}
\begin{equation}\label{eq:dstar_omega_hat}
    d^*\widehat\omega_i
    =
    -\langle F_{\nabla_i},\nabla_i\Psi_i\rangle_{\mathrm{contr}}.
\end{equation}
\end{lemma}

\begin{proof}
Since $K\Subset\mathcal R_{\mathrm{ab}}$, Lemma~\ref{lem:Li_dichotomy_away_Cinfty} gives $|\Phi_i|\geqslant m_i/2$ on $K$ for all sufficiently large $i$. Hence $\Psi_i=\Phi_i/|\Phi_i|$ is well-defined on $K$.

Let $\{e_a\}$ be a local orthonormal frame. For any vector field $Y$,
using the definition of the codifferential and the $\mathrm{Ad}$-invariance
of the inner product, we compute
\begin{align*}
    (d^*\omega_i)(Y)
    &=
    -\sum_a
    (D_{e_a}\omega_i)(e_a,Y)
    \\
    &=
    -m_i^{-1}
    \sum_a
    \left\langle
    (\nabla_i)_{e_a}F_{\nabla_i}(e_a,Y),
    \Phi_i
    \right\rangle
    \\
    &\quad
    -m_i^{-1}
    \sum_a
    \left\langle
    F_{\nabla_i}(e_a,Y),
    \nabla_{i,e_a}\Phi_i
    \right\rangle.
\end{align*}
Since
\[
    d_{\nabla_i}^*F_{\nabla_i}(Y)
    =
    -\sum_a
    (\nabla_i)_{e_a}F_{\nabla_i}(e_a,Y),
\]
the first term is
\[
    m_i^{-1}
    \left\langle
    d_{\nabla_i}^*F_{\nabla_i}(Y),
    \Phi_i
    \right\rangle.
\]
Using the second equation in \eqref{eq:YMH_second_order},
$d_{\nabla_i}^*F_{\nabla_i} = [\nabla_i\Phi_i,\Phi_i]$, we obtain
\[
    \left\langle
    d_{\nabla_i}^*F_{\nabla_i},
    \Phi_i
    \right\rangle
    =
    \left\langle
    [\nabla_i\Phi_i,\Phi_i],
    \Phi_i
    \right\rangle
    =
    0
\]
by $\mathrm{Ad}$-invariance. Therefore
\begin{equation}\label{eq:dstar_omega_basic_formula}
    d^*\omega_i
    =
    -m_i^{-1}
    \sum_a
    \left\langle
    F_{\nabla_i}(e_a,\cdot),
    \nabla_{i,e_a}\Phi_i
    \right\rangle.
\end{equation}
We now decompose this expression with respect to the splitting determined
by $\Psi_i$. Write
\[
    F_{\nabla_i}
    =
    (F_{\nabla_i})^\parallel
    +
    (F_{\nabla_i})^\perp,
    \qquad
    \nabla_i\Phi_i
    =
    (\nabla_i\Phi_i)^\parallel
    +
    (\nabla_i\Phi_i)^\perp.
\]
Since
\[
    (F_{\nabla_i})^\parallel
    =
    \widehat\omega_i\otimes\Psi_i,
    \qquad
    (\nabla_i\Phi_i)^\parallel
    =
    d|\Phi_i|\otimes\Psi_i,
\]
and the parallel and transverse summands are orthogonal, the mixed terms
vanish. Thus \eqref{eq:dstar_omega_basic_formula} becomes
\[
    d^*\omega_i
    =
    -m_i^{-1}
    \sum_a
    \widehat\omega_i(e_a,\cdot)\,e_a(|\Phi_i|)
    -
    m_i^{-1}
    \sum_a
    \left\langle
    (F_{\nabla_i})^\perp(e_a,\cdot),
    (\nabla_i\Phi_i)^\perp_{e_a}
    \right\rangle.
\]
The first summation is $(\nabla|\Phi_i|)\lrcorner\,\widehat\omega_i$. This proves the decomposition \eqref{eq:dstar_omega_decomposition}.

We next compute the codifferential of the longitudinal curvature normalized by $|\Phi_i|$
form. Again using a local orthonormal frame $\{e_a\}$,
\[
    (d^*\widehat\omega_i)(Y)
    =
    -\sum_a
    D_{e_a}\widehat\omega_i(e_a,Y).
\]
Therefore
\[
    d^*\widehat\omega_i
    =
    \langle d_{\nabla_i}^*F_{\nabla_i},\Psi_i\rangle
    -
    \langle F_{\nabla_i},\nabla_i\Psi_i\rangle_{\mathrm{contr}}.
\]
Since
$d_{\nabla_i}^*F_{\nabla_i}=[\nabla_i\Phi_i,\Phi_i]$ and $\Psi_i$ is parallel to $\Phi_i$, one has $\langle[\nabla_i\Phi_i,\Phi_i],\Psi_i\rangle=0$. Hence
\[
    d^*\widehat\omega_i
    =
    -\langle
    F_{\nabla_i},
    \nabla_i\Psi_i
    \rangle_{\mathrm{contr}},
\]
which proves \eqref{eq:dstar_omega_hat}.

We now estimate the error terms. First, the longitudinal component of
$F_{\nabla_i}$ does not enter $\zeta_i$ by construction. Let
$K'\Subset K$. By the higher order transverse decay estimates in
Lemma~\ref{lem:Li_dichotomy_away_Cinfty}, for every $\ell\geqslant 0$,
\[
    \|(F_{\nabla_i})^\perp\|_{C^\ell(K')}
    +
    \|(\nabla_i\Phi_i)^\perp\|_{C^\ell(K')}
    \leqslant
    C_{K',\ell}e^{-c_{K',\ell}m_i}.
\]
Using the product rule, we obtain, after changing constants,
\[
    \|\zeta_i\|_{C^\ell(K')}
    \leqslant
    C_{K',\ell}e^{-c_{K',\ell}m_i}.
\]
Thus $\zeta_i\to0$ locally in $C^\ell$ for every $\ell\geqslant0$.

We next estimate $d^*\widehat\omega_i$. Since $F_{\nabla_i}^{\parallel}=\widehat\omega_i\otimes\Psi_i$ and $\nabla_i\Psi_i\perp\Psi_i$, only the transverse curvature component
contributes to $\langle F_{\nabla_i},\nabla_i\Psi_i\rangle_{\mathrm{contr}}$.
Therefore
\[
    d^*\widehat\omega_i
    =
    -\langle
    (F_{\nabla_i})^\perp,
    \nabla_i\Psi_i
    \rangle_{\mathrm{contr}}.
\]
Again by Lemma~\ref{lem:Li_dichotomy_away_Cinfty}, for every $\ell\geqslant 0$,
\[
    \|(F_{\nabla_i})^\perp\|_{C^\ell(K')}
    +
    \|\nabla_i\Psi_i\|_{C^\ell(K')}
    \leqslant
    C_{K',\ell}e^{-c_{K',\ell}m_i}.
\]
The product rule then gives
\[
    \|d^*\widehat\omega_i\|_{C^\ell(K')}
    \leqslant
    C_{K',\ell}e^{-c_{K',\ell}m_i}.
\]
Since $K'\Subset K$ was arbitrary, this proves that $d^*\widehat\omega_i\to0$ locally in $C^\ell$ for every $\ell\geqslant0$.

It remains to prove that $m_i^{-1}(\nabla|\Phi_i|)\lrcorner\,\widehat\omega_i\to0$ in $L^1_{\mathrm{loc}}(K)$. Let $K'\Subset K$. Then
\[
\begin{aligned}
\left\|
m_i^{-1}(\nabla|\Phi_i|)\lrcorner\,\widehat\omega_i
\right\|_{L^1(K')}
&\leqslant
C m_i^{-1}
\|\widehat\omega_i\|_{L^2(K')}
\|d|\Phi_i|\|_{L^2(K')}.
\end{aligned}
\]
The final part of the proof of
Lemma~\ref{lem:closedness_longitudinal_curvature_large_mass} observes
that $K'\Subset\mathcal R_{\mathrm{ab}}\subset X\setminus\mathcal S$
and consequently establishes
$\|\widehat\omega_i\|_{L^2(K')}=o(m_i^{1/2})$ and
$\|d|\Phi_i|\|_{L^2(K')}=o(m_i^{1/2})$. Therefore
$m_i^{-1}\|\widehat\omega_i\|_{L^2(K')}
\|d|\Phi_i|\|_{L^2(K')}=o(1)$, and hence
$m_i^{-1}(\nabla|\Phi_i|)\lrcorner\,\widehat\omega_i\to0$ in
$L^1(K')$.

Since $K'\Subset K$ was arbitrary, and since $\zeta_i\to0$ locally in
$C^\ell$ for every $\ell\geqslant0$, the decomposition \eqref{eq:dstar_omega_decomposition} gives $d^*\omega_i\to0$ locally on $K$ in $L^1$.
\end{proof}

Together, the previous two lemmas show that the longitudinal curvature
component is closed and coclosed up to terms tending to zero on compact
subsets of $\mathcal R_{\mathrm{ab}}$.

\begin{corollary}[Hodge equations for the longitudinal component in the large mass limit]
\label{cor:Hodge_longitudinal_component_large_mass}
Let $K\Subset \mathcal R_{\mathrm{ab}}$. Then
\[
    d\omega_i\to0,
    \qquad
    d^*\omega_i\to0
\]
locally on $K$ in $L^1$, while, for every $\ell\geqslant 0$,
\[
    d\widehat\omega_i\to0,
    \qquad
    d^*\widehat\omega_i\to0
\]
locally on $K$ in $C^\ell$.

In particular, any local $L^2$ weak limit of either $\omega_i$ or
$\widehat\omega_i$ on $K$ is a harmonic $2$-form; in the
Calabi--Yau case, every such limit is primitive. Since
\[
    F_{\nabla_i}^{\parallel}
    =
    \widehat\omega_i\otimes\Psi_i,
\]
this shows that the only curvature component which can survive on
$\mathcal R_{\mathrm{ab}}$ without exponential decay is a real-valued
$2$-form satisfying the Hodge equations in the limit.
\end{corollary}

\begin{proof}
The convergence of $d\omega_i$ in $L^1_{\mathrm{loc}}$ and the
$C^\ell_{\mathrm{loc}}$ convergence of $d\widehat\omega_i$ follow from
Lemma~\ref{lem:closedness_longitudinal_curvature_large_mass}. The
corresponding statements for the codifferentials follow from
Lemma~\ref{lem:coclosedness_longitudinal_curvature_large_mass}.

Let first $\omega_i\rightharpoonup\omega_{\mathrm{har}}$ weakly in
$L^2_{\mathrm{loc}}$ along a subsequence. Since
\[
    d\omega_i\to0,
    \qquad
    d^*\omega_i\to0
\]
locally in $L^1$, we may pass to the limit in the distributional identities
and obtain
\[
    d\omega_{\mathrm{har}}=0,
    \qquad
    d^*\omega_{\mathrm{har}}=0.
\]
Hence $\omega_{\mathrm{har}}$ is harmonic.

Similarly, if
\[
    \widehat\omega_i\rightharpoonup\widehat\omega_{\mathrm{har}}
\]
weakly in $L^2_{\mathrm{loc}}$, then the convergences
\[
    d\widehat\omega_i\to0,
    \qquad
    d^*\widehat\omega_i\to0
\]
locally in $C^\ell$, in particular locally uniformly, imply
\[
    d\widehat\omega_{\mathrm{har}}=0,
    \qquad
    d^*\widehat\omega_{\mathrm{har}}=0
\]
in the distributional sense. Thus $\widehat\omega_{\mathrm{har}}$ is harmonic. In the Calabi--Yau case,
$\Lambda\omega_i=\Lambda\widehat\omega_i=0$ by
Lemma~\ref{lem:exact_scalar_abelian_equations}. Passing to the weak $L^2_{\mathrm{loc}}$ limits gives
$\Lambda\omega_{\mathrm{har}}
=\Lambda\widehat\omega_{\mathrm{har}}=0$.
\end{proof}

The Hodge estimates obtained above place the longitudinal component
within the scope of elliptic compactness theory. The key difference
between $\omega_i$ and $\widehat\omega_i$ is that closedness and
coclosedness of the mass-renormalized forms $\omega_i$ hold only in a
weak local sense, whereas for the longitudinal forms normalized by
$|\Phi_i|$, $\widehat\omega_i$, they hold up to terms tending to zero
in all local $C^\ell$ norms on $\mathcal R_{\mathrm{ab}}$.

\begin{proposition}[Elliptic compactness of the longitudinal component]
\label{prop:elliptic_compactness_longitudinal_curvature}
Let
\[
    K'\Subset K\Subset \mathcal R_{\mathrm{ab}}.
\]
Suppose first that
\[
    \sup_i\|\omega_i\|_{L^2(K)}<\infty.
\]
Then, after passing to a subsequence, there exists a harmonic $2$-form
$\omega_{\mathrm{har}}$ on a neighbourhood of $K'$ such that
\[
    \omega_i\rightharpoonup\omega_{\mathrm{har}}
    \qquad
    \text{weakly in }L^2(K').
\]
If, in addition, for some $p>1$ the sequence $\omega_i$ is locally bounded
in $L^p$ and
\[
    d\omega_i\to0,
    \qquad
    d^*\omega_i\to0
\]
in $L^p_{\mathrm{loc}}(K')$, then, after passing to a subsequence,
\[
    \omega_i\to\omega_{\mathrm{har}}
    \qquad
    \text{strongly in }W^{1,p}_{\mathrm{loc}}(K').
\]
In particular, when $p=2$ this gives strong $L^2_{\mathrm{loc}}$
convergence.

For the longitudinal forms normalized by $|\Phi_i|$, suppose that
\[
    \sup_i\|\widehat\omega_i\|_{L^2(K)}<\infty.
\]
Then, after passing to a subsequence, there exists a harmonic $2$-form
$\widehat\omega_{\mathrm{har}}$ on a neighbourhood of $K'$ such that
\[
    \widehat\omega_i\to\widehat\omega_{\mathrm{har}}
    \qquad
    \text{smoothly on }K'.
\]
In the Calabi--Yau case, each of the corresponding limiting forms is primitive.
\end{proposition}

\begin{proof}
We first consider the mass-renormalized forms $\omega_i$. Choose a
compact domain $K''$ such that
\[
    K'\Subset K''\Subset K.
\]
By the assumed $L^2(K)$ bound, after passing to a subsequence,
\[
    \omega_i\rightharpoonup\omega_{\mathrm{har}}
    \qquad
    \text{weakly in }L^2(K'').
\]
By Corollary~\ref{cor:Hodge_longitudinal_component_large_mass},
\[
    d\omega_i\to0,
    \qquad
    d^*\omega_i\to0
\]
locally in $L^1$. Passing to the limit in the sense of distributions gives
\[
    d\omega_{\mathrm{har}}=0,
    \qquad
    d^*\omega_{\mathrm{har}}=0.
\]
Thus $\omega_{\mathrm{har}}$ is harmonic on $\operatorname{int}K''$,
which is a neighbourhood of $K'$.

Assume now that, for some $p>1$, the sequence $\omega_i$ is locally bounded
in $L^p$ and that
\[
    d\omega_i,d^*\omega_i\to0
    \qquad
    \text{in }L^p_{\mathrm{loc}}(K').
\]
Let
\[
    U\Subset V\Subset W\Subset K'.
\]
The local elliptic estimate for the Hodge--Dirac operator
$\slashed D:=d+d^*$, first applied on $V$ with data on $W$, gives
\[
    \|\omega_i\|_{W^{1,p}(V)}
    \leqslant
    C
    \left(
        \|\omega_i\|_{L^p(W)}
        +
        \|d\omega_i\|_{L^p(W)}
        +
        \|d^*\omega_i\|_{L^p(W)}
    \right).
\]
Hence $\omega_i$ is uniformly bounded in $W^{1,p}(V)$. By Rellich
compactness, after passing to a subsequence,
\[
    \omega_i\longrightarrow\omega_{\mathrm{har}}
    \qquad
    \text{strongly in }L^p(V),
\]
where the limit is the same harmonic form obtained from the weak
$L^2_{\mathrm{loc}}$ convergence. We now apply the elliptic estimate to
$\omega_i-\omega_{\mathrm{har}}$ on $U$, with data on $V$:
\[
\begin{aligned}
    \|\omega_i-\omega_{\mathrm{har}}\|_{W^{1,p}(U)}
    \leqslant C\bigl(
        &\|\omega_i-\omega_{\mathrm{har}}\|_{L^p(V)}
        +
        \|d\omega_i\|_{L^p(V)}  \\
        &+
        \|d^*\omega_i\|_{L^p(V)}
    \bigr),
\end{aligned}
\]
because $d\omega_{\mathrm{har}}=d^*\omega_{\mathrm{har}}=0$. Every term on the
right-hand side tends to zero. Thus
\[
    \omega_i\longrightarrow\omega_{\mathrm{har}}
    \qquad
    \text{strongly in }W^{1,p}(U).
\]
Since $U\Subset K'$ was arbitrary, this proves the asserted
$W^{1,p}_{\mathrm{loc}}$ convergence.

We now turn to the longitudinal forms normalized by $|\Phi_i|$ $\widehat\omega_i$.
By Corollary~\ref{cor:Hodge_longitudinal_component_large_mass}, for
every $\ell\geqslant 0$,
\[
    d\widehat\omega_i\to0,
    \qquad
    d^*\widehat\omega_i\to0
\]
locally in $C^\ell$ on $K$. Equivalently, 
\[
    \slashed{D}\widehat\omega_i\to0
    \qquad
    \text{locally in }C^\ell
\]
for every $\ell\geqslant 0$. Again, since the operator $\slashed{D}$ is elliptic, standard
interior Schauder estimates give, for every
\[
    U\Subset V\Subset K,
    \qquad
    \ell\geqslant 0,
    \qquad
    \alpha\in(0,1),
\]
an estimate of the form
\[
    \|\widehat\omega_i\|_{C^{\ell+1,\alpha}(U)}
    \leqslant
    C_{\ell,\alpha,U,V}
    \left(
        \|\widehat\omega_i\|_{L^2(V)}
        +
        \|\slashed{D}\widehat\omega_i\|_{C^{\ell,\alpha}(V)}
    \right).
\]
Since $\slashed{D}\widehat\omega_i\to0$ locally in $C^j$ for every $j$, the
$C^{\ell,\alpha}$ norms of $\slashed{D}\widehat\omega_i$ are uniformly bounded, and
indeed tend to zero after increasing the integer regularity index if
necessary. Together with the assumed $L^2(K)$ bound, the Schauder estimates
therefore give uniform $C^{\ell+1,\alpha}$ bounds for
$\widehat\omega_i$ on compact subsets of $K$, for every $\ell$ and every
$\alpha\in(0,1)$.

By Arzelà--Ascoli and a diagonal argument over an exhaustion of $K'$ and
over $\ell$, after passing to a subsequence we obtain
\[
    \widehat\omega_i\to\widehat\omega_{\mathrm{har}}
\]
smoothly on $K'$. Passing to the limit in
\[
    d\widehat\omega_i\to0,
    \qquad
    d^*\widehat\omega_i\to0
\]
gives
\[
    d\widehat\omega_{\mathrm{har}}=0,
    \qquad
    d^*\widehat\omega_{\mathrm{har}}=0.
\]
Thus $\widehat\omega_{\mathrm{har}}$ is harmonic.
\end{proof}

\begin{corollary}[Equivalence of the two local longitudinal compactness regimes]
\label{cor:equivalent_longitudinal_L2_compactness}
Let $K'\Subset K\Subset\mathcal R_{\mathrm{ab}}$. The following are
equivalent:
\[
    \sup_i\|\omega_i\|_{L^2(K)}<\infty,
    \qquad
    \sup_i\|\widehat\omega_i\|_{L^2(K)}<\infty.
\]
If either condition holds, then, after passing to a subsequence, there is
a harmonic $2$-form $\omega_{\mathrm{har}}$ on a neighbourhood of $K'$
such that
\[
    \widehat\omega_i\longrightarrow\omega_{\mathrm{har}}
    \quad\text{smoothly on }K',
    \qquad
    \omega_i\longrightarrow\omega_{\mathrm{har}}
    \quad\text{strongly in }L^2(K').
\]
In the Calabi--Yau case, the common limit is primitive.

If instead
$\|\widehat\omega_i\|_{L^2(K)}\to\infty$, then
\[
    \left|
        \|\widehat\omega_i\|_{L^2(K)}
        -
        \|\omega_i\|_{L^2(K)}
    \right|
    \longrightarrow0,
\]
and, after normalizing each form by its own $L^2(K)$ norm, the two
sequences have the same local $L^2$ limits.
\end{corollary}

\begin{proof}
By Corollary~\ref{cor:subcritical_longitudinal_curvature_Rab},
\[
    \|\widehat\omega_i-\omega_i\|_{L^2(K)}\longrightarrow0.
\]
This proves the equivalence of the two boundedness conditions. Under
either condition, the forms normalized by $|\Phi_i|$ are $L^2$-bounded, so
Proposition~\ref{prop:elliptic_compactness_longitudinal_curvature} gives,
after passing to a subsequence, smooth convergence of
$\widehat\omega_i$ on $K'$ to a harmonic form. The $L^2$ comparison then
gives strong $L^2(K')$ convergence of $\omega_i$ to the same limit.
The final assertion follows from the reverse triangle inequality and,
after division by the diverging norms, from the same $L^2$ comparison.
\end{proof}

The previous proposition and corollary show that, on every compact subset
of $\mathcal R_{\mathrm{ab}}$, a uniform local $L^2$ bound for either
longitudinal normalization yields subsequential smooth convergence of the
form normalized by $|\Phi_i|$ and strong $L^2$ convergence of the mass-normalized form
to the same harmonic $2$-form. The next corollary records the complementary normalized blow-up alternative when the local $L^2$ norms have infinite limsup.

\begin{corollary}[Local harmonic compactness for the longitudinal curvature normalized by $|\Phi_i|$]
\label{cor:local_harmonic_compactness_longitudinal_curvature}
Let $K'\Subset K\Subset \mathcal R_{\mathrm{ab}}$. If
\[
    \sup_i\|\widehat\omega_i\|_{L^2(K)}<\infty,
\]
then, after passing to a subsequence, there exists a harmonic $2$-form
$\widehat\omega_{\mathrm{har}}$ on a neighbourhood of $K'$ such that
\[
    \widehat\omega_i\to\widehat\omega_{\mathrm{har}}
    \qquad
    \text{smoothly on }K'.
\]
If instead
\[
    \limsup_{i\to\infty}
    \|\widehat\omega_i\|_{L^2(K)}
    =
    +\infty,
\]
then, after passing to a subsequence and relabelling, we may assume that
\[
    \|\widehat\omega_i\|_{L^2(K)}
    \longrightarrow
    +\infty.
\]
Define
\[
    \overline\omega_i
    :=
    \frac{\widehat\omega_i}
    {\|\widehat\omega_i\|_{L^2(K)}}.
\]
Then, after passing to a subsequence, there exists a harmonic $2$-form
$\overline\omega_{\mathrm{har}}$ on a neighbourhood of $K'$ such that
\[
    \overline\omega_i\to\overline\omega_{\mathrm{har}}
    \qquad
    \text{smoothly on }K'.
\]
In the Calabi--Yau case, the harmonic limit in either alternative is primitive.
\end{corollary}
\begin{remark}
At this stage, the local compactness argument alone allows the limiting
form to vanish on $K'$ if the normalized $L^2$ mass escapes towards
$K\setminus K'$. Corollary~\ref{cor:harmonic_correction_alternatives}
shows that this does not occur for the unbounded alternative in the Main
Setting.
\end{remark}
\begin{proof}
The bounded case is the second part of
Proposition~\ref{prop:elliptic_compactness_longitudinal_curvature}. In the
unbounded case, the normalized forms satisfy
\[
    \|\overline\omega_i\|_{L^2(K)}=1.
\]
Moreover, by Corollary~\ref{cor:Hodge_longitudinal_component_large_mass},
for every $\ell\geqslant 0$,
\[
    d\widehat\omega_i\to0,
    \qquad
    d^*\widehat\omega_i\to0
\]
locally in $C^\ell$ on $K$. Since
\[
    \|\widehat\omega_i\|_{L^2(K)}\to+\infty,
\]
the same estimates hold for $\overline\omega_i$, indeed with errors divided
by $\|\widehat\omega_i\|_{L^2(K)}$. 
The proof of Proposition~\ref{prop:elliptic_compactness_longitudinal_curvature}
uses only a uniform local $L^2$ bound and the local $C^\ell$ convergence of
the Hodge errors. These hypotheses hold for $\overline\omega_i$. Hence, the same elliptic compactness argument, applied first on an intermediate
compact domain $K''$ with $K'\Subset K''\Subset K$, gives smooth
subconvergence on $K'$ to a $2$-form harmonic on a neighbourhood of $K'$.

In the Calabi--Yau case,
$\Lambda\overline\omega_i=0$ for every sufficiently large $i$, since $\Lambda\widehat\omega_i=0$. Passing to the smooth limit gives $\Lambda\overline\omega_{\mathrm{har}}=0$.
\end{proof}

\begin{remark}[Relation with Li's corrected longitudinal forms]
\label{rmk:residual_harmonic_component}
On compact subsets of $\mathcal R_{\mathrm{ab}}$, the transverse
components decay exponentially and the uncorrected longitudinal forms
$\widehat\omega_i$ satisfy the Hodge equations up to errors that vanish
in the large mass limit. Li's theorem gives
$L^2$-harmonic corrections $\sigma_i$ such that
\[
    \frac{1}{4\pi}\widehat\omega_i+\sigma_i
\]
converges strongly in $L^1_{\mathrm{loc}}$. In
Proposition~\ref{prop:smooth_convergence_Li_corrected_longitudinal}, the
Hodge equations are combined with this convergence to obtain smooth
convergence on compact subsets of $\mathcal R_{\mathrm{ab}}$.

The corrections may be unbounded in
$\mathscr H_{(2)}^2(X)\simeq H_c^2(X;\mathbb R)$, so Li's corrected
compactness does not by itself bound the uncorrected longitudinal
curvature entering the stress-energy and regularity scale analysis.
The bounded and unbounded alternatives are treated in
Section~\ref{sec:frequency_away_Cinfty}.
\end{remark}


\section{Longitudinal curvature compactness and regularity scales}
\label{sec:frequency_away_Cinfty}

We study compactness and the ordinary Yang--Mills--Higgs regularity scale
on the open set $\mathcal R_{\mathrm{ab}}$. Under the standing convention of Section~\ref{sec: abelianization}, this region is assumed nonempty, and all statements below are local on its compact subsets. The global estimate for Li's $L^2$-harmonic corrections in
Proposition~\ref{prop:subcritical_harmonic_correction_growth}, however, does not require $\mathcal R_{\mathrm{ab}}$ to be nonempty; its proof uses a compact domain in $X\setminus\mathcal S$ and Li's global corrected $L^1_{\mathrm{loc}}$ convergence.

The Higgs fields are uniformly large there on compact subsets of $\mathcal R_{\mathrm{ab}}$, the transverse curvature components and the covariant derivatives of the normalized Higgs directions decay exponentially, and the real-valued longitudinal forms
\[
    \widehat\omega_i
    :=
    \langle F_{\nabla_i},\Psi_i\rangle,
    \qquad
    \Psi_i:=\frac{\Phi_i}{|\Phi_i|},
\]
are closed and coclosed up to errors tending to zero in every local
$C^\ell$ norm. Corollary~\ref{cor:subcritical_longitudinal_curvature_Rab} already gives,
for every compact domain $K\Subset\mathcal R_{\mathrm{ab}}$,
\[
    m_i^{-1/2}
    \|\widehat\omega_i\|_{L^\infty(K)}
    \longrightarrow0,
    \qquad
    \|\widehat\omega_i\|_{L^2(K)}
    =
    o(m_i^{1/2}).
\]
Thus the basic local subcriticality of the uncorrected longitudinal
curvature is independent of Li's harmonic corrections. The rôle of the
corrections below is to provide global control of the harmonic drift,
local subcriticality for all derivatives, and a canonical description of the
normalized unbounded alternative.

The first subsection treats the locally $L^2$-bounded case for the
uncorrected forms and then combines the Hodge equations with Li's corrected
compactness theorem. We prove that the corrected forms converge smoothly
on compact subsets of $\mathcal R_{\mathrm{ab}}$. If the corresponding $L^2$-harmonic corrections are uniformly bounded, the translated configurations converge smoothly to an abelian $\Theta$-monopole. If their $L^2(X)$ norms have infinite limsup, then,
after selecting a subsequence along which the norms tend to infinity,
their normalizations converge to a nonzero global $L^2$-harmonic
$2$-form. The second subsection treats this divergent subsequence and
estimates the regularity scale away from the zero set of the limiting
harmonic form. The final subsection studies its zeros and states the conditional estimate under the quantitative order-$q$ profile hypothesis. Away from the zero set of the limiting form, the regularity scale is comparable to the inverse square-root of the local
longitudinal amplitude. At a zero, the exponent $1/(q+2)$ requires the
additional hypothesis at shrinking scales of Definition~\ref{def:quantitative_order_q_profile}.

\subsection{Local compactness and Li's corrected longitudinal forms}
\label{subsec:bounded_longitudinal_regular}

Recall the regularity scale
$\mathfrak r_i:X\to(0,r_0]$ introduced in
Subsection~\ref{subsec:calibrated_concentration_regular_scales}:
\[
    \mathfrak r_i(x)
    =
    \sup\left\{
        0<r\leqslant r_0:
        s^{4-n}
        \int_{B_s(x)}
        e_i
        <
        \varepsilon_0
        \text{ for every }0<s\leqslant r
    \right\}.
\]
The small energy condition is imposed at all smaller radii in order to
make $\mathfrak r_i(x)$ a genuine regularity scale. By the standard
codimension-four almost-monotonicity formula for the Yang--Mills--Higgs
energy, recalled in Appendix~\ref{app: B}, this definition is equivalent,
up to changing $r_0$ and $\varepsilon_0$ by fixed geometric constants, to
checking the corresponding small energy condition at the outer radius.
Thus, a positive lower bound for $\mathfrak r_i(x)$ is precisely the input
needed to apply the $\varepsilon$-regularity theorem uniformly.

We first consider the case in which the longitudinal curvature
$\widehat\omega_i$, normalized by $|\Phi_i|$, is locally bounded in $L^2$. In this regime, these Hodge estimates, together with the elliptic compactness theory of Section~\ref{sec: abelianization}, give uniform curvature
bounds on smaller compact subsets.

\begin{proposition}[$L^2$-bounded longitudinal curvature implies local regularity]
\label{prop:bounded_longitudinal_curvature_regular}
Let $K'\Subset K\Subset \mathcal R_{\mathrm{ab}}$. Assume that
\[
    \sup_i\|\widehat\omega_i\|_{L^2(K)}<\infty,
    \qquad
    \widehat\omega_i:=\langle F_{\nabla_i},\Psi_i\rangle.
\]
Then the Yang--Mills--Higgs energy densities $e_i$ are uniformly bounded
on $K'$. In particular, there exists $\mathfrak{r}_K>0$ such that
\[
    \mathfrak{r}_i(x)\geqslant \mathfrak{r}_K
\]
for every $x\in K'$ and all sufficiently large $i$.
\end{proposition}

\begin{proof}
Choose an intermediate compact domain $K'\Subset K''\Subset K$. By Corollary~\ref{cor:local_harmonic_compactness_longitudinal_curvature}, the uniform $L^2(K)$ bound for $\widehat\omega_i$ gives smooth
precompactness of the longitudinal curvature forms normalized by $|\Phi_i|$ on $K''$.
In particular, after increasing the constant if necessary, there exists
$C_K<\infty$ such that
\[
    \sup_{K''}|\widehat\omega_i|\leqslant C_K
\]
for all sufficiently large $i$.

Since
\[
    F_{\nabla_i}^{\parallel}
    =
    \widehat\omega_i\otimes\Psi_i,
\]
the longitudinal curvature components are uniformly bounded on $K''$.
On the other hand, Lemma~\ref{lem:Li_dichotomy_away_Cinfty} gives
exponential decay of the transverse curvature components
$(F_{\nabla_i})^\perp$ on compact subsets of
$\mathcal R_{\mathrm{ab}}$. Hence
\[
    \sup_{K''}|F_{\nabla_i}|\leqslant C_K
\]
for all sufficiently large $i$.

The $\Theta$-monopole equation gives
\[
    |\nabla_i\Phi_i|
    =
    |F_{\nabla_i}\wedge\Theta|
    \leqslant
    C|F_{\nabla_i}|,
\]
and therefore
\[
    \sup_{K''}|\nabla_i\Phi_i|\leqslant C_K.
\]
Thus
\[
    \sup_{K''} e_i \leqslant C_K
\]
for all sufficiently large $i$.

Let
\[
    d_K:=\operatorname{dist}(K',X\setminus K'')>0.
\]
For every $x\in K'$ and every $0<s\leqslant d_K$, one has $B_s(x)\subset K''$, and hence
\[
    s^{4-n}\int_{B_s(x)}e_i
    \leqslant
    C_K s^{4-n}\operatorname{vol}(B_s(x))
    \leqslant
    C_K s^4,
\]
after adjusting $C_K$ using the bounded geometry of $(X,g)$ on $K''$.

Choose
\[
    0<\mathfrak{r}_K\leqslant \min\{r_0,d_K\}
\]
so small that
\[
    C_K\mathfrak{r}_K^4<\varepsilon_0.
\]
Then, for every $x\in K'$, every $0<s\leqslant\mathfrak{r}_K$, and all
sufficiently large $i$,
\[
    s^{4-n}
    \int_{B_s(x)}e_i
    <
    \varepsilon_0.
\]
By the definition of the regularity scale, this gives
\[
    \mathfrak{r}_i(x)\geqslant\mathfrak{r}_K
\]
for every $x\in K'$ and all sufficiently large $i$.
\end{proof}

\begin{corollary}[Smooth compactness in the $L^2$-bounded longitudinal alternative]
\label{cor:smooth_compactness_L2_bounded_longitudinal}
Assume that, for every $K\Subset \mathcal R_{\mathrm{ab}}$, one has
\[
    \sup_i\|\widehat\omega_i\|_{L^2(K)}<\infty.
\]
Then, after passing to a subsequence and applying gauge transformations,
the translated configurations
\[
    (\nabla_i,\widetilde\Phi_i),
    \qquad
    \widetilde\Phi_i:=\Phi_i-m_i\Psi_i,
    \qquad
    \Psi_i:=\frac{\Phi_i}{|\Phi_i|},
\]
converge smoothly on compact subsets of
$\mathcal R_{\mathrm{ab}}$ to a smooth abelian
$\Theta$-monopole $(\nabla_{\mathrm{ab}},\widetilde\Phi_{\mathrm{ab}})$.
\end{corollary}

\begin{proof}
Choose a compact exhaustion
\[
    K_1\Subset K_2\Subset\cdots\Subset\mathcal R_{\mathrm{ab}},
    \qquad
    \bigcup_{j\geqslant1}K_j=\mathcal R_{\mathrm{ab}},
\]
with each $K_j$ contained in the interior of $K_{j+1}$. After passing to
a subsequence, we may also assume that, for every $j$, the Higgs field
$\Phi_i$ is nonvanishing on $K_j$ for all $i\geqslant j$. Hence
$\Psi_i$, $v_i$, and $\widetilde\Phi_i$ are defined on $K_j$ for all
$i\geqslant j$.

By the assumed local $L^2$ bound and
Proposition~\ref{prop:elliptic_compactness_longitudinal_curvature}, a
diagonal subsequence extraction gives a harmonic real-valued $2$-form
$\widehat\omega_{\mathrm{har}}$ on
$\mathcal R_{\mathrm{ab}}$ such that
\[
    \widehat\omega_i
    \longrightarrow
    \widehat\omega_{\mathrm{har}}
\]
smoothly on compact subsets of $\mathcal R_{\mathrm{ab}}$.

We next consider the scalar residual fields
\[
    v_i
    :=
    m_i-|\Phi_i|.
\]
By Proposition~\ref{prop:scalar_Green_convergence_Rab}, the functions
$v_i$ are uniformly bounded in $C^0$ on every compact subset of
$\mathcal R_{\mathrm{ab}}$. The exact scalar abelian equation
\eqref{eq:exact_scalar_abelian_equations} gives
\[
    \widehat\omega_i\wedge\Theta
    =
    -*dv_i.
\]
Since $\Theta$ is parallel in both special holonomy settings, the smooth
local bounds for $\widehat\omega_i$ imply uniform bounds for every
derivative of $dv_i$. Together with the local $C^0$ bound, this yields
uniform local $C^\ell$ bounds for $v_i$ for every $\ell\geqslant0$.
After passing to a further diagonal subsequence, $v_i$ converges
smoothly on compact subsets to a smooth function
$v_{\mathrm{ab}}$.

On the other hand,
Proposition~\ref{prop:scalar_Green_convergence_Rab} gives convergence of
the full sequence
\[
    v_i\longrightarrow v_T
    \qquad
    \text{locally uniformly on }\mathcal R_{\mathrm{ab}}.
\]
Therefore
\begin{equation}\label{eq:identification_vab_vT}
    v_{\mathrm{ab}}
    =
    v_T
    =
    4\pi
    \int_X G(\,\cdot\,,y)\,d\|T\|(y).
\end{equation}
Passing to the limit in the exact scalar equation gives
\begin{equation}
\label{eq:limiting_scalar_abelian_equation}
    \widehat\omega_{\mathrm{har}}\wedge\Theta
    =
    -*dv_T.
\end{equation}
The curvature decomposition
\[
    F_{\nabla_i}
    =
    \widehat\omega_i\otimes\Psi_i
    +
    (F_{\nabla_i})^\perp
\]
and
Corollary~\ref{cor:exp_decay_transverse_Higgs_direction} show that, on
every compact subset of $\mathcal R_{\mathrm{ab}}$, all covariant
derivatives of $F_{\nabla_i}$ are uniformly bounded. Indeed, the
longitudinal forms normalized by $|\Phi_i|$ $\widehat\omega_i$ are uniformly bounded in
all local derivatives, while
$(F_{\nabla_i})^\perp$ and all its covariant derivatives decay
exponentially. Likewise,
\[
    \widetilde\Phi_i
    =
    -v_i\Psi_i,
\]
and the uniform local derivative bounds for $v_i$, together with the
higher order exponential decay of $\nabla_i\Psi_i$, give uniform local
bounds for all covariant derivatives of $\widetilde\Phi_i$ and
$\Psi_i$.

Fix $j$. Cover a neighbourhood of $K_j$ by finitely many sufficiently
small geodesic balls
\[
    U_{j,1},\ldots,U_{j,N_j}
    \Subset
    \mathcal R_{\mathrm{ab}}.
\]
Because the curvature is uniformly bounded on a slightly larger compact
set, the balls may be chosen so that $\|F_{\nabla_i}\|_{L^{n/2}(U_{j,\alpha})}$ is uniformly smaller than the Uhlenbeck threshold for all sufficiently
large $i$. Local Uhlenbeck gauge fixing therefore gives Coulomb gauges on
the $U_{j,\alpha}$ in which the connection forms are uniformly bounded
in $W^{1,n/2}$, and the preceding curvature and field estimates, together
with the monopole equations and elliptic bootstrapping, give uniform
bounds in every $C^\ell$ norm. After passing to a subsequence, the locally
gauged triples $(\nabla_i,\widetilde\Phi_i,\Psi_i)$ converge smoothly on relatively compact subsets of each
$U_{j,\alpha}$.

The Coulomb gauges obtained on different balls cannot be chosen
independently on their overlaps. We therefore apply the standard
Uhlenbeck patching construction to the transition functions of the local
gauges; see
\cite[Lemma~4.4.5]{donaldson1990geometry} and
\cite[Chapter~7]{wehrheim2004uhlenbeck}. After modifying the local gauges
by smooth gauge transformations and passing to a further subsequence, the
transition functions converge smoothly to a fixed limiting cocycle. The patched local limits consequently define a smooth principal $G$-bundle
\[
    P_{\mathrm{ab}}
    \longrightarrow
    \mathcal R_{\mathrm{ab}},
\]
together with a smooth connection $\nabla_{\mathrm{ab}}$ and smooth
sections
$\widetilde\Phi_{\mathrm{ab}},\Psi_{\mathrm{ab}}\in
\Gamma(\mathfrak g_{P_{\mathrm{ab}}})$.

Applying the patching construction inductively over the exhaustion
$(K_j)$ and then taking a diagonal subsequence gives compatible bundle
isomorphisms
\[
    u_{i,j}
    :
    P_{\mathrm{ab}}|_{K_j}
    \longrightarrow
    P|_{K_j},
    \qquad
    i\geqslant j,
\]
such that
\[
    u_{i,j}^*
    (\nabla_i,\widetilde\Phi_i,\Psi_i)
    \longrightarrow
    (\nabla_{\mathrm{ab}},
     \widetilde\Phi_{\mathrm{ab}},
     \Psi_{\mathrm{ab}})
\]
smoothly on $K_j$. The isomorphisms can be chosen compatibly under
restriction as $j$ increases. Since the approximating connections all
live on the fixed bundle $P$ and no singular set is removed inside
$\mathcal R_{\mathrm{ab}}$, the resulting limiting bundle
$P_{\mathrm{ab}}$ is isomorphic to
$P|_{\mathcal R_{\mathrm{ab}}}$. After fixing such an identification, the
preceding convergence is precisely smooth local convergence after gauge
transformations on $\mathcal R_{\mathrm{ab}}$.

The limit is now identified directly. Since $|\Psi_i|=1$, smooth
convergence gives
\[
    |\Psi_{\mathrm{ab}}|=1.
\]
Moreover, the exponential decay of $\nabla_i\Psi_i$ implies
\[
    \nabla_{\mathrm{ab}}\Psi_{\mathrm{ab}}
    =
    0.
\]
Passing to the limit in
\[
    F_{\nabla_i}
    =
    \widehat\omega_i\otimes\Psi_i
    +
    (F_{\nabla_i})^\perp
\]
gives
\[
    F_{\nabla_{\mathrm{ab}}}
    =
    \widehat\omega_{\mathrm{har}}
    \otimes\Psi_{\mathrm{ab}},
\]
while the identity
\[
    \widetilde\Phi_i
    =
    -v_i\Psi_i
\]
gives
\[
    \widetilde\Phi_{\mathrm{ab}}
    =
    -v_T\Psi_{\mathrm{ab}}.
\]
Using
\eqref{eq:limiting_scalar_abelian_equation} and
$\nabla_{\mathrm{ab}}\Psi_{\mathrm{ab}}=0$, we obtain
\[
\begin{aligned}
    F_{\nabla_{\mathrm{ab}}}\wedge\Theta
    &=
    \bigl(
        \widehat\omega_{\mathrm{har}}\wedge\Theta
    \bigr)
    \otimes\Psi_{\mathrm{ab}}  \\
    &=
    -*dv_T
    \otimes\Psi_{\mathrm{ab}}  \\
    &=
    *
    \nabla_{\mathrm{ab}}
    \bigl(
        -v_T\Psi_{\mathrm{ab}}
    \bigr)  \\
    &=
    *
    \nabla_{\mathrm{ab}}
    \widetilde\Phi_{\mathrm{ab}}.
\end{aligned}
\]
In the Calabi--Yau case, the smooth convergence and $F_{\nabla_i}\wedge\omega^2=0$ also give $F_{\nabla_{\mathrm{ab}}}\wedge\omega^2=0$. Thus $(\nabla_{\mathrm{ab}}, \widetilde\Phi_{\mathrm{ab}})$ is a smooth $\Theta$-monopole. Finally, $\nabla_{\mathrm{ab}}\Psi_{\mathrm{ab}}=0$ shows that
$\nabla_{\mathrm{ab}}$ preserves the real line generated by
$\Psi_{\mathrm{ab}}$. Hence its holonomy reduces to the stabilizer circle
\[
    T_G
    =
    \begin{cases}
        \mathrm U(1),&G=\mathrm{SU}(2),\\
        \mathrm{SO}(2),&G=\mathrm{SO}(3),
    \end{cases}
\]
and the limiting $\Theta$-monopole is abelian.
\end{proof}

\begin{corollary}[The scalar part of the limiting abelian monopole]
\label{cor:scalar_part_abelian_limit}
In the setting of
Corollary~\ref{cor:smooth_compactness_L2_bounded_longitudinal}, the
limiting translated Higgs field is
\[
    \widetilde\Phi_{\mathrm{ab}}
    =
    -v_T\Psi_{\mathrm{ab}},
    \qquad
    \nabla_{\mathrm{ab}}\Psi_{\mathrm{ab}}=0,
\]
where
\[
    v_T
    =
    4\pi
    \int_X G(\,\cdot\,,y)\,d\|T\|(y)
    >0
    \qquad
    \text{on }\mathcal R_{\mathrm{ab}}.
\]
Moreover,
\[
    d^*dv_T=0,
    \qquad
    \widehat\omega_{\mathrm{har}}\wedge\Theta
    =
    -*dv_T
\]
on $\mathcal R_{\mathrm{ab}}$.
\end{corollary}

\begin{proof}
The identity
$\widetilde\Phi_{\mathrm{ab}}=-v_T\Psi_{\mathrm{ab}}$ follows from
\eqref{eq:identification_vab_vT}. The remaining assertions follow from
Proposition~\ref{prop:scalar_Green_convergence_Rab} and
\eqref{eq:limiting_scalar_abelian_equation}.
\end{proof}

\begin{lemma}[Type of $L^2$-harmonic $2$-forms under special holonomy]
\label{lem:L2_harmonic_two_forms_Xi_ASD}
In either special holonomy setting of the {\mainsettingref},
\[
    \mathscr H_{(2)}^2(X)
    =
    \left\{
        \sigma\in L^2\Omega^2(X):
        d\sigma=0,\quad
        *(\sigma\wedge\Xi)=-\sigma
    \right\}.
\]
Thus every $\sigma\in\mathscr H_{(2)}^2(X)$ lies in $\Lambda^2_{14}$ in the $\mathrm G_2$ case, while in the Calabi--Yau case it is a primitive real $(1,1)$-form: $\sigma\in\Lambda^{1,1}_0$ and $\Lambda\sigma=0$. In particular, in both cases,
\[
\sigma\wedge\Theta=0.
\]
\end{lemma}

\begin{proof}
The special holonomy decomposition of $\Lambda^2T^*X$ is parallel, so
the Hodge Laplacian preserves each of its irreducible summands.

In the $\mathrm G_2$ case,
\[
    \Lambda^2
    =
    \Lambda^2_7\oplus\Lambda^2_{14},
\]
and the parallel endomorphism
\[
    \alpha
    \longmapsto
    *(\alpha\wedge\varphi)
\]
has eigenvalue $2$ on $\Lambda^2_7$ and eigenvalue $-1$ on
$\Lambda^2_{14}$. The parallel identification
$\Lambda^1\simeq\Lambda^2_7$, given by contraction with $\varphi$,
intertwines the Hodge Laplacians. Hence the $\Lambda^2_7$-component of
an $L^2$-harmonic two-form corresponds to an $L^2$-harmonic one-form.
Since $(X,g)$ is complete and Ricci-flat, the Bochner formula and the
standard cutoff argument imply that every $L^2$-harmonic one-form is
parallel. As an AC manifold has infinite volume, a nonzero parallel
form cannot lie in $L^2$. The $\Lambda^2_7$-component therefore
vanishes, and every element of $\mathscr H_{(2)}^2(X)$ lies in
$\Lambda^2_{14}$. Thus
\[
    *(\sigma\wedge\varphi)=-\sigma,
    \qquad
    \sigma\wedge\psi=0.
\]
In the Calabi--Yau case, the parallel real decomposition is
\[
    \Lambda^2
    =
    \mathbb R\omega
    \oplus
    \Lambda^{1,1}_0
    \oplus
    \Lambda^2_6,
\]
where $\Lambda^2_6$ is the real $(2,0)+(0,2)$ summand. The endomorphism
\[
    \alpha
    \longmapsto
    *(\alpha\wedge\omega)
\]
has eigenvalues $2$, $-1$, and $1$ on these three summands,
respectively. The $\mathbb R\omega$-component of an $L^2$-harmonic
two-form is $f\omega$ for an $L^2$-harmonic function $f$, and hence
vanishes by the standard cutoff argument. The parallel identification
$\Lambda^1\simeq\Lambda^2_6$, given by contraction with
$\operatorname{Re}\Omega$, shows as above that the
$\Lambda^2_6$-component also vanishes. Therefore every element of
$\mathscr H_{(2)}^2(X)$ lies in $\Lambda^{1,1}_0$, and hence
\[
    *(\sigma\wedge\omega)=-\sigma,
    \qquad
    \sigma\wedge\operatorname{Re}\Omega=0.
\]
This proves that every $L^2$-harmonic two-form satisfies
$*(\sigma\wedge\Xi)=-\sigma$ and $\sigma\wedge\Theta=0$. Conversely, if
$\sigma\in L^2\Omega^2(X)$ is closed and
$*(\sigma\wedge\Xi)=-\sigma$, then
\[
    *\sigma=-\sigma\wedge\Xi.
\]
Since $d\Xi=0$, it follows that $d*\sigma=0$, and therefore
$d^*\sigma=0$. Thus $\sigma$ is $L^2$-harmonic.
\end{proof}

\begin{remark}[Global $L^2$-harmonic corrections and local harmonic limits]
\label{rmk:global_L2_harmonic_vs_local_harmonic_limit}
Lemma~\ref{lem:L2_harmonic_two_forms_Xi_ASD} applies in particular to
Li's corrections $\sigma_i$, to their limits in the bounded
alternative, and to the normalized limit in the unbounded alternative.
Thus all these forms satisfy
\[
    *(\sigma\wedge\Xi)=-\sigma,
    \qquad
    \sigma\wedge\Theta=0.
\]
No assumption that the holonomy be exactly $\mathrm G_2$ or
$\mathrm{SU}(3)$ is needed.

The same conclusion does not follow for the harmonic form
$\widehat\omega_{\mathrm{har}}$ arising as a smooth local limit of the uncorrected longitudinal curvatures on $\mathcal R_{\mathrm{ab}}$. This form is defined only on the possibly incomplete open set $\mathcal R_{\mathrm{ab}}$ and is not known to extend to an element of $\mathscr H_{(2)}^2(X)$. In the $\mathrm G_2$ case it may therefore have a $\Lambda^2_7$-component. In the Calabi--Yau case it is primitive,
but may have a nonzero real $(2,0)+(0,2)$ component. Consequently, the
limiting equation
\[
    \widehat\omega_{\mathrm{har}}\wedge\Theta
    =
    -*dv_T
\]
does not by itself imply $dv_T=0$.

If, however, $\widehat\omega_{\mathrm{har}}$ extends to a global
$L^2$-harmonic two-form on $X$, then
Lemma~\ref{lem:L2_harmonic_two_forms_Xi_ASD} gives
$\widehat\omega_{\mathrm{har}}\wedge\Theta=0$. Hence
\[
    dv_T=0,
    \qquad
    \nabla_{\mathrm{ab}}\widetilde\Phi_{\mathrm{ab}}=0,
    \qquad
    F_{\nabla_{\mathrm{ab}}}\wedge\Theta=0.
\]
Equivalently, the limiting Higgs field is parallel and the limiting connection is an abelian $\Xi$-anti-self-dual instanton on $\mathcal R_{\mathrm{ab}}$, that is, an abelian $\mathrm G_2$-instanton in the $\mathrm G_2$ case and an abelian Hermitian--Yang--Mills connection in the Calabi--Yau case.
\end{remark}

We now combine these local Hodge estimates with Li's corrected
compactness theorem.

\begin{proposition}[Smooth convergence of Li's corrected longitudinal
forms]
\label{prop:smooth_convergence_Li_corrected_longitudinal}
For each $i$, fix an $L^2$-harmonic correction $\sigma_i\in \mathscr{H}_{(2)}^2(X)$ furnished by \cite[Theorem~1.12]{li2025large}, with the $\mathrm{SO}(3)$ extension
described in Remark~\ref{rmk:Li_SO3_variant}, and set
\[
    \widehat\omega_i^{\mathrm{corr}}
    :=
    \frac{1}{4\pi}\widehat\omega_i+\sigma_i.
\]
After passing to the subsequence appearing in Li's theorem,
\[
    \widehat\omega_i^{\mathrm{corr}}
    \longrightarrow
    \widetilde F_\infty
\]
smoothly on compact subsets of $\mathcal R_{\mathrm{ab}}$, where
$\widetilde F_\infty=F_\infty/(2\pi)$ is Li's normalized limiting
abelian curvature form. In the Calabi--Yau case, one also has
$\Lambda\widetilde F_\infty=0$ in the distributional sense on $X$.
\end{proposition}

\begin{proof}
By \cite[Theorem~1.12]{li2025large},
Remark~\ref{rmk:Li_SO3_variant} in the $\mathrm{SO}(3)$ case, and the
normalization comparison in
Remark~\ref{rem:identification_with_Li_current}, Li's globally defined
corrected forms satisfy
\[
    \frac{1}{4\pi}\omega_i+\sigma_i
    \longrightarrow
    \widetilde F_\infty
\]
strongly in $L^1_{\mathrm{loc}}(X)$. On $\mathcal R_{\mathrm{ab}}$,
Corollary~\ref{cor:subcritical_longitudinal_curvature_Rab} gives
\[
    \widehat\omega_i-\omega_i
    \longrightarrow0
    \qquad
    \text{strongly in }L^2_{\mathrm{loc}}(\mathcal R_{\mathrm{ab}}),
\]
and hence strongly in $L^1_{\mathrm{loc}}$. Therefore
\[
    \widehat\omega_i^{\mathrm{corr}}
    =
    \frac{1}{4\pi}\widehat\omega_i+\sigma_i
    \longrightarrow
    \widetilde F_\infty
\]
strongly in $L^1_{\mathrm{loc}}(\mathcal R_{\mathrm{ab}})$.

In the Calabi--Yau case, the globally defined form
\[
    \omega_i=m_i^{-1}\langle F_{\nabla_i},\Phi_i\rangle
\]
satisfies
\[
    \Lambda\omega_i
    =m_i^{-1}\langle\Lambda F_{\nabla_i},\Phi_i\rangle
    =0.
\]
Moreover, Lemma~\ref{lem:L2_harmonic_two_forms_Xi_ASD} gives
$\Lambda\sigma_i=0$. Li's strong $L^1_{\mathrm{loc}}(X)$ convergence
of the globally defined corrected forms
\[
    \frac{1}{4\pi}\omega_i+\sigma_i
    \longrightarrow
    \widetilde F_\infty
\]
therefore implies
$\Lambda\widetilde F_\infty=0$ distributionally on $X$.

Since each $\sigma_i$ is harmonic,
\[
    d\widehat\omega_i^{\mathrm{corr}}
    =
    \frac{1}{4\pi}d\widehat\omega_i,
    \qquad
    d^*\widehat\omega_i^{\mathrm{corr}}
    =
    \frac{1}{4\pi}d^*\widehat\omega_i.
\]
Corollary~\ref{cor:Hodge_longitudinal_component_large_mass} therefore
gives, for every $\ell\geqslant0$,
\[
    d\widehat\omega_i^{\mathrm{corr}}
    \longrightarrow0,
    \qquad
    d^*\widehat\omega_i^{\mathrm{corr}}
    \longrightarrow0
\]
locally in $C^\ell$ on $\mathcal R_{\mathrm{ab}}$.

Because $\mathcal S\subset\mathcal C$, one has
$\mathcal R_{\mathrm{ab}}\subset X\setminus \mathcal S$, and Li's limiting form is
smooth there. Passing to the distributional limit in the preceding Hodge
equations gives
\[
    d\widetilde F_\infty=0,
    \qquad
    d^*\widetilde F_\infty=0
    \qquad
    \text{on }\mathcal R_{\mathrm{ab}}.
\]
Let $K'\Subset K\Subset\mathcal R_{\mathrm{ab}}$ and set
\[
    \eta_i
    :=
    \widehat\omega_i^{\mathrm{corr}}-\widetilde F_\infty.
\]
Choose open sets $U$ and $V$ such that
\[
    K'\subset U\Subset V\Subset K.
\]
Since
\[
    \Delta \eta_i
    =
    d d^*\widehat\omega_i^{\mathrm{corr}}
    +
    d^*d\widehat\omega_i^{\mathrm{corr}},
\]
the preceding higher order Hodge estimates imply
$\Delta \eta_i\to0$ locally in every $C^\ell$ norm. Standard interior estimates for the Hodge Laplacian give, for every
$\ell\geqslant0$, an integer $\ell'$ and a constant
$C_{K',U,V,\ell}$ such that
\[
    \|\eta_i\|_{C^\ell(K')}
    \leqslant
    C_{K',U,V,\ell}
    \left(
        \|\eta_i\|_{L^1(V)}
        +
        \|\Delta \eta_i\|_{C^{\ell'}(V)}
    \right).
\]
Both terms on the right tend to zero. Hence $\eta_i\to0$ smoothly on $K'$.
Since $K'$ was arbitrary, the convergence is smooth on compact subsets of
$\mathcal R_{\mathrm{ab}}$.
\end{proof}

\begin{proposition}[Subcritical growth of Li's harmonic corrections and longitudinal curvature]
\label{prop:subcritical_harmonic_correction_growth}
For the corrections $\sigma_i\in\mathscr H_{(2)}^2(X)$ fixed in
Proposition~\ref{prop:smooth_convergence_Li_corrected_longitudinal}, one
has
\begin{equation}
\label{eq:subcritical_harmonic_correction_growth}
    \|\sigma_i\|_{L^2(X)}
    =
    o(m_i^{1/2}).
\end{equation}
Moreover, for every compact set
$K\Subset\mathcal R_{\mathrm{ab}}$ and every integer $\ell\geqslant0$,
\begin{equation}
\label{eq:subcritical_longitudinal_curvature_growth}
    m_i^{-1/2}
    \|\widehat\omega_i\|_{C^\ell(K)}
    \longrightarrow0,
\end{equation}
and
\begin{equation}
\label{eq:subcritical_full_curvature_growth}
    m_i^{-1/2}
    \left(
        \|F_{\nabla_i}\|_{C^\ell(K)}
        +
        \|\nabla_i\Phi_i\|_{C^\ell(K)}
    \right)
    \longrightarrow0.
\end{equation}
\end{proposition}

\begin{proof}
Since $\mathcal S$ is compact and has finite $\mathcal H^{n-3}$-measure, it
has empty interior. Choose a compact domain
\[
    K_0\Subset X\setminus\mathcal S
\]
with nonempty interior. The restriction seminorm
\[
    \sigma
    \longmapsto
    \|\sigma\|_{L^1(K_0)}
\]
is a norm on the finite-dimensional space $\mathscr H_{(2)}^2(X)$.
Indeed, if it vanishes, then the smooth harmonic form $\sigma$ vanishes
on $\operatorname{int}K_0$, and unique continuation gives
$\sigma\equiv0$ on $X$. Hence there is a constant $C_{K_0}$ such that
\begin{equation}
\label{eq:harmonic_restriction_norm}
    \|\sigma\|_{L^2(X)}
    \leqslant
    C_{K_0}\|\sigma\|_{L^1(K_0)}
    \qquad
    \text{for every }\sigma\in\mathscr H_{(2)}^2(X).
\end{equation}
Choose $\chi\in C_c^\infty(X\setminus\mathcal S)$ with
$\chi\equiv1$ on $K_0$. Since
$\mu_i\stackrel{*}{\rightharpoonup}8\pi\|T\|$ and
$\operatorname{spt}\|T\|=\mathcal S$,
\[
    m_i^{-1}
    \int_{K_0}|F_{\nabla_i}|^2
    \leqslant
    \int_X\chi\,d\mu_i
    \longrightarrow0.
\]
Let
\[
    \omega_i:=m_i^{-1}\langle F_{\nabla_i},\Phi_i\rangle.
\]
Since $|\Phi_i|\leqslant m_i$, one has $|\omega_i|\leqslant |F_{\nabla_i}|$.
H\"older's inequality therefore gives
\begin{equation}
\label{eq:global_longitudinal_L1_subcritical}
    \|\omega_i\|_{L^1(K_0)}
    =
    o(m_i^{1/2}).
\end{equation}
By Li's corrected compactness theorem, with the normalization comparison
used in Proposition~\ref{prop:smooth_convergence_Li_corrected_longitudinal},
the globally defined corrected forms
\[
    \frac{1}{4\pi}\omega_i+\sigma_i
\]
converge strongly in $L^1_{\mathrm{loc}}(X)$. In particular, they are
uniformly bounded in $L^1(K_0)$. Hence
\[
    \sigma_i
    =
    \left(\frac{1}{4\pi}\omega_i+\sigma_i\right)
    -
    \frac{1}{4\pi}\omega_i,
\]
and \eqref{eq:harmonic_restriction_norm} together with
\eqref{eq:global_longitudinal_L1_subcritical} proves
\eqref{eq:subcritical_harmonic_correction_growth}.

Let now $K\Subset\mathcal R_{\mathrm{ab}}$. By finite-dimensionality and
elliptic regularity, for every $\ell\geqslant0$ there is a constant
$C_{K,\ell}$ such that
\[
    \|\sigma\|_{C^\ell(K)}
    \leqslant
    C_{K,\ell}\|\sigma\|_{L^2(X)}
    \qquad
    \text{for every }\sigma\in\mathscr H_{(2)}^2(X).
\]
Proposition~\ref{prop:smooth_convergence_Li_corrected_longitudinal}
gives a uniform $C^\ell(K)$ bound for
$\widehat\omega_i^{\mathrm{corr}}$. Therefore
\[
    \widehat\omega_i
    =
    4\pi
    \bigl(
        \widehat\omega_i^{\mathrm{corr}}-\sigma_i
    \bigr),
\]
and \eqref{eq:subcritical_harmonic_correction_growth} gives
\eqref{eq:subcritical_longitudinal_curvature_growth}. Notice that the
case $\ell=0$ is also the direct consequence recorded in
Corollary~\ref{cor:subcritical_longitudinal_curvature_Rab}.

For the full curvature, differentiate the decomposition
\[
    F_{\nabla_i}
    =
    \widehat\omega_i\otimes\Psi_i
    +
    (F_{\nabla_i})^\perp.
\]
The higher order exponential estimates of
Lemma~\ref{lem:Li_dichotomy_away_Cinfty}, together with
\eqref{eq:subcritical_longitudinal_curvature_growth} and the product
rule, give
\[
    m_i^{-1/2}
    \|F_{\nabla_i}\|_{C^\ell(K)}
    \longrightarrow0.
\]
Since $\Theta$ is parallel and
$F_{\nabla_i}\wedge\Theta=*\nabla_i\Phi_i$, the same estimate holds for
$\nabla_i\Phi_i$. This proves
\eqref{eq:subcritical_full_curvature_growth}.
\end{proof}

\begin{corollary}[Alternatives for Li's harmonic corrections]
\label{cor:harmonic_correction_alternatives}
Exactly one of the following alternatives holds.

\begin{enumerate}
\item The sequence is uniformly bounded in $L^2(X)$:
\[
    \sup_i\|\sigma_i\|_{L^2(X)}<\infty.
\]
Then, after passing to a subsequence, there exists
$\sigma_\infty\in \mathscr{H}_{(2)}^2(X)$ such that
$\sigma_i\to\sigma_\infty$ smoothly on compact subsets of $X$, and
\[
    \widehat\omega_i
    \longrightarrow
    4\pi\bigl(\widetilde F_\infty-\sigma_\infty\bigr)
\]
smoothly on compact subsets of $\mathcal R_{\mathrm{ab}}$. Consequently,
after applying gauge transformations, the translated configurations
\[
    (\nabla_i,\widetilde\Phi_i),
    \qquad
    \widetilde\Phi_i:=\Phi_i-m_i\Psi_i,
\]
converge smoothly on compact subsets of $\mathcal R_{\mathrm{ab}}$ to a
smooth abelian $\Theta$-monopole.

\item The $L^2(X)$ norms have infinite limsup:
\[
    \limsup_{i\to\infty}
    \|\sigma_i\|_{L^2(X)}
    =
    +\infty.
\]
After passing to a subsequence and relabelling, write
\[
    b_i:=\|\sigma_i\|_{L^2(X)},
\]
so that $b_i\to+\infty$. By Proposition~\ref{prop:subcritical_harmonic_correction_growth},
\[
    b_i=o(m_i^{1/2}).
\]
There exists $\sigma_\infty\in \mathscr{H}_{(2)}^2(X)$ with $\|\sigma_\infty\|_{L^2(X)}=1$ such that
\[
    \frac{\sigma_i}{b_i}
    \longrightarrow
    \sigma_\infty
\]
smoothly on compact subsets of $X$, and
\begin{equation}
\label{eq:global_harmonic_correction_limit}
    \frac{1}{4\pi b_i}\widehat\omega_i
    \longrightarrow
    -\sigma_\infty
\end{equation}
smoothly on compact subsets of $\mathcal R_{\mathrm{ab}}$. Moreover, for
every compact domain $K\Subset\mathcal R_{\mathrm{ab}}$,
\begin{equation}
\label{eq:local_amplitude_harmonic_correction}
    \frac{\|\widehat\omega_i\|_{L^2(K)}}{4\pi b_i}
    \longrightarrow
    \|\sigma_\infty\|_{L^2(K)}>0,
\end{equation}
and therefore
\begin{equation}
\label{eq:local_normalized_harmonic_correction}
    \frac{\widehat\omega_i}
    {\|\widehat\omega_i\|_{L^2(K)}}
    \longrightarrow
    -\frac{\sigma_\infty}
    {\|\sigma_\infty\|_{L^2(K)}}
\end{equation}
smoothly on compact subsets of $\operatorname{int}K$.
\end{enumerate}

In particular, the first alternative always holds if
$\mathscr{H}_{(2)}^2(X)=0$. By the $L^2$ Hodge theorem used in
\cite[Theorem~1.12]{li2025large}, in the present AC dimensions one has
$\mathscr{H}_{(2)}^2(X)\simeq H_c^2(X;\mathbb R)$, so it is enough that
$H_c^2(X;\mathbb R)=0$. If, in addition,
$\mathcal O\subset\mathcal S$, then
$\mathcal R_{\mathrm{ab}}=X\setminus \mathcal S$, and the convergence in the first
alternative is smooth on compact subsets of $X\setminus \mathcal S$.
\end{corollary}

\begin{proof}
The space $\mathscr{H}_{(2)}^2(X)$ is finite-dimensional. If
\[
    \sup_i\|\sigma_i\|_{L^2(X)}<\infty,
\]
then, after passing to a subsequence, $\sigma_i$ converges in $L^2(X)$ to
some $\sigma_\infty\in\mathscr{H}_{(2)}^2(X)$. Equivalence of norms on
this finite-dimensional space, together with elliptic regularity, gives
smooth convergence on compact subsets of $X$. The first conclusion then
follows from
Proposition~\ref{prop:smooth_convergence_Li_corrected_longitudinal}, and
the smooth convergence of the translated configurations follows from
Corollary~\ref{cor:smooth_compactness_L2_bounded_longitudinal}.

If instead
\[
    \limsup_{i\to\infty}
    \|\sigma_i\|_{L^2(X)}
    =
    +\infty,
\]
choose a subsequence, relabelled by $i$, such that
\[
    b_i:=\|\sigma_i\|_{L^2(X)}
    \longrightarrow
    +\infty.
\]
By compactness of the unit sphere
in this finite-dimensional space, there exists $\sigma_\infty\in \mathscr{H}_{(2)}^2(X)$ of unit $L^2$ norm
such that $b_i^{-1}\sigma_i\to\sigma_\infty$ smoothly on compact subsets
of $X$. Since
\[
    \frac{1}{4\pi b_i}\widehat\omega_i
    =
    \frac{1}{b_i}\widehat\omega_i^{\mathrm{corr}}
    -
    \frac{1}{b_i}\sigma_i,
\]
Proposition~\ref{prop:smooth_convergence_Li_corrected_longitudinal} gives
\eqref{eq:global_harmonic_correction_limit}.

Let $K\Subset\mathcal R_{\mathrm{ab}}$ be a compact domain. If
$\|\sigma_\infty\|_{L^2(K)}=0$, then $\sigma_\infty$ vanishes on the
nonempty open set $\operatorname{int}K$. Unique continuation for harmonic
forms on the connected manifold $X$ would imply
$\sigma_\infty\equiv0$, contradicting its unit $L^2$ norm. Hence
$\|\sigma_\infty\|_{L^2(K)}>0$. Taking $L^2(K)$ norms in
\eqref{eq:global_harmonic_correction_limit} proves
\eqref{eq:local_amplitude_harmonic_correction}, and division by this
limit gives
\eqref{eq:local_normalized_harmonic_correction}.
\end{proof}

\begin{remark}[Comparison with Li's singular abelian limit]
\label{rmk:comparison_with_Li_singular_limit}
Assume that the first alternative of
Corollary~\ref{cor:harmonic_correction_alternatives} holds. Let
$(\nabla_{\mathrm{ab}},\widetilde\Phi_{\mathrm{ab}})$ be the smooth
abelian $\Theta$-monopole obtained as the limit of the translated
configurations, and write
\[
    \widetilde\Phi_{\mathrm{ab}}
    =
    -v_T\Psi_{\mathrm{ab}},
    \qquad
    \nabla_{\mathrm{ab}}\Psi_{\mathrm{ab}}=0.
\]
If $\sigma_i\to\sigma_\infty$ is the limiting harmonic correction, then
Proposition~\ref{prop:smooth_convergence_Li_corrected_longitudinal} and
the smooth convergence of the translated configurations give
\begin{equation}
\label{eq:comparison_smooth_and_Li_curvatures}
    \widetilde F_\infty
    =
    \frac{1}{4\pi}
    \left\langle
        F_{\nabla_{\mathrm{ab}}},
        \Psi_{\mathrm{ab}}
    \right\rangle
    +
    \sigma_\infty
    \qquad
    \text{on }\mathcal R_{\mathrm{ab}}.
\end{equation}
Equivalently, since
$F_\infty=2\pi\widetilde F_\infty$,
\begin{equation}
\label{eq:comparison_smooth_and_Li_curvatures_unnormalized}
    F_\infty
    =
    \frac12
    \left\langle
        F_{\nabla_{\mathrm{ab}}},
        \Psi_{\mathrm{ab}}
    \right\rangle
    +
    2\pi\sigma_\infty
    \qquad
    \text{on }\mathcal R_{\mathrm{ab}}.
\end{equation}
The scalar component is determined by
Proposition~\ref{prop:scalar_Green_convergence_Rab} and
\eqref{eq:identification_vab_vT}:
\[
    v_i
    \longrightarrow
    v_T
    =
    4\pi
    \int_X G(\,\cdot\,,y)\,d\|T\|(y)
\]
locally uniformly on $\mathcal R_{\mathrm{ab}}$, while
Lemma~\ref{lem:exact_scalar_abelian_equations} identifies
$\Phi_\infty^{\mathrm{Li}}$ with the negative of one half of this
limit. Thus Li's limiting scalar Higgs field is represented on
$\mathcal R_{\mathrm{ab}}$ by
\[
    -\frac12v_T
    =
    -2\pi
    \int_X G(\,\cdot\,,y)\,d\|T\|(y)
\]
with the normalization of this paper. Li's scalar variable is
$(|\Phi_i|^2-m_i^2)/(2m_i)$, while the nonnegative quadratic Higgs defect
used above is $(m_i^2-|\Phi_i|^2)/(2m_i)$.

Consequently, on $\mathcal R_{\mathrm{ab}}$, Li's singular abelian
monopole is represented in real-valued abelian variables by
\[
    \left(
        \frac12
        \left\langle
            F_{\nabla_{\mathrm{ab}}},
            \Psi_{\mathrm{ab}}
        \right\rangle
        +
        2\pi\sigma_\infty,
        -\frac12v_T
    \right).
\]
Since $\sigma_\infty\wedge\Theta=0$, these variables satisfy the
abelian monopole equation
\[
    F_\infty\wedge\Theta
    =
    -\frac12*dv_T
    =
    *d\left(-\frac12v_T\right).
\]
By
Lemma~\ref{lem:L2_harmonic_two_forms_Xi_ASD}, the limiting correction
satisfies
\[
    *(\sigma_\infty\wedge\Xi)=-\sigma_\infty,
    \qquad
    \sigma_\infty\wedge\Theta=0
\]
on $X$. In particular, if $H_c^2(X;\mathbb R)=0$, then
$\sigma_\infty=0$, and Li's limit is precisely the scalar abelian
reduction of the smooth limit, after the normalization conversion in
\eqref{eq:comparison_smooth_and_Li_curvatures_unnormalized}.
\end{remark}

\begin{remark}[Cohomology classes of the corrections in the bounded alternative]
\label{rmk:lattice_corrections_bounded_alternative}
Under the identification
$\mathscr H_{(2)}^2(X)\simeq H_c^2(X;\mathbb R)$, Li's corrections may be
chosen so that their cohomology classes lie in the integral lattice when
$G=\mathrm{SU}(2)$, with the corresponding half-integral lattice when
$G=\mathrm{SO}(3)$. For such a choice, the bounded alternative has the
following additional consequence. A bounded subset of this lattice is
finite, so, after passing to a subsequence, the classes $[\sigma_i]$ are
constant. Since each class has a unique $L^2$-harmonic representative,
there is a form $\sigma_\infty\in\mathscr H_{(2)}^2(X)$ such that $\sigma_i=\sigma_\infty$ for all sufficiently large $i$.
\end{remark}

\begin{remark}[The cohomogeneity-one examples]
\label{rmk:homogeneous_Hc2_vanishing}
Each of the three ambient manifolds in
Theorem~\ref{thm:intro_homogeneous_examples} is the total space of an
oriented real vector bundle of rank three over a compact base:
$\Lambda_-^2(S^4)$, $\Lambda_-^2(\mathbb{CP}^2)$, and $T^*S^3$.
If $X\to B$ denotes any of these rank-three bundles, the Thom
isomorphism gives
$H_c^k(X;\mathbb R)\simeq H^{k-3}(B;\mathbb R)$, and hence
$H_c^2(X;\mathbb R)=0$. Moreover,
Corollary~\ref{cor:homogeneous_examples_dream_case} gives
$\mathcal O=\varnothing$ for these families. Therefore
Corollary~\ref{cor:harmonic_correction_alternatives} yields smooth
convergence, after gauge, of the translated monopoles on compact subsets
of the complement of the calibrated zero section.

For the same cohomogeneity-one families, smooth convergence on compact
subsets of the complement of the zero section was proved directly by
the second author from the explicit solutions; see
\cites{oliveira2014monopoles,oliveira2016calabi}. Thus the argument above,
which combines the vanishing of $H_c^2(X;\mathbb R)$ with the absence of
the obstruction established in Section~\ref{sec:homogeneous_examples},
gives a second proof of this smooth local convergence.
\end{remark}

Proposition~\ref{prop:smooth_convergence_Li_corrected_longitudinal} and
Corollary~\ref{cor:harmonic_correction_alternatives} reduce the remaining
analysis to the second alternative of that corollary. After selecting
and relabelling the divergent subsequence, we estimate the regularity
scale in this case.

\subsection{Unbounded harmonic corrections}
\label{subsec:unbounded_longitudinal_alternative}

Assume that the second alternative of
Corollary~\ref{cor:harmonic_correction_alternatives} occurs, and write
\[
    b_i:=\|\sigma_i\|_{L^2(X)}\longrightarrow+\infty,
    \qquad
    \frac{\sigma_i}{b_i}\longrightarrow\sigma_\infty,
\]
where $\sigma_\infty\in \mathscr{H}_{(2)}^2(X)$ has unit $L^2$ norm. Fix
compact domains
\[
    K'\Subset K\Subset\mathcal R_{\mathrm{ab}}
\]
and set
\[
    \mathfrak a_i
    :=
    \|\widehat\omega_i\|_{L^2(K)}.
\]
By \eqref{eq:local_amplitude_harmonic_correction},
$\mathfrak a_i\to+\infty$, and
\[
    \frac{\mathfrak a_i}{4\pi b_i}
    \longrightarrow
    \|\sigma_\infty\|_{L^2(K)}>0.
\]
Independently of this comparison with the harmonic corrections,
Corollary~\ref{cor:subcritical_longitudinal_curvature_Rab} gives
\begin{equation}
\label{eq:local_longitudinal_amplitude_subcritical}
    \mathfrak a_i=o(m_i^{1/2}).
\end{equation}
The same conclusion is also consistent with
$b_i=o(m_i^{1/2})$ and the preceding limit relation.
Consequently,
\[
    \overline\omega_i
    :=
    \mathfrak a_i^{-1}\widehat\omega_i
    \longrightarrow
    \overline\omega_{\mathrm{har}}
    :=
    -\frac{\sigma_\infty}
    {\|\sigma_\infty\|_{L^2(K)}}
\]
smoothly on $K'$. In particular,
$\overline\omega_{\mathrm{har}}$ cannot vanish identically on any
nonempty open subset. The next proposition applies on compact subsets
where this limiting harmonic form is bounded away from zero. Near its
zero set, no estimate at shrinking scales is inferred without the additional
quantitative hypothesis introduced later.

\begin{proposition}[Regularity scale where the harmonic limit is nonzero]
\label{prop:regularity_scale_nonzero_harmonic_profile}
Let $U'\Subset U\Subset K'$ be open sets such that
\[
    |\overline\omega_{\mathrm{har}}|\geqslant c_U>0
    \qquad\text{on }U.
\]
Then there exist constants $0<c<C<\infty$ such that
\[
    c\mathfrak{a}_i^{-1/2}
    \leqslant
    \mathfrak{r}_i(x)
    \leqslant
    C\mathfrak{a}_i^{-1/2}
\]
for every $x\in U'$ and all sufficiently large $i$.
\end{proposition}

\begin{proof}
Since
\[
    \overline\omega_i\to\overline\omega_{\mathrm{har}}
\]
smoothly on $K'$ and
\[
    |\overline\omega_{\mathrm{har}}|\geqslant c_U
    \qquad\text{on }U,
\]
after decreasing $c_U$ and increasing a constant $C_U$ if necessary, we
have
\[
    c_U
    \leqslant
    |\overline\omega_i|
    \leqslant
    C_U
\]
on $U$, for all sufficiently large $i$. Hence
\[
    c_U\mathfrak{a}_i
    \leqslant
    |\widehat\omega_i|
    \leqslant
    C_U\mathfrak{a}_i
\]
on $U$.

Since
\[
    F_{\nabla_i}^{\parallel}
    =
    \widehat\omega_i\otimes\Psi_i,
\]
we have
\[
    |F_{\nabla_i}^{\parallel}|
    =
    |\widehat\omega_i|.
\]
The transverse curvature components decay exponentially on compact
subsets of $\mathcal R_{\mathrm{ab}}$ by
Lemma~\ref{lem:Li_dichotomy_away_Cinfty}. Therefore, on $U$ and for all
sufficiently large $i$,
\[
    |F_{\nabla_i}|
    \leqslant
    C_U\mathfrak{a}_i.
\]
The $\Theta$-monopole equation gives
\[
    |\nabla_i\Phi_i|
    =
    |F_{\nabla_i}\wedge\Theta|
    \leqslant
    C|F_{\nabla_i}|,
\]
and hence
\[
    e_i
    \leqslant
    C_U\mathfrak{a}_i^2
    \qquad
    \text{on }U.
\]
Let
\[
    d_U:=\operatorname{dist}(U',X\setminus U)>0.
\]
Set
\[
    \rho_U
    :=
    \min\{r_0,d_U\}.
\]
All radii used below will be chosen smaller than $\rho_U$, so that the
corresponding balls are contained in $U$ and lie within the range of the
regularity scale definition and the local volume estimates.

For every $x\in U'$ and every $0<r<\rho_U$, the ball $B_r(x)$ is contained in
$U$. Using the local volume upper bound on $U$, we obtain
\[
    r^{4-n}
    \int_{B_r(x)}
    e_i
    \leqslant
    C_U\mathfrak{a}_i^2r^{4-n}\operatorname{vol}(B_r(x))
    \leqslant
    C_U\mathfrak{a}_i^2r^4.
\]
Choose $c>0$ so small that
\[
    C_Uc^4<\varepsilon_0.
\]
Since $\mathfrak{a}_i\to+\infty$, for all sufficiently large $i$ one has
\[
    c\mathfrak{a}_i^{-1/2}<\rho_U.
\]
Thus, for every $x\in U'$ and every
\[
    0<s\leqslant c\mathfrak{a}_i^{-1/2},
\]
we have
\[
    s^{4-n}
    \int_{B_s(x)}
    e_i
    \leqslant
    C_U\mathfrak{a}_i^2s^4
    \leqslant
    C_Uc^4
    <
    \varepsilon_0.
\]
By the definition of $\mathfrak{r}_i(x)$,
\[
    \mathfrak{r}_i(x)\geqslant c\mathfrak{a}_i^{-1/2}.
\]
For the opposite inequality, we use the lower bound for the longitudinal
curvature. Since
\[
    |F_{\nabla_i}^{\parallel}|
    =
    |\widehat\omega_i|
    \geqslant
    c_U\mathfrak{a}_i
\]
on $U$, we have, after changing $c_U>0$,
\[
    e_i
    \geqslant
    |F_{\nabla_i}|^2
    \geqslant
    |F_{\nabla_i}^{\parallel}|^2
    \geqslant
    c_U\mathfrak{a}_i^2
    \qquad
    \text{on }U.
\]
Using the local volume lower bound on $U$, for every
\[
    0<r<\rho_U
\]
we obtain
\[
    r^{4-n}
    \int_{B_r(x)}
    e_i
    \geqslant
    c_U\mathfrak{a}_i^2r^{4-n}\operatorname{vol}(B_r(x))
    \geqslant
    c_U\mathfrak{a}_i^2r^4,
\]
after again changing $c_U>0$.

Choose $C>0$ so large that
\[
    c_UC^4>\varepsilon_0.
\]
Since $\mathfrak{a}_i\to+\infty$, for all sufficiently large $i$ one has
\[
    C\mathfrak{a}_i^{-1/2}<\rho_U.
\]
Therefore, for every $x\in U'$,
\[
    (C\mathfrak{a}_i^{-1/2})^{4-n}
    \int_{B_{C\mathfrak{a}_i^{-1/2}}(x)}
    e_i
    \geqslant
    c_U\mathfrak{a}_i^2(C\mathfrak{a}_i^{-1/2})^4
    =
    c_UC^4
    >
    \varepsilon_0.
\]
Thus the small energy condition fails at the radius
\[
    C\mathfrak{a}_i^{-1/2}.
\]
By the definition of the regularity scale,
\[
    \mathfrak{r}_i(x)\leqslant C\mathfrak{a}_i^{-1/2}.
\]
This proves the proposition.
\end{proof}

Thus, on compact subsets where
$\overline\omega_{\mathrm{har}}$ is bounded away from zero, the regularity
scale is comparable to $\mathfrak{a}_i^{-1/2}$. This conclusion does not
apply at zeros of the harmonic limit. The next subsection records the
fixed scale vanishing information for the limit and isolates the
additional hypothesis needed to obtain estimates at shrinking scales for the
sequence.

\subsection{Vanishing order and conditional regularity scale estimates}
\label{subsec:harmonic_profile_frequency}

Let $x$ be a zero of a nontrivial harmonic limit
$\overline\omega_{\mathrm{har}}$. Standard unique continuation gives a
finite vanishing order and corresponding fixed scale estimates for the
limit. These estimates alone do not control
$\overline\omega_i$ at radii tending to zero with $i$. The latter control
is imposed separately in
Definition~\ref{def:quantitative_order_q_profile}; the resulting
regularity scale estimate is conditional on that hypothesis.

We first recall the fixed scale statement for harmonic forms.

\begin{lemma}[Finite vanishing order for harmonic forms]
\label{lem:finite_vanishing_order_harmonic_forms}
Let $\eta$ be a harmonic $2$-form on a geodesic ball
$B_{r_0}(x)$:
\[
    d\eta=0,
    \qquad
    d^*\eta=0.
\]
If $\eta$ is not identically zero in a neighbourhood of $x$, then there
is an integer $q\geqslant0$ and a nonzero homogeneous harmonic
$2$-form $\eta_q$ of degree $q$ on
$T_xX\simeq\mathbb R^n$ such that, in normal coordinates centred at
$x$,
\[
    \eta(y)
    =
    \eta_q(y)
    +
    O(|y|^{q+1}).
\]
Consequently, for all sufficiently small $r>0$,
\begin{equation}
\label{eq:harmonic_form_vanishing_order_bounds}
    c r^{n+2q}
    \leqslant
    \int_{B_r(x)}|\eta|^2
    \leqslant
    C r^{n+2q},
    \qquad
    \sup_{B_r(x)}|\eta|
    \leqslant
    C r^q,
\end{equation}
for constants $0<c<C<\infty$.

Moreover, harmonic forms satisfy strong unique continuation. Thus, if
$\eta$ is harmonic on a connected open set $\Omega$ and vanishes on a
nonempty open subset, then $\eta\equiv0$ on $\Omega$. In particular, a
nontrivial harmonic form has zero set with empty interior.
\end{lemma}

\begin{proof}
By the Bochner--Weitzenb\"ock formula and $d\eta=0=d^*\eta$, we have $\nabla^*\nabla\eta+\mathcal R\eta = \Delta\eta = 0$. Thus, in a local trivialization, $\eta$ solves a second-order elliptic system with scalar principal symbol and smooth coefficients. The standard Almgren--Garofalo--Lin frequency theory for such systems gives an almost-monotone frequency at sufficiently small scales, finite vanishing
order, and strong unique continuation; see
\cites{AronszajnKrzywickiSzarski1962,GarofaloLin1986,
GarofaloLin1987,HanLin2011}.

The finite vanishing order is an integer because it is the degree of the
first nonzero Taylor term. Freezing the principal part at $x$ shows that
this first term is a homogeneous harmonic form on $T_xX$. The expansion
and the estimates
\eqref{eq:harmonic_form_vanishing_order_bounds} follow from Taylor's
theorem and local elliptic estimates. Strong unique continuation gives the
final assertion.
\end{proof}

Applied to $\overline\omega_{\mathrm{har}}$, the lemma shows that its zero
set has empty interior and that every zero has finite order. Proposition~\ref{prop:regularity_scale_nonzero_harmonic_profile}
applies on compact subsets of its nonvanishing set. Near a zero, however, smooth
convergence on fixed compact subsets does not transfer the vanishing
estimates uniformly to the radii relevant to the regularity scale. We
therefore state the required control at shrinking scales as a separate
hypothesis.

\begin{definition}[Quantitative order-$q$ profile]
\label{def:quantitative_order_q_profile}
Let $x\in\operatorname{int}K'$ and let $q\geqslant 0$ be an integer. We say that the normalized
longitudinal curvature forms $\overline\omega_i$ have a \textbf{quantitative
order-$q$ profile at $x$} if there exist constants $r_x>0$, with $B_{r_x}(x)\Subset K'$, and $0<c<C<\infty$ such that, for all sufficiently large $i$ and all
$0<r<r_x$,
\begin{equation}\label{eq:quantitative_profile_upper}
    \sup_{B_r(x)}
    |\overline\omega_i|
    \leqslant
    C r^q,
\end{equation}
and
\begin{equation}\label{eq:quantitative_profile_lower}
    \int_{B_r(x)}
    |\overline\omega_i|^2
    \geqslant
    c r^{n+2q}.
\end{equation}
\end{definition}
\begin{remark}
The quantitative profile condition is genuinely stronger than smooth
convergence
\[
    \overline\omega_i
    \longrightarrow
    \overline\omega_{\mathrm{har}}
\]
on fixed compact subsets. In particular, when $q>0$,
\eqref{eq:quantitative_profile_upper} implies
\[
    \overline\omega_i(x)=0
\]
for every sufficiently large $i$, and more generally imposes a uniform
order-$q$ upper vanishing estimate at the fixed center $x$. Neither this
exact centering nor the control down to arbitrarily small radii follows
from the fact that the limiting harmonic form has order $q$ at $x$.
\end{remark}
\begin{remark}
If $\overline\omega_{\mathrm{har}}$ is a nonzero harmonic form with vanishing
order $q$ at $x$, then the estimates
\[
    \sup_{B_r(x)}
    |\overline\omega_{\mathrm{har}}|
    \leqslant
    C r^q,
    \qquad
    \int_{B_r(x)}
    |\overline\omega_{\mathrm{har}}|^2
    \geqslant
    c r^{n+2q}
\]
hold for all sufficiently small $r$. Smooth convergence transfers these
bounds to $\overline\omega_i$ on fixed scales. The additional content of
Definition~\ref{def:quantitative_order_q_profile} is that the same control
is assumed uniformly at the shrinking scales relevant to the curvature
problem. Establishing such a condition for the almost harmonic sequence
$\overline\omega_i$ would require a separate quantitative frequency
analysis controlling the Hodge errors and the lower-order jets at those
shrinking scales. That analysis is not carried out in the present paper;
the regularity scale result below is conditional on the stated profile
hypothesis.
\end{remark}

Under this quantitative profile hypothesis, one obtains the following
two-sided estimate for the regularity scale.

\begin{proposition}[Quantitative vanishing order and regularity scale]
\label{prop:vanishing_order_regular_scale}
Assume that
\[
    \mathfrak{a}_i:=\|\widehat\omega_i\|_{L^2(K)}
    \longrightarrow+\infty
\]
and that
\[
    \overline\omega_i:=\mathfrak{a}_i^{-1}\widehat\omega_i
    \longrightarrow
    \overline\omega_{\mathrm{har}}
\]
smoothly on $K'$. Suppose that $\overline\omega_i$ has a quantitative
order-$q$ profile at $x\in K'$ in the sense of
Definition~\ref{def:quantitative_order_q_profile}. Then there exist
constants $0<c<C<\infty$ such that
\[
    c \mathfrak{a}_i^{-1/(q+2)}
    \leqslant
    \mathfrak{r}_i(x)
    \leqslant
    C \mathfrak{a}_i^{-1/(q+2)}
\]
for all sufficiently large $i$.
\end{proposition}

\begin{proof}
Let $r_x>0$ be as in
Definition~\ref{def:quantitative_order_q_profile}. Set
\[
    \rho_x
    :=
    \min\left\{
        r_x,
        r_0,
        \operatorname{dist}(x,X\setminus K')
    \right\}.
\]
All radii below are chosen smaller than $\rho_x$. 

Since
$\widehat\omega_i
    =
    \mathfrak{a}_i\overline\omega_i$,
and $|F_{\nabla_i}^{\parallel}|
    =
    |\widehat\omega_i|$, 
the quantitative profile bounds give direct control of the longitudinal
curvature at small scales around $x$.

We first prove the lower bound for $\mathfrak{r}_i(x)$. Let
\[
    r_i:=c_0\mathfrak{a}_i^{-1/(q+2)}
\]
with $c_0>0$ to be chosen. Since $\mathfrak{a}_i\to+\infty$, for all sufficiently
large $i$ one has
\[
    r_i<\rho_x.
\]
For every $0<s\leqslant r_i$, the upper profile estimate gives
\[
    |\widehat\omega_i|
    \leqslant
    C \mathfrak{a}_i s^q
    \qquad
    \text{on }B_s(x).
\]
The transverse curvature components decay exponentially on compact subsets
of $\mathcal R_{\mathrm{ab}}$, and the $\Theta$-monopole equation gives
\[
    |\nabla_i\Phi_i|
    =
    |F_{\nabla_i}\wedge\Theta|
    \leqslant
    C|F_{\nabla_i}|.
\]
Therefore, uniformly for $0<s\leqslant r_i$,
\[
    s^{4-n}
    \int_{B_s(x)}
    e_i
    \leqslant
    C \mathfrak{a}_i^2 s^{4+2q}
    +
    C e^{-c m_i}s^4.
\]
In particular,
\[
    s^{4-n}
    \int_{B_s(x)}
    e_i
    \leqslant
    C c_0^{4+2q}
    +
    o_i(1)
\]
for every $0<s\leqslant r_i$. Choosing $c_0>0$ sufficiently small and then
taking $i$ sufficiently large, we obtain
\[
    s^{4-n}
    \int_{B_s(x)}
    e_i
    <
    \varepsilon_0
\]
for every $0<s\leqslant r_i$. Hence
\[
    \mathfrak{r}_i(x)\geqslant c\mathfrak{a}_i^{-1/(q+2)}
\]
for some $c>0$.

We now prove the upper bound. Let
\[
    R_i:=C_0\mathfrak{a}_i^{-1/(q+2)}
\]
with $C_0>0$ to be chosen. For all sufficiently large $i$,
\[
    R_i<\rho_x.
\]
The lower profile estimate gives
\[
    \int_{B_{R_i}(x)}
    |\widehat\omega_i|^2
    =
    \mathfrak{a}_i^2
    \int_{B_{R_i}(x)}
    |\overline\omega_i|^2
    \geqslant
    c \mathfrak{a}_i^2 R_i^{n+2q}.
\]
Since $|F_{\nabla_i}^{\parallel}|=|\widehat\omega_i|$, we have
\[
\begin{aligned}
    R_i^{4-n}
    \int_{B_{R_i}(x)}
    e_i
    &\geqslant
    R_i^{4-n}
    \int_{B_{R_i}(x)}
    |F_{\nabla_i}^{\parallel}|^2        \\
    &=
    R_i^{4-n}
    \int_{B_{R_i}(x)}
    |\widehat\omega_i|^2                \\
    &\geqslant
    c \mathfrak{a}_i^2 R_i^{4+2q}
    =
    c C_0^{4+2q}.
\end{aligned}
\]
Choosing $C_0>0$ sufficiently large gives
\[
    R_i^{4-n}
    \int_{B_{R_i}(x)}
    e_i
    >
    \varepsilon_0.
\]
Thus the small energy condition fails at the radius $R_i$, and by the
definition of the regularity scale,
\[
    \mathfrak{r}_i(x)
    \leqslant
    C \mathfrak{a}_i^{-1/(q+2)}
\]
for some $C>0$.
\end{proof}

\begin{corollary}[Regularity scale relative to the mass]
\label{cor:regularity_scale_mass_comparison}
For every compact set $K_0\Subset\mathcal R_{\mathrm{ab}}$,
\[
    \inf_{x\in K_0}
    m_i^{1/4}\mathfrak r_i(x)
    \longrightarrow
    +\infty.
\]
This conclusion is independent of the alternative for the harmonic corrections.

Assume in addition that the second alternative of
Corollary~\ref{cor:harmonic_correction_alternatives} holds, fix compact
domains $K'\Subset K\Subset\mathcal R_{\mathrm{ab}}$, and set
\[
    \mathfrak a_i
    :=
    \|\widehat\omega_i\|_{L^2(K)}.
\]
If the normalized longitudinal curvature has a quantitative order-$q$
profile at $x\in K'$ in the sense of
Definition~\ref{def:quantitative_order_q_profile}, then
\[
    m_i^{1/(2(q+2))}\mathfrak r_i(x)
    \longrightarrow
    +\infty.
\]
\end{corollary}

\begin{proof}
The first assertion is
Corollary~\ref{cor:regularity_scale_mass_separation_Rab}. Under the
quantitative order-$q$ profile hypothesis,
Proposition~\ref{prop:vanishing_order_regular_scale} and
\eqref{eq:local_longitudinal_amplitude_subcritical} give
\[
    m_i^{1/(2(q+2))}\mathfrak r_i(x)
    \geqslant
    c
    \left(
        \frac{m_i^{1/2}}{\mathfrak a_i}
    \right)^{1/(q+2)}
    \longrightarrow+\infty.
\]
\end{proof}

\begin{remark}
When $q=0$, Proposition~\ref{prop:vanishing_order_regular_scale} gives the
same exponent $1/2$ as
Proposition~\ref{prop:regularity_scale_nonzero_harmonic_profile}. The two
results have different hypotheses: the latter applies where the harmonic
limit is bounded away from zero, whereas the former assumes the
quantitative order-$0$ bounds of
Definition~\ref{def:quantitative_order_q_profile} at the chosen point.
\end{remark}

The geometric identification
\[
    \mathcal C
    =
    \mathcal S\cup\mathcal O,
    \qquad
    \mathcal R_{\mathrm{ab}}
    =
    X\setminus(\mathcal S\cup\mathcal O)
\]
was established in
Corollary~\ref{cor:nonabelian_locus_obstruction_identity}, independently
of the harmonic correction analysis. The new global estimate in this
section is
\[
    \|\sigma_i\|_{L^2(X)}=o(m_i^{1/2}),
\]
while Proposition~\ref{prop:subcritical_harmonic_correction_growth}
upgrades the basic local $C^0$ subcriticality from
Corollary~\ref{cor:subcritical_longitudinal_curvature_Rab} to all local
derivatives of the longitudinal curvature and the full monopole fields.

If Li's harmonic corrections are uniformly bounded in $L^2(X)$, the
translated configurations converge smoothly after gauge on compact
subsets of $\mathcal R_{\mathrm{ab}}$. If their $L^2(X)$ norms have
infinite limsup, then, after passing to a subsequence along which these
norms tend to infinity, the normalized uncorrected longitudinal
curvatures converge smoothly on compact subsets of
$\mathcal R_{\mathrm{ab}}$ to the restriction of a nonzero global
$L^2$-harmonic $2$-form. Independently of this alternative,
$m_i^{1/4}\mathfrak r_i\to+\infty$ locally uniformly throughout
$\mathcal R_{\mathrm{ab}}$. Where the limiting harmonic form is bounded
away from zero, the sharper comparison
$\mathfrak r_i\asymp\mathfrak a_i^{-1/2}$ holds. At a point satisfying
the additional quantitative order-$q$ hypothesis, one has
$\mathfrak r_i\asymp\mathfrak a_i^{-1/(q+2)}$ and
$m_i^{1/(2(q+2))}\mathfrak r_i(x)\to+\infty$. At zeros of the limiting harmonic form, no comparison with
$\mathfrak a_i$ is asserted without a quantitative profile hypothesis valid
at shrinking scales.

The section does not address the size of $\mathcal O$, the
nonemptiness or density of $\mathcal R_{\mathrm{ab}}$, or the escape of the loci $C_{\Lambda_0,i}$ to infinity. These questions are discussed in
Section~\ref{sec:resulting_picture}.


\section{The resulting picture in the large mass regime and open problems}
\label{sec:resulting_picture}

The main results give the hierarchy
\[
    \mathcal S
    \subset
    \mathcal Z
    \subset
    \mathcal C
    =
    \mathcal S\cup\mathcal O.
\]
Here $\mathcal S$ is the compact support of the $\Theta$-calibrated
$(n-3)$-cycle $T:=T_{\mathrm{PPS}}=T_{\mathrm{Li}}$, $\mathcal Z$ is
the Kuratowski upper limit of the zero sets of the Higgs fields,
$\mathcal C$ is the limiting nonabelian locus, namely the Kuratowski
upper limit of Li's curvature concentration loci, and $\mathcal O$ is the
closed set where effective
codimension-three monotonicity of the mass-renormalized energy fails locally. Thus
\[
    \mathcal C\setminus\mathcal S
    =
    \mathcal O\setminus\mathcal S,
    \qquad
    \mathcal R_{\mathrm{ab}}
    =
    X\setminus(\mathcal S\cup\mathcal O).
\]
When $n=3$, the identity $r^{3-n}=1$ makes effective
codimension-three monotonicity automatic. Together with the
three-dimensional monopole theory,
Remark~\ref{rmk:three_dimensional_benchmark_C_O} gives equality of the
three concentration and nonabelian loci, an empty obstruction locus, and
a residual Dirac monopole on the whole complement of the concentration
set. By contrast, the possible excess $\mathcal C\setminus\mathcal S$
and the harmonic corrections appearing below are higher-dimensional
phenomena.

The inclusion $\mathcal S\subset\mathcal Z$ is quantitative. If
$Z_i=\Phi_i^{-1}(0)$, then Proposition~\ref{prop:S_subset_Z} gives
\[
    \sup_{x\in\mathcal S} d(x,Z_i)\longrightarrow0.
\]
The zeros used in this argument are forced by transverse charge. On the
normal balls provided by Li's analysis, one passes to a local
$\mathrm{SU}(2)$ lift; the first Chern number of the eigenline on an
outer linking sphere equals the local multiplicity of $T$. Equivalently,
in the $\mathrm{SO}(3)$ description, the Euler number of the induced
$\mathrm{SO}(2)$ reduction is twice this charge. The enclosed ball must
therefore contain a zero. This is a statement about the total charge in
the normal ball and does not assign a canonical charge to each individual
zero.

Li's limiting singular abelian monopole from \cite[Theorem~1.12]{li2025large} satisfies
\[
    dF_\infty=2\pi T
\]
in the sense of distributions. In the Calabi--Yau case,
Proposition~\ref{prop:smooth_convergence_Li_corrected_longitudinal} also gives $\Lambda F_\infty=0$ distributionally. In the alternative with bounded harmonic corrections, its restriction to $\mathcal R_{\mathrm{ab}}$ is related to the smooth abelian limit of the translated configurations by Remark~\ref{rmk:comparison_with_Li_singular_limit}. Thus $\mathcal S$ is the support of the limiting Dirac source. Points of
$\mathcal Z\setminus\mathcal S$, if they occur, are limits of Higgs
zeros which carry no codimension-three charge detected by $T$ and lie in
$\mathcal O\setminus\mathcal S$. At every point of
$\mathcal C\setminus\mathcal S=\mathcal O\setminus\mathcal S$,
Proposition~\ref{prop:excess_points_determine_characteristic_cores}
produces points $p_i$ and radii $R_i\downarrow0$ for which Li's weighted
curvature potential on $B_{R_i}(p_i)$ is bounded below while the
mass-renormalized Yang--Mills--Higgs energy tends to zero. By Proposition~\ref{prop:regularity_scale_inside_characteristic_core}, the
ordinary regularity scale at the selected centre is at most the
characteristic radius. If its ratio to $m_i^{-1}$ tends to infinity, then
Proposition~\ref{prop:characteristic_core_above_mass_scale} shows that
the limit at the mass scale is flat with parallel Higgs field, the weighted
potential remains outside every fixed multiple of $m_i^{-1}$, and the
codimension-four scale-invariant curvature energy becomes unbounded on
suitable intervening annuli.

If $\mathcal O\subset\mathcal S$, then
\[
    \mathcal S=\mathcal Z=\mathcal C,
    \qquad
    \mathcal R_{\mathrm{ab}}=X\setminus\mathcal S,
\]
and Corollary~\ref{cor:Hausdorff_convergence_unobstructed} gives local
Hausdorff convergence of the zero sets and Li loci on every fixed
compact neighbourhood of $\mathcal S$. This does not imply tightness of
the curvature concentration loci on the AC end. Independently of the harmonic correction
alternative, Corollary~\ref{cor:subcritical_longitudinal_curvature_Rab}
gives
\[
    m_i^{-1/2}\widehat\omega_i
    \longrightarrow0
    \qquad
    \text{in }L^\infty_{\mathrm{loc}}(\mathcal R_{\mathrm{ab}}),
\]
and shows that the difference between the mass-normalized and
$|\Phi_i|$-normalized longitudinal forms converges to zero in
$L^p_{\mathrm{loc}}$ for every
$1\leqslant p<2n/(n-2)$. Moreover,
Corollary~\ref{cor:regularity_scale_mass_separation_Rab} gives
\[
    \inf_{x\in K}m_i^{1/4}\mathfrak r_i(x)
    \longrightarrow+\infty
\]
for every compact set $K\Subset\mathcal R_{\mathrm{ab}}$.

On $\mathcal R_{\mathrm{ab}}$,
Proposition~\ref{prop:smooth_convergence_Li_corrected_longitudinal} gives
smooth convergence of Li's corrected longitudinal forms, while
Proposition~\ref{prop:subcritical_harmonic_correction_growth} gives
\[
    \|\sigma_i\|_{L^2(X)}=o(m_i^{1/2}).
\]
By Corollary~\ref{cor:harmonic_correction_alternatives}, smooth
compactness of the translated configurations follows when the
corresponding sequence in
$\mathscr H_{(2)}^2(X)\simeq H_c^2(X;\mathbb R)$ is uniformly bounded.
If its $L^2(X)$ norms have infinite limsup, then, after passing to a
subsequence along which the norms tend to infinity, the normalized
corrections converge to a nonzero global $L^2$-harmonic form.

\subsection{Regimes at characteristic scales at excess points}
\label{subsec:characteristic_scale_regimes}

Let $x\in\mathcal C\setminus\mathcal S$, and choose points
$p_i\to x$, characteristic radii
$R_i=r_{2\Lambda_0,i}(p_i)\to0$, and ordinary regularity scales
\[
    \lambda_i
    :=
    \mathfrak r_i(p_i)
    \leqslant
    R_i
\]
as above.  After passage to a subsequence, the dimensionless quantities
\[
    \tau_i
    :=
    m_i\lambda_i
\]
fall into one of the usual three regimes: $\tau_i\to0$,
$\tau_i\to\tau\in(0,+\infty)$, or $\tau_i\to+\infty$.

If $\tau_i\to0$, then the same rescaling argument as in
Proposition~\ref{prop:zero_centered_scale_regimes}\textnormal{(a)} shows
that every smooth pointed $\lambda_i$-scale limit has vanishing Higgs
field and satisfies the tangent space instanton equation
\[
    F_{\nabla_\infty}\wedge\Theta_x=0,
\]
together with
$\Lambda_{\omega_x}F_{\nabla_\infty}=0$ in the Calabi--Yau case.  The
regime $\tau_i\to\tau\in(0,+\infty)$ is the one compatible with a
finite nonzero rescaled mass, but no transverse dimension reduction or
finite energy conclusion at an arbitrary point of
$\mathcal C\setminus\mathcal S$ follows from the present estimates.

The regime $\tau_i\to+\infty$ is constrained by
Proposition~\ref{prop:characteristic_core_above_mass_scale}.  The
limit at the mass scale is flat with parallel Higgs field, while, for every
fixed $L>1$, the characteristic potential remains in
$B_{R_i}(p_i)\setminus B_{L/m_i}(p_i)$.  Moreover,
\eqref{ineq:large_annular_codim4_energy_characteristic_core} implies
that the codimension-four scale-invariant curvature energy is not
uniformly bounded on all annuli between these scales.  This conclusion
is specific to Li's weighted potential: it is compatible with the local
$o(m_i)$ energy bound
\eqref{ineq:core_energy_vanishing} because the two quantities have
different scaling.

At a point of $\mathcal C\setminus\mathcal Z$, the centres $p_i$ and
balls $B_{R_i}(p_i)$ may be chosen in a fixed neighbourhood which is
eventually disjoint from the zero sets of the Higgs fields. The preceding
alternatives constrain such a point but do not exclude it. In the regime
$\tau_i\to0$, one must understand approximations to tangent space
instantons by configurations whose Higgs fields have no zeros in that
neighbourhood.  In the finite regime one
needs a transverse bubble extraction theorem.  In the regime
$\tau_i\to+\infty$, one needs an estimate excluding the annular
concentration quantified by
\eqref{eq:divergent_annular_codim4_energy_characteristic_core}.

The contrast with dimension three is important.  For large mass
Bogomolny monopoles, the energy identity and moving centre compactness
force every nontrivial finite mass profile to occur at scale comparable
to $m_i^{-1}$, and no secondary bubbling occurs inside a mass one
profile
\cite[Lemma~8.7 and Proposition~8.3]{fadeloliveira2026limitv5}.  For the
three-dimensional Yang--Mills--Higgs functional with Higgs
self-interaction, the selected bubbling scales are likewise comparable
to the distinguished parameter, and the complete analysis includes
bubble extraction and vanishing neck energy
\cite[Lemma~5.17 and Propositions~5.19 and~5.27]{cheng2025su2}.  In the
present higher-dimensional problem, the instanton regime below the mass
scale and the annular alternative above it are not removed by the
available identities.

\subsection{Blow-up scales centred at Higgs zeros}
\label{subsec:zero_centered_scale_regimes}

The scales furnished by Proposition~\ref{prop:zero_centered_TU_concentration},
with centres chosen in the zero sets of the Higgs fields, admit the
following elementary limiting descriptions. They do not by themselves
determine which regime occurs.

\begin{proposition}[Scale regimes centred at Higgs zeros]
\label{prop:zero_centered_scale_regimes}
Let $x\in\mathcal Z$, choose a subsequence along which the zero is realized and points
$x_i\in Z_i$ as in
Proposition~\ref{prop:zero_centered_TU_concentration}, and set
\[
    \lambda_i
    :=
    \mathfrak r_i(x_i),
    \qquad
    \tau_i
    :=
    m_i\lambda_i.
\]
Choose an oriented orthonormal frame
$(E_1,\ldots,E_n)$ on a neighbourhood of $x$, and let
\[
    \iota_i:\mathbb R^n\longrightarrow T_{x_i}X,
    \qquad
    \iota_\infty:\mathbb R^n\longrightarrow T_xX
\]
be the oriented linear isometries determined by
$\iota_i(e_a)=E_a(x_i)$ and
$\iota_\infty(e_a)=E_a(x)$.
For a scale $\rho_i\to0$, let
\[
    \delta_{i,\rho_i}(y)
    :=
    \exp_{x_i}\bigl(\rho_i\iota_i(y)\bigr),
\]
and define on the expanding domains in $\mathbb R^n$
\[
    g_{i,\rho_i}
    :=
    \rho_i^{-2}\delta_{i,\rho_i}^*g,
    \qquad
    \nabla_{i,\rho_i}
    :=
    \delta_{i,\rho_i}^*\nabla_i,
    \qquad
    \Phi_{i,\rho_i}
    :=
    \rho_i\delta_{i,\rho_i}^*\Phi_i.
\]
Set
\[
    \Theta_x
    :=
    \iota_\infty^*(\Theta|_x),
\]
and, in the Calabi--Yau case, set likewise
$\omega_x:=\iota_\infty^*(\omega|_x)$.
Then the following hold.

\begin{enumerate}[label=\textnormal{(\alph*)}]
\item If $\tau_i\to0$, then every smooth pointed limit of the
$\lambda_i$-rescaled configurations has vanishing Higgs field and
satisfies
\[
    F_{\nabla_\infty}\wedge\Theta_x=0.
\]
In the Calabi--Yau case, it also satisfies
\[
    \Lambda_{\omega_x}F_{\nabla_\infty}=0.
\]
Thus every such limit is a tangent space instanton in the corresponding
special holonomy sense.

\item If $\tau_i\to+\infty$, then, after passing to a further subsequence
and applying gauge transformations, the rescalings at the mass scale
\[
    (\nabla_{i,m_i^{-1}},\Phi_{i,m_i^{-1}})
\]
converge smoothly on compact subsets of $\mathbb R^n$ to a flat
connection with identically zero Higgs field.

\item For every sufficiently small fixed $\delta>0$,
\begin{equation}
\label{eq:zero_centered_Green_localization}
    2m_i^{-2}
    \int_{B_\delta(x_i)}
    G(x_i,y)|\nabla_i\Phi_i|^2(y)\,\vol(y)
    \longrightarrow1.
\end{equation}
If, in addition, $\tau_i\to+\infty$, then
\begin{equation}
\label{eq:zero_centered_Green_neck}
    \lim_{R\to\infty}
    \lim_{i\to\infty}
    2m_i^{-2}
    \int_{B_\delta(x_i)\setminus B_{R/m_i}(x_i)}
    G(x_i,y)|\nabla_i\Phi_i|^2(y)\,\vol(y)
    =
    1.
\end{equation}
\end{enumerate}
\end{proposition}

\begin{proof}
By the maximum principle,
$|\Phi_i|\leqslant m_i$. Therefore, in the $\lambda_i$-rescaling,
\[
    |\Phi_{i,\lambda_i}|
    \leqslant
    m_i\lambda_i
    =
    \tau_i.
\]
If $\tau_i\to0$, every smooth pointed limit has zero Higgs field. Under
the same rescaling, the parallel forms
\[
    \Theta_{i,\lambda_i}
    :=
    \lambda_i^{3-n}\delta_{i,\lambda_i}^*\Theta
\]
converge smoothly on compact subsets to the constant tangent space form
$\Theta_x$, and the rescaled monopole equation passes to
\[
    F_{\nabla_\infty}\wedge\Theta_x=0.
\]
In the Calabi--Yau case, the rescaled K\"ahler forms
\[
    \omega_{i,\lambda_i}
    :=
    \lambda_i^{-2}\delta_{i,\lambda_i}^*\omega
\]
converge to $\omega_x$, and the primitivity equations pass to
$\Lambda_{\omega_x}F_{\nabla_\infty}=0$. This proves
part~\textnormal{(a)}.

Assume next that $\tau_i\to+\infty$. Fix $R>0$. For all sufficiently
large $i$,
\[
    \frac{R}{m_i}
    <
    \lambda_i.
\]
The definition of $\lambda_i$ gives
\[
    \left(\frac{R}{m_i}\right)^{4-n}
    \int_{B_{R/m_i}(x_i)}
    e_i
    <
    \varepsilon_0.
\]
Equivalently, the rescalings at the mass scale satisfy
\[
    R^{4-n}
    \int_{B_R(0)}
    e(\nabla_{i,m_i^{-1}},\Phi_{i,m_i^{-1}})
    <
    \varepsilon_0.
\]
Theorem~\ref{thm: total_epsilon_regularity}, local Uhlenbeck gauge
fixing, and a diagonal argument therefore give smooth convergence modulo
gauge on compact subsets of $\mathbb R^n$ to a Yang--Mills--Higgs
configuration $(\nabla_\infty,\Phi_\infty)$. For every $R>0$, the
interior estimate gives
\[
    \sup_{B_{R/2}(0)}
    e(\nabla_\infty,\Phi_\infty)
    \leqslant
    C\varepsilon_0R^{-4}.
\]
Letting $R\to\infty$ shows that
$F_{\nabla_\infty}=0$ and
$\nabla_\infty\Phi_\infty=0$. Since
$\Phi_{i,m_i^{-1}}(0)=0$, one has
$\Phi_\infty(0)=0$, and hence $\Phi_\infty\equiv0$. This proves
part~\textnormal{(b)}.

At the zero $x_i$, Theorem~\ref{thm: alternative_finite_mass} gives
\[
    m_i^2
    =
    2\int_X
    G(x_i,y)|\nabla_i\Phi_i|^2(y)\,\vol(y).
\]
Choose $\delta>0$ so small that the points $x_i$ and the balls
$B_{2\delta}(x_i)$ lie in a fixed compact set for all sufficiently large
$i$. The Green kernel is uniformly bounded on
\[
    \left\{
        (x_i,y):
        y\in X\setminus B_\delta(x_i)
    \right\},
\]
and the fixed class energy identity gives
\[
    \int_X|\nabla_i\Phi_i|^2
    =
    4\pi k m_i.
\]
Consequently,
\[
\begin{aligned}
    0
    &\leqslant
    2m_i^{-2}
    \int_{X\setminus B_\delta(x_i)}
    G(x_i,y)|\nabla_i\Phi_i|^2(y)\,\vol(y)
    \\
    &\leqslant
    \frac{C_\delta}{m_i}
    \longrightarrow0,
\end{aligned}
\]
which proves
\eqref{eq:zero_centered_Green_localization}.

Finally, assume $\tau_i\to+\infty$ and fix $R>0$. Define the rescaled
Green kernels on $B_R(0)$ by
\[
    G_i^{(m)}(z)
    :=
    m_i^{2-n}
    G\bigl(x_i,\delta_{i,m_i^{-1}}(z)\bigr).
\]
The local singular estimate for $G$, uniformly for the centers $x_i$ in
a fixed compact set, gives
\[
    0\leqslant G_i^{(m)}(z)
    \leqslant
    C|z|^{2-n}
\]
on $B_R(0)\setminus\{0\}$. By change of variables,
\begin{align*}
    &2m_i^{-2}
    \int_{B_{R/m_i}(x_i)}
    G(x_i,y)|\nabla_i\Phi_i|^2(y)\,\vol(y)
    \\
    &\qquad=
    2\int_{B_R(0)}
    G_i^{(m)}
    |\nabla_{i,m_i^{-1}}\Phi_{i,m_i^{-1}}|^2
    \,\operatorname{vol}_{g_{i,m_i^{-1}}}.
\end{align*}
The smooth convergence at the mass scale from part~\textnormal{(b)} makes the
second factor converge uniformly to zero on $B_R(0)$, while
$|z|^{2-n}$ is locally integrable. Dominated convergence therefore gives
\[
    2m_i^{-2}
    \int_{B_{R/m_i}(x_i)}
    G(x_i,y)|\nabla_i\Phi_i|^2(y)\,\vol(y)
    \longrightarrow0.
\]
Subtracting this from \eqref{eq:zero_centered_Green_localization}, then letting
$R\to\infty$, proves \eqref{eq:zero_centered_Green_neck}.
\end{proof}

\subsection{Size, tightness, excess loci, and bubbling}
\label{subsec:concentration_excess_bubbling_questions}

The questions below concern the obstruction locus, points of $\mathcal C\setminus\mathcal S$,
bubbling centred at Higgs zeros, and escape of the loci $C_{\Lambda_0,i}$ to infinity.

\begin{enumerate}

\item
\textbf{Effective codimension-three monotonicity.}
The principal analytic problem is to find natural hypotheses implying
\[
    \mathcal O=\varnothing.
\]
The stress-energy identity shows that the relevant error is the positive
imbalance
\[
    q_i^+
    =
    \bigl(|F_{\nabla_i}|^2-|\nabla_i\Phi_i|^2\bigr)^+.
\]
Proposition~\ref{prop:codim2_Morrey_criterion_positive_monotonicity}
shows that a codimension-two Morrey estimate
\[
    m_i^{-1}s^{2-n}
    \int_{B_s(y)} q_i^+
    \leqslant
    \Lambda_U,
\]
If this estimate holds uniformly for every ball $B_s(y)\subset U$
and all sufficiently large indices, then effective codimension-three
monotonicity holds on $U$ with exponent $\gamma_U=1$. If one constant
and one threshold index work for every ball in $X$, effective
monotonicity holds globally. Hence $\mathcal O=\varnothing$ and
$\mathcal S=\mathcal Z=\mathcal C$. To exclude excess points and obtain
local Hausdorff convergence, the weaker conclusion
$\mathcal O\subset\mathcal S$ is sufficient.

For $\Theta$-monopoles, the signed imbalance is the Chern--Weil density
\[
    \bigl(|F_{\nabla_i}|^2-|\nabla_i\Phi_i|^2\bigr)
    \operatorname{vol}
    =
    -\langle F_{\nabla_i}\wedge F_{\nabla_i}\rangle\wedge\Xi.
\]
Locally this form is exact, but the resulting boundary term is not
controlled directly by the nonnegative radial term in the stress-energy
identity. The problem is to identify assumptions under which its positive
part satisfies a codimension-two Morrey estimate, preferably with a power gain.

The cohomogeneity-one families in
Section~\ref{sec:homogeneous_examples} show that such a Morrey estimate
may hold even when the imbalance is not uniformly bounded pointwise; see
Remark~\ref{rmk:no_uniform_pointwise_imbalance_homogeneous}. It remains
to determine whether comparable estimates follow from hypotheses such as
stability, local minimizing properties, or a quantitative transverse
Bogomolny approximation.\\

\item
\textbf{Characteristic radii at points of $\mathcal C\setminus\mathcal S$.}
At each point
$x\in\mathcal C\setminus\mathcal S
=\mathcal O\setminus\mathcal S$,
Proposition~\ref{prop:excess_points_determine_characteristic_cores}
produces, after passing to a subsequence, points
$p_i\in C_{\Lambda_0,i}$ and radii
$R_i=r_{2\Lambda_0,i}(p_i)\to0$ such that
\[
    m_i^{-2}
    \int_{B_{R_i}(p_i)}
    \frac{|F_{\nabla_i}|^2(y)}
    {\max\{d(y,p_i),m_i^{-1}\}^{n-2}}
    \,\operatorname{vol}(y)
    \geqslant
    \frac{1}{2\Lambda_0},
\]
while the mass-renormalized energy of these balls tends to zero.
Proposition~\ref{prop:regularity_scale_inside_characteristic_core}
shows, in addition, that
\[
    \lambda_i
    :=
    \mathfrak r_i(p_i)
    \leqslant
    R_i
    \longrightarrow0.
\]
Thus every point of $\mathcal C\setminus\mathcal S$ carries ordinary codimension-four loss of
regularity inside the selected ball $B_{R_i}(p_i)$.

If $m_i\lambda_i\to+\infty$, then
Proposition~\ref{prop:characteristic_core_above_mass_scale} gives three
further conclusions: the limit at the mass scale is flat with parallel Higgs
field; the characteristic potential remains outside every ball
$B_{L/m_i}(p_i)$ with fixed $L$; and there are scales
$m_i^{-1}\ll s_i<R_i$ such that
\[
    s_i^{4-n}
    \int_{B_{2s_i}(p_i)\setminus B_{s_i}(p_i)}
    |F_{\nabla_i}|^2\,\vol
    \longrightarrow+\infty.
\]
Consequently, the regime above the inverse mass scale cannot be excluded
by extracting only finitely many annular regions with uniformly bounded
codimension-four scale-invariant curvature energy.

The remaining problem is to rule out all three possible behaviours at a
point of $\mathcal C\setminus\mathcal Z$.  In the regime
$m_i\lambda_i\to0$, one needs a classification or exclusion of
zero-free approximations to tangent space instantons.  When
$m_i\lambda_i$ stays bounded above and below, one needs transverse
dimension reduction and a finite energy bubble theorem capable of
recovering nonzero three-dimensional charge.  When
$m_i\lambda_i\to+\infty$, one needs a new estimate excluding the annular
concentration in
\eqref{eq:divergent_annular_codim4_energy_characteristic_core}; a
codimension-two Morrey bound for $q_i^+$ is one sufficient mechanism,
since it excludes the point of $\mathcal C\setminus\mathcal S$ altogether.  Proving such an
exclusion for every ball eventually containing no Higgs zeros would give
$\mathcal C\subset\mathcal Z$.  Together with the already known reverse
inclusion, this would yield $\mathcal C=\mathcal Z$.\\

\item
\textbf{Size and tightness of Li's curvature concentration loci.}
The set $\mathcal S$ is the compact support of a calibrated integral
cycle, whereas $\mathcal C$ is defined as a closed upper limit of
the sets $C_{\Lambda_0,i}$. The identity
\[
    \mathcal C = \mathcal S\cup\mathcal O
\]
reduces the local size problem for $\mathcal C$ to that for
$\mathcal O$. In particular, for every compact set $K\subset X$ and
every $n-3\leqslant d\leqslant n$,
\[
    \mathcal H^{d}(K\cap\mathcal C)
    \leqslant
    \mathcal H^{d}(K\cap\mathcal S)
    +
    \mathcal H^{d}(K\cap\mathcal O).
\]
The first term is controlled by $C_n k/\ell_T$ when $d=n-3$ and is zero when $d>n-3$; see Corollary~\ref{cor: size_S}. Since $\mathcal O$ is closed and $\mathcal S$ is nowhere dense,
Corollary~\ref{cor:nonabelian_locus_obstruction_identity} also gives
\[
    \operatorname{int}(\mathcal C)
    =
    \operatorname{int}(\mathcal O).
\]
Thus the question whether $\mathcal C$ has empty
interior is exactly the question whether the effective monotonicity
obstruction locus has empty interior. The present arguments do not give
rectifiability, Hausdorff measure bounds, or capacity estimates for
$\mathcal O$. In particular, they do not rule out
\[
    \operatorname{int}(\mathcal C)\neq\varnothing
    \qquad\text{or}\qquad
    \mathcal C=X.
\]
Since $\mathcal S$ has empty interior, every nonempty open subset of
$\mathcal C$ contains a point of
$\mathcal O\setminus\mathcal S$. A natural question is whether the
definition of $\mathcal O$, together with the energy and imbalance
bounds available here, implies Minkowski, Hausdorff, or capacity bounds
for this set. At present, the available quantitative substitute is Li's mass-dependent Minkowski measure estimate for tubular neighbourhoods of the loci $C_{\Lambda_0,i}$, valid in the scale range of \cite[Proposition~2.20]{li2025large}.

There is a separate problem at infinity. The definition of
$\mathcal C$ records convergent realizations in $X$ and does not
exclude sequences
\[
    p_i\in C_{\Lambda_0,i},
    \qquad
    \rho(p_i)\longrightarrow\infty.
\]
The large radius estimates do not presently exclude shrinking balls at the
mass or characteristic scales centred at such points; see
Remark~\ref{rmk:Cinfty_escape_infinity}. A tightness estimate uniform in $i$
would be needed for global Hausdorff convergence of the loci
$C_{\Lambda_0,i}$. The identity $\mathcal C=\mathcal S\cup\mathcal O$ concerns convergent realizations in $X$ and does not address this escape phenomenon.

In dimension three the corresponding loss of mass-renormalized energy is
quantized exactly. If $K_a$ are the total charges of the complete finite
clusters, then
\[
    k-\sum_aK_a
\]
is the charge escaping through the AC end, and the energy measures are
tight if and only if this integer vanishes
\cite[Theorem~1.1]{fadeloliveira2026limitv5}
\cite[Corollary~5.2]{fadel2026abelian}. This concerns the physical energy
measures and monopole cores. The higher-dimensional question above is
stronger in a different direction: it asks for tightness of Li's
auxiliary curvature concentration loci themselves.

A weaker intermediate problem is to prove that $\mathcal O$ is
bounded. Since $\mathcal S$ is compact,
$\mathcal O$ and
$\mathcal C=\mathcal S\cup\mathcal O$ are bounded
simultaneously. Both sets are closed, and the complete Riemannian
manifold $X$ is proper; hence the following conditions are equivalent:
\[
    \mathcal O\text{ is bounded},
    \qquad
    \mathcal O\text{ is compact},
    \qquad
    \mathcal C\text{ is bounded},
    \qquad
    \mathcal C\text{ is compact}.
\]
Under these conditions,
$\mathcal Z\subset\mathcal C$ is compact and there exists
$R<\infty$ such that
\[
    X\setminus B_R
    \subset
    \mathcal R_{\mathrm{ab}}.
\]
Thus the entire AC end lies in the open set $\mathcal R_{\mathrm{ab}}$, although
components of the curvature concentration loci may still escape to infinity and therefore
fail to contribute points to the Kuratowski upper limit. An endwise Morrey estimate
for $q_i^+$, uniform in the mass, is one possible route to this
boundedness conclusion.

It is also unknown whether either strict inclusion
\[
    \mathcal S\subsetneq\mathcal Z
    \qquad\text{or}\qquad
    \mathcal Z\subsetneq\mathcal C
\]
can occur. The identity
\[
    \mathcal C\setminus\mathcal S
    =
    \mathcal O\setminus\mathcal S
\]
shows more precisely that
\[
    \mathcal C=\mathcal S
    \quad\Longleftrightarrow\quad
    \mathcal O\subset\mathcal S.
\]
This does not determine
whether the intermediate zero locus $\mathcal Z$ coincides with either
endpoint or lies strictly between them. Examples with
$\mathcal O\not\subset\mathcal S$, or with nonabelian concentration
escaping to infinity, would show that the corresponding alternatives in the
compactness theory are sharp.\\

\item
\textbf{Higgs zeros away from $\mathcal S$, scales centred at zeros, and instanton bubbling.}
The inclusion $\mathcal S\subset\mathcal Z$ follows from transverse
linking charge. The reverse inclusion asks whether every compact set
$K\Subset X\setminus\mathcal S$ is disjoint from $Z_i$ for all
sufficiently large $i$. This holds whenever
$\mathcal O\subset\mathcal S$, by
Corollary~\ref{cor:Hausdorff_convergence_unobstructed}. A point of
$\mathcal Z\setminus\mathcal S$, if nonempty, consists of limits of
Higgs zeros carrying no codimension-three charge detected by $T$ and is
contained in $\mathcal O\setminus\mathcal S$; see
Remark~\ref{rmk:charged_and_neutral_higgs_degenerations}.

The unconditional estimates for the varying loci obtained from Li's covers give a
partial size statement. For every $K\Subset X\setminus\mathcal S$,
Corollary~\ref{cor:Higgs_zero_Hausdorff_content_off_S} gives
\[
    \mathcal H^{n-3}_\infty(Z_i\cap K)\longrightarrow0,
\]
and Proposition~\ref{prop:Higgs_zero_density_rate_obstruction} shows
that, if the zero sets nevertheless become dense in a fixed ball, the
maximal distance to the zero set satisfies
\[
    h_i^{\mathcal Z}\gg m_i^{-1/4}.
\]
For the curvature concentration loci, the corresponding consequence of Li's
Vitali covers is
\[
    h_i^{\mathcal C}\gg m_i^{-\frac{n-3}{n}}.
\]
These conclusions constrain the curvature concentration loci but do not bound the
Kuratowski differences $\mathcal Z\setminus\mathcal S$ or
$\mathcal C\setminus\mathcal S$ without a common quantitative
approximation rate; see
Proposition~\ref{prop:Li_cover_rate_dependent_transfer} and
Remark~\ref{rmk:Li_content_does_not_transfer_to_Cinfty}.

Proposition~\ref{prop:zero_centered_TU_concentration} shows that every
$x\in\mathcal Z$, after passing to a subsequence realizing $x$ by zeros
$x_i\in Z_i$, carries a collapsing ordinary regularity scale
\[
    \lambda_i
    :=
    \mathfrak r_i(x_i)
    \longrightarrow0
\]
with
\[
    \lambda_i^{4-n}
    \int_{B_{\lambda_i}(x_i)}
    e_i
    =
    \varepsilon_0.
\]
Thus, after relabelling the subsequence along which the zero is realized, $x$ belongs to the
corresponding codimension-four concentration set $\mathcal B_{\mathrm{TU}}$. In particular, $\mathcal Z\setminus\mathcal S$ is naturally associated with concentration centred at Higgs zeros for the ordinary codimension-four
scale-invariant energy which is invisible to the limiting mass-renormalized measure.

The corresponding analysis centred at the characteristic radii is given in
Subsection~\ref{subsec:characteristic_scale_regimes}.  In particular,
Proposition~\ref{prop:regularity_scale_inside_characteristic_core}
shows that the ordinary regularity scale is at most Li's characteristic radius,
while
Proposition~\ref{prop:characteristic_core_above_mass_scale} shows that a
regularity scale above $m_i^{-1}$ forces annular concentration of the
codimension-four scale-invariant curvature energy.  Thus the
Green-weighted neck problem centred at Higgs zeros below and the
weighted curvature problem centred at the characteristic radii are distinct: the
former starts from the exact identity at a Higgs zero, whereas the latter
starts only from Li's curvature potential lower bound.

The remaining scale question is governed by
\[
    \tau_i
    :=
    m_i\lambda_i.
\]
If $\tau_i\to0$, Proposition~\ref{prop:zero_centered_scale_regimes}
shows that every smooth pointed limit has zero Higgs field and satisfies
the tangent space instanton equations. If a transverse dimension reduction
argument produces a four-dimensional normal limit and the tangential
directions decouple, the induced equation is anti-self-duality in the
standard associative or complex normal models. To identify
points of $\mathcal Z\setminus\mathcal S$ with the regime
$m_i\mathfrak r_i(x_i)\to0$, one would need to prove
\[
    x\in\mathcal Z\setminus\mathcal S
    \quad\Longrightarrow\quad
    m_i\mathfrak r_i(x_i)\longrightarrow0
\]
for every relabelled sequence of zeros $x_i\in Z_i$ with $x_i\to x$.

If this conclusion fails, one may pass to a further subsequence for which
\[
    \liminf_{i\to\infty}
    m_i\lambda_i
    >
    0.
\]
The relevant analytic tools are already available at the level of
stress-energy identities. By Corollary~\ref{cor:Euclidean_YMH_stress_energy_rigidity}, a nontrivial finite energy full-dimensional Yang--Mills--Higgs bubble on $\mathbb R^n$, $n=6,7$, cannot occur. A finite energy four-dimensional bubble centred at a Higgs zero has identically zero
Higgs field, while Corollary~\ref{cor:four_dimensional_Higgs_no_neck} gives the corresponding Higgs no-neck estimate on a genuine four-dimensional neck once the
energies on its two endpoint annuli tend to zero. Lemma~\ref{lem:annular_YMH_no_neck} gives the corresponding coercive ambient no-neck estimate in dimensions greater than four.

The Green representation isolates the remaining issue. At every
sequence centred at Higgs zeros,
\eqref{eq:zero_centered_Green_localization} shows that the Green-weighted
Higgs energy near the zero, normalized by the factor $2m_i^{-2}$, tends
to one. If
$m_i\lambda_i\to+\infty$, the blow-up at the mass scale is trivial and
\eqref{eq:zero_centered_Green_neck} places this entire weighted energy in
the expanding annulus between the mass scale and a fixed macroscopic
scale. A sufficient no-neck statement would be
\begin{equation}
\label{eq:weighted_Higgs_no_neck_target}
    \lim_{R\to\infty}
    \limsup_{i\to\infty}
    2m_i^{-2}
    \int_{B_\delta(x_i)\setminus B_{R/m_i}(x_i)}
    G(x_i,y)|\nabla_i\Phi_i|^2(y)\,\vol(y)
    =
    0
\end{equation}
whenever $x\notin\mathcal S$ and no nontrivial monopole bubble at the mass scale
occurs. This would contradict
\eqref{eq:zero_centered_Green_neck}.

More generally, the problem of separating the relevant scales would
follow from a
transverse bubble-tree alternative centred at Higgs zeros with the following
conclusions whenever
$\liminf_i m_i\mathfrak r_i(x_i)>0$:
\begin{enumerate}[label=\textnormal{(\roman*)}]
\item a three-dimensional monopole tangent forms and persists in the
tangential directions, forcing $x\in\mathcal S$;

\item a finite energy full-dimensional Yang--Mills--Higgs bubble forms;

\item a finite four-dimensional Yang--Mills bubble tree forms, the Higgs
energy vanishes on every neck, and the Higgs values match across
successive bubble regions.
\end{enumerate}
The first alternative would contradict $x\notin\mathcal S$, the second
is excluded by
Corollary~\ref{cor:Euclidean_YMH_stress_energy_rigidity}, and the third
would have to be combined with the local $L^2$ recovery of
$|\Phi_i|/m_i$ from Lemma~\ref{lem:largeness}, a quantitative selection of outer spheres with the required estimates,
and value matching across the necks. These inputs
would be incompatible with an inner bubble centred at a Higgs zero.

This discussion should be read together with Li's quantitative energy
identity \cite[Theorem~5.5]{li2025large}. Stated in the normalization
of the present paper, that result says that, on every fixed large
region, after discarding a subset of arbitrarily small
$\mathcal H^{n-3}$-measure in the smooth locus of the limiting
calibrated cycle, finitely many three-dimensional monopole bubbles at the mass scale
on the normal fibres account, up to an arbitrarily small error,
for the mass-renormalized curvature energy
$m_i^{-1}\int |F_{\nabla_i}|^2$. Thus any order-$m_i$ local curvature
energy in a fixed compact region is exhausted by monopole bubbling at the mass scale
along $\mathcal S$.

Consequently, at a point $x\in X\setminus\mathcal S$ the local bound
\begin{equation}\label{eq:local_L2_energy_bound}
    \int_{B_r(x)}e_i=o(m_i)
\end{equation}
rules out residual curvature concentration carrying energy of order
$m_i$. Any remaining codimension-four phenomenon centred at a Higgs zero at a point of
$\mathcal Z\setminus\mathcal S$ would therefore have to be sublinear in
$m_i$ and hence invisible to the mass-renormalized measures. The
remaining missing ingredients are more precise:
\begin{itemize}
\item a transverse dimension reduction and bubble extraction theorem under the local bound \eqref{eq:local_L2_energy_bound} available at points of $X\setminus\mathcal S$;

\item control of the tangential and mixed curvature components needed to
pass from the ambient equations to approximate Yang--Mills--Higgs
equations on normal slices;

\item a weighted Higgs no-neck estimate such as
\eqref{eq:weighted_Higgs_no_neck_target};

\item exclusion of collections of four-dimensional bubble regions whose
total energy is $o(m_i)$ and is therefore negligible in the
mass-renormalized energy; the order-$m_i$ part is already accounted for by
Li's monopole bubbling at the mass scale along $\mathcal S$;

\item tangential persistence of a three-dimensional bubble, sufficient
to produce positive mass-renormalized codimension-three density and
hence membership in $\mathcal S$.
\end{itemize}

The local charge on an outer normal sphere records the total transverse
charge in the enclosed ball and may split among several clusters. In the
three-dimensional theory, the complete moving centre decomposition shows
that the concentration multiplicity is the sum of the charges of all
mass one Euclidean profiles in the cluster and that no mass-renormalized
energy remains between those profiles and the ambient scale
\cite[Proposition~8.3]{fadeloliveira2026limitv5}. The following questions
ask for a transverse counterpart of that identity along the
positive-dimensional calibrated cycle:
\begin{itemize}
\item Can the local multiplicity of $T$ be recovered from a quantitative
charge decomposition of the nearby zero clusters?

\item Is the local multiplicity equal to the sum of the charges of all
finite energy Bogomolny monopoles on $\mathbb R^3$ arising from
transverse blow-up limits?

\item Under what hypotheses does nonzero linking charge persist
quantitatively over an $(n-3)$-dimensional family of nearby normal
slices?
\end{itemize}
A single charged linking sphere need not yield positive codimension-three
density for the mass-renormalized energy. The exact calibrated pairing
in Remark~\ref{rmk:charged_and_neutral_higgs_degenerations} shows,
however, that an affirmative answer to the last question, combined with
the topological energy lower bound on each charged fibre, would force the
limiting point to lie in $\mathcal S$. By contrast, if $x\in\mathcal C$
but $x\notin\mathcal Z$, then some neighbourhood of $x$ is disjoint from
the Higgs zero sets for all
sufficiently large indices. This is distinct from instanton bubbling
centred at Higgs zeros, which is relevant to
$\mathcal Z\setminus\mathcal S$. Global Hausdorff convergence requires, in addition,
a tightness estimate excluding zeros escaping along the AC end.

\begin{remark}[A codimension-four example without Higgs zeros for a larger structure group]
\label{rmk:larger_group_zero_free_instanton_core}
The simplest exact model of a codimension-four bubble whose energy
vanishes after normalization by $m_i^{-1}$ appears once the rank-one
structure group hypothesis is dropped.
It illustrates the analogue of
$\mathcal C\setminus\mathcal Z$, rather than
$\mathcal Z\setminus\mathcal S$.

Indeed, let $X=T^3\times\mathbb R^4$ carry the standard flat torsion-free $\mathrm G_2$-structure for which $T^3\times\{0\}$ is associative, and take
$G=\mathrm{SU}(2)\times\mathrm U(1)$. Let $\nabla_i$ be a charge-one BPST
instanton of scale $\lambda_i$ on the $\mathbb R^4$ factor, pulled back
to $X$, and set
\[
    \Phi_i=m_i\zeta,
\]
where $\zeta\neq0$ is constant in the commuting
$\mathfrak u(1)$ factor. Then
\[
    \nabla_i\Phi_i=0,
    \qquad
    F_{\nabla_i}\wedge\psi=0,
\]
so $(\nabla_i,\Phi_i)$ is an exact $\mathrm G_2$-monopole and
$Z_i=\varnothing$. If $m_i\to\infty$ and
$m_i\lambda_i\to0$, then on every compact subset the
mass-renormalized Yang--Mills--Higgs measures converge to zero. On the
other hand, for every $p\in T^3\times\{0\}$,
\[
    \int_{B_{m_i^{-1}}(p)}|F_{\nabla_i}|^2\,\operatorname{vol}
    \asymp
    m_i^{-3},
\]
and hence the characteristic potential
\eqref{eq:Li_potential_our_convention} remains bounded below by a
positive constant at $p$. Thus the curvature concentration along $T^3\times\{0\}$ remains
visible to the weighted potential, although its energy vanishes after
normalization by $m_i^{-1}$ and the Higgs field has no zeros.

This example is reducible, has larger structure group, and lies on a manifold which is not AC, so it is outside the {\mainsettingref}. Its significance is
structural. For $\mathrm{SU}(2)$ or $\mathrm{SO}(3)$, the centralizer of
a nonzero Higgs field is the abelian rank-one Cartan subalgebra, leaving
no commuting nonabelian factor in which to insert such an instanton.
For a larger group, a nonzero Higgs field may have nonabelian
centralizer, and the preceding construction becomes possible; compare
\cite[Question~4]{li2025large}.
\end{remark}
\end{enumerate}

\subsection{The abelian region and harmonic corrections}
\label{subsec:clearing_harmonic_questions}

\begin{enumerate}[resume]

\item
\textbf{The open set $\mathcal R_{\mathrm{ab}}$ and Li's harmonic corrections.}
The local smooth analysis of
Sections~\ref{sec: abelianization} and~\ref{sec:frequency_away_Cinfty}
applies on
\[
    \mathcal R_{\mathrm{ab}}
    =
    X\setminus(\mathcal S\cup\mathcal O),
\]
the maximal clearing region characterized by
Lemma~\ref{lem:maximal_abelian_clearing_region}. The present theory does not imply that this set is nonempty or dense. We know that it equals $X\setminus\mathcal S$ precisely when $\mathcal O\subset\mathcal S$. Thus, determining whether
$\mathcal R_{\mathrm{ab}}$ is nonempty or dense reduces to understanding
the size and location of $\mathcal O$.

Proposition~\ref{prop:smooth_convergence_Li_corrected_longitudinal}
shows that
\[
    \frac{1}{4\pi}\widehat\omega_i+\sigma_i
    \longrightarrow
    \widetilde F_\infty
\]
smoothly on compact subsets of $\mathcal R_{\mathrm{ab}}$, and
Proposition~\ref{prop:subcritical_harmonic_correction_growth} gives the
global estimate
\[
    \|\sigma_i\|_{L^2(X)}=o(m_i^{1/2}).
\]
Smooth compactness of the translated configurations follows if the
sequence
\[
    \sigma_i
    \in
    \mathscr H_{(2)}^2(X)
    \simeq
    H_c^2(X;\mathbb R)
\]
is uniformly bounded in $L^2(X)$. In that case,
Corollary~\ref{cor:harmonic_correction_alternatives} gives smooth
convergence on compact subsets of $\mathcal R_{\mathrm{ab}}$ to an
abelian $\Theta$-monopole. If instead
$\limsup_{i\to\infty}\|\sigma_i\|_{L^2(X)}=+\infty$, pass to the
subsequence selected in the second alternative of that corollary and
write
\[
    b_i:=\|\sigma_i\|_{L^2(X)}\longrightarrow+\infty.
\]
Then $b_i=o(m_i^{1/2})$, and the forms
$(4\pi b_i)^{-1}\widehat\omega_i$ converge smoothly on compact subsets
of $\mathcal R_{\mathrm{ab}}$ to the restriction of the nonzero global
$L^2$-harmonic form $-\sigma_\infty$.

Under the identification
$\mathscr H_{(2)}^2(X)\simeq H_c^2(X;\mathbb R)$, Li's corrections may be
chosen so that their cohomology classes lie in the integral lattice in
the $\mathrm{SU}(2)$ case, with the corresponding half-integral
modification for $\mathrm{SO}(3)$. In the bounded alternative, such a
choice is eventually constant after passing to a subsequence; see
Remark~\ref{rmk:lattice_corrections_bounded_alternative}. It remains to
determine whether the fixed principal bundle, its second characteristic class, the fixed asymptotic monopole class, or additional geometric hypotheses force these lattice classes to remain bounded. Equivalently, can
\[
    \limsup_{i\to\infty}
    \|\sigma_i\|_{L^2(X)}
    =
    +\infty
\]
for a large mass sequence with fixed topological data? Li notes that no
example is known. The condition $H_c^2(X;\mathbb R)=0$ rules this out,
but no general boundedness criterion is currently available.

\item
\textbf{Quantitative behaviour near the zero set of the limiting
$L^2$-harmonic form.}
If the second alternative of
Corollary~\ref{cor:harmonic_correction_alternatives} occurs, pass to the
subsequence selected there and write
$b_i:=\|\sigma_i\|_{L^2(X)}\to+\infty$. One has $b_i=o(m_i^{1/2})$.
The corollary gives a nonzero form
$\sigma_\infty\in\mathscr H_{(2)}^2(X)$ such that the forms
$(4\pi b_i)^{-1}\widehat\omega_i$ converge smoothly on compact subsets
of $\mathcal R_{\mathrm{ab}}$ to $-\sigma_\infty$. Thus the normalized
local limits in Section~\ref{sec:frequency_away_Cinfty} are restrictions
of the same global $L^2$-harmonic form and cannot vanish identically on a
compact domain with nonempty interior.

Where this harmonic form is bounded away from zero,
Proposition~\ref{prop:regularity_scale_nonzero_harmonic_profile} and
Corollary~\ref{cor:regularity_scale_mass_comparison} give
\[
    \mathfrak r_i\asymp b_i^{-1/2},
    \qquad
     m_i^{1/4}\mathfrak r_i\longrightarrow+\infty,
\]
up to fixed factors determined by the local $L^2$-norm of
$\sigma_\infty$. Near a zero,
Proposition~\ref{prop:vanishing_order_regular_scale} gives the exponent
$1/(q+2)$ only under the quantitative order-$q$
profile condition of
Definition~\ref{def:quantitative_order_q_profile}. Under this
condition,
\[
    m_i^{1/(2(q+2))}\mathfrak r_i(x)
    \longrightarrow
    +\infty.
\]
Although unique continuation excludes infinite-order vanishing for a
nonzero harmonic form, finite vanishing order on fixed scales does not by
itself imply the required profile at shrinking scales.

It remains to derive the quantitative profile condition from the smooth
convergence of the corrected forms and additional estimates for the error
$\widehat\omega_i/(4\pi)+\sigma_i-\widetilde F_\infty$, or to identify
geometric hypotheses under which the required control holds uniformly
near the nodal set of $\sigma_\infty$.
\end{enumerate}


\appendix

\section{Mean value inequalities}\label{app: A}

The estimates in this appendix are purely local. They depend only on bounded geometry of the background Riemannian manifold. In what follows, we say a Riemannian manifold $(X^n,g)$ is of \textbf{bounded geometry} if the following conditions hold:
\begin{itemize}
	\item The global injectivity radius of $(X,g)$ is positive: $\inj_g(X)>0$ (in particular, $(X^n,g)$ is complete);
	\item The Riemann curvature tensor $\mathrm{Rm}$ of $(X,g)$, together with its covariant derivatives satisfy uniform $L^{\infty}$-bounds on $X$: for all $j\geqslant 0$,
	\[
	\| D^j\mathrm{Rm}\|_{L^{\infty}(X)} \lesssim_j 1.
	\] 
\end{itemize}Here and henceforth, we use the notation $a\lesssim b$ to mean $a\leqslant cb$, for some constant $c>0$ depending only on the geometry of $(X^n,g)$; when we add a subscript like in $a\lesssim_j b$, we mean that the hidden constant $c$ depends, furthermore, on the subscript $j$.

To begin, we need the following standard result.
\begin{lemma}\label{lem: mean_value}
	Let $(X^n,g)$ be a Riemannian $n$-manifold of bounded geometry. Then there is $0<\delta<\inj_g(X)$ with the following significance. If $p\in X$, $r\in (0,\delta]$ and $C\geqslant 0$, then every $f\in C^2(B_r(p),[0,\infty))$ satisfies:
	\begin{equation}\label{eq: first_mean_value}
		\Delta f\leqslant C\quad\Longrightarrow\quad f(p)\lesssim Cr^2 + r^{-n}\int_{B_r(p)} f.
	\end{equation}
\end{lemma}
\begin{proof}
    For a proof adapted from \cite[Step 2 in the Proof of Theorem
B.1]{hohloch2009hypercontact}, see \cite[\S 3.3]{fadel2026ebook}.
\end{proof}
The following general nonlinear mean value theorem combines various results of this kind appearing in \cite{wehrheim2005energy}, \cite[Appendix B]{hohloch2009hypercontact}, and \cite[Appendix A]{walpuski2017compactness}.
\begin{theorem}[Mean value inequalities]\label{thm: gen_mean_value}
	Under the hypotheses of Lemma~\ref{lem: mean_value}, let $x\in X$, $r\in (0,\delta]$ and suppose $f\in C^2(B_r(x),[0,\infty))$ satisfies the following conditions:
	\begin{itemize}
		\item[(C1)] There is $d\in\mathbb{N}$ such that if $B_s(y)\subseteq B_{r/2}(x)$ then
		\begin{equation}\label{eq: mon_ass}
			s^{d-n}\int_{B_s(y)} f\lesssim r^{d-n}\int_{B_r(x)}f + \tau,
		\end{equation} for some $\tau=\tau(r)\geqslant 0$.
		\item[(C2)] There are constants $a_0,a_1,a\geqslant 0$ and $\alpha\in [1,(d+2)/d]$ such that
		\begin{equation}\label{eq: boch_ass}
			\Delta f\leqslant a_0 + a_1 f + af^{\alpha}.
		\end{equation}
	\end{itemize} Setting
	\[
	\varepsilon := r^{d-n}\int_{B_r(x)}f,
	\] we have:
	\begin{itemize}
		\item[(i)] If $\alpha<(d+2)/d$ then
		\begin{equation}\label{eq: subcritical}
			\sup_{B_{\frac{r}{4}}(x)} f \lesssim_d a_0r^2 + (a_1^{d/2} + r^{-d})(\varepsilon + \tau) + \left(a^{d/2}(\varepsilon + \tau)\right)^{\gamma_{\alpha,d}}, 
		\end{equation} where $\gamma_{\alpha,d}:= 2/(2+d-\alpha d)$.
		\item[(ii)] If $\alpha = (d+2)/d$, then there is a constant $\hbar>0$ (depending only on the geometry and $d$) such that
		\begin{equation}\label{eq: critical}
			a^{d/2}(\varepsilon+\tau)<\hbar\quad\Longrightarrow\quad\sup_{B_{\frac{r}{4}}(x)} f\lesssim_d a_0 r^2 + (a_1^{d/2} + r^{-d})(\varepsilon + \tau)
		\end{equation}
	\end{itemize}
\end{theorem}
\begin{proof}[Proof of Theorem~\ref{thm: gen_mean_value}]
The proof is based on the so-called `Heinz trick'. In order for the
reader to see precisely where each relevant constant appears in the
estimates, in this proof we mostly avoid the notation `$\lesssim$' and
show the constants explicitly. Thus, choose $\delta>0$ as in
Lemma~\ref{lem: mean_value}, let $c_1>0$ be the hidden constant in
\eqref{eq: first_mean_value}, and let $c_2>0$ be the hidden constant in
\eqref{eq: mon_ass}. Both $\hbar$ and the final estimates depend only on
$c_1$, $c_2$, and $d$.

Fix $y\in B_{r/4}(x)$. If $f(y)=0$, there is nothing to prove.
Otherwise, define $h:[0,r/4]\to[0,\infty)$ by
\[
    h(s)
    :=
    \left(\frac{r/4-s}{r/4}\right)^d
    \max_{\overline{B_s(y)}}f.
\]
Since $h(0)=f(y)>0$ and $h(r/4)=0$, there exist
$s^\ast\in[0,r/4)$ and
$y^\ast\in\overline{B_{s^\ast}(y)}$ such that
\[
    h(s^\ast)=\max_{0\leqslant s\leqslant r/4}h(s),
    \qquad
    F:=f(y^\ast)=\max_{\overline{B_{s^\ast}(y)}}f.
\]
Set
\[
    s_0:=\frac{r/4-s^\ast}{2}>0.
\]
Then
\[
    \max_{\overline{B_{s_0}(y^\ast)}}f
    \leqslant
    \max_{\overline{B_{s^\ast+s_0}(y)}}f
    \leqslant
    2^dF.
\]
Hence, by \eqref{eq: boch_ass}, on $B_{s_0}(y^\ast)$ one has
\[
    \Delta f
    \leqslant
    a_0+a_1(2^dF)+a(2^dF)^\alpha.
\]
Lemma~\ref{lem: mean_value} therefore gives
\begin{equation}\label{eq: mv_0}
    F
    \leqslant
    c_1\left(
    (a_0+a_1(2^dF)+a(2^dF)^\alpha)s^2
    +s^{-n}\int_{B_s(y^\ast)}f
    \right)
\end{equation}
for every $0<s\leqslant s_0$. Moreover,
\[
    d(x,y^\ast)+s
    \leqslant
    \frac r4+s^\ast+\frac{r/4-s^\ast}{2}
    =
    \frac{3r}{8}+\frac{s^\ast}{2}
    <
    \frac r2
\]
for every $0<s\leqslant s_0$. Thus
$B_s(y^\ast)\subset B_{r/2}(x)$, and condition~\textup{(C1)} applies.
Consequently,
\begin{equation}\label{eq: mv_1}
    F
    \leqslant
    c_1a_0r^2/4
    +c_1\left(a_1(2^dF)+a(2^dF)^\alpha\right)s^2
    +c_1c_2s^{-d}(\varepsilon+\tau),
    \qquad 0<s\leqslant s_0.
\end{equation}
We now make a case by case distinction.\\
	
	\textbf{Case 1.} $F\leqslant c_1a_0r^2$\\
	
	In this case $f(y)\leqslant F\leqslant c_1a_0 r^2\lesssim a_0r^2$, which already proves the assertion.\\
	
	\textbf{Case 2.} $F\geqslant c_1a_0r^2$ and $a_1 (2^d F)\geqslant a(2^d F)^{\alpha}$.\\
	
	From \eqref{eq: mv_1} we derive
	\begin{equation}\label{eq: mv_2}
		F\leqslant F/4 + 2c_1 a_1 2^d F s^2 + c_1c_2s^{-d}(\varepsilon+\tau),\quad\forall\,0<s\leqslant s_0.
	\end{equation} so that we have two possibilities:
	\begin{itemize}
		\item If $c_1 a_1 2^d s_0^2<1/8$, then by \eqref{eq: mv_2} we get
		\[
		F\leqslant 2c_1c_2s_0^{-d}(\varepsilon+\tau).
		\]
Since $r/4-s^{\ast}=2s_0$, this gives
\[
f(y)=h(0)\leqslant h(s^{\ast})
=
\left(\frac{r/4-s^{\ast}}{r/4}\right)^d F
=
\frac{2^{3d}s_0^d}{r^d}F\leqslant 2^{3d+1}c_1c_2r^{-d}(\varepsilon+\tau)\lesssim_d r^{-d}(\varepsilon+\tau),
		\]  and so the desired estimate holds.
		\item Otherwise, we choose $s\leqslant s_0$ such that $c_1 a_1 2^d s^2 = 1/8$. Then by \eqref{eq: mv_2} we have
		\[
		f(y)\leqslant F\leqslant 2c_1c_2s^{-d}(\varepsilon+\tau) = 2c_1c_2(8c_1 a_1 2^d)^{d/2}(\varepsilon+\tau)\lesssim_d a_1^{d/2}(\varepsilon+\tau),
		\] as we wanted.
	\end{itemize}
	
	\textbf{Case 3.} $F\geqslant c_1a_0r^2$ and $a_1 (2^d F)\leqslant a(2^d F)^{\alpha}$.\\	
	
	Here from \eqref{eq: mv_1} we get
	\begin{equation}\label{eq: mv_3}
		F\leqslant F/4 + 2c_1 a(2^d F)^{\alpha} s^2 + c_1c_2s^{-d}(\varepsilon+\tau),\quad\forall\,0<s\leqslant s_0.
	\end{equation}
Hence we have the following possibilities:
	\begin{itemize}
		\item If $c_1 a 2^{\alpha d} F^{\alpha-1} s_0^2<1/8$ then by \eqref{eq: mv_3}, 
		\[
		F\leqslant 2c_1c_2 s_0^{-d}(\varepsilon+\tau),
		\] and again by the same argument as above we get
		\[
		f(y)\leqslant h(s^{\ast})\leqslant 2^{3d+1}c_1c_2r^{-d}(\varepsilon+\tau)\lesssim_d r^{-d}(\varepsilon+\tau),
		\] as required. 
		\item Otherwise, we choose $s\leqslant s_0$ such that $c_1 a 2^{\alpha d} F^{\alpha-1} s^2 = 1/8$ so that
		\begin{equation}\label{eq: mv_3'}
			F\leqslant 2c_1c_2 s^{-d}(\varepsilon+\tau) = 2c_1c_2 (c_1 2^{\alpha d + 3})^{d/2} a^{d/2} F^{(\alpha d - d)/2}(\varepsilon+\tau).
		\end{equation}
If $\alpha<(d+2)/d$, then $2+d-\alpha d>0$, and hence
		\[
		f(y)\leqslant F\leqslant \left(2c_1c_2 (c_1 2^{\alpha d + 3})^{d/2}a^{d/2}(\varepsilon+\tau)\right)^{\gamma_{\alpha,d}}\lesssim_d \left(a^{d/2}(\varepsilon+\tau)\right)^{\gamma_{\alpha,d}},
		\]
		as required. For the critical exponent $\alpha=(d+2)/d$, one has
		$(\alpha d-d)/2=1$, so \eqref{eq: mv_3'} implies
		$1\leqslant c a^{d/2}(\varepsilon+\tau)$ for some $c>0$ depending
		only on $c_1,c_2$, and $d$. Thus, if
		$a^{d/2}(\varepsilon+\tau)<\hbar:=c^{-1}$, this case is excluded.
		We are done.
	\end{itemize}
\end{proof}


\section{Stress-energy identities and monotonicity}\label{app: B}

In this appendix, we collect the local stress-energy identities used in
the paper and several consequences relevant to the bubbling questions in
Section~\ref{sec:resulting_picture}. We begin with the monotonicity
formula for the full Yang--Mills--Higgs energy of a solution of
\eqref{eq:YMH_second_order}. It is a convenient reformulation of
\cite[Theorem~2.1]{afuni2019regularity}, adapted to the local geometric
framework used throughout the paper. We then record its standard
codimension-four almost-monotonicity consequence, annular estimates, and
the Euclidean stress-energy identity and finite energy rigidity statement.
Finally, we derive a localized monotonicity identity for the Higgs energy
component alone.

In what follows, recall the definition \eqref{eq: YMH_density} of the
energy density $e(\nabla,\Phi)$.

\subsection{Full stress-energy identity and codimension-four monotonicity}
\label{subsec:full_stress_energy_monotonicity}

\begin{theorem}[Monotonicity formula with error terms]\label{thm: monotonicity_crit_pt}
Let $(X^n,g)$ be an oriented Riemannian $n$-manifold, $n\geqslant4$, and let
$P$ be a principal $G$-bundle over $X$, where $G$ is compact. Let
$x\in X$, set $r_x:=d(x,\cdot)$, and assume that on $B_{r_0}(x)$ one has
\[
-C_H r_x^2g
\leqslant
\mathrm{Hess}\,\!\left(\frac12r_x^2\right)-g
\leqslant
C_H r_x^2g
\]
for some $r_0\in(0,\inj_g(x))$ and $C_H\geqslant0$. If
$(\nabla,\Phi)\in\mathscr{A}(P)\times\Gamma(\mathfrak{g}_P)$ solves \eqref{eq:YMH_second_order}, then
there exists a constant $c\geqslant0$, depending only on $n$ and $C_H$, such
that
\[
\theta_x(r)
:=
e^{cr^2}r^{3-n}
\int_{B_r(x)}e(\nabla,\Phi),
\]
satisfies, for every $0<r<r_0$,
\begin{align}
\frac{d}{dr}\theta_x(r)
&\geqslant
2e^{cr^2}r^{3-n}
\int_{\partial B_r(x)}
\left(
|\partial_{r_x}\lrcorner F_\nabla|^2
+
|\nabla_{\partial_{r_x}}\Phi|^2
\right)
\nonumber\\
&\quad
-
e^{cr^2}r^{2-n}
\int_{B_r(x)}
\bigl(|F_\nabla|^2-|\nabla\Phi|^2\bigr).
\label{eq: monotonicity}
\end{align}
\end{theorem}

\begin{proof}
Consider the stress-energy tensor
\[
\mathbb{T}(v,w)
=
g(v,w)e(\nabla,\Phi)
-
2\langle v\lrcorner F_\nabla,w\lrcorner F_\nabla\rangle
-
2\langle\nabla_v\Phi,\nabla_w\Phi\rangle.
\]
Its trace is
\[
\langle \mathbb{T},g\rangle
=
(n-3)e(\nabla,\Phi)
-
\bigl(|F_\nabla|^2-|\nabla\Phi|^2\bigr).
\]
Moreover, there is a dimensional constant $C_n>0$ such that
\[
    |\mathbb{T}|\leqslant C_ne(\nabla,\Phi).
\]
Since $(\nabla,\Phi)$ solves \eqref{eq:YMH_second_order}, the stress-energy
tensor is divergence-free, $D^*\mathbb{T}=0$. Hence, for every precompact domain
$\Omega\subset X$ and every vector field $Y$ on $\overline\Omega$,
\begin{equation}\label{eq: div_thm_T}
\int_{\partial\Omega}\mathbb{T}(Y,\mathbf n)
=
\int_\Omega\langle \mathbb{T},DY\rangle,
\end{equation} where $\mathbf n$ denotes the outward unit normal to $\partial\Omega$.

Apply this to $\Omega=B_r(x)$ and
\[
    Y=D\left(\frac12r_x^2\right)=r_x\partial_{r_x}.
\]
Since $\mathbf n=\partial_{r_x}$ on $\partial B_r(x)$, we obtain
\begin{align}
\int_{\partial B_r(x)}e(\nabla,\Phi)
&=
2\int_{\partial B_r(x)}
\bigl(
|\partial_{r_x}\lrcorner F_\nabla|^2
+
|\nabla_{\partial_{r_x}}\Phi|^2
\bigr)
\nonumber\\
&\quad
+
r^{-1}
\int_{B_r(x)}
\left\langle
\mathbb{T},
\mathrm{Hess}\,\!\left(\frac12r_x^2\right)-g
\right\rangle
\nonumber\\
&\quad
+
r^{-1}
\int_{B_r(x)}
\left[
(n-3)e(\nabla,\Phi)
-
\bigl(|F_\nabla|^2-|\nabla\Phi|^2\bigr)
\right].
\end{align}
By the Hessian assumption and the bound $|\mathbb{T}|\leqslant C_ne(\nabla,\Phi)$,
there is a constant $C_1=C_1(n)C_H$ such that
\[
\left\langle
\mathbb{T},
\mathrm{Hess}\,\!\left(\frac12r_x^2\right)-g
\right\rangle
\geqslant
-C_1r_x^2e(\nabla,\Phi).
\]
Therefore
\begin{align}
\int_{\partial B_r(x)}e(\nabla,\Phi)
&\geqslant
2\int_{\partial B_r(x)}
\bigl(
|\partial_{r_x}\lrcorner F_\nabla|^2
+
|\nabla_{\partial_{r_x}}\Phi|^2
\bigr)
\nonumber\\
&\quad
+
r^{-1}
\int_{B_r(x)}
\bigl((n-3)-C_1r^2\bigr)e(\nabla,\Phi)
\nonumber\\
&\quad
-
r^{-1}
\int_{B_r(x)}
\bigl(|F_\nabla|^2-|\nabla\Phi|^2\bigr).
\label{ineq: der_energy}
\end{align}
By the coarea formula,
\[
\frac{d}{dr}
\int_{B_r(x)}e(\nabla,\Phi)
=
\int_{\partial B_r(x)}e(\nabla,\Phi).
\]
Thus, if
\[
    E_x(r):=\int_{B_r(x)}e(\nabla,\Phi),
\]
then
\[
\frac{d}{dr}
\left(
e^{cr^2}r^{3-n}E_x(r)
\right)
=
\bigl(2cr+(3-n)r^{-1}\bigr)
e^{cr^2}r^{3-n}E_x(r)
+
e^{cr^2}r^{3-n}E_x'(r).
\]
Combining this identity with \eqref{ineq: der_energy} gives
\begin{align*}
\frac{d}{dr}
\left(
e^{cr^2}r^{3-n}E_x(r)
\right)
&\geqslant
2e^{cr^2}r^{3-n}
\int_{\partial B_r(x)}
\bigl(
|\partial_{r_x}\lrcorner F_\nabla|^2
+
|\nabla_{\partial_{r_x}}\Phi|^2
\bigr)
\\
&\quad
+
(2c-C_1)r\,e^{cr^2}r^{3-n}E_x(r)
\\
&\quad
-
e^{cr^2}r^{2-n}
\int_{B_r(x)}
\bigl(|F_\nabla|^2-|\nabla\Phi|^2\bigr).
\end{align*}
Choosing
\[
    c\geqslant \frac{C_1}{2}
\]
and discarding the nonnegative middle term yields
\eqref{eq: monotonicity}.
\end{proof}

\begin{corollary}[Codimension-four almost-monotonicity]
\label{cor:codim4_almost_monotonicity}
Under the hypotheses of Theorem~\ref{thm: monotonicity_crit_pt}, with the
same constant $c$ as in that theorem, the function
\[
    r\longmapsto
    e^{cr^2}r^{4-n}
    \int_{B_r(x)}
    e(\nabla,\Phi)
\]
is nondecreasing on $(0,r_0)$. More precisely,
\begin{align}
\frac{d}{dr}
\left(
e^{cr^2}r^{4-n}
\int_{B_r(x)}
e(\nabla,\Phi)
\right)
&\geqslant
2e^{cr^2}r^{4-n}
\int_{\partial B_r(x)}
\left(
|\partial_{r_x}\lrcorner F_\nabla|^2
+
|\nabla_{\partial_{r_x}}\Phi|^2
\right)
\nonumber\\
&\quad
+
2e^{cr^2}r^{3-n}
\int_{B_r(x)}
|\nabla\Phi|^2
\geqslant 0.
\label{eq:codim4_almost_monotonicity}
\end{align}
\end{corollary}

\begin{proof}
Let
\[
    E_x(r):=
    \int_{B_r(x)}
    e(\nabla,\Phi).
\]
From the proof of Theorem~\ref{thm: monotonicity_crit_pt}, we have
\begin{align*}
E_x'(r)
&\geqslant
2\int_{\partial B_r(x)}
\left(
|\partial_{r_x}\lrcorner F_\nabla|^2
+
|\nabla_{\partial_{r_x}}\Phi|^2
\right)
\\
&\quad
+
r^{-1}
\int_{B_r(x)}
\bigl((n-3)-C_1r^2\bigr)e(\nabla,\Phi)
\\
&\quad
-
r^{-1}
\int_{B_r(x)}
\bigl(|F_\nabla|^2-|\nabla\Phi|^2\bigr),
\end{align*}
where $C_1$ is the geometric constant appearing in
\eqref{ineq: der_energy}. By the choice of the constant $c$ in
Theorem~\ref{thm: monotonicity_crit_pt}, we have $2c-C_1\geqslant0$. Equivalently,
\begin{align*}
E_x'(r)
&\geqslant
2\int_{\partial B_r(x)}
\left(
|\partial_{r_x}\lrcorner F_\nabla|^2
+
|\nabla_{\partial_{r_x}}\Phi|^2
\right)
\\
&\quad
+
r^{-1}
\int_{B_r(x)}
\bigl((n-4)-C_1r^2\bigr)|F_\nabla|^2
+
\bigl((n-2)-C_1r^2\bigr)|\nabla\Phi|^2.
\end{align*}

Now differentiate
\[
    e^{cr^2}r^{4-n}E_x(r).
\]
We obtain
\[
\frac{d}{dr}
\left(
e^{cr^2}r^{4-n}E_x(r)
\right)
=
e^{cr^2}r^{4-n}E_x'(r)
+
\left(2cr+(4-n)r^{-1}\right)
e^{cr^2}r^{4-n}E_x(r).
\]
Substituting the previous lower bound for $E_x'(r)$ gives
\begin{align*}
\frac{d}{dr}
\left(
e^{cr^2}r^{4-n}E_x(r)
\right)
&\geqslant
2e^{cr^2}r^{4-n}
\int_{\partial B_r(x)}
\left(
|\partial_{r_x}\lrcorner F_\nabla|^2
+
|\nabla_{\partial_{r_x}}\Phi|^2
\right)
\\
&\quad
+
e^{cr^2}r^{3-n}
\int_{B_r(x)}
(2c-C_1)r^2|F_\nabla|^2
\\
&\quad
+
e^{cr^2}r^{3-n}
\int_{B_r(x)}
\left(
2+(2c-C_1)r^2
\right)|\nabla\Phi|^2.
\end{align*}
Since $c\geqslant C_1/2$, the curvature contribution is nonnegative and
\[
    2+(2c-C_1)r^2\geqslant 2.
\]
Therefore
\begin{align*}
\frac{d}{dr}
\left(
e^{cr^2}r^{4-n}E_x(r)
\right)
&\geqslant
2e^{cr^2}r^{4-n}
\int_{\partial B_r(x)}
\left(
|\partial_{r_x}\lrcorner F_\nabla|^2
+
|\nabla_{\partial_{r_x}}\Phi|^2
\right)
\\
&\quad
+
2e^{cr^2}r^{3-n}
\int_{B_r(x)}
|\nabla\Phi|^2 ,
\end{align*}
which is \eqref{eq:codim4_almost_monotonicity}.
\end{proof}

\subsection{Annular estimates and Euclidean finite energy rigidity}
\label{subsec:annular_stress_energy_rigidity}

\begin{lemma}[Annular stress-energy estimate]
\label{lem:annular_YMH_no_neck}
Let $(X^n,g)$ be an oriented Riemannian $n$-manifold, $n>4$, and let
$(\nabla,\Phi)$ satisfy the hypotheses of
Theorem~\ref{thm: monotonicity_crit_pt}. After replacing $r_0$ by a
smaller positive constant, depending only on $n$, $C_H$, and the original
$r_0$, the following holds. If
\[
    0<4\sigma<\rho\leqslant r_0,
\]
then there exist radii
\[
    \widehat\sigma\in[\sigma,2\sigma],
    \qquad
    \widehat\rho\in[\rho/2,\rho]
\]
such that
\begin{equation}
\label{eq:annular_YMH_no_neck}
\begin{aligned}
    \int_{B_{\widehat\rho}(x)\setminus B_{\widehat\sigma}(x)}
    e(\nabla,\Phi)
    \leqslant
    C\Bigg(
        &\int_{B_{2\sigma}(x)\setminus B_\sigma(x)}
        e(\nabla,\Phi)
        \\
        &+
        \int_{B_\rho(x)\setminus B_{\rho/2}(x)}
        e(\nabla,\Phi)
    \Bigg),
\end{aligned}
\end{equation}
where $C<\infty$ depends only on $n$ and $C_H$.

Consequently, suppose that $x_i\in X$ are points for which the Hessian hypothesis of Theorem~\ref{thm: monotonicity_crit_pt} holds on $B_{r_0}(x_i)$ with the same constants $C_H$ and $r_0$. If $(\nabla_i,\Phi_i)$ is a sequence of solutions, set
$e_i:=e(\nabla_i,\Phi_i)$. If $0<4\sigma_i<\rho_i\leqslant r_0$ and $a_i>0$ satisfy
\[
    a_i
    \int_{B_{2\sigma_i}(x_i)\setminus B_{\sigma_i}(x_i)}
    e_i
    \longrightarrow0,
    \qquad
    a_i
    \int_{B_{\rho_i}(x_i)\setminus B_{\rho_i/2}(x_i)}
    e_i
    \longrightarrow0,
\]
then the good radii may be chosen so that
\[
    a_i
    \int_{B_{\widehat\rho_i}(x_i)\setminus
    B_{\widehat\sigma_i}(x_i)}
    e_i
    \longrightarrow0.
\]
\end{lemma}

\begin{proof}
Set
\[
    Y
    :=
    D\left(\frac12r_x^2\right)
    =
    r_x\partial_{r_x}
\]
and define
\[
    Q(s)
    :=
    s\int_{\partial B_s(x)}
    \left(
        e(\nabla,\Phi)
        -
        2|\partial_{r_x}\lrcorner F_\nabla|^2
        -
        2|\nabla_{\partial_{r_x}}\Phi|^2
    \right).
\]
Applying \eqref{eq: div_thm_T} to the annulus
$B_b(x)\setminus\overline{B_a(x)}$ gives
\begin{align}
    Q(b)-Q(a)
    &=
    \int_{B_b(x)\setminus B_a(x)}
    \left(
        (n-4)|F_\nabla|^2
        +
        (n-2)|\nabla\Phi|^2
    \right)
    \nonumber\\
    &\quad
    +
    \int_{B_b(x)\setminus B_a(x)}
    \left\langle
        \mathbb T,
        \mathrm{Hess}\,\!\left(\frac12r_x^2\right)-g
    \right\rangle.
\label{eq:annular_stress_energy_identity}
\end{align}
Since $|\mathbb T|\leqslant C_ne(\nabla,\Phi)$ and
$r_x\leqslant\rho$ on the annulus under consideration,
\[
    \left|
    \int_{B_b(x)\setminus B_a(x)}
    \left\langle
        \mathbb T,
        \mathrm{Hess}\,\!\left(\frac12r_x^2\right)-g
    \right\rangle
    \right|
    \leqslant
    C_nC_H\rho^2
    \int_{B_b(x)\setminus B_a(x)}
    e(\nabla,\Phi).
\]
After decreasing $r_0$ so that
\[
    C_nC_Hr_0^2
    \leqslant
    \frac{n-4}{2},
\]
equation \eqref{eq:annular_stress_energy_identity} gives
\[
    \int_{B_b(x)\setminus B_a(x)}
    e(\nabla,\Phi)
    \leqslant
    C\bigl(|Q(b)|+|Q(a)|\bigr).
\]
Moreover,
\[
    |Q(s)|
    \leqslant
    Cs\int_{\partial B_s(x)}e(\nabla,\Phi).
\]
By the coarea formula, there exists
$\widehat\rho\in[\rho/2,\rho]$ such that
\[
    \widehat\rho
    \int_{\partial B_{\widehat\rho}(x)}
    e(\nabla,\Phi)
    \leqslant
    2
    \int_{B_\rho(x)\setminus B_{\rho/2}(x)}
    e(\nabla,\Phi).
\]
Similarly, there exists
$\widehat\sigma\in[\sigma,2\sigma]$ such that
\[
    \widehat\sigma
    \int_{\partial B_{\widehat\sigma}(x)}
    e(\nabla,\Phi)
    \leqslant
    2
    \int_{B_{2\sigma}(x)\setminus B_\sigma(x)}
    e(\nabla,\Phi).
\]
Substituting these two estimates proves
\eqref{eq:annular_YMH_no_neck}. The sequential conclusion follows by
multiplying this estimate by $a_i$.
\end{proof}

\begin{corollary}[Four-dimensional Higgs no-neck estimate]
\label{cor:four_dimensional_Higgs_no_neck}
Let $(A,\phi)$ be a smooth solution of \eqref{eq:YMH_second_order}
on a Euclidean ball $B_\rho(0)\subset\mathbb R^4$, and suppose that
$0<4\sigma<\rho$. Then there exist
\[
    \widehat\sigma\in[\sigma,2\sigma],
    \qquad
    \widehat\rho\in[\rho/2,\rho]
\]
such that
\begin{equation}
\label{eq:four_dimensional_Higgs_no_neck}
\begin{aligned}
    \int_{B_{\widehat\rho}(0)\setminus B_{\widehat\sigma}(0)}
    |\nabla_A\phi|^2
    \leqslant
    C\Bigg(
        &\int_{B_{2\sigma}(0)\setminus B_\sigma(0)}
        e(A,\phi)
        \\
        &+
        \int_{B_\rho(0)\setminus B_{\rho/2}(0)}
        e(A,\phi)
    \Bigg).
\end{aligned}
\end{equation}
\end{corollary}

\begin{proof}
In Euclidean dimension four, the geometric error in
\eqref{eq:annular_stress_energy_identity} vanishes and its bulk term is
\[
    2\int_{B_b(0)\setminus B_a(0)}
    |\nabla_A\phi|^2.
\]
Thus
\[
    2\int_{B_b(0)\setminus B_a(0)}
    |\nabla_A\phi|^2
    =
    Q(b)-Q(a).
\]
The same radius selection argument used in
Lemma~\ref{lem:annular_YMH_no_neck} gives
\eqref{eq:four_dimensional_Higgs_no_neck}.
\end{proof}

\begin{corollary}[Euclidean stress-energy identity and finite energy rigidity]
\label{cor:Euclidean_YMH_stress_energy_rigidity}
Let $(A,\phi)$ be a smooth solution of \eqref{eq:YMH_second_order}
on $\mathbb R^d$, $d\geqslant3$, such that
\[
    \int_{\mathbb R^d}
    \left(
        |F_A|^2+|\nabla_A\phi|^2
    \right)
    <
    \infty.
\]
Then
\begin{equation}
\label{eq:Euclidean_YMH_stress_energy}
    (d-4)
    \int_{\mathbb R^d}|F_A|^2
    +
    (d-2)
    \int_{\mathbb R^d}|\nabla_A\phi|^2
    =
    0.
\end{equation}
In dimension $d=3$, this is the equipartition identity
\[
    \int_{\mathbb R^3}|F_A|^2
    =
    \int_{\mathbb R^3}|\nabla_A\phi|^2.
\]
If $d=4$, then $\nabla_A\phi=0$ and $A$ is Yang--Mills; if moreover
$\phi$ vanishes at one point, then $\phi\equiv0$. If $d>4$, then
$F_A=0$ and $\nabla_A\phi=0$.
\end{corollary}

\begin{proof}
The stress-energy computation leading to \eqref{eq:annular_stress_energy_identity} is valid in every dimension. In Euclidean space, where the Hessian error vanishes, its ball version gives
\begin{align*}
    R\int_{\partial B_R(0)}
    \left(
        e(A,\phi)
        -
        2|\partial_r\lrcorner F_A|^2
        -
        2|\nabla_{A,\partial_r}\phi|^2
    \right)
    &=
    \int_{B_R(0)}
    \left(
        (d-4)|F_A|^2
        +
        (d-2)|\nabla_A\phi|^2
    \right).
\end{align*}
Finite total energy implies that there exists a sequence
$R_j\to\infty$ such that
\[
    R_j
    \int_{\partial B_{R_j}(0)}
    e(A,\phi)
    \longrightarrow0.
\]
Letting $R=R_j$ and then $j\to\infty$ proves
\eqref{eq:Euclidean_YMH_stress_energy}. For $d=3$, this is the displayed
equipartition identity. For $d=4$, it gives
$\nabla_A\phi=0$, and the connection equation then reduces to the
Yang--Mills equation. If $\phi$ vanishes at one point, the fact that $\phi$
is parallel gives $\phi\equiv0$. For $d>4$, both coefficients in
\eqref{eq:Euclidean_YMH_stress_energy} are positive, so
$F_A=0$ and $\nabla_A\phi=0$.
\end{proof}

\begin{remark}[Terminology]
The identity \eqref{eq:Euclidean_YMH_stress_energy} is obtained from the
stress-energy identity by testing against the radial dilation vector
field. Equivalently, it is the scaling identity for the Euclidean
Yang--Mills--Higgs energy. Its three-dimensional specialization is the
equipartition theorem of Jaffe--Taubes; see
\cite[Corollary~II.2.2]{Jaffe1980}.
\end{remark}

\subsection{Localized Higgs energy identity}
\label{subsec:localized_Higgs_energy_identity}

Corollary~\ref{cor:codim4_almost_monotonicity} is the usual
codimension-four almost-monotonicity for the full Yang--Mills--Higgs
energy. In Section~\ref{sec: Li_concentration}, we also need a localized
monotonicity identity for the Higgs energy component alone. Although its
derivation is parallel to the proof of
Theorem~\ref{thm: monotonicity_crit_pt}, the resulting formula has
codimension-two scaling and a different error term. This error measures
the failure of the Higgs stress-energy tensor to be divergence-free in
the presence of curvature.

\begin{lemma}[Localized Higgs energy monotonicity]
\label{lem:localized_Higgs_energy_identity}
Let $(X^n,g)$ be an oriented Riemannian $n$-manifold, $n\geqslant4$, and let
$P$ be a principal $G$-bundle over $X$, where $G$ is compact. Let
$x\in X$, set $r_x:=d(x,\cdot)$, and assume that on $B_{r_0}(x)$ one has
\[
-C_Hr_x^2g
\leqslant
\mathrm{Hess}\,\!\left(\frac12r_x^2\right)-g
\leqslant
C_Hr_x^2g
\]
for some $r_0\in(0,\inj_g(x))$ and $C_H\geqslant0$. If
$(\nabla,\Phi)\in\mathscr{A}(P)\times\Gamma(\mathfrak{g}_P)$ solves \eqref{eq:YMH_second_order}, define
\[
    H_x(r)
    :=
    \int_{B_r(x)}
    |\nabla\Phi|^2.
\]
Then there exists a constant $c'\geqslant0$, depending only on $n$ and
$C_H$, such that, for every $0<r<r_0$,
\begin{align}
\frac{d}{dr}
\left(
e^{c'r^2}r^{2-n}H_x(r)
\right)
&\geqslant
2e^{c'r^2}r^{2-n}
\int_{\partial B_r(x)}
|\nabla_{\partial_{r_x}}\Phi|^2
\nonumber\\
&\quad
+
2e^{c'r^2}r^{1-n}
\int_{B_r(x)}
r_x
\left\langle
[\partial_{r_x}\lrcorner F_\nabla,\Phi],
\nabla\Phi
\right\rangle.
\label{eq:localized_Higgs_energy_monotonicity}
\end{align}
Consequently,
\begin{align}
\frac{d}{dr}
\left(
e^{c'r^2}r^{2-n}H_x(r)
\right)
&\geqslant
-
2e^{c'r^2}r^{1-n}
\left|
\int_{B_r(x)}
r_x
\left\langle
[\partial_{r_x}\lrcorner F_\nabla,\Phi],
\nabla\Phi
\right\rangle
\right|.
\label{ineq:localized_Higgs_energy_weighted}
\end{align}
\end{lemma}

\begin{proof}
Consider the stress-energy tensor of the Higgs field,
\[
    \widehat{\mathbb{T}}(v,w)
    :=
    \langle\nabla_v\Phi,\nabla_w\Phi\rangle
    -
    \frac12|\nabla\Phi|^2g(v,w).
\]
Its trace satisfies
\[
    2\langle\widehat{\mathbb{T}},g\rangle
    =
    (2-n)|\nabla\Phi|^2.
\]
Moreover, there exists a dimensional constant $C_2(n)>0$ such that
\[
    |\widehat{\mathbb{T}}|\leqslant C_2(n)|\nabla\Phi|^2.
\]
Using $\nabla^*\nabla\Phi=0$, one computes
\[
    D^*\widehat{\mathbb{T}}(v)
    =
    \left\langle
    [v\lrcorner F_\nabla,\Phi],
    \nabla\Phi
    \right\rangle.
\]
Apply the divergence theorem to the vector field
\[
    Y
    =
    r_x\partial_{r_x}
    =
    D\left(\frac12r_x^2\right)
\]
on $B_r(x)$. Since $Y=r\partial_{r_x}$ on $\partial B_r(x)$, the boundary
term is
\[
    \int_{\partial B_r(x)}
    \widehat{\mathbb{T}}(Y,\partial_{r_x})
    =
    r\int_{\partial B_r(x)}
    |\nabla_{\partial_{r_x}}\Phi|^2
    -
    \frac r2
    \int_{\partial B_r(x)}
    |\nabla\Phi|^2.
\]
On the other hand,
\[
    DY
    =
    \mathrm{Hess}\,\!\left(\frac12r_x^2\right),
\]
and hence
\[
    \langle\widehat{\mathbb{T}},DY\rangle
    =
    \langle\widehat{\mathbb{T}},g\rangle
    +
    \left\langle
    \widehat{\mathbb{T}},
    \mathrm{Hess}\,\!\left(\frac12r_x^2\right)-g
    \right\rangle.
\]
Thus, we obtain
\begin{align*}
r\int_{\partial B_r(x)}
|\nabla_{\partial_{r_x}}\Phi|^2
-
\frac r2 H_x'(r)
&=
\frac{2-n}{2}H_x(r)
+
\int_{B_r(x)}
\left\langle
\widehat{\mathbb{T}},
\mathrm{Hess}\,\!\left(\frac12r_x^2\right)-g
\right\rangle
\\
&\quad
-
\int_{B_r(x)}
r_x
\left\langle
[\partial_{r_x}\lrcorner F_\nabla,\Phi],
\nabla\Phi
\right\rangle.
\end{align*}
Rearranging gives the exact identity
\begin{align}
H_x'(r)
&=
2\int_{\partial B_r(x)}
|\nabla_{\partial_{r_x}}\Phi|^2
+
\frac{n-2}{r}H_x(r)
\nonumber\\
&\quad
-
\frac{2}{r}
\int_{B_r(x)}
\left\langle
\widehat{\mathbb{T}},
\mathrm{Hess}\,\!\left(\frac12r_x^2\right)-g
\right\rangle
\nonumber\\
&\quad
+
\frac{2}{r}
\int_{B_r(x)}
r_x
\left\langle
[\partial_{r_x}\lrcorner F_\nabla,\Phi],
\nabla\Phi
\right\rangle.
\label{eq:localized_Higgs_energy_exact}
\end{align}
By the Hessian assumption and the bound for $\widehat{\mathbb{T}}$, after increasing
$C_2(n)$ if necessary,
\[
\left|
\left\langle
\widehat{\mathbb{T}},
\mathrm{Hess}\,\!\left(\frac12r_x^2\right)-g
\right\rangle
\right|
\leqslant
C_2(n)C_Hr_x^2|\nabla\Phi|^2.
\]
Consequently,
\begin{align}
H_x'(r)
&\geqslant
2\int_{\partial B_r(x)}
|\nabla_{\partial_{r_x}}\Phi|^2
+
r^{-1}
\bigl(n-2-C_2(n)C_Hr^2\bigr)H_x(r)
\nonumber\\
&\quad
+
\frac{2}{r}
\int_{B_r(x)}
r_x
\left\langle
[\partial_{r_x}\lrcorner F_\nabla,\Phi],
\nabla\Phi
\right\rangle.
\label{ineq:Higgs_energy_unweighted}
\end{align}
Now differentiate
\[
    e^{c'r^2}r^{2-n}H_x(r).
\]
Using \eqref{ineq:Higgs_energy_unweighted}, we get
\begin{align*}
\frac{d}{dr}
\left(
e^{c'r^2}r^{2-n}H_x(r)
\right)
&\geqslant
2e^{c'r^2}r^{2-n}
\int_{\partial B_r(x)}
|\nabla_{\partial_{r_x}}\Phi|^2
\\
&\quad
+
\bigl(2c'-C_2(n)C_H\bigr)
e^{c'r^2}r^{3-n}H_x(r)
\\
&\quad
+
2e^{c'r^2}r^{1-n}
\int_{B_r(x)}
r_x
\left\langle
[\partial_{r_x}\lrcorner F_\nabla,\Phi],
\nabla\Phi
\right\rangle.
\end{align*}
Choose
\[
    c'\geqslant \frac12 C_2(n)C_H.
\]
Discarding the nonnegative middle term gives
\eqref{eq:localized_Higgs_energy_monotonicity}. Finally, discarding the
nonnegative boundary term gives
\eqref{ineq:localized_Higgs_energy_weighted}.
\end{proof}


\section{Bochner--Weitzenb\"ock formulae and
\texorpdfstring{$\varepsilon$}{epsilon}-regularity}\label{app: C}

Let $d_\nabla$ denote the exterior covariant derivative and
$\Delta_\nabla=d_\nabla d_\nabla^*+d_\nabla^*d_\nabla$ the associated Hodge
Laplacian.  We use the convention
$d_\nabla^*=(-1)^{n(p+1)+1}*d_\nabla*$ on $p$-forms, and denote by
$\nabla^*\nabla=-\sum_i\nabla^2_{e_i,e_i}$ the rough Laplacian.  For an
$\mathfrak{g}_P$-valued tensor $S$,
\[
\frac{1}{2}\Delta|S|^2
= \langle \nabla^*\nabla S, S\rangle
-|\nabla S|^2.
\]
With respect to a local orthonormal frame, the standard Weitzenb\"ock formulae
for $\mathfrak{g}_P$-valued $1$-forms and $2$-forms are
\[
(\Delta_\nabla a)_i
=
(\nabla^*\nabla a)_i
+
[(F_\nabla)_{ki},a_k]
+
\Ric_{ki}a_k,
\]
and
\[
(\Delta_\nabla\varphi)_{ij}
=
(\nabla^*\nabla\varphi)_{ij}
+
[F_{ki},\varphi_{kj}]
-
[F_{kj},\varphi_{ki}]
+
[\mathcal R_2(\varphi)]_{ij},
\]
where
\[
[\mathcal R_2(\varphi)]_{ij}
=
\varphi_{\Ric(e_i),e_j}
+
\varphi_{e_i,\Ric(e_j)}
+
\varphi_{e_k,R_{e_i,e_j}e_k}.
\]
We now record the Bochner identities and local estimates used in the paper.

\begin{lemma}\label{lem: bochner_formula}
Let $(\nabla,\Phi)\in\mathscr A(P)\times\Gamma(\mathfrak g_P)$ be a solution
of \eqref{eq:YMH_second_order}. Then
\begin{align}
\frac12\Delta|\nabla\Phi|^2
=&
-\Ric^g(\nabla\Phi,\nabla\Phi)
-
2\langle *[*F_\nabla,\nabla\Phi],\nabla\Phi\rangle
\nonumber\\
&-
|[\nabla\Phi,\Phi]|^2
-
|\nabla(\nabla\Phi)|^2 ,
\label{eq: gen_bochner_dphi}
\end{align}
and
\begin{align}
\frac12\Delta|F_\nabla|^2
=&
-\langle F_\nabla\circ(\Ric^g\wedge\mathrm{Id}+2R^g),F_\nabla\rangle
-
\langle[\nabla\Phi,\nabla\Phi],F_\nabla\rangle
\nonumber\\
&-
\sum_{i,j,k}
\langle[F_{ij},F_{jk}],F_{ki}\rangle
-
|[F_\nabla,\Phi]|^2
-
|\nabla F_\nabla|^2.
\label{eq: gen_bochner_curv}
\end{align}
Here $F_{ij}:=F_\nabla(e_i,e_j)$ with respect to a local orthonormal frame.
\end{lemma}

\begin{proof}
The identities follow from the standard Bochner--Weitzenb\"ock formulas,
combined with \eqref{eq:YMH_second_order}.

For the Higgs term, the Weitzenb\"ock formula for
$\mathfrak g_P$-valued $1$-forms gives
\[
\Delta_\nabla(\nabla\Phi)
=
\nabla^*\nabla(\nabla\Phi)
+
(\nabla\Phi)\circ\Ric^g
+
*[*F_\nabla,\nabla\Phi].
\]
Since $(\nabla,\Phi)$ solves \eqref{eq:YMH_second_order},
\[
\Delta_\nabla(\nabla\Phi)
=
[[\nabla\Phi,\Phi],\Phi]
-
*[*F_\nabla,\nabla\Phi].
\]
Taking the inner product with $\nabla\Phi$ and using Ad-invariance gives
\eqref{eq: gen_bochner_dphi}.

For the curvature term, the Bianchi identity and
\eqref{eq:YMH_second_order} imply
\[
\Delta_\nabla F_\nabla
=
d_\nabla d_\nabla^*F_\nabla
=
d_\nabla[\nabla\Phi,\Phi]
=
[[F_\nabla,\Phi],\Phi]
-
[\nabla\Phi,\nabla\Phi].
\]
Combining this with the Weitzenb\"ock formula for $\mathfrak g_P$-valued
$2$-forms gives \eqref{eq: gen_bochner_curv}.
\end{proof}

\begin{corollary}[Bochner-type estimates]\label{cor: bochner_estimate}
Let $(\nabla,\Phi)\in\mathscr A(P)\times\Gamma(\mathfrak g_P)$ solve
\eqref{eq:YMH_second_order}. Then:
\begin{itemize}
\item[(a)] If $\Ric_g\geqslant-\kappa_U g$ on a domain $U\subset X$ and
$\|F_\nabla\|_{L^\infty(U)}<\infty$, then
\begin{equation}\label{ineq: Bochner_ineq_nablaPhi}
\Delta|\nabla\Phi|^2
\lesssim
(\kappa_U+\|F_\nabla\|_{L^\infty(U)})|\nabla\Phi|^2
\quad\text{on }U.
\end{equation}

\item[(b)] If $\|R^g\|_{L^\infty(U)}\leqslant c_U$ on $U$ and
$e=e(\nabla,\Phi)$ given by \eqref{eq: YMH_density}, then
\begin{equation}\label{eq: bochner-estimate}
\Delta e
\lesssim
c_U e+e^{3/2}
\quad\text{on }U.
\end{equation}
\end{itemize}
\end{corollary}

We now record two standard consequences of these inequalities. The first is a
local estimate for the Higgs derivative.

\begin{lemma}\label{lem:bounded_curvature_moser}
Let $(\nabla,\Phi)$ be a solution of \eqref{eq:YMH_second_order}. There
exists $r_1>0$, depending only on the bounded geometry constants of
$(X,g)$, such that, for every $x\in X$ and $0<r\leqslant r_1$,
\begin{equation}\label{ineq:moser}
\sup_{B_{r/2}(x)}|\nabla\Phi|^2
\lesssim
\left(
\|F_\nabla\|_{L^\infty(B_r(x))}^{n/2}
+
\|\Ric\|_{L^\infty(B_r(x))}^{n/2}
+
 r^{-n}
\right)
\int_{B_r(x)}|\nabla\Phi|^2.
\end{equation}
\end{lemma}

\begin{proof}
Set $f:=|\nabla\Phi|^2$ and
\[
    A
    :=
    \|F_\nabla\|_{L^\infty(B_r(x))}
    +
    \|\Ric\|_{L^\infty(B_r(x))}.
\]
By \eqref{ineq: Bochner_ineq_nablaPhi}, one has
$\Delta f\lesssim Af$ on $B_r(x)$. Fix $y\in B_{r/2}(x)$ and put
$R:=r/2$. Then $B_R(y)\subset B_r(x)$. Apply the linear case of
Theorem~\ref{thm: gen_mean_value} on $B_R(y)$ with $d=n$; condition
\textup{(C1)} follows simply from inclusion of balls. Since the centre
$y$ belongs to $B_{R/4}(y)$, the theorem gives
\[
    f(y)
    \lesssim
    \left(A^{n/2}+R^{-n}\right)
    \int_{B_R(y)}f
    \lesssim
    \left(A^{n/2}+r^{-n}\right)
    \int_{B_r(x)}f.
\]
Taking the supremum over $y\in B_{r/2}(x)$ proves
\eqref{ineq:moser}.
\end{proof}

Combining the stress-energy monotonicity formula
Theorem~\ref{thm: monotonicity_crit_pt} with the Bochner estimate
\eqref{eq: bochner-estimate} gives the following classical
$\varepsilon$-regularity theorem for the total Yang--Mills--Higgs energy
density.

\begin{theorem}[See, e.g., {\cite[Theorem~B]{afuni2019regularity}}]
\label{thm: total_epsilon_regularity}
Let $(X^n,g)$ be an oriented Riemannian $n$-manifold of bounded geometry,
with $n\geqslant4$, and let $P\to X$ be a principal $G$-bundle, where
$G$ is compact. Then there are
constants $r_0>0$ and $\varepsilon_0>0$ such that the following holds. If
$(\nabla,\Phi)\in\mathscr{A}(P)\times\Gamma(\mathfrak{g}_P)$ solves \eqref{eq:YMH_second_order} and
\[
\varepsilon
:=
r^{4-n}
\mathcal E_{B_r(x)}(\nabla,\Phi)
<
\varepsilon_0
\]
for some $x\in X$ and $0<r\leqslant r_0$, then, for every integer
$j\geqslant0$,
\begin{equation}\label{ineq:estimates_coulomb}
\sup_{B_{r/2}(x)}
\left(
|\nabla^jF_\nabla|^2
+
|\nabla^{j+1}\Phi|^2
\right)
\lesssim_j
r^{-4-2j}\varepsilon
=
r^{-n-2j}\frac12\int_{B_r(x)}e(\nabla,\Phi).
\end{equation}
The constants depend only on $j$, $n$, the bounded geometry constants of
$(X,g)$ and the structure constants of $\mathfrak g$.
\end{theorem}

\begin{proof}
We include the proof, taking the opportunity to correct a minor radius
mistake in the proof of the corresponding estimate for higher derivatives in
\cite[Proposition~3.5]{fadel2020asymptotic}: smallness obtained on an interior ball yields a
Coulomb gauge only on a still smaller ball, rather than on the original
ball.

Fix $y\in B_{r/2}(x)$ and set $R:=r/2$. Then
$B_R(y)\subset B_r(x)$. We apply
Theorem~\ref{thm: gen_mean_value} to $f=e(\nabla,\Phi)$ on $B_R(y)$,
with $d=4$, $\tau=0$, $a_0=0$, and $a_1,a\lesssim1$.
If $B_s(z)\subset B_{R/2}(y)$, then
$B_{R/2}(z)\subset B_R(y)$; codimension-four almost-monotonicity,
applied at $z$ between the radii $s$ and $R/2$, therefore gives
\[
    s^{4-n}\int_{B_s(z)}e(\nabla,\Phi)
    \lesssim
    R^{4-n}\int_{B_R(y)}e(\nabla,\Phi),
\]
which is condition~\textup{(C1)}. The Bochner estimate
\eqref{eq: bochner-estimate} gives condition~\textup{(C2)} with the
critical exponent $\alpha=3/2=(d+2)/d$.

Moreover,
\[
    R^{4-n}\int_{B_R(y)}e(\nabla,\Phi)
    \leqslant
    2^{n-3}\varepsilon,
\]
because
$\mathcal E_{B_r}(\nabla,\Phi)=\frac12\int_{B_r(x)}e(\nabla,\Phi)$.
Thus, after choosing
\[
    \varepsilon_0
    \leqslant
    2^{3-n}\hbar a^{-2},
\]
the critical smallness condition in
Theorem~\ref{thm: gen_mean_value} is satisfied. Since
$y\in B_{R/4}(y)$, the theorem yields
\[
    e(\nabla,\Phi)(y)
    \lesssim
    \left(a_1^2+R^{-4}\right)
    R^{4-n}\int_{B_R(y)}e(\nabla,\Phi)
    \lesssim
    r^{-4}\varepsilon,
\]
where $r_0$ is chosen sufficiently small, depending only on bounded
geometry. Taking the supremum over $y\in B_{r/2}(x)$ proves
\eqref{ineq:estimates_coulomb} for $j=0$. The constants depend only on
$n$, the bounded geometry constants of $(X,g)$, and the structure
constants of $\mathfrak g$.

It remains to prove the estimates for $j\geqslant1$. We first obtain the
zeroth-order estimate on a slightly larger interior ball. Let
$y\in B_{3r/4}(x)$. Then $B_{r/4}(y)\subset B_r(x)$, and hence
\[
\left(\frac r4\right)^{4-n}
\mathcal E_{B_{r/4}(y)}(\nabla,\Phi)
\lesssim_n
\varepsilon.
\]
After decreasing $\varepsilon_0$ by a fixed dimensional factor, the case
$j=0$, applied on $B_{r/4}(y)$, gives
\[
e(\nabla,\Phi)(y)
\leqslant
\sup_{B_{r/8}(y)}e(\nabla,\Phi)
\lesssim
r^{-4}\varepsilon.
\]
Since $y\in B_{3r/4}(x)$ was arbitrary,
\begin{equation}\label{ineq:interior_e_bound_for_bootstrap}
\sup_{B_{3r/4}(x)}
\left(
|F_\nabla|^2+|\nabla\Phi|^2
\right)
\lesssim
r^{-4}\varepsilon.
\end{equation}
Fix $p>n/2$. The bounded geometry volume estimates and
\eqref{ineq:interior_e_bound_for_bootstrap} imply
\[
r^{2-n/p}
\|F_\nabla\|_{L^p(B_{3r/4}(x))}
\lesssim_p
\varepsilon^{1/2},
\qquad
\|F_\nabla\|_{L^{n/2}(B_{3r/4}(x))}
\lesssim
\varepsilon^{1/2}.
\]
After decreasing $\varepsilon_0$ once more, the second quantity is below
Uhlenbeck's constant. Upon choosing normal coordinates and rescaling to
unit size, \cite[Theorem~1.3]{Uhlenbeck1982a} therefore gives a Coulomb
gauge on $B_{2r/3}(x)$ in which the connection form $A$ satisfies
\[
d^*A=0,
\qquad
r^{1-n/p}\|A\|_{L^p(B_{2r/3}(x))}
+
r^{2-n/p}\|\nabla A\|_{L^p(B_{2r/3}(x))}
\lesssim_p
r^{2-n/p}\|F_\nabla\|_{L^p(B_{3r/4}(x))}
\lesssim_p
\varepsilon^{1/2}.
\]
In this gauge, the Yang--Mills--Higgs system
\eqref{eq:YMH_second_order} is an elliptic
system. For the uniformity of the estimates with respect to the size of
$\Phi$, it is useful to write the corresponding covariant equations for
$F_\nabla$ and $a:=\nabla\Phi$ schematically as
\[
    \nabla^*\nabla F_\nabla+\mathcal A_\Phi F_\nabla
    =R^g*F_\nabla+F_\nabla*F_\nabla+a*a,
\]
\[
    \nabla^*\nabla a+\mathcal A_\Phi a
    =R^g*a+F_\nabla*a,
\]
where
\[
    \mathcal A_\Phi U:=-[[U,\Phi],\Phi],
    \qquad
    \langle\mathcal A_\Phi U,U\rangle=|[U,\Phi]|^2\geqslant0.
\]
After differentiating these equations, the terms containing an
undifferentiated $\Phi$ are absorbed, using Ad-invariance and Young's
inequality, by the coercive commutator terms associated with
$\mathcal A_\Phi$. The remaining terms involve only the ambient curvature
and lower covariant derivatives of $F_\nabla$ and $a$. The fixed inclusion
\[
    B_{r/2}(x)\Subset B_{2r/3}(x),
\]
together with the estimate in Coulomb gauge and
\eqref{ineq:interior_e_bound_for_bootstrap}, therefore allows the standard
interior covariant elliptic estimates to be iterated. Consequently, for
every integer $j\geqslant1$,
\[
\sup_{B_{r/2}(x)}
\left(
|\nabla^jF_\nabla|^2
+
|\nabla^{j+1}\Phi|^2
\right)
\lesssim_j
r^{-4-2j}\varepsilon.
\]
Together with the case $j=0$, this proves
\eqref{ineq:estimates_coulomb} for all $j\geqslant0$.
\end{proof}


\bibliography{references}

@article{fadeloliveira2026limitv5,
  author = {Fadel, D. and Oliveira, G.},
  title = {The limit of large mass monopoles},
  year = {2026},
  eprint = {1803.04117v5},
  archivePrefix = {arXiv},
  primaryClass = {math.DG},
  note = {Self-contained corrected and expanded version of the published article}
}

@article{fadel2026abelian,
  author = {Fadel, D.},
  title = {The large mass limit of monopoles: abelian limits and {D}irac singularities},
  year = {2026},
  eprint = {2607.29667v1},
  archivePrefix = {arXiv},
  primaryClass = {math.DG}
}

@unpublished{cheng2025su2,
  author = {Cheng, D. R. and Fadel, D. and Lara, L.},
  title = {{SU(2)} Yang--Mills--Higgs functional with {H}iggs self-interaction on 3-manifolds},
  eprint = {2505.08076},
  archivePrefix = {arXiv},
  primaryClass = {math.DG},
  year = {2025}
}

@article{stein2023invariant,
  title={$\mathrm{SU}(2)^2$-invariant Gauge Theory on Asymptotically Conical {C}alabi--{Y}au 3-Folds},
  author={Stein, J.},
  journal={The Journal of Geometric Analysis},
  volume={33},
  pages={121},
  year={2023}
}

@article{AronszajnKrzywickiSzarski1962,
    author  = {Aronszajn, N. and Krzywicki, A. and Szarski, J.},
    title   = {A unique continuation theorem for exterior differential forms on {R}iemannian manifolds},
    journal = {Arkiv f{\"o}r Matematik},
    volume  = {4},
    pages   = {417--453},
    year    = {1962},
    doi     = {10.1007/BF02591624}
}

@article{GarofaloLin1986,
    author  = {Garofalo, N. and Lin, F.-H.},
    title   = {Monotonicity properties of variational integrals, {$A_p$} weights and unique continuation},
    journal = {Indiana University Mathematics Journal},
    volume  = {35},
    number  = {2},
    pages   = {245--268},
    year    = {1986},
    doi     = {10.1512/iumj.1986.35.35015}
}

@article{GarofaloLin1987,
    author  = {Garofalo, N. and Lin, F.-H.},
    title   = {Unique continuation for elliptic operators: a geometric-variational approach},
    journal = {Communications on Pure and Applied Mathematics},
    volume  = {40},
    number  = {3},
    pages   = {347--366},
    year    = {1987},
    doi     = {10.1002/cpa.3160400305}
}

@book{HanLin2011,
    author    = {Han, Qing and Lin, Fanghua},
    title     = {Elliptic Partial Differential Equations},
    edition   = {2},
    series    = {Courant Lecture Notes in Mathematics},
    volume    = {1},
    publisher = {Courant Institute of Mathematical Sciences, New York University; American Mathematical Society},
    address   = {New York; Providence, RI},
    year      = {2011}
}

@article{li2025large,
  author        = {Li, Y.},
  title         = {The large mass limit of {$\mathrm G_2$} and {C}alabi--{Y}au monopoles},
  year          = {2025},
  eprint        = {2503.12075v1},
  archivePrefix = {arXiv},
  primaryClass  = {math.DG},
  note          = {Preprint}
}

@article{parise2025nonabelian,
  author  = {Parise, D. and Pigati, A. and Stern, D.},
  title   = {Nonabelian {Y}ang--{M}ills--{H}iggs and {P}lateau's problem in codimension three},
  journal = {J. Reine Angew. Math.},
  volume  = {837},
  year    = {2026},
  pages   = {253--288},
  doi     = {10.1515/crelle-2026-0038}
}

@misc{fadel2026ebook,
  author       = {Fadel, D. and S{\'a} Earp, H. N.},
  title        = {Introduction to {Y}ang--{M}ills theory in higher dimensions},
  year         = {2026},
  note = {Revised and expanded version of the 32$^\mathrm{o}$ Col{\'o}quio Brasileiro de Matem{\'a}tica lecture notes ``Gauge Theory in Higher Dimensions'', available from the author's webpage \url{https://drive.google.com/file/d/1R5CYfxZm8lR39807gGwgFZJMvFJZa7cw/view}}
}

@article{hohloch2009hypercontact,
	title="{Hypercontact structures and Floer homology}",
	author={Hohloch, S. and Noetzel, G. and Salamon, D. A.},
	journal={Geometry \& Topology},
	volume={13},
	number={5},
	pages={2543--2617},
	year={2009},
	publisher={Mathematical Sciences Publishers}
}

@article{fadel2023asymptotics,
  author  = {Fadel, D.},
  title   = {Asymptotics of finite energy monopoles on {AC} $3$-manifolds},
  journal = {J. Geom. Anal.},
  volume  = {33},
  number  = {1},
  year    = {2023},
  pages   = {Paper No.~17, 70 pp.},
  doi     = {10.1007/s12220-022-01095-8}
}

@article{fadel2020asymptotic,
  author    = {Fadel, D. and Nagy, {\'A}. and Oliveira, G.},
  title     = {The asymptotic geometry of {$\mathrm G_2$}-monopoles},
  journal   = {Mem. Amer. Math. Soc.},
  volume    = {303},
  number    = {1521},
  year      = {2024},
  pages     = {v+85},
  doi       = {10.1090/memo/1521},
  publisher = {American Mathematical Society}
}

@article{oliveira2016calabi,
	title={{C}alabi--{Y}au monopoles for the {S}tenzel metric},
	author={Oliveira, G.},
	journal={Communications in Mathematical Physics},
	volume={341},
	number={2},
	pages={699--728},
	year={2016},
	publisher={Springer}
}

@book{wehrheim2004uhlenbeck,
	title={Uhlenbeck compactness},
	author={Wehrheim, K.},
	volume={1},
	year={2004},
	publisher={European Mathematical Society}
}

@article{wehrheim2005energy,
	title={Energy quantization and mean value inequalities for nonlinear boundary value problems},
	author={Wehrheim, K.},
	journal={Journal of the European Mathematical Society},
	volume={7},
	number={3},
	pages={305--318},
	year={2005},
	publisher={European Matemathical Society, EMS}
}

@article{afuni2019regularity,
	title={Regularity and vanishing theorems for {Y}ang--{M}ills--{H}iggs pairs},
	author={Afuni, A.},
	journal={Archiv der Mathematik},
	volume={112},
	number={5},
	pages={547--558},
	year={2019},
	publisher={Springer}
}

@article{walpuski2017g2,
	title={$\rm{G}_2$-instantons, associative submanifolds and {F}ueter sections},
	author={Walpuski, T.},
	journal={Communications in Analysis and Geometry},
	volume={25},
	number={4},
	pages={847--893},
	year={2017},
	publisher={International Press of Boston}
}

@article{fadel2019limit,
	title={The limit of large mass monopoles},
	author={Fadel, D. and Oliveira, G.},
	journal={Proceedings of the London Mathematical Society},
	volume={119},
	number={6},
	pages={1531--1559},
	year={2019},
	publisher={Wiley Online Library}
}

@article{van2009regularity,
  title         = {Regularity of asymptotically conical {R}icci-flat {K}\"ahler metrics},
  author        = {van Coevering, C.},
  year          = {2010},
  eprint        = {0912.3946},
  archivePrefix = {arXiv},
  primaryClass  = {math.DG},
  doi           = {10.48550/arXiv.0912.3946}
}

@phdthesis{oliveira2014thesis,
	title={Monopoles in higher dimensions},
	author={Oliveira, G.},
	year={2014},
	school={Imperial College London}
}

@article{walpuski2017compactness,
	title={A compactness theorem for {F}ueter sections},
	author={Walpuski, T.},
	journal={Commentarii Mathematici Helvetici},
	volume={92},
	number={4},
	pages={751--776},
	year={2017}
}

@article{oliveira2014monopoles,
	title={Monopoles on the {B}ryant--{S}alamon $\rm{G}_2$-manifolds},
	author={Oliveira, G.},
	journal={Journal of Geometry and Physics},
	volume={86},
	pages={599--632},
	year={2014},
	publisher={Elsevier}
}

@article{tian2000gauge,
	title={Gauge theory and calibrated geometry, {I}},
	author={Tian, G.},
	journal={Annals of Mathematics},
	volume={151},
	number={1},
	pages={193--268},
	year={2000},
	publisher={Princeton University}
}

@book{donaldson1990geometry,
	title={The geometry of four--manifolds},
	author={Donaldson, S. K. and Kronheimer, P. B.},
	year={1990},
	publisher={Clarendon Press}
}

@incollection{donaldson2009gauge,
  title         = {Gauge theory in higher dimensions, {II}},
  author        = {Donaldson, S. K. and Segal, E.},
  booktitle     = {Geometry of Special Holonomy and Related Topics},
  series        = {Surveys in Differential Geometry},
  volume        = {16},
  pages         = {1--42},
  publisher     = {International Press},
  address       = {Somerville, MA},
  year          = {2011},
  eprint        = {0902.3239},
  archivePrefix = {arXiv},
  primaryClass  = {math.DG}
}

@article{harvey1982calibrated,
	title={Calibrated geometries},
	author={Harvey, R. and Lawson, H. B.},
	journal={Acta Mathematica},
	volume={148},
	number={1},
	pages={47--157},
	year={1982},
	publisher={Springer}
}

@article{walpuski2014spin,
  title   = {$\mathrm{Spin}(7)$-instantons, {C}ayley submanifolds and {F}ueter sections},
  author  = {Walpuski, T.},
  journal = {Communications in Mathematical Physics},
  volume  = {352},
  number  = {1},
  pages   = {1--36},
  year    = {2017},
  doi     = {10.1007/s00220-016-2735-0},
  eprint  = {1409.6705},
  archivePrefix = {arXiv},
  primaryClass  = {math.DG}
}

@book{Jaffe1980,
	Address = {Mass.},
	Author = {Jaffe, A. and Taubes, C.~H.},
	Isbn = {3-7643-3025-2},
	Mrclass = {81E10 (53C80 81-02)},
	Mrnumber = {MR614447 (82m:81051)},
	Mrreviewer = {Masatsugu Minami},
	Note = {Structure of static gauge theories},
	Pages = {v+287},
	Publisher = {Birkh\"auser Boston},
	Series = {Progress in Physics},
	Title = {Vortices and monopoles},
	Volume = {2},
	Year = {1980}}

@article{Taubes1999,
	Author = {Taubes, C.~H.},
	Coden = {JDGEAS},
	Fjournal = {Journal of Differential Geometry},
	Issn = {0022-040X},
	Journal = {J. Differential Geom.},
	Mrclass = {53D45 (53D35 57R17 57R57)},
	Mrnumber = {MR1761081 (2002i:53119)},
	Mrreviewer = {Yi-Jen Lee},
	Number = {3},
	Pages = {453--609},
	Title = {{${\rm GR}={\rm SW}$}: counting curves and connections},
	Url = {http://projecteuclid.org/getRecord?id=euclid.jdg/1214425348},
	Volume = {52},
	Year = {1999}}

@article{Taubes1999a,
	Author = {Taubes, C.~H.},
	Coden = {JDGEAS},
	Fjournal = {Journal of Differential Geometry},
	Issn = {0022-040X},
	Journal = {J. Differential Geom.},
	Mrclass = {53D45 (57R57)},
	Mrnumber = {MR1728301 (2000i:53123)},
	Mrreviewer = {Ignasi Mundet-Riera},
	Number = {2},
	Pages = {203--334},
	Title = {{${\rm Gr}\Rightarrow{\rm SW}$}: from pseudo-holomorphic curves to {S}eiberg-{W}itten solutions},
	Url = {http://projecteuclid.org/getRecord?id=euclid.jdg/1214425068},
	Volume = {51},
	Year = {1999}}

@article{Taubes1996,
	Author = {Taubes, C.~H.},
	Doi = {10.1090/S0894-0347-96-00211-1},
	Fjournal = {Journal of the American Mathematical Society},
	Issn = {0894-0347},
	Journal = {J. Amer. Math. Soc.},
	Mrclass = {57R57 (53C15 58D10 58D27 58G30)},
	Mrnumber = {MR1362874 (97a:57033)},
	Mrreviewer = {Dietmar A. Salamon},
	Number = {3},
	Pages = {845--918},
	Title = {{${\rm SW}\Rightarrow{\rm Gr}$}: from the {S}eiberg--{W}itten equations to pseudo-holomorphic curves},
	Url = {http://dx.doi.org/10.1090/S0894-0347-96-00211-1},
	Volume = {9},
	Year = {1996}}

@article{Uhlenbeck1982a,
	Author = {Uhlenbeck, K.~K.},
	Coden = {CMPHAY},
	Fjournal = {Communications in Mathematical Physics},
	Issn = {0010-3616},
	Journal = {Comm. Math. Phys.},
	Mrclass = {53C05 (49F10 58E20 81E10)},
	Mrnumber = {MR648356 (83e:53035)},
	Mrreviewer = {Wolfgang L{\"u}cke},
	Number = {1},
	Pages = {31--42},
	Title = {{Connections with {$L^p$} bounds on curvature}},
	Url = {http://projecteuclid.org/getRecord?id=euclid.cmp/1103920743},
	Volume = {83},
	Year = {1982}}

\end{document}